\documentclass[12pt,twoside,leqno,openany]{amsart}

\pdfoutput=1

\usepackage{amsbsy,amscd,amsfonts,amsmath,amssymb,amsthm,color,
fancybox,fancyhdr,footmisc,graphics,graphicx,ifthen,mathrsfs,
multicol,pdfpages,rotating,times,wasysym,pifont,blkarray}
\usepackage[dvipsnames,svgnames,x11names]{xcolor}

\usepackage[all]{xy}
\usepackage[utf8]{inputenc}
\usepackage[T1]{fontenc}
\usepackage{mathtools}
\newtagform{EngelLie}[\scriptstyle]{$}{$}
\makeatletter\newcommand{\leqnomode}{\tagsleft@true}
\newcommand{\reqnomode}{\tagsleft@false}\makeatother

\newtheorem{Theorem}[equation]{Theorem}

\newtheorem{Proposition}[equation]{Proposition}

\newtheorem{Lemma}[equation]{Lemma}

\newtheorem{Corollary}[equation]{Corollary}
\newtheorem{Assertion}[equation]{Assertion}

\newtheorem{Observation}[equation]{Observation}

\theoremstyle{definition}

\newtheorem{Definition}[equation]{Definition}

\newtheorem{Principle}[equation]{Principle}
\newtheorem{Question}[equation]{Question}

\newcommand{\C}{\mathbb{C}}

\newcommand{\N}{\mathbb{N}}

\newcommand{\R}{\mathbb{R}}

\newcommand{\KK}{\text{\sc k}}

\newcommand{\NN}{\text{\sc n}}

\newcommand{\ZZ}{\text{\sc z}}

\newcommand{\kaux}{{\text{\usefont{T1}{qcs}{m}{sl}k}}}

\newcommand{\maux}{{\text{\usefont{T1}{qcs}{m}{sl}m}}}

\newcommand{\Aaux}{{\text{\usefont{T1}{qcs}{m}{sl}A}}}
\newcommand{\Baux}{{\text{\usefont{T1}{qcs}{m}{sl}B}}}

\newcommand{\Faux}{{\text{\usefont{T1}{qcs}{m}{sl}F}}}

\newcommand{\Iaux}{{\text{\usefont{T1}{qcs}{m}{sl}I}}}
\newcommand{\Jaux}{{\text{\usefont{T1}{qcs}{m}{sl}J}}}

\newcommand{\Paux}{{\text{\usefont{T1}{qcs}{m}{sl}P}}}
\newcommand{\Qaux}{{\text{\usefont{T1}{qcs}{m}{sl}Q}}}
\newcommand{\Raux}{{\text{\usefont{T1}{qcs}{m}{sl}R}}}

\newcommand{\Vaux}{{\text{\usefont{T1}{qcs}{m}{sl}V}}}
\newcommand{\Waux}{{\text{\usefont{T1}{qcs}{m}{sl}W}}}

\definecolor{blue}{cmyk}{1.,1.,0.,0.63}
\definecolor{red}{cmyk}{0.,1.,1.,0.63}
\definecolor{green}{cmyk}{1.,0.,1.,0.63}
\definecolor{black}{cmyk}{1.,1.,1.,1.}

\newcommand{\red}{\textcolor{red}}

\makeatletter
\renewcommand{\@fnsymbol}[1]
{\ensuremath{\ifcase#1\or $*$\or $**$\or $***$\or $****$\or $*****$
\else\@ctrerr\fi}}
\makeatother

\newcommand{\HEAD}[2]{%
\pagestyle{fancy}
\fancyhead[RO]{\tiny\sf\thepage}
\fancyhead[CO]{{\tiny\sf #1}}
\fancyhead[LE]{\tiny\sf\thepage}
\fancyhead[CE]{{\tiny\sf #2}}
\fancyfoot{}}

\numberwithin{equation}{section}

\newcommand{\Section}[1]{
\renewcommand{\thesection}{\bf\arabic{section}}
\section{#1}
\renewcommand{\thesection}{\arabic{section}}}

\newcommand{\Subsection}[1]{
\refstepcounter{equation}
\medskip\noindent{\bf\arabic{section}.\arabic{equation}.~#1.}}

\newcommand{\style}[1]{\text{\footnotesize{\sf #1}}}

\newcommand{\stylesmall}[1]{{\sf #1}}

\renewcommand{\arg}{\style{arg}}

\newcommand{\Aut}{\style{Aut}}

\renewcommand{\cos}{\style{cos}}

\renewcommand{\deg}{\style{deg}}

\renewcommand{\dim}{\style{dim}}
\newcommand{\dimsmall}{\stylesmall{dim}}

\renewcommand{\exp}{\style{exp}}

\renewcommand{\Im}{\style{Im}}

\newcommand{\Levi}{\style{Levi}}

\renewcommand{\lim}{\style{lim}}

\renewcommand{\log}{\style{log}}

\newcommand{\ord}{\style{ord}}

\renewcommand{\Re}{\style{Re}}

\renewcommand{\sin}{\style{sin}}

\newcommand{\Span}{\style{Span}}

\newcommand{\bdot}{{\boldsymbol{\cdot}}}

\newcommand{\centersmallbullet}{{}_{{}^{{}^{
\scriptscriptstyle{\bullet\!}}}}}

\newcommand{\Hall}{\Hall}

\newcommand{\medcup}{\mathbin{\scalebox{1.5}{\ensuremath{\cup}}}}

\newcommand{\smallbullet}{{\scriptscriptstyle{\bullet}}}

\newcommand{\smallsum}[1]{
\underset{#1}{\raisebox{1pt}{$\sum$\,}}
}

\newcommand{\vf}{\vfill\end{document}}

\newcommand{\zero}[1]{\underline{#1}_{\red{\circ}}}

\newcommand{\Lbar}{\overline{\mathcal{L}_1}}
\newcommand{\Kbar}{\overline{\mathcal{K}}}
\newcommand{\Lbark}{\overline{\mathcal{L}_1}(\kaux)}
\newcommand{\LbarLbark}{\overline{\mathcal{L}_1} \, \overline{\mathcal{L}_1}(\kaux)}
\newcommand{\KLbark}{\mathcal{K} \, \overline{\mathcal{L}_1}(\kaux)}
\newcommand{\LLbark}{\mathcal{L}_1 \, \overline{\mathcal{L}_1}(\kaux)}
\newcommand{\LbarLk}{ \overline{\mathcal{L}_1} \, \mathcal{L}_1(\kaux)}
\newcommand{\LbarLbarLbark}{\overline{\mathcal{L}_1} \, \overline{\mathcal{L}_1} \, \overline{\mathcal{L}_1}(\kaux)}
\newcommand{\KLbar}{\mathcal{K} \, \overline{\mathcal{L}_1}}
\newcommand{\KLbarLbark}{\mathcal{K} \, \overline{\mathcal{L}_1} \, \overline{\mathcal{L}_1}(\kaux)}
\newcommand{\LbarLbarLk}{\overline{\mathcal{L}_1} \, \overline{\mathcal{L}_1} \, \mathcal{L}_1(\kaux)}
\newcommand{\Lkbar}{\mathcal{L}_1(\overline{\kaux})}
\newcommand{\LLkbar}{\mathcal{L}_1 \, \mathcal{L}_1(\overline{\kaux})}
\newcommand{\barP}{\overline{{\text{\usefont{T1}{qcs}{m}{sl}P}}}}

\newcommand{\bc}{\mathbb{C}}
\newcommand{\ch}{\mathscr{H}}
\newcommand{\cn}{\mathscr{N}}
\newcommand{\bz}{\mathbb{Z}}
\newcommand{\br}{\mathbb{R}}
\newcommand{\rw}{\longrightarrow}

\begin{document}

\setcounter{section}{0}

$\:$

\bigskip\bigskip\bigskip\bigskip\bigskip

\begin{center}

{\large\bf Normal Forms\footnotemark[1]
for Rigid $\mathfrak{C}_{2,1}$ 
Hypersurfaces $M^5 \subset \C^3$}
\label{normal-forms-rigid-M5-C3}

\footnotetext[1]{\,
This work was supported
in part by the Polish National Science Centre (NCN) 
via the grant number 2018/29/B/ST1/02583.}

\medskip

\bigskip\bigskip

Zhangchi {\sc Chen}\footnotemark[2],
Wei Guo {\sc Foo}\footnotemark[3],
Joël~{\sc Merker}\footnotemark[2],
The Anh {\sc Ta}\footnotemark[2]

\footnotetext[2]{\,\,Laboratoire de Mathématiques d'Orsay,
CNRS, Université Paris-Saclay, 91405 Orsay Cedex,
France.}

\footnotetext[3]{\,\,Hua Loo-Keng Center for Mathematical
Sciences, AMSS, Chinese Academy of Sciences, Beijing.}

\medskip

\hfill
{\footnotesize\sf\em
Dedicated to the memory of Alexander 
Isaev\,\dots}\footnotemark[4]

\footnotetext[4]{\,
\dots who, in February 2019, visited Orsay University and
with his energetic character, gave impetus,
and fostered with 
breadth 
exciting exchanges about 
relationships between CR geometry and Affine geometry.}

\end{center}\bigskip

\begin{center}
\begin{minipage}[t]{12.5cm}
\parindent 0.53cm
\footnotesize
\noindent
{\sc Abstract}.
Consider a $2$-nondegenerate constant
Levi rank $1$ rigid $\mathcal{C}^\omega$
hypersurface $M^5 \subset \mathbb{C}^3$
in coordinates $(z, \zeta, w = u + iv)$:
\[
u
\,=\,
F\big(z,\zeta,\overline{z},\overline{\zeta}\big).
\]
The Gaussier-Merker model $u = \frac{z \overline{z} +
\frac{1}{2} z^2 \overline{\zeta} + \frac{1}{2}
\overline{z}^2 \zeta}{1 - \zeta \overline{\zeta}}$
was shown by Fels-Kaup 2007 to be locally CR-equivalent to
the light cone $\{ x_1^2 + x_2^2 - x_3^2 = 0\}$.
Another representation is the tube $u = \frac{x^2}{1-y}$.

Inspired by Alexander Isaev, we study {\sl rigid} biholomorphisms:
\[
(z,\zeta,w)
\,\,\,\longmapsto\,\,\,
\big(
f(z,\zeta),\,
g(z,\zeta),\,
\rho\,w+h(z,\zeta)
\big)
\,\,=:\,\,
(z',\zeta',w').
\]
The G-M model has $7$-dimensional rigid automorphisms group.

A Cartan-type reduction to an $\{e\}$-structure was
done by Foo-Merker-Ta in arxiv.org/abs/1904.02562/.
Three relative invariants appeared: $\Vaux_0$, $\Iaux_0$
(primary) and $\Qaux_0$ (derived).
In Pocchiola's formalism,
Section~8 provides a finalized expression for $\Qaux_0$.

The goal is to establish the Poincar\'e-Moser complete normal form:
\[
u
\,=\,
\frac{z\overline{z}+\frac{1}{2}\,z^2\overline{\zeta}
+\frac{1}{2}\,\overline{z}^2\zeta}{
1-\zeta\overline{\zeta}}
+
\sum_{a,b,c,d\in\N
\atop
a+c\geqslant 3}\,
G_{a,b,c,d}\,
z^a\zeta^b\overline{z}^c\overline{\zeta}^d,
\]
with $0 = G_{a,b,0,0} = G_{a,b,1,0} = G_{a,b,2,0}$ and
$0 = G_{3,0,0,1} = {\rm Im}\, G_{3,0,1,1}$.

We apply the method of Chen-Merker
arxiv.org/abs/1908.07867 to catch (relative) invariants at
{\em every} point, not only at the central point,
as the coefficients $G_{0,1,4,0}$, $G_{0, 2, 3, 0}$, 
${\rm Re}\, G_{3,0,1,1}$.
With this, a complete {\sl brige}
Poincar\'e $\longleftrightarrow$ Cartan is constructed.

In terms of $F$, the numerators of $\Vaux_0$, $\Iaux_0$,
$\Qaux_0$ incorporate $11$, $52$, $824$ differential monomials.
\hfill
{\scriptsize
[Message~to~the~busy~reader:~Section~{\ref{introduction-normal-form-rigid}}~explains~and~summarizes~all~the~ideas.]}
\end{minipage}
\end{center}

\Section{\bf Introduction}
\label{introduction-normal-form-rigid}
\HEAD{{\ref{introduction-normal-form-rigid}}.~{\sf Introduction}
}{
Zhangchi {\sc Chen}, Wei Guo {\sc Foo}, Joël {\sc Merker}, 
The Anh {\sc Ta}}

The problem of equivalence for CR manifolds was begun by Poincaré in
1907, who, by a plain counting argument, pointed out that real
hypersurfaces $M^3 \subset \C^2$ must {\em a priori} possess
infinitely many {\em invariants} under biholomorphic transformations.

\smallskip

\hfill
\begin{minipage}[t]{14.25cm}
\baselineskip=0.37cm\parindent=0.37cm
{\scriptsize{
Nous pourrons [\dots] supposer que $F$ est de la forme
\[
F
\,=\,
X
-
\Phi(Y,X,X'),
\]
et il y a alors
\[
N'
\,=\,
{\textstyle{\frac{(n+1)(n+2)(n+3)}{6}}}
-
1
\]
coefficients arbitraires r\'eels [\dots].
Enfin, les \'equations de la transformation peuvent s'\'ecrire
\def\theequation{3}\begin{equation}
Z
\,=\,
\psi(z,z'),
\ \ \ \ \ \ \ \ \ \ \ \ \ \ \ \ \ \ \ \
Z'
\,=\,
\psi_1(z,z'),
\end{equation}
$\psi$ et $\psi_1$ \'etant deux fonctions analytiques 
complexes d\'eveloppables suivant les puissances de $z$ 
et de $z'$: nous avons besoin des termes jusqu'au
$n^{\rm e}$ ordre, ce qui fait
\[
2\,
\big[
{\textstyle{
\frac{(n+1)(n+2)}{2}}}
-
1
\big]
\]
coefficients arbitraires complexes, ou, ce qui revient au m\^eme,
\[
N''
\,=\,
2\,n^2
+
6\,n
\]
coefficients arbitraires r\'eels que nous appellerons les coefficients
$C$.
\hfill{\cite[pp.~194--195]{Poincare-1907}}
}}
\end{minipage}

\smallskip

Thus in $\C^2$, there are more
hypersurfaces, namely $\sim \frac{n^3}{6}$, 
than there are biholomorphisms, namely $\sim 2\, n^2$,
did argue Poincaré.

As in the theory that Lie erected in the end
of the XIX\textsuperscript{th} Century with his
students Engel, Scheffers, Kowalevski and others, 
the existence of (local) invariants creates a (local) classification
problem, not even terminated nowadays for hypersurfaces in $\C^3$.

Analogously, given the action of a finite-dimensional Lie group on a
manifold $M$ which induces an action on (local) graphs embedded in
$M$, Lie discovered that prolongations of the $G$-action to jet
bundles of sufficiently high order automatically create infinitely
many differential invariants~{\cite{Lie-Merker-2015, Olver-1995}},
hence various classification problems can be undertaken.

Throughout all of this memoir, concentrated on CR geometry, all CR
manifolds will be assumed real analytic ($\mathcal{C}^\omega$).  An
elementary complex Frobenius theorem proved {\em e.g.} by Paulette
Libermann in~{\cite{Libermann-1955}}, guarantees embedabbility in some
$\C^\NN$.  We will restrict ourselves to the definite class of {\em
hypersurfaces} $M^{2n+1} \subset \C^{n+1}$, which are automatically
CR. Results for {\em embedded} hypersurfaces $M^{2n+1} \subset
\C^{n+1}$ of class $\mathcal{C}^\infty$ or $\mathcal{C}^\KK$ with $\KK
\gg 1$ sufficiently high can be formulated, and proofs easily adapted.
In fact, only $\mathcal{C}^\omega$ hypersurfaces $M^3 \subset \C^2$
and $M^5 \subset \C^3$ will be studied here.

The interest of studying {\sl rigidly equivalent}\,\,---\,\,in 
Alexander Isaev's terminology\,\,---\,\,{\sl rigid} hypersurfaces 
was pointed out to us during his February 2019 stay in Orsay.
In recent publications~{\cite{Isaev-2016, 
Isaev-2016-bis,
Isaev-2018, 
Isaev-2019}}, Alexander tackled
to {\em integrate} Pocchiola's zero
CR curvature equations $\Waux = 0 = \Jaux$ 
of tube and rigid
$2$-nondegenerate constant Levi rank $1$ 
hypersurfaces 
$M^5 \subset \C^3$ (more will be said later).

A local hypersurface $M^{2n+1} 
\subset \C^{n+1}$ with coordinates $\ZZ = (\ZZ_1, \dots,
\ZZ_{n+1})$
is said to be {\sl rigid} if 
there exists an infinitesimal CR automorphism, namely a vector
field $T$ tangent to $M$ of the form $T = X + \overline{X}$ with a
nonzero holomorphic vector field $X = \sum_{i=1}^{n+1}\, a_i(\ZZ)\,
\partial_{\ZZ_i}$, which is {\em transversal} to the complex tangent
space $T^cM$ in the sense that $TM = T^cM \oplus \R T$.  After a local
biholomorphic straightening, one makes $X = i\,
\frac{\partial}{\partial w}$ with $w = \ZZ_{n+1}$, and tangency of $X
+ \overline{X} = 2\, \frac{\partial}{ \partial v}$ to $M$ shows that,
restricting considerations to dimensions $n+1 = 2, 3$, writing
coordinates $\C^2 \ni (z, w)$ and $\C^3 \ni (z,\zeta, w)$, the
right-hand side $\mathcal{C}^\omega$ graphing functions:
\[
M^3
\colon\ \ \
u
\,=\,
F(z,\overline{z}),
\ \ \ \ \ \ \ \ \ \ \ \ \ \ \ \ \ \ \ \ \ \ \ \ \ \
M^5
\colon\ \ \
u
\,=\,
F(z,\zeta,\overline{z},\overline{\zeta}),
\] 
are independent of $v$, where $w = u + i\, v$:

Alexander Isaev's concept of {\sl rigid biholomorphic transformation}
is less popular or widespread. In $\C^2$ and in $\C^3$, such are
biholomorphisms of the form:
\[
(z,w)
\,\,\longmapsto\,\,
\big(
f(z),\,\rho\,w+g(z)
\big),
\ \ \ \ \ \ \ \ \ \ \ 
(z,\zeta,w)
\,\,\longmapsto\,\,
\big(
f(z,\zeta),\,
g(z,\zeta),\,
\rho\,w+h(z,\zeta)
\big),
\]
where $f$, $g$, $h$ are holomorphic of their arguments, {\em
independently of $w$}, and where $\rho \in \R^\ast$. 
The interest is that rigid biholomorphisms trivially send rigid
hypersurfaces to rigid hypersurfaces: they respect the pre-given CR
symmetry, and much more will be explained later.

As Poincaré did, but without assuming that the origin is left fixed,
for any integer $d \geqslant 1$, writing $f(z) = \sum_{0 \leqslant k
\leqslant d}\, f_k\, z^k$ with $f_k \in \C$ and similarly $g(z) =
\sum\, g_k\, z^k$, the (rough) ``number'' of rigid biholomorphisms of
degree $\leqslant d$ is the number of incoming 
{\em real} parameters, namely
$2\,(d+1)+1+2\,(d+1) = 4\,d+5 \sim 4\,d$, while the (rough)
``number'' of rigid hypersurfaces $\big\{ u = \sum_{j+k\leqslant d}\,
F_{j,k}\, x^j y^k \big\}$ of degree $\leqslant d$ too, with $F_{j,k}
\in \R$, is equal to $\binom{d+2}{2} \sim \frac{1}{2}\, d^2$, hence
much larger as $d \longrightarrow \infty$.

Similarly in $\C^3$, the (rough) ``space'' of rigid biholomorphisms of
degree $\leqslant d$ is of real dimension:
\[
2\,
{\textstyle{\binom{d+2}{2}}}
+
2\,
{\textstyle{\binom{d+2}{2}}}
+
1
+
2\,
{\textstyle{\binom{d+2}{2}}}
\,=\,
3\,(d+2)(d+1)
+
1
\,\sim\,
3\,d^2,
\]
much smaller than the dimension of the ``space'' of hypersurfaces of
degree $\leqslant d$ too:
\[
{\textstyle{\binom{d+4}{4}}}
\,\sim\,
{\textstyle{\frac{1}{24}}}\,
d^4.
\]

\smallskip

To classify CR manifolds, two methods exist in the supermarket: 
that of Cartan, and that of Moser.

\smallskip

Cartan devised a quite sophisticated and proteiform
{\sl method of equivalence}. Given
a manifold $M$ equipped with a certain class of geometric,
say CR here, structures, Cartan's 
{\sl method of equivalence} consists in
constructing a bundle $\pi \colon P \longrightarrow M$
together with an {\sl absolute (co)parallelism}
on $P$, namely a coframe of
everywhere linearly independent $1$-forms
$\theta^1, \dots, \theta^{\dimsmall\, P}$
on $P$ such that:
\[
\xymatrix{
P
\ar[rr]^\Pi
\ar[d]_\pi
&
&
P'
\ar[d]^{\pi'}
\\
M
\ar[rr]_\Phi
&
&
M'
}
\]

\smallskip\noindent$\bullet$\,
every local CR diffeomorphism $\Phi \colon M \longrightarrow M'$
between two CR manifolds
lifts uniquely as a diffeomorphism 
$\Pi \colon P \longrightarrow P'$ satisfying
$\Pi^\ast {\theta'}^i = \theta^i$ for
$1 \leqslant i \leqslant \dim\, P$, with
$P'$ and the ${\theta'}^i$ similarly constructed;

\smallskip\noindent$\bullet$\,
conversely, every diffeomorphism $\Pi \colon
P \longrightarrow P'$ commuting with projections $\pi$, $\pi'$
whose horizontal part is a diffeomorphims $M \longrightarrow M'$
and which satisfies $\Pi^\ast {\theta'}^i = \theta^i$ for
$1 \leqslant i \leqslant \dim\, P$, has a horizontal
part which is {\em Cauchy-Riemann} diffeomorphism
(or, more generally, a diffeomorphism respecting the
considered geometric structure).

\smallskip

[Beyond, there can exist {\sl Cartan connections} associated
to (modifications of) $P \longrightarrow M$, but we
will not need this concept.]

Rexpressing the exterior differentials $d\theta^i$ and
$d{\theta'}^i$ from both sides in terms
of the basic $2$-forms provided by the two ambient coframes:
\[
d\theta^i
\,=\,
\sum_{j<k}\,
T_{j,k}^i(p)\,
\theta^j\wedge\theta^k
\ \ \ \ \ \ \ \ \ \ \ \ \ \ \ \ \ \ \ \
\text{and}
\ \ \ \ \ \ \ \ \ \ \ \ \ \ \ \ \ \ \ \
d{\theta'}^i
\,=\,
\sum_{j<k}\,
{T'}_{j,k}^i(p')\,
{\theta'}^j\wedge{\theta'}^k,
\]
certain {\sl structure functions} appear, defined for $p \in P$ and
for $p' \in P'$, and the 
exact pullback relations $\Pi^\ast {\theta'}^i = 
\theta^i$ force {\em individual invariancy} of all them:
\[
{T'}_{j,k}^i
\big(\Phi(p)\big)
\,=\,
T_{j,k}^i(p)
\eqno
{\scriptstyle{(\forall\,p\,\in\,P)}}.
\]

As is known, Cartan's method is computationally {\em extremely
intensive}, especially in CR geometry, where several
normalizations and prolongations are required.
Explicit expressions of intermediate torsion
coefficients which conduct to the final $T_{j,k}^i(p)$
grow dramatically in complexity. 

One reason for such a complexity is the presence of large isotropy
groups for the CR automorphisms groups of (standard) models, which
imposes a great number of steps.  Another reason is the {\em nonlinear}
character of differential algebraic polynomial expressions that must
be handled progressively.  The last reason is that Cartan's method
studies geometric structures {\em at every point} of
the base manifold, and there is a price to pay for this generality.

In most existing references ({\em cf.} the bibliography), the trick
that Cartan himself devised to avoid nonlinear complications while
retaining anyway some essential information, is the so-called {\sl
Cartan Lemma}.  It is explicit only at the level of linear
algebra. Even admitting to only deal with linear algebra computations,
as Chern always did, Cartan's method is often long and demanding.

\smallskip

\hfill
\begin{minipage}[t]{14.25cm}
\baselineskip=0.37cm\parindent=0.37cm
{\scriptsize{
In his works, Moser usually searched for wisdom rather
than simply knowledge, and thus he strongly emphasized
developments of methods and insights over pushing
a specific result to the limit. Accordingly, he sometimes
described the outcome of his own work as methods
rather than theorems.
\hfill{\cite[p.~1348]{Katok-Hasselblatt-2002}}
}}
\end{minipage}

\medskip

Moser's method is more `down to Earth', computationally speaking,
since it usually proceeds at only {\em one} point, often the origin,
of a manifold, manipulating power series expanded at that point.
Hence it needs geometric objects of class $\mathcal{C}^\omega$, while
adaptations to the $\mathcal{C}^\infty$ or $\mathcal{C}^{\KK \gg
1}$ classes can concern only formal Taylor expansions at the point.

Coming from problems and techniques in Dynamical Systems and Celestial
Mechanics, Moser's method consists in constructing certain {\sl normal
forms} for the objects studied, in order to simplify them and hence
to enable one to rapidly determine whether two given objects are {\em
the same}, up to equivalence.

For instance, for our rigid toy hypersurfaces $\{ u = F(z,
\overline{z}) \}$ in $\C^2$, assuming that they are Levi nondegenerate
at the origin:
\[
u
\,=\,
z\overline{z}
+
{\rm O}_{z,\overline{z}}(3)
\,\,=\,\,
z\overline{z}
+
\sum_{j+k\geqslant 3}\,
F_{j,k}\,
z^j\overline{z}^k,
\]
Moser's game consists in applying several local
rigid biholomorphisms
in order to obtain a simpler graphing
function $F(z, \overline{z})$, {\em e.g.}
with as many as possible coefficients $F_{j,k} = 0$ disappearing,
so that the equation becomes closest as possible to the
model Heisenberg sphere $\{ u = z\overline{z} \}$.

It is not difficult to realize that the isotropy subgroup of 
the origin,
namely the group of {\em rigid} biholomorphisms fixing $(0,0) \in
\C^2$, is $2$-dimensional, and consists of 
weighted scalings coupled with
`horizontal rotations':
\leqnomode\usetagform{default}
\begin{align}
\label{2D-isotropy-introduction}
z'
\,=\,
\rho^{1/2}\,
e^{i\varphi}\,
z,
\ \ \ \ \ \ \ \ \ \ \ \ \ \ \ \ \ \ \ \
w'
\,=\,
\rho\,w,
\end{align}
with $\rho \in \R^\ast$ and $\varphi \in \R$.
Then Section~{\ref{rigid-toy-C-2}}
will elementarily show that one 
can annihilate all $F_{j,0} = 0 = F_{0,k}$ 
and all $F_{j,1} = 0 = F_{1,k}$ as well,
except of course $F_{1,1} = 1$, 
bringing any two rigid hypersurfaces in $M \subset
\C^2$ and $M' \subset {\C'}^2$ to the {\em normalized} forms:
\[
u
\,=\,
z\overline{z}
+
\smallsum{j,k\geqslant 2}\,
F_{j,k}\,
z^j\overline{z}^k
\ \ \ \ \ \ \ \ \ \ \ \ \ \ \ \ \ \ \ \
\text{and}
\ \ \ \ \ \ \ \ \ \ \ \ \ \ \ \ \ \ \ \
u'
\,=\,
z'\overline{z}'
+
\smallsum{j,k\geqslant 2}\,
F_{j,k}'\,
{z'}^j{\overline{z}'}^k,
\]
and then an analysis of what freedom remains in the group of rigid
biholomorphisms will (easily) show that only {\em two}
real parameters
remain free
to send $M$ in normal form to $M'$ also in normal form, namely
$(\rho, \varphi)$ above.  Moreover, it will follows that $M$ and $M'$
are rigidly biholomorphically equivalent if and only if they exchange
through such a trivial scaling-rotation transformation, hence if and
only if there exist $\rho \in \R_+^\ast$ and $\varphi \in \R$ such
that:
\[
F_{j,k}
\,=\,
\rho^{\frac{j+k-2}{2}}\,\,
e^{i\,\varphi\,(j-k)}\,\,
F_{j,k}'
\eqno
{\scriptstyle{(j\,\geqslant\,2,\,\,k\,\geqslant\,2)}}.
\]
Thus, once two normal forms are constructed, whether
$M \sim M'$ or not can be straightforwardly seen.

What is true of the toy will be true of higher dimensional CR objects.
In particular, crude normal forms cannot be made unique, they are
defined only up to the action of a
certain finite-dimensional Lie group, namely the
isotropy sugroup of the (always transitive) model.

Beyond, in most circumstances, {\em e.g.} when $F_{2, 2} \neq 0$
above, one can push further Moser's method, and obtain normal forms
for which {\em all} remaining coefficients $F_{j,k}$ are uniquely
defined, so that $F_{j,k} = F_{j,k}'$
exactly, with no isotropy ambiguity.
This is analog to what one can do in Cartan's method
when some curvature torsion coefficients are nonvanishing:
one can indeed normalize some group parameters present 
in some $T_{j,k}^i$ 
further and further, 
and thereby
decrease the dimension of the bundle $P \longrightarrow M$,
reducing it to smaller subbundles
$P \supsetneqq P_1 \subsetneqq P_2 \supsetneqq \cdots$.

In comparison to Cartan's method, we repeat that one drawback of
Moser's method is that it seems to capture invariants only at one
point.  Fortunately, Moser's method can be applied simultaneously to
all nearby points, 
especially to determine all homogeneous models of a given
class of geometries,
and in a CR context, 
this was done {\em e.g.} in 
Loboda's works~{\cite{Loboda-2001, 
Loboda-2003, Loboda-2013}}.

Recently, Chen-Merker~{\cite{Chen-Merker-2019}} found an
alternative (probably known) method
to capture differential invariants at {\em all points} while working
{\em only at one point}. This method avoids then to 
move the origin everywhere nearby by translations, and it works
most of the times, namely when
the group of transformations is 
only assumed transitive, either finite or infinite
dimensional, {\em see} especially~{\cite[Sec.~12]{Chen-Merker-2019}}.
Hence this method clearly applies to the group 
of rigid biholomorphisms.
Chen-Merker studied mainly parabolic (real) surfaces $S^2 \subset
\R^3$ under the group of special affine transformations of $\R^3$, and
developed an analog of Moser's method in this context.

Links between Affine Geometry and CR geometry have been studied in
depth by Alexander Isaev in his monograph~{\cite{Isaev-2011}}.  Here,
to a given a parabolic surface $\{ u = F(x,y) \}$, namely a surface
whose graphing function $F$ satisfies everywhere:
\[
F_{xx}
\,\neq\,
0
\,\equiv\,
\left\vert\!
\begin{array}{cc}
F_{xx} & F_{xy}
\\
F_{yx} & F_{yy}
\end{array}
\!\right\vert,
\]
one can associate the {\sl tube} hypersurface $M^5 \subset \C^3$
defined as $M^5 := S^2 \times (i\,\R)^3$.  The
paper~{\cite{Merker-2019}} shows that Pocchiola's invariant $\Waux$
associated to $M^5$ produces a seemingly new affine invariant
$\Waux_{\sf aff}$ for parabolic $S^2 \subset \R^3$.  During Alexander
Isaev's stay in Orsay, and after fruitful exchanges with
Peter Olver, it became clear that an independent
study of affine differential invariants of
parabolic surfaces $S^2 \subset \R^3$ should be 
endeavoured, and this was pushed to an end
in~{\cite{Chen-Merker-2019}}.

There, by keeping memory of all terms in the power series that lie
above those coefficients that are progressively normalized,
Chen-Merker obtained certain (complicated) differential-algebraic
expressions made from Taylor coefficients at the origin, from which
one can straightforwardly recover differential invariants {\em at
every point}.  But traditionally instead, people only look at lowest
order currently normalized coefficients in each step, so that
computations remain simple.

Since the technique of~{\cite{Chen-Merker-2019}} seems not to have
been well developed or understood by CR geometers up to now, we
decided to write up the present memoir.  Its main goal is to construct
a {\em bridge:}
\[
\text{\rm Cartan's method}
\ \ \ \ \ \ \ \ \ \ \ \ \ \ \ \ \ \ \ \
\includegraphics[angle=1,scale=0.15]{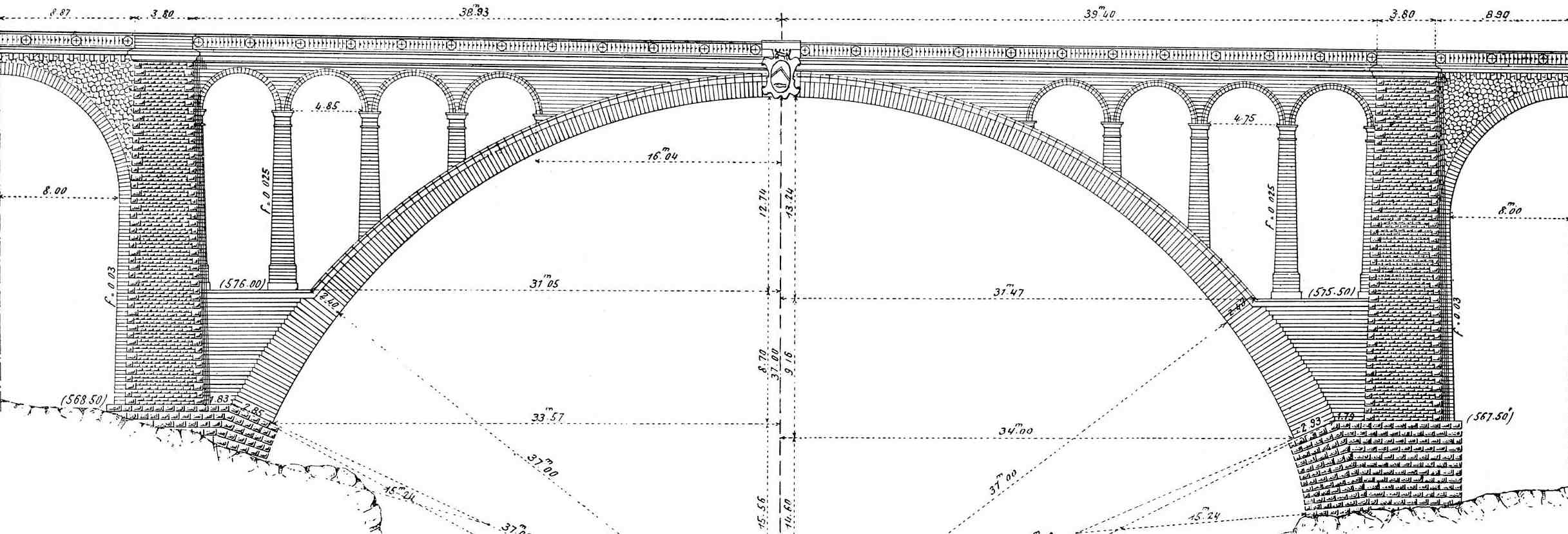}
\ \ \ \ \ \ \ \ \ \ \ \ \ \ \ \ \ \ \ \
\text{\rm Moser's method},
\]
and exhibit how differential invariants pass from one side 
of the river to the
other side, {\em computationally}.  Reading the toy
Section~{\ref{rigid-toy-C-2}} below is enough to understand the key
arch-ideas of such a bridge.  We indeed first focus on the toy case of
rigid equivalences of rigid hypersurfaces in $\C^2$ (easily
reached results),
before passing to the not so simple case of rigid equivalences in the
rigid class denoted $\mathfrak{C}_{2,1}$ by Alexander Isaev which
consists, as written above, of $2$-nondegenerate constant Levi rank
$1$ hypersurfaces $M^5 \subset \C^3$ with $0 \in M$.

In $\C^2$, on the Cartan side of the bridge, we construct in
Section~{\ref{rigid-toy-C-2}} an absolute parallelism
on $P^5 := M^3 \times \C$ equipped with coordinates
$(z, \overline{z}, v, {\sf c}, \overline{\sf c})$
consisting of $5$ differential $1$-forms:
\[
\big\{
\rho,\ 
\zeta,\
\overline{\zeta},\
\pi,\ 
\overline{\pi}
\big\}
\eqno
{\scriptstyle{(\overline{\rho}\,=\,\rho)}},
\]
which satisfy invariant structure equations of the shape:
\[
\aligned
d\rho 
&
\,=\, 
(\pi+\overline{\pi})
\wedge\rho 
+ 
i\,
\zeta\wedge\overline{\zeta},
\\
d\zeta 
&
\,=\,
\pi\wedge\zeta,
&
\ \ \ \ \ \ \ \ \ \ \ \ \ \ \ \ \ \ \ \
\ \ \ \ \ \ \ \ \ \ \ \ \ \ \ \ \ \ \ \
d\overline{\zeta}
&
\,=\,
\overline{\pi}\wedge\overline{\zeta},
\\
d\pi
&
\,=\,
{\textstyle{\frac{1}{{\sf c}\overline{\sf c}}}}\,
\Raux\,
\zeta
\wedge
\overline{\zeta},
&
\ \ \ \ \ \ \ \ \ \ \ \ \ \ \ \ \ \ \ \
\ \ \ \ \ \ \ \ \ \ \ \ \ \ \ \ \ \ \ \
d\overline{\pi}
&
\,=\,
-\,
{\textstyle{\frac{1}{{\sf c}\overline{\sf c}}}}\,
\overline{\Raux}\,
\zeta
\wedge
\overline{\zeta},
\endaligned
\]
where there is only invariant function:
\[
\Raux 
:=
\frac{
F_{zz\overline{z}\overline{z}}\,F_{z\overline{z}} 
-
F_{z z\overline{z}}\, 
F_{z\overline{z}\overline{z}} 
}{
(F_{z\overline{z}})^2}.
\]
We show that $M$ is rigidly equivalent to
$\{ u = z\overline{z} \}$ if and only if $\Raux(F) \equiv 0$.

On the Moser side of the bridge, 
starting from a given $u = \sum_{j+k \geqslant 1}\,
F_{j,k}\, z^j \overline{z}^k$ passing by the origin,
we perform as said above
a few normalizing biholomorphisms in order to reach:
\reqnomode\usetagform{EngelLie}
\begin{align}
0
&
\,=\,
F_{j,0}
\,=\,
F_{0,k}
\tag{(j\,\geqslant\,1,\,\,k\,\geqslant\,1),}
\\
1
&
\,=\,
F_{1,1},
\notag
\\
0
&
\,=\,
F_{j,1}
\,=\,
F_{1,k}
\tag{(j\,\geqslant\,2,\,\,k\,\geqslant\,2),}
\end{align}
and the key feature of the method is to
{\em keep track} of all performed rigid biholomorphic transformations,
which will give us at the end:
\[
u
\,=\,
z\overline{z}
+
\Big[
\frac{F_{2,2}\,F_{1,1}-F_{2,1}\,F_{1,2}}{F_{1,1}^3}
\Big]\,
z^2\overline{z}^2
+
z^2\overline{z}^3
\big(\cdots\big)
+
z^3\overline{z}^2
\big(\cdots\big),
\]
and from this rational expression of 
the final $F_{2,2}'$ coefficient at the origin, 
it is easy to recognize\big/reconstitute\big/translate
Cartan's invariant $\Raux(F)$ {\em at every point}
(up to a nowhere vanishing factor $
{\sf const} \cdot F_{z\overline{z}}$).
Why this is so has already been explained
in~{\cite[Sec.~12]{Chen-Merker-2019}}
and will not be repeated here.

\begin{Principle}
{\sl 
In all CR equivalence problems (and outside CR geometry too),
there exists a way of computing with power series at only 
{\rm one} point which generates
all Cartan-like invariants together with
their syzygies.}
\end{Principle}

In fact, relations (syzygies) require the theory
of {\sl recurrence relations}, developed
for infinite-dimensional Lie groups by 
Olver-Pohjanpelto~{\cite{Olver-Pohjanpelto-2008,
Olver-Pohjanpelto-2009}},
but we will not touch this aspect here.

Because such a `{\sl bridge-principle}' has neither been
constructed nor really noticed in CR geometry,
a joint forthcoming publication will tackle to build it
also for {\em nonrigid} $M^5 \subset \C^3$
that are $2$-nondegenerate and
have constant rank $1$ Levi form,
thereby recovering the full explicit expressions of 
Pocchiola's invariants
$\Waux$ and $\Jaux$ at every point, not only
as number-coefficients 
at one given point as in~{\cite[Thm.~2]{Kolar-Kossovskiy-2019}}.

\smallskip

The first question is: {\sl what is the appropriate local graphed
model for $2$-nondegenerate constant Levi rank $1$ hypersurfaces
$M^5 \subset \C^3$?}  Of course,
it is known from the recent Cartan-theoretic
achievements in~{\cite{Isaev-Zaitsev-2013, Medori-Spiro-2014,
Merker-Pocchiola-2018}} that 
the local model is any neighborhood of any smooth
point of the tube in $\C^3$ over the light cone in $\R^3$ having
equation $x_2^2 - x_3^2 = x_1^2$.
But it is not graphed! We claim that in
different notations, this cone has local graphed equation:
\[
u
\,=\,
\frac{x^2}{1-y},
\]
with $x$, $y$, $u$ being the real parts of three complex coordinates
on $\C^3 \ni (z, \zeta, w)$. 
As we agreed orally
with Alexander Isaev, this is the best, most
compact existing graphed equation.  It happens to also be the central
model 
of parabolic surface $S^2 \subset \R^3$ occurring
in~{\cite{Chen-Merker-2019}}.

The claim is easy. By CR-homogeneity, one can recenter at
any smooth point, {\em e.g.}
at $(0, 1, 1)$, write $(1+x_2)^2 - (1+x_3)^2
= x_1^2$, factor, divide, get $x_2 - x_3 = \frac{x_1^2}{2+x_2+x_3}$,
and linearly change coordinates.

However, this tube graphed equation contains many 
pluriharmonic terms:
\[
\frac{w+\overline{w}}{2}
\,=\,
\frac{(z+\overline{z})^2}{4-2\,\zeta-2\,\overline{\zeta}}
\,=\,
\frac{1}{8}\,
z^2\zeta
+
\frac{1}{8}\,
\overline{z}^2\overline{\zeta}
+
\cdots,
\]
that Moser's method would compulsorily kill at the very beginning.
Thus, $u = \frac{x^2}{1-y}$ is not the right start.
Similarly, $u = x^2 = \frac{1}{2}\, z^2 + 
\frac{1}{2}\, \overline{z}^2 + \cdots$ in $\C^2$ is not the right start
from Moser's point of view.

The right graphed equation for the model light cone $M_{\sf LC}
\subset \C^3$ in $\mathfrak{C}_{2,1}$ was discovered by
Gaussier-Merker in~{\cite{Gaussier-Merker-2003}}: 
\[
M_{\sf LC}
\colon
\ \ \ \ \ \ \ \ \ \ \ \ \ \ \ \ \ \ \ \
u
\,=\,
\frac{z\overline{z}+\frac{1}{2}\,z^2\overline{\zeta}
+\frac{1}{2}\,\overline{z}^2\zeta}{
1-\zeta\overline{\zeta}}
\,\,=:\,
\maux\big(z,\zeta,\overline{z},\overline{\zeta}\big),
\]
and before commenting about very
funny zig-zag errors made in the field at
that time, we review the naive reasoning.  Here, the letter $\maux$ is
from $\maux$odel. By luck, $M_{\sf LC}$ is rigid!

Start with $M^5 \subset \C^3$, with $0 \in M$, rigid,
graphed as:
\[
u
\,=\,
F(z,\zeta,\overline{z},\overline{\zeta}).
\]
Constant Levi rank $1$ means, possibly after a linear
transformation in $\C_{z,\zeta}^2$, that:
\leqnomode\usetagform{default}
\begin{align}
\label{1-lrk-intro}
F_{z\overline{z}}
\,\neq\,
0
\,\equiv\,
\left\vert\!
\begin{array}{cc}
F_{z\overline{z}} & F_{\zeta\overline{z}}
\\
F_{z\overline{\zeta}} & F_{\zeta\overline{\zeta}}
\end{array}
\!\right\vert
\,\,=:\,
\Levi(F),
\end{align}
while $2$-nondegeneracy means that:
\leqnomode\usetagform{default}
\begin{align}
\label{2-ndg-intro}
0
\,\neq\,
\left\vert\!
\begin{array}{cc}
F_{z\overline{z}} & F_{z\overline{z}}
\\
F_{zz\overline{\zeta}} & F_{zz\overline{\zeta}}
\end{array}
\!\right\vert.
\end{align}
By direct symbolic computations, 
Propositions~{\ref{Prp-Levi-determinant}}
and~{\ref{Prp-2ndg-determinant}} 
will establish {\em invariancy} of these vanishing/nonvanishing
properties under rigid changes of holomorphic coordinates.

At the origin, $M_{\sf LC}$ of equation:
\[
u
\,=\,
z\overline{z}
+
{\textstyle{\frac{1}{2}}}\,
z^2\overline{\zeta}
+
{\textstyle{\frac{1}{2}}}\,
\overline{z}^2\zeta
+
{\rm O}_{z,\zeta,\overline{z},\overline{\zeta}}(4),
\]
is obviously $2$-nondegenerate, thanks to the cubic monomial
$\frac{1}{2}\, z^2\overline{\zeta}$ which gives
that~({\ref{2-ndg-intro}}) at $(z,\zeta) = (0, 0)$ becomes
$\left\vert \begin{smallmatrix} 1 & 0 \\ \ast & 1
\end{smallmatrix} \right\vert = 1$.
As for constant Levi rank $1$, order two terms $u = z\overline{z} +
\cdots$ show that this condition is true at the origin, and simple
computations show that~({\ref{1-lrk-intro}}) is identically zero:
\[
\left\vert\!
\begin{array}{cc}
\frac{1}{1-\zeta\overline{\zeta}} 
&
\frac{\overline{z}+z\overline{\zeta}}{(1-\zeta\overline{\zeta})^2}
\\
\frac{z+\overline{z}\zeta}{(1-\zeta\overline{\zeta})^2}
&
\frac{(\overline{z}+z\overline{\zeta})(z+\overline{z}\zeta)}{
(1-\zeta\overline{\zeta})^3}
\end{array}
\!\right\vert
\,\,\equiv\,\,
0
\eqno
{\scriptstyle{(\text{--\,\,indeed!})}}.
\]

So how to easily produce one simple example? 
How $M_{\sf LC}$ was born?

Normalizing the Levi form at the origin,
one can assume $F = z\overline{z} + \cdots$. 
Hence the $2$-nondegeneracy determinant~({\ref{2-ndg-intro}})
becomes at the origin
$\big\vert \begin{smallmatrix} 1 & 0 \\ \ast & 
F_{zz\overline{\zeta}}(0)
\end{smallmatrix} \big\vert = 1$.
Thus, a monomial like $\frac{1}{2}\, z^2\overline{\zeta}$
must be present. Since $F$ is real, its conjugate
$\frac{1}{2}\, \overline{z}^2 \zeta$
also comes:
\[
u
\,=\,
F
\,=\,
z\overline{z}
+
{\textstyle{\frac{1}{2}}}\,
z^2\overline{\zeta}
+
{\textstyle{\frac{1}{2}}}\,
\overline{z}^2
\zeta
+
\sum_{k\geqslant 4}\,
F^k\big(z,\zeta,\overline{z},\overline{\zeta}\big);
\]
here of course, the $F^k$ are homogeneous polynomials of degree $k$.
Without remainders, {\em i.e.} with all $F^k = 0$,
the cubic equation is {\em not} of constant Levi rank $1$
(exercise). 

The idea of Gaussier-Merker
was to take the simplest possible successive $F^4$,
$F^5$, $F^6, \dots$ in order to guarantee $\Levi(F) \equiv 0$.
Thus, plug all this in:
\[
0
\overset{\text{\bf ?}}{\,\,\equiv\,\,}
\left\vert\!
\begin{array}{cc}
1+F_{z\overline{z}}^4+F_{z\overline{z}}^5+F_{z\overline{z}}^5+\cdots
&
\overline{z}+F_{\zeta\overline{z}}^4+F_{\zeta\overline{z}}^5
+F_{\zeta\overline{z}}^6+\cdots
\\
z+F_{z\overline{\zeta}}^4+F_{z\overline{\zeta}}^5
+F_{z\overline{\zeta}}^6+\cdots
&
\ \ \ \ \ \ \
F_{\zeta\overline{\zeta}}^4+F_{\zeta\overline{\zeta}}^5
+F_{\zeta\overline{\zeta}}^6+\cdots
\end{array}
\!\right\vert.
\]
At first, look at terms of order $2$, get $0 = F_{\zeta
\overline{\zeta}}^4 - z\overline{z}$, 
integrate as the simplest possible $F^4 := z \overline{z} 
\zeta \overline{\zeta}$. Next, plug this $F^4$ in, chase only
homogeneous terms of degree $3$, get $F_{\zeta \overline{\zeta}}^5 =
z^2 \overline{\zeta} + \overline{z}^2 \zeta$, and integrate most
simply as $F^5 := \frac{1}{2}\, z^2 \overline{\zeta} \big( \zeta
\overline{\zeta} \big) + \frac{1}{2}\, \overline{z}^2\zeta \big( \zeta
\overline{\zeta} \big)$.  Next, plug this $F^5$ in, get $F_{\zeta
\overline{\zeta}}^6 = 4\, z \overline{z} \zeta \overline{\zeta}$,
integrate $F^6 := z \overline{z} \big( \zeta \overline{\zeta}
\big)^2$, and so on.

An easy induction then shows that powers $\big( \zeta \overline{\zeta}
\big)^k$ appear, and a geometric summation reconstitutes the
denominator $\frac{1}{1 - \zeta \overline{\zeta}}$ in the
Gaussier-Merker model.  \hfill$\bigtriangleup$

\smallskip

Gaussier-Merker made an error when computing (by hand) the Lie algebra
of infinitesimal CR automorphisms of $M_{\sf LC}$, and found a
$7$-dimensional Lie algebra.  This looked `coherent' with a paper
published by Ebenfelt in the Duke Mathematical Journal (year 2000),
which pretended to bound by $7$ the dimension of the CR automorphism
group of any $\mathfrak{C}_{2,1}$ hypersurface $M^5 \subset
\C^3$\,\,---\,\,but due to an incorrect expression of the initial
$G$-structure, Ebenfelt's paper appeared later to be wrong.  Experts
of Cartan theory know how sensitive can be any little error in
normalizations\big/reductions of $G$-structures.

Then the masters Fels-Kaup of Lie transformation groups cleaned up the
subject, showing in~{\cite{Fels-Kaup-2007}}, {\em inter alia}, that
the Gaussier-Merker model is locally biholomorphically equivalent to
the tube over the light cone, so that everybody was wrong before.
They proceeded as follows.

Let $S_{2\times 2} \equiv \R^3 \subset \R^{2\times 2}$ be the space of
all real symmetric $2 \times 2$ matrices.  The open set $\Omega 
{\sf C}^+ \subset S_{2\times 2}$ 
consisting of positive definite matrices has boundary the
future light cone, which may be represented as:
\[
{\sf LC}^+
\,=\,
\left\{
\left(\!
\begin{array}{cc}
t+x_1 & x_2
\\
x_2 & t-x_1
\end{array}
\!\right)
\in
S_{2\times 2}
\colon\,\,
t^2=x_1^2+x_2^2,\,\,
t>0
\right\}.
\]
The objects of study are the following tube 
domain\,\,---\,\,Siegel's upper
half plane up to the factor $i$\,\,---\,\,and 
its boundary hypersurface:
\[
{\sf H}
\,:=\,
\Omega{\sf C}^+
\times
iS_{2\times 2}
\ \ \ \ \ \ \ \ \ \ \ \ \ \ \ \ \ \ \ \
\text{and}
\ \ \ \ \ \ \ \ \ \ \ \ \ \ \ \ \ \ \ \
{\sf T}
\,:=\,
{\sf LC}^+
\times
iS_{2\times 2}.
\]
The global CR automorphism group of ${\sf T}$ consists of just affine
transformations, while the global biholomorphic transformation group
$\Aut ({\sf H})$ of the domain ${\sf H}$ is
known for a long time to consist of the
$10$-dimensional group of all biholomorphic transformations $z
\longmapsto (az + ib) (icz+d)^{-1}$, where $z =
\left( \begin{smallmatrix} w & z_1 \\ z_1 & z_2
\end{smallmatrix} \right)$ with $(z_1, z_2, w) \in \C^3$,
and where $\left( \begin{smallmatrix} a & b \\ c & d \end{smallmatrix}
\right)$ belongs to the real symplectic subgroup ${\sf SP}_2(\R)
\subset {\sf SL}_4(\R)$.

Differentiating this action yields that the algebra of infinitesimal
automorphisms $\mathfrak{aut} ({\sf H})$ 
of the {\em domain}
is equal to $\mathfrak{sp}(2,
\R) \cong \mathfrak{so}_{2,3}(\R)$, also $10$-dimensional.

Fels-Kaup then asked how such automorphisms could be inherited by 
(transmitted to) the
boundary ${\sf T} = \partial {\sf H}$. 

They chose a Cartan subalgebra
of $\mathfrak{so}_{2,3}(\R)$ represented by $\R \zeta_1 \oplus \R
\zeta_2$, where:
\[
\zeta_1
\,:=\,
2w\,\partial_w
\ \ \ \ \ \ \ \ \ \ \ \ \ \ \ \ \ \ \ \
\text{and}
\ \ \ \ \ \ \ \ \ \ \ \ \ \ \ \ \ \ \ \
\zeta_2
\,:=\,
z_1\,\partial_{z_1}
+
2z_2\,\partial_{z_2},
\]
and they showed that any hypersurface $M^5 \subset \C^3$ whose
graphing function starts as $w + \overline{w} = 2z_1 \overline{z}_1 +
z_1^2 \overline{z}_2 + \overline{z}_1^2 z_2 + {\rm O}(4)$ such that
$\mathfrak{hol}(M, 0)$ {\em includes} $\zeta_1$ and $i\zeta_2$ is
locally homogeneous if and only if $\mathfrak{hol}(M, 0)$ {\em also
contains} the two further infinitesimal transformations:
\[
(1-z_2)\,
\partial_{z_1}
+
2z_1\,\partial_w
\ \ \ \ \ \ \ \ \ \ \ \ \ \ \ \ \ \ \ \
\text{and}
\ \ \ \ \ \ \ \ \ \ \ \ \ \ \ \ \ \ \ \
-\,z_1z_2\,\partial_{z_1}
+
(1-z_2^2)\,\partial_w.
\]
Analyzing further structure-theoretic features of the
simple Lie algebra
$\mathfrak{so}_{2,3}(\R)$, they showed that this holds if and only if
the graphed equation reads as the Gaussier-Merker model (up to a
factor $2$):
\leqnomode\usetagform{default}
\begin{align}
\label{Fels-Kaup-model}
w+\overline{w}
\,=\,
\frac{2\,z_1\overline{z}_1+z_1^2\overline{z}_2+\overline{z}_1^2z_2}{
1-z_2\overline{z}_2},
\end{align}
thus giving another natural way to produce this model. 
The main thing was that $\mathfrak{aut}_{\rm CR}$ 
is $10$-dimensional, not $7$!

Fels-Kaup also deduced an explicit rational biholomorphism from this
model~({\ref{Fels-Kaup-model}}) onto a subdomain of ${\sf T}$:
\[
(z_1,z_2,w)
\,\,\,\longmapsto\,\,\,
\frac{1}{1+z_2}\,
\left(\!
\begin{array}{cc}
w+wz_2+z_1^2 & \sqrt{2}\,z_1
\\
\sqrt{2}\,z_1 & 1 - z_2
\end{array}
\!\right).
\]

At about the same time, Fels-Kaup in Acta Mathematica made the
breakthrough of classifying all homogeneous models $M \in
\mathfrak{C}_{2,1}$. They showed that, excepting the light cone, all
such $M$ are in fact {\em simply homogeneous}\,\,---\,\,isotropy Lie
subgroup reduced to identity\,\,---\,\,and necessarily {\em tube},
namely biholomorphically equivalent to $S^2 + (i\R)^3$, for some
surface $S^2 \subset \R^3$ which is simply homogeneous with respect to
the affine group ${\sf A}_3(\R)$.  Fels-Kaup's complete classification
is:

\medskip\noindent{\bf (1)}\,
$S = \{ x_1^2 + x_2^2 = x_3^2, \, x_3 >0\}$ the future light cone;

\medskip\noindent{\bf (2a)}\,
$S = \{ r ( \cos t, \sin t, e^{ \omega t})\in \mathbb{ R}^3: \, r
\in \mathbb{ R}^+ \ \text{\rm and}\ t \in \mathbb{ R}\}$ with
$\omega > 0$ arbitrary;

\medskip\noindent{\bf (2b)}\,
$S = \{ r ( 1, t, e^t) \in \mathbb{ R}^3: \, r \in \mathbb{ R}^+
\ \text{\rm and}\ t \in \mathbb{ R}\}$;

\medskip\noindent{\bf (2c)}\,
$S = \{ r ( 1, e^t, e^{ \theta t}) \in \mathbb{ R}^3 : \, r \in
\mathbb{ R}^+ \ \text{\rm and}\ t \in \mathbb{ R}\}$ with $\theta > 2$
arbitrary;

\medskip\noindent{\bf (3)}\,
$S = \{ c ( t) + r c' ( t) \in \mathbb{ R}^3 : \, r \in \mathbb{
R}^+ \ \text{\rm and} \ t \in \mathbb{ R}\}$, where $c ( t) := (t,
t^2, t^3)$ parametrizes the {\em twisted cubic} $\{ (t, t^2, t^3): \,
t \in \mathbb{ R}\}$ in $\mathbb{ R}^3$ and $c'(t) = (1, 2t, 3t^2)$.

The limit case $\omega = 0$ in {\small\bf (2a)} 
regives the future light cone
{\small\bf (1)}, 
while the limit case $\theta = 2$ in 
{\small\bf (2c)} gives $\{ x \in
\mathbb{ R}^3: \, x_1 x_3 = x_2^2 \ \text{\rm and}\ x_1, x_2 >0 \}$
which is locally linearly (but not globally) equivalent to 
{\small\bf (1)}.  These
five (families of) surfaces are known to be pairwise locally
inequivalent under affine transformations
({\cite{Doubrov-Komrakov-Rabinovich-1996,
Eastwood-Ezhov-1999}}).

As spectacular as they were, the Fels-Kaup articles did not treat the
equivalence problem for {\em all} hypersurfaces $M^5 \subset \C^3$ in
the class $\mathfrak{C}_{2,1}$.  Indeed, like in Riemannian geometry,
it is well known that homogeneous CR manifolds are rather rare in the
set of all CR manifolds.  Although Lie-theoretic methods seem to be
undoubtedly the best to determine homogeneous structures, they lose
their power when dealing with {\sl generic}, non-homogeneous,
structures.  Only Cartan's and Moser's methods of equivalence are able
to handle {\em all} geometric objects of a given kind.

Thus, it was only in the years $2010$'s that the three
papers~{\cite{Isaev-Zaitsev-2013, Medori-Spiro-2014,
Merker-Pocchiola-2018}} achieved the construction of
$10$-dimensional $\{e\}$-structure bundles (or Cartan connections)
$P^{10} \longrightarrow M^5$.

Among these, only Pocchiola's Ph.D.~{\cite{Pocchiola-2013}}, published
as~{\cite{Merker-Pocchiola-2018}}, really performed sufficiently
advanced computations to determine what are the primary curvature
invariants, he called $\Waux$ and $\Iaux$.  Let us review
Pocchiola's results. 
We also follow the article~{\cite{Foo-Merker-2019}},
written because Alexander Isaev insisted that {\em all} details be
made public, while Pocchiola intensively used his computer.

Recall that we denote the class of 
(local) hypersurfaces $M^5 \subset \C^3$
passing by the origin $0 \in M$ 
that are $2$-nondegenerate and whose Levi form
has constant rank $1$ as:
\[
\mathfrak{C}_{2,1}.
\]
Consider therefore a not necessarily rigid
hypersurface $M^5 \subset \C^3$ which
belongs to this class $\mathfrak{C}_{2,1}$, and
which is graphed as:
\[
u
\,=\,
F\big(z_1,z_2,\overline{z}_1,\overline{z}_2,v\big).
\]
The two natural generators of $T^{1,0}M$ and
$T^{0,1}M$ are:
\[
\mathcal{L}_1
\,:=\,
\frac{\partial}{\partial z_1}
-
i\,
\frac{F_{z_1}}{1+i\,F_v}\,
\frac{\partial}{\partial v}
\ \ \ \ \ \ \ \ \ \ \ \ \ \ \ \ \ \ \ \
\text{and}
\ \ \ \ \ \ \ \ \ \ \ \ \ \ \ \ \ \ \ \
\mathcal{L}_2
\,:=\,
\frac{\partial}{\partial z_2}
-
i\,
\frac{F_{z_2}}{1+i\,F_v}\,
\frac{\partial}{\partial v},
\]
in the intrinsic coordinates $(z_1, z_2, \overline{z}_1, 
\overline{z}_2, v)$ on $M$.
We will use the abbreviations:
\[
\Aaux^1
\,:=\,
-\,i\,
\frac{\Faux_{z_1}}{
1+i\,\Faux_v}
\ \ \ \ \ \ \ \ \ \ \ \ \ \ \ \ \ \ \ \
\text{and}
\ \ \ \ \ \ \ \ \ \ \ \ \ \ \ \ \ \ \ \
\Aaux^2
\,:=\,
-\,i\,
\frac{\Faux_{z_2}}{
1+i\,\Faux_v}.
\]

Clearly, the real differential $1$-form:
\[
\varrho_0
\,:=\,
dv
-
\Aaux^1\,dz_1
-
\Aaux^2\,dz_2
-
\overline{\Aaux}^1\,
d\overline{z}_1
-
\overline{\Aaux}^2\,
d\overline{z}_2
\]
has kernel:
\[
\big\{
\varrho_0
=
0
\big\}
\,=\,
T^{1,0}M
\oplus
T^{0,1}M.
\]
At various points:
\[
p
\,=\,
\big(z_1, z_2, \overline{z}_1, \overline{z}_2, v\big)
\,\in\,
M,
\] 
and in terms of $\varrho_0$, the hypothesis
that $M$ has everywhere
degenerate Levi form writes as:
\[
0
\,\equiv\,
\,=\,
\left\vert\!
\begin{array}{cc}
\varrho_0\big(i\,[\mathcal{L}_1,\overline{\mathcal{L}}_1]\big)
&
\varrho_0\big(i\,[\mathcal{L}_2,\overline{\mathcal{L}}_1]\big)
\\
\varrho_0\big(i\,[\mathcal{L}_1,\overline{\mathcal{L}}_2]\big)
&
\varrho_0\big(i\,[\mathcal{L}_2,\overline{\mathcal{L}}_2]\big)
\end{array}
\!\right\vert
(p).
\]

The hypothesis that the Levi form has 
constant rank equal to $1$\,\,---\,\,not to 
$0$!\,\,---\,\,expresses as the fact that the real
CR-transversal vector field:
\[
\mathcal{T}
\,:=\,
i\,
\big[\mathcal{L}_1,\overline{\mathcal{L}}_1\big]
\,=\,
i\,
\Big(
\mathcal{L}_1\big(\overline{\Aaux}^1\big)
-
\overline{\mathcal{L}}_1\big(\Aaux^1\big)
\Big)
\frac{\partial}{\partial v}
\,=:\,
\ell\,
\frac{\partial}{\partial v},
\]
has nowhere vanishing real coefficient:
\[
\ell
\,:=\,
i\,
\Big(
\overline{\Aaux}_{z_1}^1
+
\Aaux^1\,\overline{\Aaux}_v^1
-
\Aaux_{\overline{z}_1}^1
-
\overline{\Aaux}^1\,
\Aaux_v^1
\Big)
\,\,\neq\,\,
0.
\]

The Levi kernel bundle $K^{1,0}M \subset T^{1,0}M$
is then generated by:
\[
\mathcal{K}
\,:=\,
\kaux\,\mathcal{L}_1
+
\mathcal{L}_2,
\]
where:
\[
\kaux
\,:=\,
-\,
\frac{
\mathcal{L}_2\big(\overline{\Aaux}^1\big)
-
\overline{\mathcal{L}}_1\big(\Aaux^2\big)}{
\mathcal{L}_1\big(\overline{\Aaux}^1\big)
-
\overline{\mathcal{L}}_1\big(\Aaux^1\big)}
\]
is the fundamental {\sl slant function}.
As is known from~{\cite{Merker-Pocchiola-Sabzevari-2013-5-CR-II,
Pocchiola-2013, Merker-Pocchiola-2018}},
the hypothesis of $2$-nondegeneracy is then equivalent to
the nonvanishing:
\[
0
\,\neq\,
\overline{\mathcal{L}}_1(\kaux).
\]

Also, the conjugate field $\overline{\mathcal{K}}$ generates
the conjugate Levi kernel bundle
$K^{0,1}M \subset T^{0,1}M$.
There also is a second fundamental function:
\[
\Paux
\,:=\,
\frac{\ell_{z_1}+\Aaux^1\,\ell_v-\ell\,\Aaux_v^1}{\ell}.
\]

Pocchiola conducted in~{\cite{Pocchiola-2013}} the Cartan equivalence
method for such $M^5 \in \mathfrak{C}_{2,1}$ under general (local)
biholomorphic transformations. Reduction to an explicit
$\{e\}$-structure was later done in~{\cite{Foo-Merker-2019}}, after
Alexander Isaev insisted through e-mail exchanges to do this as was
done in~{\cite{Isaev-Zaitsev-2013}}, though in a non-explicit
way. However, such a task is not essential from the point of view of
Cartan's theory, as was well understood by Pocchiola, and as we will
explain in a while.

For now, introducing the five $1$-forms:
\[
\aligned
\rho_0
&
=
\frac{dv-\Aaux^1dz_1-\Aaux^2dz_2
-\overline{\Aaux}^1d\overline{z}_1
-\overline{\Aaux}^2d\overline{z}_2}{\ell},
\\
\kappa_0
&
=
dz_1-\kaux\,dz_2,
\\
\zeta_0
&
=
dz_2,
\\
\overline{\kappa}_0
&
=
d\overline{z}_1
-
\overline{\kaux}\,d\overline{z}_2,
\\
\overline{\zeta}_0
&
=
d\overline{z}_2,
\endaligned
\]
after very, very intensive computations, redone manually by Foo-Merker
in~{\cite{Foo-Merker-2019}} all along $\sim 50$ pages, Pocchiola
obtained modifications $\big\{ \rho, \kappa, \zeta, \overline{\kappa},
\overline{\zeta} \big\}$ of these $1$-forms $\big\{ \rho_0, \kappa_0,
\zeta_0, \overline{\kappa}_0, \overline{\zeta}_0 \big\}$, together
with four complicated $1$-forms $\pi^1$, $\pi^2$, $\overline{\pi}^1$,
$\overline{\pi}^2$ which satisfy structure equations of the specific
concise shape:
\leqnomode\usetagform{default}
\begin{align}
\label{concise-final-drho-dkappa-dzeta}
d\rho
&
\,=\,
\big(
\pi^1
+
\overline{\pi}^1
\big)
\wedge\rho
+
i\,\kappa\wedge\overline{\kappa},
\notag
\\
d\kappa
&
\,=\,
\pi^2\wedge\rho
+
\pi^1\wedge\kappa
+
\zeta\wedge\overline{\kappa},
\notag
\\
d\zeta
&
\,=\,
\big(\pi^1-\overline{\pi}^1\big)
\wedge\zeta
+
i\,\pi^2\wedge\kappa
\,+
\\
&
\ \ \ \ \
+
\Raux\,
\rho\wedge\zeta
+
i\,
\frac{1}{\overline{\sf c}^3}\,
\overline{\Jaux}_0\,\rho\wedge\overline{\kappa}
+
\frac{1}{{\sf c}}\,
\Waux_0\,
\kappa\wedge\zeta,
\notag
\end{align}
in which $\Raux$ is a secondary invariant:
\[
\Raux
\,:=\,
\Re\,
\left[
i\,
\frac{{\sf e}}{{\sf c}{\sf c}}\,
\Waux_0
+
\frac{1}{{\sf c}\overline{\sf c}}
\bigg(
-\,\frac{i}{2}\,
\overline{\mathcal{L}}_1\big(\Waux_0\big)
+
\frac{i}{2}\,
\bigg(
-\,\frac{1}{3}\,
\frac{\overline{\mathcal{L}}_1\big(
\overline{\mathcal{L}}_1(\kaux)\big)}{
\overline{\mathcal{L}}_1(\kaux)}
+
\frac{1}{3}\,
\overline{\Paux}
\bigg)\,
\Waux_0
\bigg)
\right],
\]
expressed in terms of
Pocchiola's two primary invariants
whose explicit expressions have been
confirmed in~{\cite{Foo-Merker-2019}}
(and also after~{\cite{Pocchiola-2013}}
by Alexander Isaev in~{\cite{Isaev-2016}}
assuming $M$ is rigid):
\[
\aligned
\Waux_0
&
\,:=\,
-\,\frac{1}{3}\,
\frac{\mathcal{K}\big(\overline{\mathcal{L}}_1\big(
\overline{\mathcal{L}}_1(\kaux)\big)\big)}{
\overline{\mathcal{L}}_1(\kaux)^2}
+
\frac{1}{3}\,
\frac{\mathcal{K}\big(\overline{\mathcal{L}}_1(\kaux)\big)\,\,
\overline{\mathcal{L}}_1\big(\overline{\mathcal{L}}_1(\kaux)\big)}{
\overline{\mathcal{L}}_1(\kaux)^3}
\,+
\\
&
\ \ \ \ \
+
\frac{2}{3}\,
\frac{\mathcal{L}_1\big(\mathcal{L}_1(\overline{\kaux})\big)}{
\mathcal{L}_1(\overline{\kaux})}
+
\frac{2}{3}\,
\frac{\mathcal{L}_1\big(\overline{\mathcal{L}}_1(\kaux)\big)}{
\overline{\mathcal{L}}_1(\kaux)}
+
\frac{i}{3}\,
\frac{\mathcal{T}(\kaux)}{\overline{\mathcal{L}}_1(\kaux)},
\\
\overline{\Jaux}_0
&
\,:=\,
\frac{1}{6}\,
\frac{\overline{\mathcal{L}}_1\big(
\overline{\mathcal{L}}_1\big(
\overline{\mathcal{L}}_1\big(
\overline{\mathcal{L}}_1(\kaux)\big)\big)\big)}{
\overline{\mathcal{L}}_1(\kaux)}
-
\frac{5}{6}\,
\frac{\overline{\mathcal{L}}_1\big(
\overline{\mathcal{L}}_1\big(
\overline{\mathcal{L}}_1(\kaux)\big)\big)\,\,
\overline{\mathcal{L}}_1\big(
\overline{\mathcal{L}}_1(\kaux)\big)
}{
\overline{\mathcal{L}}_1(\kaux)^2}
-
\frac{1}{6}\,
\frac{\overline{\mathcal{L}}_1\big(
\overline{\mathcal{L}}_1\big(
\overline{\mathcal{L}}_1(\kaux)\big)\big)
}{
\overline{\mathcal{L}}_1(\kaux)}\,
\overline{\Paux}
\,+
\\
&
\ \ \ \ \
+
\frac{20}{27}\,
\frac{\overline{\mathcal{L}}_1\big(\overline{\mathcal{L}}_1
(\kaux)\big)^3}{
\overline{\mathcal{L}}_1(\kaux)^3}
+
\frac{5}{18}\,
\frac{\overline{\mathcal{L}}_1\big(
\overline{\mathcal{L}}_1(\kaux)\big)^2}{
\overline{\mathcal{L}}_1(\kaux)^2}\,
\overline{\Paux}
+
\frac{1}{6}\,
\frac{\overline{\mathcal{L}}_1\big(
\overline{\mathcal{L}}_1(\kaux)\big)\,\,
\overline{\mathcal{L}}_1\big(\overline{\Paux}\big)}{
\overline{\mathcal{L}}_1(\kaux)}
-
\frac{1}{9}\,
\frac{\overline{\mathcal{L}}_1\big(
\overline{\mathcal{L}}_1(\kaux)\big)}{
\overline{\mathcal{L}}_1(\kaux)}\,\,
\overline{\Paux}\,\overline{\Paux}
\,-
\\
&
\ \ \ \ \
-
\frac{1}{6}\,
\overline{\mathcal{L}}_1\big(
\overline{\mathcal{L}}_1\big(
\overline{\Paux}\big)\big)
+
\frac{1}{3}\,
\overline{\mathcal{L}}_1\big(\overline{\Paux}\big)\,
\overline{\Paux}
-
\frac{2}{27}\,
\overline{\Paux}\,
\overline{\Paux}\,
\overline{\Paux}.
\endaligned
\]
When $M$ is assumed to be rigid for simplicity, the numerator of
$\Waux_0$ contains $52$ differential monomials.  When $M$ is not
assumed rigid, it contains hundreds of
thousands of differential monomials instead!
Furthermore, the numerator of $\Jaux_0$ is even huger!

Thus, as is known, the complexity increases spectacularly from rigid
to nonrigid CR manifolds.  This justifies, in a way, to devote some
mathematical works to {\em rigid} CR manifolds, as Alexander Isaev
did, and as we do in the present memoir.

The full $\{e\}$-structure obtained by Foo-Merker
in~{\cite{Foo-Merker-2019}} for nonrigid $M^5 \subset
\C^3$ 
shows that a unique prolongation of
$G$-structure is needed, introducing one further parameter ${\sf t}
\in \R$, together with a (very complicated) real modified
Maurer-Cartan form $\Lambda = d{\sf t} + \cdots$ and that all
appearing torsion coefficients are {\em secondary invariants}.  The
constructed
bundle $P^{10} \longrightarrow M^5$ is equipped with ten coordinates:
\[
\big(
z_1,z_2,\overline{z}_1,\overline{z}_2,v,\,
{\sf c},\overline{\sf c},
{\sf e},\overline{\sf e},{\sf t}
\big),
\]
with ${\sf c} \in \C^\ast$, ${\sf e} \in \C$, ${\sf t} \in \R$, 
together with a collection of ten 
complex-valued $1$-form
which make a frame for $TP^{10}$, denoted:
\[
\big\{
\rho,\,
\kappa,\,
\zeta,\,
\overline{\kappa},\,
\overline{\zeta},\,
\pi^1,\,
\overline{\pi}^1,\,
\pi^2,\,
\overline{\pi}^2,\,
\Lambda
\big\}
\eqno
{\scriptstyle{(\overline{\rho}\,=\,\rho,\,\,
\overline{\Lambda}\,=\,\Lambda)}},
\] 
and which satisfy $10$ invariant structure equations;
however, we will not write the structure equations for
$d\pi^1$, $d\overline{\pi}^1$,
$d\pi^2$, $d\overline{\pi}^2$, $d\Lambda$,
because they are not simple,
and anyway, they incorporate only
secondary invariants.

Thus quite unexpectedly, Pocchiola discovered that 
all primary invariants appear {\em before} 
prolongation of the equivalence problem,
that is to say,
they already appear 
at the beginning of the story,
in the structure 
equations~({\ref{concise-final-drho-dkappa-dzeta}}).

This phenomenon is in some sense
`counter-intuitive' to CR geometers, 
since for Levi nondegenerate CR structures
$M^{2n+1} \subset \C^{n+1}$,
and for the corresponding second order {\sc pde} systems,
{\em no}
curvatures appear after absorption before prolongation
(summation convention holds):
\[
\aligned
d\omega
&
\,=\,
\omega^\alpha
\wedge
\omega_\alpha
+
\omega\wedge\varphi,
\\
d\omega^\alpha
&
\,=\,
\omega^\beta\wedge\varphi_\beta^\alpha
+
\omega\wedge\varphi^\alpha,
\\
d\omega_\alpha
&
\,=\,
\varphi_\alpha^\beta
\wedge
\omega_\beta
+
\omega_\alpha\wedge\varphi
+
\omega\wedge\varphi_\alpha,
\endaligned
\]
while primary and secondary invariants appear afterwards, {\em e.g.}
like $S_{\beta\rho}^{\alpha\sigma}$ and $R_{\beta\gamma}^\alpha$,
$T_\beta^{\alpha\gamma}$ in:
\[
\aligned
d\varphi_\beta^\alpha
&
\,=\,
{\textstyle{\frac{1}{2}}}\,
\delta_\beta^\alpha\,
\psi\wedge\omega
-
\varphi_\beta^\gamma
\wedge
\varphi_\gamma^\alpha
-
\varphi_\beta\wedge\omega^\alpha
-
\varphi^\alpha\wedge\omega_\beta
+
\delta_\beta^\alpha\,
\omega^\gamma\wedge\varphi_\gamma
\,+
\\
&
\ \ \ \ \
+
S_{\beta\rho}^{\alpha\sigma}\,
\omega^\rho
\wedge\omega_\sigma
+
R_{\beta\gamma}^\alpha\,
\omega^\gamma\wedge\omega
+
T_\beta^{\alpha\gamma}\,
\omega_\gamma\wedge\omega.
\endaligned
\]

Next, in the `flat case' where both $\Jaux_0 \equiv 0 \equiv \Waux_0$
vanish identically, which implies $\Raux \equiv 0$ too, Pocchiola's
structure equations reduce to constant coefficients:
\begin{align}
\label{structure-constant-coefficients}
d\rho
&
\,=\,
\big(\pi^1+\overline{\pi}^1\big)
\wedge\rho
+
i\,\kappa\wedge\overline{\kappa},
\notag
\\
d\kappa
&
\,=\,
\pi^2
\wedge\rho
+
\pi^1\wedge\kappa
+
\zeta\wedge\overline{\kappa},
\\
d\zeta
&
\,=\,
\big(
\pi^1
-
\overline{\pi}^1
\big)
\wedge\zeta
+
i\,
\pi^2\wedge\kappa.
\notag
\end{align}
Then a key point is to show that after prolongation, precisely the
structure equations of the Gaussier-Merker model
pop up, namely
(conjugate equations are unwritten):
\[
\aligned
d\rho
&
\,=\,
\pi^1\wedge\rho
+
\overline{\pi}^1\wedge\rho
+
i\,\kappa\wedge\overline{\kappa},
\\
d\kappa
&
\,=\,
\pi^1\wedge\kappa
+
\pi^2\wedge\rho
+
\zeta\wedge\overline{\kappa},
\\
d\zeta
&
\,=\,
i\,\pi^2\wedge\kappa
+
\pi^1\wedge\zeta
-
\overline{\pi}^1\wedge\zeta,
\\
d\pi^1
&
\,=\,
i\,\kappa\wedge\overline{\pi}^2
+
\zeta\wedge\overline{\zeta}
+
\Lambda\wedge\rho,
\\
d\pi^2
&
\,=\,
\pi^2\wedge\overline{\pi}^1
+
\zeta\wedge\overline{\pi}^2
+
\Lambda\wedge\kappa,
\\
d\Lambda
&
\,=\,
i\,\pi^2\wedge\overline{\pi}^2
+
\Lambda\wedge\pi^1
+
\Lambda\wedge\overline{\pi}^1,
\endaligned
\]
and not the structure equations of any other
kind of hypersurface $M^5 \subset \C^3$. This was done by
Pocchiola at the very end of~{\cite{Pocchiola-2013}}, not published
in~{\cite{Merker-Pocchiola-2018}} for reasons of space. 

In the meanwhile, Wei Guo Foo found that Pocchiola missed the presence
of a purely imaginary function $h = i\, H$ with $\overline{H} = H$ in
computations starting from~({\ref{structure-constant-coefficients}}),
{\em which could have destroyed 
Pocchiola's main result} ({\bf !}), because
some (phantom) primary invariants could have then existed in the
structure equations for $d\pi^1$, $d\overline{\pi}^1$, $d\pi^2$,
$d\overline{\pi}^2$, $d\Lambda$,
exactly as in Cartan-Chern-Moser's
computations{\bf !}

Fortunately, this function $h = i\, H$ could be shown to vanish, 
hence phantoms remained phantoms, and the
correction to the (unpublished) end of~{\cite{Pocchiola-2013}} will
appear as~{\cite{Merker-Pocchiola-2019}}, 
prepublished at the end of~{\cite{Foo-Merker-2019}}. 
Maybe Pocchiola just did not type a proper presentation, and
was anyway right
in his manuscripts.

Lastly, we recall that Cartan adopted Lie's principle of thought
({\cite[Chap.~1]{Lie-Merker-2015}), as we do too, which admits that
either a given differential invariant, call it $\Paux$, is
identically zero, or is assumed to be nowhere zero, after
restriction to an appropriate open subset:
\[
\xymatrix{
&&
\Paux\,\equiv\,0,
\\
\Paux
\ar[urr]
\ar[drr]
&&
\\
&&
\Paux\,\not\equiv\,0.
}
\]
Mixed cases where some
invariant is nonzero on some nonempty open subset and
vanishes on a nonempty closed subset are excluded from exploration.

{\em Therefore there is essentially no necessity to set
up an $\{e\}$-structure when $\Waux_0 \equiv 0 \equiv
\Jaux_0$}, because when either $\Waux_0 \not\equiv 0$,
hence $\Waux_0 \neq 0$ after restriction,
or $\Jaux_0 \not\equiv 0$, hence $\Jaux_0 \neq 0$
after restriction, {\em Cartan's method commands to
continue the group parameter normalizations}!

Pocchiola indeed listened to captain Cartan, and was able to prove the

\begin{Theorem}
{\rm {\cite{Pocchiola-2013, Merker-Pocchiola-2018,
Foo-Merker-2019, Merker-Pocchiola-2019}}}
Only two primary invariants, $\Waux_0$ and $\Jaux_0$, occur for
biholomorphic equivalences of $\mathfrak{C}_{2,1}$ real analytic
hypersurfaces $M^5 \subset \C^3$, and:
\[
0
\,\equiv\,
\Waux_0
\,\equiv\,
\Jaux_0
\ \ \ \ \ \ \ 
\Longleftrightarrow
\ \ \ \ \ \ \ 
M\,\,
\text{is equivalent to the
Gaussier-Merker model}.
\]
Furthermore, when either $\Waux_0 \neq 0$ or
$\Jaux_0 \neq 0$, the equivalence problem
reduces to a $5$-dimensional $\{e\}$-structure on $M^5$.
\end{Theorem}

As a corollary known from general Cartan theory, every non-flat $M^5
\in \mathfrak{C}_{2,1}$ has CR automorphisms group of dimension
$\leqslant 5$.  This confirmed the same dimensional gap estimate $10
\downarrow 5$ obtained by Fels-Kaup in~{\cite{Fels-Kaup-2008}}, 
who assumed $M$ to be homogeneous from the
beginning.

\smallskip

Now, as said,
we will work with {\em rigid} hypersurfaces, which is
easier. Only in a future publication will we complete
the views of~{\cite{Kolar-Kossovskiy-2019}} by comparing
them with Pocchiola's results in a deeper way, 
inspired by the present article.

We start by presenting the Moser side of the river.
But before we really treat $\mathfrak{C}_{2,1}$
hypersurfaces $M^5 \subset \C^3$, 
let us explain first how we can get rid of 
{\sl infinity} in the local Lie group of rigid biholomorphisms
by performing what we will call as 
in~{\cite{Kolar-Kossovskiy-2019}}
a {\sl prenormalization}, which is here,
as we already saw, to reach:
\leqnomode\usetagform{default}
\begin{align}
\label{prenormalized-M3-C2-intro}
u
\,=\,
z\overline{z}
+
\sum_{j,k\geqslant 2}\,
F_{j,k}\,
z^j\overline{z}^k,
\end{align}
with $\overline{F_{k,j}} = F_{j,k}$.

How can we do this? Simple! First, starting from a general
$u = \sum_{j+k \geqslant 1}\, F_{j,k}\, z^j \overline{z}^k$,
we get rid of all harmonic terms $F_{j,0}\, z^j$,
$F_{0,k}\, \overline{z}^k$
in the graphing function
by setting:
\[
z'
\,:=\,
z,
\ \ \ \ \ \ \ \ \ \ \ \ \ \ \ \ \ \ \ \ \ \ \ \ \ \
w'
\,:=\,
w
-
2\,\sum_{j\geqslant 1}\,
F_{j,0}\,z^j,
\]
and we get a new graphed equation of the form (dropping primes):
\[
u
\,=\,
\sum_{j\geqslant 1
\atop
k\geqslant 1}\,
F_{j,k}\,
z^j\overline{z}^k.
\]
By this, we have erased an {\em infinite} number of coefficients
$F_{j,0}$, $F_{0,k}$, which was possible thanks to 
the infinite dimensionality of the group of
rigid biholomorphisms.
More precisely, we have consumed $1$ function of $1$ complex
variable.

Next, assuming Levi nondegeneracy at the origin,
making an elementary linear transformation
(exercise), we can assume:
\[
\aligned
u
&
\,=\,
z\overline{z}
+
\sum_{j+k\geqslant 3
\atop
j,\,k\geqslant 1}\,
F_{j,k}\,
z^j\overline{z}^k
\\
&
\,=\,
z\overline{z}
+
\overline{z}\,
\Big(
\sum_{j\geqslant 2}\,
F_{j,0}\,z^j
\Big)
+
z\,
\Big(
\sum_{k\geqslant 2}\,
F_{0,k}\,
\overline{z}^k
\Big)
+
\sum_{j\geqslant 2
\atop
k\geqslant 2}\,
F_{j,k}\,
z^j\overline{z}^k.
\endaligned
\]
Here, the presence of the monomial $z\overline{z}$ is
very advantageous in that it enables to {\em capture}
all monomials $\overline{z}\, z^j$ and their
conjugates $z\, \overline{z}^k$ in a tricky but simple 
{\em factorization}, in which we abbreviate
$\Lambda(z) := \sum_{j\geqslant 2}\, F_{j,0}\, z^j$:
\[
u
\,=\,
\Big(
z
+
\Lambda(z)
\Big)\,
\Big(
\overline{z}
+
\overline{\Lambda}(\overline{z})
\Big)
-
\Lambda(z)\,
\overline{\Lambda}(\overline{z})
+
\sum_{j\geqslant 2
\atop
k\geqslant 2}\,
F_{j,k}\,
z^j\overline{z}^k.
\]
The same factorization idea will work soon for 
$M^5 \in \mathfrak{C}_{2,1}$. Then by making the biholomorphism:
\[
z'
\,:=\,
z
+
\Lambda(z)
\,=\,
z
+
{\rm O}_z(2),
\ \ \ \ \ \ \ \ \ \ \ \ \ \ \ \ \ \ \ \ \ \ \ \ \ \
w'
\,:=\,
w,
\]
it is not difficult to see 
(details in Section~{\ref{rigid-toy-C-2}}) that 
we come
to the prenormalized form~({\ref{prenormalized-M3-C2-intro}}).  Observe
that we have consumed a second infinity, again $1$ function of $1$
complex variable.

Why do we call this {\sl prenormal} form?  Firstly, because it is
in a sense easily and almost freely got from the
assumptions. Secondly, because one key aspect of power series normal
forms is the progressive reduction of {\sl stability groups}, not well
emphasized in~{\cite{Jacobowitz-1990, Kolar-Kossovskiy-2019}}.  The
reader is referred to Sections~13 and~16 of
Chen-Merker~{\cite{Chen-Merker-2019}} to see examples 
of curves $C^1 \subset \R^2$ and surfaces
$S^2 \subset \R^3$ modulo the group of special affine
transformations for which successive stability groups are 
explicitly described.

The presence of group structure reduction
{\sl also} in Moser's theory of normal forms
is in surprising {\sl homology}, not to say {\sl harmony},
with Cartan's method of equivalence,
whose main gist {\em is} group structure reduction.

Plato's Philosophy states that 
Mathematical objects are one and the same in their World.
Various theories elaborate different concept to
grasp these Ideas. The more adequate the concepts are,
the more unitary they are.
What we are claiming is again a
good sign of {\sl Unity} in Mathematics.

Indeed, once a prenormalization is obtained,
in order to normalize $F(z, \overline{z})$
further, it is natural to assume
that the next rigid biholomorphic transformations
$(z, w) \longmapsto (z', w')$ 
to be used should keep unchanged the `{\sl shape}'
of the prenormalization, namely send:
\[
u
\,=\,
z\overline{z}
+
\sum_{j,k\geqslant 2}\,
F_{j,k}\,
z^j\overline{z}^k
\ \ \ \ \ \ \ \ \ \ \ \ \ \ \ \ \ \ \ \
\text{to}
\ \ \ \ \ \ \ \ \ \ \ \ \ \ \ \ \ \ \ \
u'
\,=\,
z'\overline{z}'
+
\sum_{j,k\geqslant 2}\,
F_{j,k}'\,
{z'}^j{\overline{z}'}^k.
\]
This of course imposes many contraints on the map $(z, w) \longmapsto
(z', w')$. And in the rigid context, it is easy to see (in
Section~{\ref{rigid-toy-C-2}}), that only a {\em finite-dimensional}
Lie group remains. Thus, after prenormalization is performed,
one is led back to Lie's original theory~{\cite{Lie-Merker-2015,
Olver-1995}}
in jet spaces for finite-dimensional
continuous groups, which can be safely and naturally applied, 
to finish.

\smallskip

Next, what about $\mathfrak{C}_{2,1}$
rigid hypersurface $M^5 \subset \C^3$? Quite the same!

In coordinates $(z, \zeta, \overline{z}, \overline{\zeta}) 
\in \C^3$, we start at the origin with:
\[
u
\,=\,
\sum_{a+b+c+d\geqslant 1}\,
F_{a,b,c,d}\,
z^a\zeta^b\overline{z}^c\overline{\zeta}^d.
\]
Abbreviating $\chi (z, \zeta) := \sum_{a+b\geqslant 1}\, 
F_{a,b,0,0}\, 
z^a \zeta^b$, we similarly get rid of {\em pluri}harmonic
terms thanks to $z' := z$, $\zeta' := \zeta$,
$w' := w - 2\, \chi(z, \zeta)$, receiving, after dropping
primes, a right-hand side graphing function
$F$ which satisfies:
\[
0
\,=\,
F_{a,b,0,0}
\,=\,
F_{0,0,c,d}.
\]

Next, since $M$ is $2$-nondegenerate and has 
Levi form of rank $1$ at the origin, it is not difficult
({\em see} Section~{\ref{prenormalization}})
to bring its cubic approximation to:
\[
u
\,=\,
z\overline{z}
+
{\textstyle{\frac{1}{2}}}\,
z^2\overline{\zeta}
+
{\textstyle{\frac{1}{2}}}\,
\overline{z}^2\zeta
+
\sum_{
\substack{
a+b+c+d\geqslant 4
\\
a+b\geqslant 1
\\
c+d\geqslant 1
}}\,
F_{a,b,c,d}\,
z^a\zeta^b\overline{z}^c\overline{\zeta}^d.
\]
And now, the same idea of {\sl absorption} by factorization
pops up. But compared to
$M^3 \subset \C^2$, there is a difference:
{\em two} nontrivial monomials
$z\overline{z}$ (self-conjugate) and 
$\frac{1}{2}\, \overline{z}^2\zeta$ 
(with its equivalent conjugate) can be used to absorb
infinities. Writing them as
$\overline{z} \big( z \big)$ and
$\overline{z}^2 \big( \frac{1}{2}\, \zeta \big)$,
we may therefore {\em capture} all holomorphic
monomials behind $\overline{z} \big( \cdots \big)$
and behind $\overline{z}^2 \big( \cdots \big)$,
by making the rigid biholomorphism:
\[
\aligned
z
+
\sum_{a+b\geqslant 1}\,
F_{a,b,1,0}\,
z^a\zeta^b
&
\,=:\,
z',
\\
{\textstyle{\frac{1}{2}}}\,\zeta
+
\sum_{a+b\geqslant 2}\,
F_{a,b,2,0}\,
z^a\zeta^b
&
\,=:\,
\zeta',
\endaligned
\]
with unchanged $w' := w$.
The true story is a little more subtle, requires more care,
and will be told with rigorous details in
Section~{\ref{prenormalization}}.

Therefore, after having consumed {\em three} holomorphic
functions of the two complex variables
$(z, \zeta)$, we end up with a graph $u = F(z, \zeta, \overline{z}, 
\overline{\zeta})$ which is {\sl prenormalized}
in the sense that:
\[
\aligned
0
&
\,=\,
F_{a,b,0,0}
\,=\,
F_{0,0,c,d},
\\
0
&
\,=\,
F_{a,b,1,0}
\,=\,
F_{1,0,c,d},
\\
0
&
\,=\,
F_{a,b,2,0}
\,=\,
F_{2,0,c,d},
\endaligned
\]
except of course $F_{1,0,1,0} = 1$ and
$F_{2,0,0,1} = \frac{1}{2} = F_{0,1,2,0}$.
An equivalent way to express prenormalization
is to write that (exercise):
\[
u
\,=\,
F
\,=\,
z\overline{z}
+
{\textstyle{\frac{1}{2}}}\,
\overline{z}^2\zeta
+
{\rm O}_{\overline{z}}(3)
+
{\rm O}_{\overline{\zeta}}(1).
\]
The next task is to normalize $F$ beyond prenormalization.

Because in
$\C^2$ a general rigid hypersurface 
$u = F = z\overline{z} + {\rm O}_{z, \overline{z}}(3)$ 
is naturally represented as a perturbation
of the (flat) model $u = z\overline{z}$, 
we represent a general rigid $M \in \mathfrak{C}_{2,1}$
as a perturbation of the Gaussier-Merker model:
\[
u
\,=\,
F\big(z,\zeta,\overline{z},\overline{\zeta}\big)
\,=\,
\maux(z,\zeta,\overline{z},\overline{\zeta})
+
G\big(z,\zeta,\overline{z},\overline{\zeta}\big),
\]
but\,\,---\,\,warning!\,\,---, the remainder function
$G$ here cannot be arbitrary, it must be so that
$\Levi(\maux + G) \equiv 0$.

Next, inspired by~{\cite{Kolar-Kossovskiy-2019}},
we show in the key
Proposition~{\ref{Prp-G-O-z-3}} that 
in prenormalized coordinates, one necessarily has:
\[
G
\,=\,
{\rm O}_{z,\overline{z}}(3).
\]
Since the Gaussier-Merker function:
\[
\maux(z,\zeta,\overline{z},\overline{\zeta})
\,=\,
\frac{z\overline{z}+\frac{1}{2}\,z^2\overline{\zeta}
+\overline{z}^2\zeta}{1-\zeta\overline{\zeta}}
\]
is homogeneous of degree $2$ in $(z, \overline{z})$,
this conducts us, as in~{\cite{Kolar-Kossovskiy-2019}},
to assign the following weights to the coordinate
variables:
\[
[z]
\,:=\,
1
\,=:\,
[\overline{z}],
\ \ \ \ \ \ \ \ \ \ \ \ \ \ \ \ \ \ \ \
[\zeta]
\,:=\,
0
\,=:\,
\big[\overline{\zeta}\big],
\ \ \ \ \ \ \ \ \ \ \ \ \ \ \ \ \ \ \ \
[w]
\,:=\,
2
\,=:\,
[\overline{w}].
\]

Similarly as for rigid $M^3 \subset \C^2$, we next ask:
{\sl which rigid transformations stabilize prenormalization?},
and we will again realize that only a {\em finite-dimensional}
Lie group remains.

Thus we take $M$ in $\C^3 \ni (z_1, z_2, w)$
graphed as $u = F = \maux + G$
and $M'$ in $\C^3 \ni (z_1', z_2', w')$ 
graphed as $u' = F' = \maux' + G'$,
with $G$ prenormalized:
\leqnomode\usetagform{default}
\begin{align}
\label{simultaneous-normalizations-G}
G
\,=\,
{\rm O}_{\overline{z}}(3)
+
{\rm O}_{\overline{\zeta}}(1)
\,=\,
{\rm O}_{z,\overline{z}}(3),
\end{align}
(none condition implies the other),
and the same about $G'$.
The goal is to {\em normalize further} $G'$. 

Without waiting, we expand $G$ in weighted homogeneous parts:
\[
G
\,=\,
\sum_{\nu\geqslant 3}\,
G_\nu,
\ \ \ \ \ \ \ \ \ \ \ \ \ \ \ \ \ \ \ \ \ \ \ \ \ \
G_\nu
\,=\,
\sum_{a+c=\nu}\,
z^a\overline{z}^c\,
G_{a,c}(\zeta,\overline{\zeta}),
\]
and the same for $G'$, with, unlike in Moser's theory for Levi
nondegenerate hypersurfaces in $\C^{n+1}$, coefficient-functions
$G_{a,c}$ which are {\em analytic}, not polynomial.

The elementary Proposition~{\ref{Prp-initial-f-g-h}} shows that,
composing in advance with some element of the $2$-dimensional isotropy
group~({\ref{2D-isotropy-introduction}}) of the origin for the
Gaussier-Merker model, we can assume that the normalizing map has
weighted expansion of the form:
\leqnomode\usetagform{default}
\begin{align}
\label{intro-f-g-h-nu}
f
\,=\,
z
+f_2+f_3+\cdots,
\ \ \ \ \ \ \ \ \ \ \ \ \ \ \ \ \ \ \ \
g
\,=\,
\zeta
+
g_1+g_2+\cdots,
\ \ \ \ \ \ \ \ \ \ \ \ \ \ \ \ \ \ \ \
h
\,=\,
w
+h_3+h_4+\cdots,
\end{align}
where, for $\nu = 3, 4, 5, \dots$,
the appearing holomorphic functions $f_{\nu-1}$,
$g_{\nu-2}$, $h_\nu$ are weighted homogeneous.
Keeping good memory of this pre-composition,
there will remain at the end a $2$-dimensional
ambiguity in the obtained normal form.

As in Jacobowitz's~{\cite[Ch.~3]{Jacobowitz-1990}}
presentation of Moser's method,
with increasing weights $\nu = 3, 4, 5, \dots$, we shall
perform successive holomorphic rigid transformations
of the shape:
\[
z'
\,:=\,
z+f_{\nu-1},
\ \ \ \ \ \ \ \ \ \ \ \ \ \ \ \ \ \ \ \
\zeta'
\,:=\,
\zeta+g_{\nu-2},
\ \ \ \ \ \ \ \ \ \ \ \ \ \ \ \ \ \ \ \
w'
\,:=\,
w+h_\nu.
\]

Then in the main Proposition~{\ref{Prp-5-2}}, we 
will show that through any such
biholomorphism~({\ref{intro-f-g-h-nu}}) which
transforms:
\[
u
\,=\,
\maux
+
G_3+\cdots+G_{\nu-1}
+
G_\nu
+
{\rm O}(\nu+1)
\ \ \ \ \ \ \ 
\text{into}
\ \ \ \ \ \ \
u'
\,=\,
\maux
+
G_3'+\cdots+G_{\nu-1}'
+
G_\nu'
+
{\rm O}'(\nu+1),
\]
homogeneous terms are kept untouched up to order $\leqslant \nu-1$:
\[
G_\mu'\big(z,\zeta,\overline{z},\overline{\zeta}\big)
\,=\,
G_\mu\big(z,\zeta,\overline{z},\overline{\zeta}\big)
\eqno
{\scriptstyle{(3\,\leqslant\,\mu\,\leqslant\,\nu-1)}},
\]
while:
\[
G_\nu'\big(z,\zeta,\overline{z},\overline{\zeta}\big)
\,=\,
G_\nu\big(z,\zeta,\overline{z},\overline{\zeta}\big)
-
2\,\Re\,
\Big\{
{\textstyle{
\frac{\overline{z}+z\overline{\zeta}}{1-\zeta\overline{\zeta}}}}\,
f_{\nu-1}(z,\zeta)
+
{\textstyle{
\frac{(\overline{z}+z\overline{\zeta})^2}{
2(1-\zeta\overline{\zeta})^2}}}\,
g_{\nu-2}(z,\zeta)
-
{\textstyle{\frac{1}{2}}}\,
h_\nu(z,\zeta)
\Big\}.
\]
Here, the {\sl freedom},
which consists of a triple $\{f_{\nu-1}, g_{\nu-2}, h_\nu \big\}$
of holomorphic functions of the two complex
variables $(z, \zeta)$, can be used to simplify\big/normalize
$G_\nu'$ in comparison with $G_\nu$.

It is important to point out that in this paper, we dispense ourselves
completely of making a {\sl formal} theory of normal form {\em before}
conducting a {\sl geometric} reduction to normal form, we come
directly to (geometric) heart.

Then we study the initial weights
$\nu = 3, 4, 5$, even
restricting our attention firstly to total
degree $a + b + c + d \leqslant 5$.
In Section~{\ref{normal-form}},
we show that only two monomials (up to conjugation) 
remain after prenormalization in:
\[
G_3
\,=\,
2\,\Re\,
\Big\{
z^3\overline{\zeta}\,
G_{3,0,0,1}
+
z^3\overline{\zeta}^2\,
G_{3,0,0,2}
\Big\}
+
{\rm O}_{z,\zeta,\overline{z},\overline{\zeta}}(6).
\]
Using the freedom~({\ref{intro-f-g-h-nu}}) and taking account of
preservation of prenormalization, similarly as
in~{\cite{Chen-Merker-2019}}, we show that we can annihilate
$G_{3,0,0,1}' := 0$.  And then, we show that no other Taylor
coefficient of $G_3$ can be normalized, if one requires preservation
of $G_{3, 0, 0, 1} = 0 = G_{3, 0, 0, 1}'$

In particular, this implies that there is no invariant of
(differential) order $4$, and this confirms the results
of~{\cite{Foo-Merker-Ta-2019}}, to be reviewed and compared in a
while.

Next, we study $\nu = 4$, still with $a + b + c + d \leqslant 5$,
and there are again only two monomials:
\[
G_4
\,=\,
2\,\Re\,
\Big\{
z^4\overline{\zeta}\,
G_{4,0,0,1}
+
z^3\overline{z}\overline{\zeta}\,
G_{3,0,1,1}
\Big\}
+
{\rm O}_{z,\zeta,\overline{z},\overline{\zeta}}(6).
\]
Using the freedom~({\ref{intro-f-g-h-nu}}) and taking account of
preservation of all preceding normalizations, 
we show that we can annihilate
$\Im\, G_{3,0,1,1}' := 0$.  And then, we show that no other Taylor
coefficient of $G_4$ can be normalized.

Lastly, for every remaining $\nu \geqslant 5$, 
we verify that only the identity tranformation
$z' = z$, $\zeta' = \zeta$, $w' = w$,
stabilizes prenormalization {\em and}:
\[
0
\,=\,
G_{3,0,0,1}
\,=\,
G_{3,0,0,1}',
\ \ \ \ \ \ \ \ \ \ \ \ \ \ \ \ \ \ \ \ \ \ \ \ \ \
0
\,=\,
\Im\,
G_{3,0,1,1}
\,=\,
\Im\,
G_{3,0,1,1}'.
\]
namely we show that $0 = f_{\nu-1} = g_{\nu-2} = h_\nu$,
necessarily.

Moser's algorithm therefore terminates, and we may 
at last state our main 

\begin{Theorem}
Every hypersurface $M^5 \in \mathfrak{C}_{2,1}$ 
is equivalent, through a local rigid biholomorphism,
to a rigid $\mathcal{C}^\omega$ hypersurface ${M'}^5 \subset {\C'}^3$
which, dropping primes for target coordinates, is a perturbation of
the Gaussier-Merker model:
\[
u
\,=\,
\frac{z\overline{z}+\frac{1}{2}\,z^2\overline{\zeta}
+\frac{1}{2}\,\overline{z}^2\zeta}{
1-\zeta\overline{\zeta}}
+
\sum_{a,b,c,d\in\N
\atop
a+c\geqslant 3}\,
G_{a,b,c,d}\,
z^a\zeta^b\overline{z}^c\overline{\zeta}^d,
\]
with a simplified remainder $G$ which:

\smallskip\noindent{\bf (1)}\,
is normalized to be an ${\rm O}_{z, \overline{z}}(3)$;

\smallskip\noindent{\bf (2)}\,
satisfies the prenormalization conditions $G = 
{\rm O}_{\overline{z}}(3) + {\rm O}_{\overline{\zeta}}(1) = 
{\rm O}_z(3) + {\rm O}_\zeta(1)$:
\[
\aligned
G_{a,b,0,0}
&
\,=\,
0
\,=\,
G_{0,0,c,d},
\\
G_{a,b,1,0}
&
\,=\,
0
\,=\,
G_{1,0,c,d},
\\
G_{a,b,2,0}
&
\,=\,
0
\,=\,
G_{2,0,c,d};
\endaligned
\]

\smallskip\noindent{\bf (3)}\,
satisfies in addition the sporadic normalization conditions:
\[
\aligned
G_{3,0,0,1}
&
\,=\,
0
\,=\,
G_{0,1,3,0},
\notag
\\
\Im\,G_{3,0,1,1}
&
\,=\,
0
\,=\,
\Im\,G_{1,1,3,0}.
\endaligned
\]

Furthermore, two such rigid $\mathcal{C}^\omega$ hypersurfaces $M^5
\subset \C^3$ and ${M'}^5 \subset {\C'}^3$, both brought into such a
normal form, are rigidly biholomorphically equivalent if and only if
there exist two constants $\rho \in \R_+^\ast$, $\varphi \in \R$, such
that for all $a$, $b$, $c$, $d$:
\[
G_{a,b,c,d}
\,=\,
G_{a,b,c,d}'\,
\rho^{\frac{a+c-2}{2}}\,
e^{i\varphi(a+2b-c-2d)}.
\]
\end{Theorem}

Now, before talking about any bridge, 
we must survey the results of the
article~{\cite{Foo-Merker-Ta-2019}}, from
Cartan's side of the river. These results were finalized after the
stay in Orsay of Alexander Isaev, who raised the problem.  The reader
is referred to the introduction of~{\cite{Foo-Merker-Ta-2019}} for
more extensive information.

Consider as before a rigid $M^5 \subset \C^3$ with $0 \in M$, which is
$2$-nondegenerate and has Levi form of constant rank $1$, {\em i.e.}
belongs to the class $\mathfrak{C}_{2,1}$, and which is graphed as:
\[
u
\,=\,
F\big(z_1,z_2,\overline{z}_1,\overline{z}_2\big).
\]
The letter $\zeta$ is protected, hence not used instead of $z_2$,
since $\zeta$ will denote a $1$-form.  The two natural generators of
$T^{1,0}M$ and $T^{0,1}M$ are:
\[
\mathcal{L}_1
\,:=\,
\partial_{z_1}
-
i\,F_{z_1}\,\partial_v
\ \ \ \ \ \ \ \ \ \ \ \ \ \ \ \ \ \ \ \
\text{and}
\ \ \ \ \ \ \ \ \ \ \ \ \ \ \ \ \ \ \ \
\mathcal{L}_2
\,:=\,
\partial_{z_2}
-
i\,F_{z_2}\,\partial_v,
\]
in the intrinsic coordinates $(z_1, z_2, \overline{z}_1,
\overline{z}_2, v)$ on $M$.  The Levi kernel bundle $K^{1,0}M \subset
T^{1,0}M$ is generated by:
\[
\mathcal{K}
\,:=\,
\kaux\,\mathcal{L}_1
+
\mathcal{L}_2,
\ \ \ \ \ \ \ \ \ \ \ \ \ \ \ \ \ \ \ \
\text{where}
\ \ \ \ \ \ \ \ \ \ \ \ \ \ \ \ \ \ \ \
\kaux
\,:=\,
-\,
\frac{F_{z_2\overline{z}_1}}{
F_{z_1\overline{z}_1}},
\]
is the slant function.
The hypothesis of $2$-nondegeneracy is equivalent to
the nonvanishing:
\[
0
\,\neq\,
\overline{\mathcal{L}}_1(\kaux).
\]
Also, the conjugate $\overline{\mathcal{K}}$ generates
the conjugate Levi kernel bundle
$K^{0,1} \subset T^{0,1}M$.

There is a second fundamental function, and no more:
\[
\Paux
\,:=\,
\frac{F_{z_1z_1\overline{z}_1}}{F_{z_1\overline{z}_1}}.
\]
In the rigid case, it looks so simple{\bf !}
But in the nonrigid case, $\Paux$ has
a numerator involving {\bf 69} differential monomials{\bf !}

Foo-Merker-Ta produced in~{\cite{Foo-Merker-Ta-2019}} 
reduction to an $\{e\}$-structure
for the equivalence problem, under {\em rigid}
(local) biholomorphic transformations, 
of such rigid $M^5 \in \mathfrak{C}_{2,1}$.
They constructed an invariant $7$-dimensional bundle
$P^7 \longrightarrow M^5$ equipped with coordinates:
\[
\big(
z_1,z_2,\overline{z}_1,\overline{z}_2,v,\,
{\sf c},\overline{\sf c}
\big),
\]
with ${\sf c} \in \C$,
together with a collection of seven 
complex-valued $1$-form
which make a frame for $TP^7$, denoted:
\[
\big\{
\rho,\,
\kappa,\,
\zeta,\,
\overline{\kappa},\,
\overline{\zeta},\,
\alpha,\,
\overline{\alpha}
\big\}
\eqno
{\scriptstyle{(\overline{\rho}\,=\,\rho)}},
\] 
which satisfy $7$ invariant structure equations of the form:
\[
\aligned
d\rho 
&
\,=\,
\big(\alpha+\overline{\alpha}\big)
\wedge
\rho 
+ 
i\,
\kappa
\wedge
\overline{\kappa},
\\
d\kappa 
&
\,=\,
\alpha
\wedge
\kappa 
+ 
\zeta
\wedge
\overline{\kappa},
\\
d\zeta 
&
\,=\,
\big(
\alpha-\overline{\alpha}
\big)
\wedge
\zeta
+
\frac{1}{\sf c}\,\Iaux_{0}\,
\kappa\wedge\zeta
+
\frac{1}{\overline{\sf c}\overline{\sf c}}\,
\Vaux_{0}\,
\kappa\wedge\overline{\kappa},
\\
d\alpha 
&
\,=\, 
\zeta\wedge\overline{\zeta}
-
\frac{1}{\sf c}\,
\Iaux_{0}\,
\zeta\wedge\overline{\kappa}
+
\frac{1}{{\sf c}\overline{\sf c}}\,\Qaux_{0}\,
\kappa
\wedge
\overline{\kappa}
+
\frac{1}{\overline{\sf c}}\,
\overline{\Iaux}_{0}\, 
\overline{\zeta}
\wedge
\kappa,
\endaligned
\]
conjugate structure equations for 
$d\overline{\kappa}$, $d\overline{\zeta}$, $d\overline{\alpha}$
being easily deduced.

Here, as in Pocchiola's Ph.D., there are exactly {\em two}
primary Cartan-curvature invariants:
\[
\aligned
\Iaux_0
&
\,:=\,
-\,\frac{1}{3}\,
\frac{\mathcal{K}\big(\overline{\mathcal{L}}_1
\big(\overline{\mathcal{L}}_1(\kaux)\big)\big)}{
\overline{\mathcal{L}}_1(\kaux)^2}
+
\frac{1}{3}\,
\frac{\mathcal{K}\big(\overline{\mathcal{L}}_1(\kaux)\big)\,
\overline{\mathcal{L}}_1\big(\overline{\mathcal{L}}_1(\kaux)\big)}{
\overline{\mathcal{L}}_1(\kaux)^3}
\,+
\\
&
\ \ \ \ \
+
\frac{2}{3}\,
\frac{\mathcal{L}_1\big(\mathcal{L}_1(\overline{\kaux})\big)}{
\mathcal{L}_1(\overline{\kaux})}
+
\frac{2}{3}\,
\frac{\mathcal{L}_1\big(\overline{\mathcal{L}}_1(\kaux)\big)}{
\overline{\mathcal{L}}_1(\kaux)},
\\
\Vaux_0
&
\,:=\,
-\,\frac{1}{3}\,
\frac{\overline{\mathcal{L}}_1\big(\overline{\mathcal{L}}_1\big(
\overline{\mathcal{L}}_1(\kaux)\big)\big)}{
\overline{\mathcal{L}}_1(\kaux)}
+
\frac{5}{9}\,
\bigg(
\frac{\overline{\mathcal{L}}_1\big(\overline{\mathcal{L}}_1
(\kaux)\big)}{
\overline{\mathcal{L}}_1(\kaux)}
\bigg)^2
\,-
\\
&
\ \ \ \ \
-\,
\frac{1}{9}\,
\frac{\overline{\mathcal{L}}_1\big(\overline{\mathcal{L}}_1
(\kaux)\big)\,\,\overline{\Paux}}{
\overline{\mathcal{L}}_1(\kaux)}
+
\frac{1}{3}\,
\overline{\mathcal{L}}_1(\overline{\Paux})
-
\frac{1}{9}\,
\overline{\Paux}\,
\overline{\Paux}.
\endaligned
\]
One can check that Pocchiola's $\Waux_0$
which occurs under {\em general} biholomorphic transformations
of $\C^3$ (not necessarily rigid!),
when written for a {\em rigid} $M^5 \subset \C^3$, 
identifies with:
\[
\Iaux_0
\big(
F(z_1,z_2,\overline{z}_1,\overline{z}_2)
\big)
\,\equiv\,
\Waux_0
\big(
F(z_1,z_2,\overline{z}_1,\overline{z}_2)
\big).
\]

Furthermore, there is {\em one} secondary invariant 
whose unpolished expression
is:
\[
\aligned
\Qaux_0
\,:=\,
\frac{1}{2}\,
\overline{\mathcal{L}}_1
\big(
\Iaux_0
\big)
-
\frac{1}{3}\,
\bigg(
\Paux
-
\frac{\mathcal{L}_1\big(\mathcal{L}_1(\overline{\kaux})\big)}{
\mathcal{L}_1(\overline{\kaux})}
\bigg)\,
\overline{\Iaux}_0
-
\frac{1}{6}\,
\bigg(
\overline{\Paux}
-
\frac{\overline{\mathcal{L}}_1\big(\overline{\mathcal{L}}_1
(\kaux)\big)}{
\overline{\mathcal{L}}_1(\kaux)}
\bigg)\,
\Iaux_0
-
\frac{1}{2}\,
\frac{\mathcal{K}(\Vaux_0)}{
\overline{\mathcal{L}}_1(\kaux)}.
\endaligned
\]

Visibly indeed, the vanishing of $\Iaux_0$ and $\Vaux_0$ implies the
vanishing of $\Qaux_0$. In fact, a consequence of Cartan's general
theory is:
\[
0
\,\equiv\,
\Iaux_0
\,\equiv\,
\Vaux_0
\ \ \ \ \ \ \ 
\Longleftrightarrow
\ \ \ \ \ \ \ 
M\,\,
\text{is rigidly equivalent to the
Gaussier-Merker model}.
\]

In~{\cite{Foo-Merker-Ta-2019}}, by deducing new
relations from the structure equations above,
it was proved that $\Qaux_0$ is real-valued,
but a finalized expression was missing there. 
A clean finalized expression of
$\Qaux_0$, in terms of only the two fundamental functions $\kaux$,
$\Paux$ (and their conjugates), from which 
one immediately sees real-valuedness, is:
\[
\aligned
\Qaux_0
&
\,:=\,
2\,\Re\,
\bigg\{
\frac{1}{9}\,
\frac{\mathcal{K}\big(\overline{\mathcal{L}}_1(\kaux)\big)\,
\overline{\mathcal{L}}_1\big(\overline{\mathcal{L}}_1(\kaux)\big)^2}{
\overline{\mathcal{L}}_1(\kaux)^4}
\,-
\\
&
\ \ \ \ \ \ \ \ \ \ \ \ \ \ \
-\,
\frac{1}{9}\,
\frac{\mathcal{K}\big(\overline{\mathcal{L}}_1
\big(\overline{\mathcal{L}}_1(\kaux)\big)\big)\,
\overline{\mathcal{L}}_1\big(\overline{\mathcal{L}}_1(\kaux)\big)}{
\overline{\mathcal{L}}_1(\kaux)^3}
-
\frac{1}{9}\,
\frac{\mathcal{K}
\big(\overline{\mathcal{L}}_1(\kaux)\big)\,
\overline{\mathcal{L}}_1\big(\overline{\mathcal{L}}_1(\kaux)\big)\,
\overline{\Paux}}{
\overline{\mathcal{L}}_1(\kaux)^3}
\,-
\\
&
\ \ \ \ \ \ \ \ \ \ \ \ \ \ \
-\,
\frac{1}{9}\,
\frac{
\mathcal{L}_1\big(
\overline{\mathcal{L}}_1(\kaux)\big)\,
\overline{\mathcal{L}}_1\big(\overline{\mathcal{L}}_1(\kaux)\big)}{
\overline{\mathcal{L}}_1(\kaux)^2}
+
\frac{1}{9}\,
\frac{\mathcal{K}\big(\overline{\mathcal{L}}_1
\big(\overline{\mathcal{L}}_1(\kaux)\big)\big)\,
\overline{\Paux}}{
\overline{\mathcal{L}}_1(\kaux)^2}
\,-
\\
&
\ \ \ \ \ \ \ \ \ \ \ \ \ \ \
-\,
\frac{2}{9}\,
\frac{\mathcal{L}_1\big(\overline{\mathcal{L}}_1(\kaux)\big)\,
\overline{\Paux}}{
\overline{\mathcal{L}}_1(\kaux)}
-
\frac{1}{9}\,
\frac{\overline{\mathcal{L}}_1\big(
\overline{\mathcal{L}}_1(\kaux)\big)\,
\Paux}{
\overline{\mathcal{L}}_1(\kaux)}
+
\frac{1}{3}\,
\frac{
\mathcal{L}_1\big(
\overline{\mathcal{L}}_1\big(
\overline{\mathcal{L}}_1(\kaux)\big)\big)}{
\overline{\mathcal{L}}_1(\kaux)}
+
\frac{1}{6}\,
\overline{\mathcal{L}}_1(\Paux)
\bigg\}
\\
&
\ \ \ \ \ \ \ \ \ \ \ \ \ \ \
-
\frac{1}{9}\,
\big\vert
\overline{\Paux}
\big\vert^2
+
\frac{1}{3}\,
\bigg\vert
\frac{\overline{\mathcal{L}}_1\big(
\overline{\mathcal{L}}_1(\kaux)\big)}{
\overline{\mathcal{L}}_1(\kaux)}
\bigg\vert^2.
\endaligned
\]
Section~{\ref{finalized-expression-Q0}} 
is devoted to provide the details
of the necessary, nontrivial computations.
Having $\Qaux_0$ in finalized form is required to
compare with what Moser's method 
gives on the other side of the bridge.

\smallskip

Indeed, to finish this introduction, we can at last
say that the key idea of the bridge
is presented in Sections~{\ref{caves-beneath-waterfall}}
and~{\ref{invariants-I-0-V-0-Q-0}}.

\medskip\noindent
{\bf Acknowledgments.}
The realization of this research work in Cauchy-Riemann (CR) geometry
has received generous financial support from the scientific grant
2018/29/B/ST1/02583 originating from the
Polish National Science Center (NCN).

During fall 2019,
in October in Warsaw, 
then in November and December in Paris, 
the authors benefited from countless oral exchanges
with Pawe{\l} Nurowski (Center For Theoretical Physics),
who explained, developed, 
and even {\em taught in the latest details} 
his deep knowledge of
{\sl Cartan's method of equivalence} 
(production of homogeneous models),
whose tenuous relationships
with the theory of normal forms 
{\sl à la Poincar\'e} and {\sl à la Moser}
will continue to be explored and unveiled 
in several upcoming mathematical memoirs.


\Section{\bf Rigid Equivalences of Rigid
Hypersurfaces in $\C^2$: A Toy Study}
\label{rigid-toy-C-2}
\HEAD{{\ref{rigid-toy-C-2}}.~{\sf Rigid Equivalences of Rigid
Hypersurfaces in $\C^2$: A Toy Study}
}{
Zhangchi {\sc Chen}, Wei Guo {\sc Foo}, Joël {\sc Merker}, 
The Anh {\sc Ta}}

We first consider the equivalence problem of rigid hypersurfaces in
$\C^2$ under the action of rigid biholomorphic transformations. We
will solve this problem with both Cartan's method of equivalence and
Moser's 
method of normal forms. The calculations here are simple,
and they will serve as a toy model for our more substantial problem in
$\C^3$ later. Throughout this section, we use the complex coordinates
$(z,w)$ on $\C^2$ with $w = u + iv, $ where $u, v \in \R$.

We recall that 
a real analytic hypersurface in $\C^2$ is called {\sl rigid} 
if it can
be written $\big\{ u = F(z,\overline{z}) \big\}$, where
$F$ is a converging power series in $z, \overline{z}$.
A local biholomorphic map of $\C^2$ of the form:
\leqnomode\usetagform{default}
\begin{align}
\label{Def-rigid-hypersurface-in-C2}
(z,w) 
\,\longmapsto\,
\big( 
f(z),\, 
a\,w
+
g(z)
\big),
\end{align}
with $a \in \R^{*}$, $c \in \R$, will be called called {\sl rigid}.
Most of the times, we will assume that the origin is fixed,
whence $0 = f(0) = g(0)$.

Since rigid transformations send rigid hypersurfaces to hypersurfaces
which are again rigid, it then makes sense to consider rigid
equivalences of rigid hypersurfaces in $\C^2$, as we do here. The
homogeneous model here is (still) the Heisenberg sphere $\{ u =
z\overline{z} \}$, whose rigid automorphisms 
fixing the origin can be extracted from the
set of general automorphisms of the sphere (exercise).

As a starter, consider a rigid biholomorphic map $(z,w) \longmapsto
\big(f(z),\, a\,w + 
g(z) \big) =: (z', w')$ between
two hypersurfaces
$\{ u = F(z,\overline{z}) \}$ in $\C^2$ and 
$\{ u' = F'(z', \overline{z}' )\}$ in $\C^2$ too. From:
\[
F'
\big(
f(z),
\overline{f}(\overline{z})
\big)
\,=\,
F'\big(z',\overline{z}'\big)
\,=\,
u'
\,=\,
a\,u
+
\Re\,g(z)
\,=\,
a\,F(z,\overline{z})
+
{\textstyle{\frac{1}{2}}}\,
g(z)
+
{\textstyle{\frac{1}{2}}}\,
\overline{g}(\overline{z}),
\]
it comes the {\sl fundamental equation}, identically
satisfied:
\leqnomode\usetagform{default}
\begin{align}
\label{C2-fundamental-identity}
F'
\big(
f(z),
\overline{f}(\overline{z})
\big)
\,\equiv\,
a\,F(z,\overline{z})
+
{\textstyle{\frac{1}{2}}}\,
g(z)
+
{\textstyle{\frac{1}{2}}}\,
\overline{g}(\overline{z}).
\end{align}

\begin{Lemma}
\label{M-3-transfer-LNDG}
Through a rigid biholomorphism between
two rigid
hypersurfaces $\{ u = F\}$ and $\{u' = F'\}$
in $\C^2$, it holds:
\[
F_{z\overline{z}}
\,=\,
{\textstyle{\frac{1}{a}}}\,
\big\vert
f_z
\big\vert^2\,
F_{z'\overline{z}'}'.
\]
\end{Lemma}

\proof
Applying $\partial_z \partial_{\overline{z}}$ eliminates
$g$ and $\overline{g}$ above and yields the result.
\endproof

Thus, $F_{z\overline{z}}$ is a {\sl relative invariant:}
it is nonvanishing in one system of coordinates
if and only if it is nonvanishing in any other system
of coordinates. Of course,
$M$ is {\sl Levi nondegenerate} in the classical sense
if and only if $F_{z\overline{z}} \neq 0$.
We will constantly assume that this holds at {\em every} point.

\Subsection{Cartan's method of equivalence}
Consider a real analytic 
graphed hypersurface $M^3 = \{ u =
F(z,\overline{z}) \}$ 
passing through the origin in $\C^2$.
Its holomorphic
tangent space $T^{1,0} M := (\C \otimes TM) \cap T^{1,0} \C$ is a
$1$-dimensional 
complex vector bundle on $M$. 
One can check directly that
the vector field $\mathcal{L} := \frac{\partial}{\partial z} - i
F_{z}\frac{\partial}{\partial v}$ generates $T^{1,0} M$,
in the intrinsic coordinates $(z, \overline{z}, v)$ on $M$. 
We abbreviate $A :=
-i\, F_{z}$ so that $\mathcal{L} = \frac{\partial}{\partial z} + A\,
\frac{\partial}{\partial v}$ and $\overline{\mathcal{L}}=
\frac{\partial}{\partial \overline{z}} +
\overline{A}\frac{\partial}{\partial v}$. 

Assume that $M$ is everywhere Levi nondegenerate,
namely $F_{z\overline{z}} \neq 0$.
Next, define the real 
vector field $\mathcal{T}$ on $M$ by $\mathcal{T}
:= -i\, [\mathcal{L},\overline{\mathcal{L}}] = \ell\, 
\frac{\partial}{\partial v},$ where $\ell := -2
F_{z\overline{z}}$. 
As in~{\cite{Foo-Merker-Ta-2019}},
introduce also the auxiliary function on $M$:
\[
\Paux
\,:=\, 
\frac{\ell_{z}}{\ell} 
\,=\, 
\frac{F_{zz\overline{z}}}{F_{z\overline{z}}}.
\]

\begin{Lemma}
\label{lm-frame-C2}
The vector fields 
$\mathcal{T}, \mathcal{L}, \overline{\mathcal{L}}$ 
constitute a frame on $\C \otimes TM$, with Lie brackets: 
\[
\big[\mathcal{T},\mathcal{L}\big]
\,=\,
-\,\Paux\,\mathcal{T},
\ \ \ \ \ \ \ \ \ \ \ \ \ \ \ \ \ \ \ \
\big[\mathcal{T},\overline{\mathcal{L}}\big] 
\,=\,
-\,\overline{\Paux}\,\mathcal{T},
\ \ \ \ \ \ \ \ \ \ \ \ \ \ \ \ \ \ \ \
\big[\mathcal{L},\overline{\mathcal{L}}\big] 
\,=\,
-\,i\,\mathcal{T}.
\eqno\qed
\]
\end{Lemma}

Next, denote by $\rho_0,
\zeta_0,\overline{\zeta}_0$ the (complex)
$1$-forms on $M$ which are dual to the (complex) vector fields
$\mathcal{T},
\mathcal{L}, \overline{\mathcal{L}}$, respectively. More
precisely, the expressions of 
$\rho_0, \zeta_0,\overline{\zeta}_0$ in
terms of $dv, dz, d\overline{z}$ are:
\[
\rho_0
\,:=\,
{\textstyle{\frac{1}{\ell}}}\,
\big(
dv
-
A\,dz
-
\overline{A}\,d\overline{z}
\big),
\ \ \ \ \ \ \ \ \ \ \ \ \ \ \ \ \ \ \ \
\zeta_0 
\,:=\,
dz,
\ \ \ \ \ \ \ \ \ \ \ \ \ \ \ \ \ \ \ \
\overline{\zeta}_0
\,=\,
d\overline{z}.
\]
This gives us an initial coframe for $\C \otimes
TM$ having structure equations:
\[
\aligned
d\rho_0
&
\,=\,
\Paux\,
\rho_0
\wedge
\zeta_0 
+ 
\overline{\Paux}\,
\rho_0
\wedge
\overline{\zeta}_0
+
i\,\zeta_0 
\wedge 
\overline{\zeta_0},
\\
d\zeta_0
&
\,=\,
d\overline{\zeta}_0
\,=\,
0.
\endaligned
\]

We now look at the action of rigid transformations on $M$ in order to
setup an initial $G$-structure.  Observe that if a rigid
biholomorphism $h \colon (z, w) \longmapsto \big( f(z),\, aw + g(z)
\big) =: (z',w')$ fixing the origin maps a rigid hypersurface $M
\subset \C^2$ to another rigid hypersurface $M' \subset {\C'}^2$, then
$h$ sends $T^{1,0}M$ to $T^{1,0}M'$, {\em i.e.}  $h_*(T^{1,0}M)
=T^{1,0}M'$.  Without loss of generality, it can be assumed that the
target $M' = \{ u' = F'(z', \overline{z}') \}$ is also graphed, and is
equipped with a similar frame $\{ \mathcal{T}', \mathcal{L}',
\overline{\mathcal{L}}' \}$.  It follows that there exists a uniquely
defined 
nowhere vanishing function $c' \colon M' 
\longrightarrow \C^\ast$ so that $h_*(\mathcal{L}) = c'
\mathcal{L}'$.

Similary, $h_*(\mathcal{T}) = a' \mathcal{T} + b' \mathcal{L} +
\overline{b}' \overline{\mathcal{L}}'$.  From
Definition~{\ref{Def-rigid-hypersurface-in-C2}}, it is clear that
$h_\ast (\partial_v) = a\, \partial_{v'}$.  Since $\mathcal{T} =
\ell\, \partial_v$ and $\mathcal{T}' = \ell'\, \partial_{v'}$, it
comes $h_\ast ( \mathcal{T} ) = a\, \frac{\ell}{\ell'}\,
\mathcal{T}'$. Hence $b'= 0$.
Furthermore:
\[
h_*(\mathcal{T})
\,=\,
h_*\big(-i\,[\mathcal{L},\overline{\mathcal{L}}]\big) 
\,=\,
-\,i\,
\big[
h_*(\mathcal{L}),\,
h_*(\overline{\mathcal{L}})
\big]
\,=\,
-\,i\,
\big[
c'\mathcal{L}',\,
\overline{c}'\overline{\mathcal{L}}'
\big]
\,=\,
c'\overline{c}'\,
\mathcal{T}',
\]
with necessarily $0 \equiv \mathcal{L}'(\overline{c}')$ 
while expanding the bracket
thanks to $b' = 0$, and we conclude that the function 
$a' = c'\overline{c}'$ is determined. 

Consequently, 
under the action of $h$, 
the frame $\{ \mathcal{T}, \mathcal{L}, 
\overline{\mathcal{L}}\}$ changes as:
\[
h_* 
\begin{pmatrix}
\mathcal{T}
\\
\mathcal{L}
\\
\overline{\mathcal{L}}
\end{pmatrix}
=
\begin{pmatrix}
c'\overline{c}' & 0 & 0 
\\
0 & c' & 0 
\\
0 & 0 & \overline{c}' 
\end{pmatrix} 
\begin{pmatrix}
\mathcal{T}'
\\
\mathcal{L}'
\\
\overline{\mathcal{L}}'
\end{pmatrix}
\eqno
{\scriptstyle{(c'\,\neq\,0)}}.
\]
This gives us the transfer relation between the two
{\em dual} coframes,
in terms of a nowhere vanishing function $c \colon M
\longrightarrow \C^\ast$:
\[
h^\ast
\begin{pmatrix}
\rho_0'
\\
\zeta_0'
\\
\overline{\zeta}_0'
\end{pmatrix}
=
\begin{pmatrix}
c\overline{c} & 0 & 0 
\\
0 & c & 0 
\\
0 & 0 & \overline{c} 
\end{pmatrix} 
\begin{pmatrix}
\rho_0 
\\
\zeta_0 
\\
\overline{\zeta}_0
\end{pmatrix}.
\]

The initial $G$-structure is now obtained as follows. 
Such a function $c$ is replaced by a free variable
${\sf c} \in \C^\ast$, an unknown of the problem.
The structure group is the 
$2$-dimensional Lie group of matrices of the form:
\[
g
=
\begin{pmatrix}
{\sf c}\overline{\sf c} & 0 & 0 
\\
0 & {\sf c} & 0 
\\
0 & 0 & \overline{\sf c} 
\end{pmatrix}
\eqno
{\scriptstyle{({\sf c}\,\neq\,0)}},
\]
and we introduce the {\sl lifted coframe:}
\[
\begin{pmatrix}
\rho
\\
\zeta
\\
\overline{\zeta} 
\end{pmatrix}
\,:=\,
g 
\cdot
\begin{pmatrix}
\rho_0
\\
\zeta_0
\\
\overline{\zeta}_0
\end{pmatrix}.
\]

We are now in the position to apply Cartan's method of equivalence to
the $G$-structure just obtained. First, we compute the Maurer-Cartan
matrix as:
\[
dg\cdot g^{-1}
\,=\,
\begin{pmatrix}
\frac{d{\sf c}}{{\sf c}} 
+
\frac{d\overline{\sf c}}{\overline{\sf c}}
& 0 & 0 
\\
0 & \frac{d{\sf c}}{{\sf c}}  & 0 
\\
0 & 0 & \frac{d\overline{\sf c}}{\overline{\sf c}}
\end{pmatrix},
\]
and there is only one (complex-valued) Maurer-Cartan form $\alpha :=
\frac{d{\sf c}}{{\sf c}}$.  
The structure equations are the following:
\[
\aligned
d\rho
&
\,=\,
\big(\alpha+\overline{\alpha}\big) 
\wedge 
\rho 
+
\frac{1}{{\sf c}}\,
\Paux\,
\rho\wedge\zeta
+
\frac{1}{\overline{\sf c}}\,
\overline{\Paux}\,\,
\rho\wedge\overline{\zeta}
+ 
i\,
\zeta\wedge\overline{\zeta},
\\
d\zeta
&
\,=\,
\alpha\wedge\zeta,
\\
d\overline{\zeta}
&
\,=\,
\overline{\alpha}\wedge\overline{\zeta}.
\endaligned
\]

We proceed to absorption of torsion by introducing the 
{\sl modified Maurer-Cartan form}:
\[
\pi
\,:=\,
\alpha 
- 
{\textstyle{\frac{1}{{\sf c}}}}\,
\Paux\,
\zeta, 
\]
in terms of which the structure equations contract as:
\[
\aligned
d\rho
&
\,=\,
(\pi+\overline{\pi}) 
\wedge 
\rho
+ 
i\,\zeta\wedge\overline{\zeta},
\\
d\zeta
&
\,=\,
\pi \wedge \zeta,
\ \ \ \ \ \ \ \ \ \ \ \ \ \ \ \ \ \ \ \
\ \ \ \ \ \ \ \ \ \ \ \ \ \ \ \ \ \ \ \
d\overline{\zeta}
\,=\,
\overline{\pi} \wedge \overline{\zeta}.
\endaligned
\]

At this point, no more absorption can be performed, because if one
modifies the $1$-form $\pi$ as $\tilde{\pi} := \pi - A\, \rho -
B\,\zeta - C\,\overline{\zeta}$, which transforms the structure
equations into:
\[
\aligned
d\rho
&
\,=\,
\big(\tilde{\pi}+\overline{\tilde{\pi}}\big) 
\wedge \rho
-
(B+\overline{C})\,
\rho\wedge\zeta
-
(\overline{B}+C)\,
\rho\wedge\overline{\zeta}
+ 
i\,\zeta \wedge \overline{\zeta} ,
\\
d\zeta
&
\,=\,
\tilde{\pi} \wedge \zeta 
+
A\,\rho\wedge\zeta
-
C\,\zeta \wedge \overline{\zeta},
\endaligned
\]
all the functions $A$, $B$, $C$ must be zero to conserve the same
shape. In other words, the prolongation reduces to identity, and $\pi$
is uniquely defined.

Therefore, Cartan's process stops, and to finish, it remains to
finalize the expression of:
\[
\aligned
d\pi
&
\,=\, 
\zero{d\alpha}
+
{\textstyle{\frac{1}{{\sf c}}}}\,
{\textstyle{\frac{d{\sf c}}{{\sf c}}}}\,
\Paux\,
\wedge\zeta 
- 
{\textstyle{\frac{1}{{\sf c}}}}\,
d\Paux
\wedge
\zeta 
-
{\textstyle{\frac{1}{{\sf c}}}}\,
\Paux\,
d\zeta
\\
&
\,=\,
0
+
{\textstyle{\frac{1}{{\sf c}}}}\,
\big(
\pi
+
{\textstyle{\frac{1}{{\sf c}}}}\,
\Paux\,
\zeta
\big)\,
\Paux
\wedge
\zeta
-
{\textstyle{\frac{1}{{\sf c}}}}\,
\big(
\Paux_z\,
dz
+
\Paux_{\overline{z}}\,
d\overline{z}
\big)
\wedge\zeta
-
{\textstyle{\frac{1}{{\sf c}}}}\,
\Paux\,
\pi\wedge\zeta
\\
&
\,=\,
-\,
{\textstyle{\frac{1}{{\sf c}}}}\,
\big(
\Paux_z\,
{\textstyle{\frac{1}{{\sf c}}}}\,
\zeta
+
\Paux_{\overline{z}}\,
{\textstyle{\frac{1}{\overline{\sf c}}}}\,
\overline{\zeta}
\big)
\wedge\zeta,
\endaligned
\]
where we need to know\big/abbreviate just:
\[
\Paux_{\overline{z}}
\,=\,
{\textstyle{
\frac{F_{zz\overline{z}\overline{z}} F_{z\overline{z}} 
-F_{zz\overline{z}}\,F_{z\overline{z}\overline{z}}}{
(F_{z\overline{z}})^2}
}}
\,\,=:\,\,
\Raux,
\]
whence:
\[
d\pi
\,=\,
{\textstyle{\frac{1}{{\sf c}\overline{\sf c}}}}\,
\Raux\,
\zeta
\wedge
\overline{\zeta}.
\]
Visibly, $\overline{\Raux} = \Raux$ is real, because 
$\overline{F} = F$ is, whence $\overline{F_{z^a\overline{z}^c}}
= F_{\overline{z}^a z^c}$.

\begin{Theorem}
The equivalence problem under local rigid biholomorphisms of
$\mathcal{C}^{\omega}$ rigid real hypersurfaces
$\{u=F(z,\overline{z})\}$ in $\C^2$ whose Levi form 
is everywhere nondegenerate reduces
to classifying $\{e\}$-structures on the $5$-dimensional bundle
$M^3 \times \C$ equipped with coordinates
$(z, \overline{z}, v, {\sf c}, \overline{\sf c})$ 
together with a coframe of $5$ differential $1$-forms:
\[
\big\{
\rho,\ 
\zeta,\
\overline{\zeta},\
\pi,\ 
\overline{\pi}
\big\}
\eqno
{\scriptstyle{(\overline{\rho}\,=\,\rho)}},
\]
which satisfy invariant structure equations of the shape:
\[
\aligned
d\rho 
&
\,=\, 
(\pi+\overline{\pi})
\wedge\rho 
+ 
i\,
\zeta\wedge\overline{\zeta},
\\
d\zeta 
&
\,=\,
\pi\wedge\zeta,
&
\ \ \ \ \ \ \ \ \ \ \ \ \ \ \ \ \ \ \ \
\ \ \ \ \ \ \ \ \ \ \ \ \ \ \ \ \ \ \ \
d\overline{\zeta}
&
\,=\,
\overline{\pi}\wedge\overline{\zeta},
\\
d\pi
&
\,=\,
{\textstyle{\frac{1}{{\sf c}\overline{\sf c}}}}\,
\Raux\,
\zeta
\wedge
\overline{\zeta},
&
\ \ \ \ \ \ \ \ \ \ \ \ \ \ \ \ \ \ \ \
\ \ \ \ \ \ \ \ \ \ \ \ \ \ \ \ \ \ \ \
d\overline{\pi}
&
\,=\,
-\,
{\textstyle{\frac{1}{{\sf c}\overline{\sf c}}}}\,
\overline{\Raux}\,
\zeta
\wedge
\overline{\zeta}.
\endaligned
\]
\end{Theorem}

Another way to see that $\overline{\Raux} = \Raux$ is real
from the structure equations is as follows,
using Poincar\'e's relation:
\[
\footnotesize
\aligned
0
\,=\,
d\circ d\rho
&
\,=\,
\big(
d\pi
+
d\overline{\pi}
\big)
\wedge
\rho
-
\big(
\pi
+
\overline{\pi}
\big)
\wedge
d\rho
+
i\,d\zeta
\wedge
\overline{\zeta}
-
i\,\zeta
\wedge
d\overline{\zeta}
\\
&
\,=\,
\frac{1}{{\sf c}\overline{\sf c}}\,
\Raux\,
\zeta
\wedge
\overline{\zeta}
\wedge
\rho
+
\frac{1}{{\sf c}\overline{\sf c}}\,
\overline{\Raux}\,
\overline{\zeta}
\wedge
\zeta
\wedge
\rho
-
\big(
\pi+\overline{\pi}
\big)
\Big[
\zero{
\big(
\pi+\overline{\pi}
\big)}
\wedge\rho
+
i\,\zeta\wedge\overline{\zeta}
\Big]
+
i\,
\pi\wedge\zeta\wedge\overline{\zeta}
-
i\,\zeta\wedge\overline{\pi}\wedge\overline{\zeta}
\\
&
\,=\,
\frac{1}{{\sf c}\overline{\sf c}}\,
\big(
\Raux
-
\overline{\Raux}
\big)\,
\rho\wedge\zeta\wedge\overline{\zeta}.
\endaligned
\]

Thus, the only invariant here is:
\leqnomode\usetagform{default}
\begin{align}
\label{R-invariant}
\Raux 
:=
\frac{
F_{zz\overline{z}\overline{z}}\,F_{z\overline{z}} 
-
F_{z z\overline{z}}\, 
F_{z\overline{z}\overline{z}} 
}{
(F_{z\overline{z}})^2}.
\end{align}
When $\Raux \equiv 0$, the structure equations have constants
coefficients, which shows,
by Cartan's theory, that all rigid hypersurfaces with $\Raux
\equiv 0$ are rigidly equivalent to each other, and equivalent to the
model $\{ u = z \overline{z} \}$.
There also are straightforward arguments to get this.

\begin{Proposition}
\label{Assertion-M-flat-R-zero}
A rigid $M = \{ u = F(z, \overline{z}) \}$ 
in $\C^2$ is
rigidly biholomorphically equivalent to
the Heisenberg sphere $\{ u' = z' \overline{z} ' \}$
if and only if:
\[
0
\,\equiv\,
\Raux(F)
\,\equiv\,
F_{zz\overline{z}\overline{z}}\,
F_{z\overline{z}}
-
F_{zz\overline{z}}\,
F_{z\overline{z}\overline{z}}.
\]
\end{Proposition}

\proof
Recall that the condition $\Raux(F) \equiv 0$ is invariant
under rigid biholomorphisms.

Trivially, $F := z \overline{z}$ implies $\Raux(F) \equiv 0$.

For the converse,
Lemma~{\ref{M-3-transfer-LNDG}}
guarantees that $M$ is of course Levi-nondegenerate too, and
by invariancy of $\Raux = 0$, we can assume that 
$F = z\overline{z} + {\rm O}_{z,\overline{z}}(3)$.

Set $G := F_{z\overline{z}}$, a function which is also real-valued,
with $G(0) = 1$. Thus:
\[
0
\,\equiv\,
G_{z\overline{z}}\,G
-
G_z\,G_{\overline{z}}
\ \ \ \ \ \ \ \ \ \ \ \ \ \ \ \ \ \ \ \
\Longleftrightarrow
\ \ \ \ \ \ \ \ \ \ \ \ \ \ \ \ \ \ \ \
\big(\log\,G\big)_{z\overline{z}}
\,\equiv\,
0.
\]
Consequently $\log\, G(z, \overline{z}) = \varphi(z) + 
\overline{\varphi} ( \overline{z} )$ for some holomorphic
function with $\varphi (0) = 0$, whence 
$G(z,\overline{z}) = \psi(z)\cdot \overline{\psi}(\overline{z})$
with $\psi(0) = 1$, and 
\[
F(z,\overline{z})
\,=\,
\int_0^z\,
\psi(\zeta)\,d\zeta
\cdot
\int_0^{\overline{z}}\,
\overline{\psi}(\overline{\zeta})\,
d\overline{\zeta}
\,=:\,
f(z)
\cdot
\overline{f}(\overline{z}),
\]
with $f(z) = z + {\rm O}_z(2)$. Thus $u = f(z)\, \overline{f}
(\overline{z})$, and the rigid biholomorphism
$z' := f(z)$ terminates.
\endproof

We know from Lemma~{\ref{M-3-transfer-LNDG}}
that $F_{z\overline{z}}$ is a relative invariant.
What about $\Raux$? It suffices to examine
how the numerator of $\Raux$ behaves
under transformations.

\begin{Lemma}
Through a rigid biholomorphism 
$(z, w) \longmapsto \big(f(z),\, a\,w + g(z) \big) =:
(z', w')$
between two rigid hypersurfaces $\{ u =
F\}$ and $\{u' = F'\}$ in $\C^2$, it holds:
\[
F_{zz\overline{z}\overline{z}}\,
F_{z\overline{z}}
-
F_{zz\overline{z}}\,
F_{z\overline{z}\overline{z}}
\,\equiv\,
{\textstyle{\frac{1}{a^2}}}\,
\big(
f_z\,\overline{f}_{\overline{z}}
\big)^3\,
\Big[
F_{z'z'\overline{z}'\overline{z}'}'\,
F_{z'\overline{z}'}'
-
F_{z'z'\overline{z}'}'\,
F_{z'\overline{z}'\overline{z}'}'
\Big].
\]
\end{Lemma}

\proof
Differentiate the fundamental 
identity~({\ref{C2-fundamental-identity}})
four appropriate times:
\[
\aligned
a\,F_{z\overline{z}}
&
\,\equiv\,
f_z\,\overline{f}_{\overline{z}}\,
F_{z'\overline{z}'}',
\\
a\,
F_{zz\overline{z}}
&
\,\equiv\,
f_{zz}\,\overline{f}_{\overline{z}}\,
F_{z'\overline{z}'}'
+
f_z\overline{f}_{\overline{z}}f_z\,
F_{z'z'\overline{z}'}',
\\
a\,
F_{z\overline{z}\overline{z}}
&
\,\equiv\,
f_z\overline{f}_{\overline{z}\overline{z}}\,
F_{z'\overline{z}'}'
+
f_z\overline{f}_{\overline{z}}\overline{f}_{\overline{z}}\,
F_{z'\overline{z}'\overline{z}'},
\\
a\,
F_{zz\overline{z}\overline{z}}
&
\,\equiv\,
f_{zz}\overline{f}_{\overline{z}\overline{z}}\,
F_{z'\overline{z}'}'
+
f_{zz}\overline{f}_{\overline{z}}\overline{f}_{\overline{z}}\,
F_{z'\overline{z}'\overline{z}'}'
+
f_z\overline{f}_{\overline{z}\overline{z}}f_z\,
F_{z'z'\overline{z}'}'
+
f_z\overline{f}_{\overline{z}}f_z\overline{f}_{\overline{z}}\,
F_{z'z'\overline{z}'\overline{z}'}',
\endaligned
\]
perform the necessary products, substract, and get the result.
\endproof

\Subsection{Method of normal forms of Moser}
In this subsection, following the method of Moser, we will approach
the equivalence problem for rigid hypersurfaces in $\C^2$ under rigid
biholomorphisms by constructing a normal form.  Notice that although
the problem is (much) simpler than that considered by Moser for
general hypersurfaces in $\C^2$, our problem here is not a special
case of what is already known. 

The goal is to simplify the defining function $u = F(z, \overline{z})$
of a given hypersurface $M^3 \subset \C^2$ as much as possible by
applying rigid holomorphic changes of variables $(z,w) \mapsto \big(
f(z),\, \rho\, w+ g(z) \big) =: (z', w')$,
with $\rho \in \R^\ast$.  We will find step by step
changes, so that the transformed graphing functions $F'$ for
successive $M' = \big\{ u' = F'(z', \overline{z'}) \big\}$ will
contain more and more zero coefficients.

Take a real analytic hypersurface
$M = \{ u = F(z, \overline{z}) \}$ 
passing through the
origin in $\C^2$, and expand:
\[
u
\,=\,
{\textstyle{\frac{1}{2}}}\,
\big(
w+\overline{w}
\big)
\,=\,
\sum_{j+k\geqslant 1}\,
F_{j,k}\,
z^j\overline{z}^k,
\]
with $F_{j,k} = \overline{F}_{k,j}$.
At first, set $z' := z$ and:
\[
w'
\,:=\,
w
-
2\,
\smallsum{j\geqslant 1}\,
F_{j,0}\,z^j,
\]
in order to subtract all harmonic 
monomials $F_{j,0}\, z^j$ and $F_{0,k}\, \overline{z}^k$
to obtain:
\[
u'
\,=\,
\sum_{j\geqslant 1
\atop
k\geqslant 1}\,
F_{j,k}\,
z^j\overline{z}^k
\,\,=\,\,
F_{1,1}\,
z\overline{z}
+
\sum_{j+k\geqslant 3
\atop
j\geqslant 1\,\text{and}\,k\geqslant 1}\,\,
F_{j,k}\,
z^j\overline{z}^k.
\]

The invariant property $F_{1, 1} \neq 0$ characterizes Levi
nondegeneracy of $M$ at the origin (hence in a neighborhood).
Switching $u \longmapsto -\, u$ if necessary, we may assume $F_{1,1} >
0$.

Next, make the rigid biholomorphism $z' := \sqrt{F_{1,1}}\, z$ with
$w' := w$, drop the prime, single out monomials of degree $1$ in
either $z$ or $\overline{z}$, factorize, and point out remainders:
\[
\!\!\!\!\!\!\!\!\!\!\!\!\!\!\!\!\!\!\!\!
\footnotesize
\aligned
u
&
\,=\,
z\overline{z}
+
\sum_{j+k\geqslant 3
\atop
j\geqslant 1\,\text{and}\,k\geqslant 1}\,\,
\frac{F_{j,k}}{
\sqrt{F_{1,1}}^{\,j+k}}\,
z^j\overline{z}^k
\\
&
\,=\,
z\overline{z}
+
\overline{z}\,
\bigg(
\frac{F_{2,1}}{F_{1,1}^{3/2}}\,
z^2
+
\sum_{j\geqslant 3}\,
\frac{F_{j,1}}{F_{1,1}^{(j+1)/2}}\,
z^j
\bigg)
+
z\,
\bigg(
\frac{F_{1,2}}{F_{1,1}^{3/2}}\,
\overline{z}^2
+
\sum_{k\geqslant 3}\,
\frac{F_{1,k}}{F_{1,1}^{(1+k)/2}}\,
\overline{z}^k
\bigg)
+
\frac{F_{2,2}}{F_{1,1}^2}\,
z^2\overline{z}^2
+
\sum_{j+k\geqslant 5
\atop
j\geqslant 2\,\text{\rm and}\,k\geqslant 2}\,
\frac{F_{j,k}}{F_{1,1}^{(j+k)/2}}\,\,
z^j\overline{z}^k
\\
&
\,=\,
\bigg(
z
+
\frac{F_{2,1}}{F_{1,1}^{3/2}}\,
z^2
+
\sum_{j\geqslant 3}\,
\frac{F_{j,1}}{F_{1,1}^{(j+1)/2}}\,
z^j
\bigg)\,
\bigg(
\overline{z}
+
\frac{F_{1,2}}{F_{1,1}^{3/2}}\,
\overline{z}^2
+
\sum_{k\geqslant 3}\,
\frac{F_{1,k}}{F_{1,1}^{(1+k)/2}}\,
\overline{z}^k
\bigg)
-
\frac{F_{2,1}\,F_{1,2}}{F_{1,1}^3}\,
z^2\overline{z}^2
-
z^2\overline{z}^3\big(\cdots\big)
-
z^3\overline{z}^2\big(\cdots\big)
\,+
\\
&
\ \ \ \ \
+
\frac{F_{2,2}}{F_{1,1}^2}\,
z^2\overline{z}^2
+
z^2\overline{z}^3\big(\cdots\big)
+
z^3\overline{z}^2\big(\cdots\big).
\endaligned
\]

Such a factorization suggests to 
perform the rigid biholomorphism:
\[
z'
\,:=\,
z
+
\frac{F_{2,1}}{F_{1,1}^{3/2}}\,
z^2
+
\sum_{j\geqslant 3}\,
\frac{F_{j,1}}{F_{1,1}^{(j+1)/2}}\,
z^j,
\]
again with untouched $w' := w$.
Its inverse is of the form $z = z' \big( 1 + {z'}^2(\cdots) \big)$,
so ${\rm O} \big( z^l \overline{z}^m \big) = 
{\rm O} \big( {z'}^l {\overline{z}'}^m \big)$, and finally,
dropping primes, we have proved the

\begin{Proposition}
Any rigid $M = \big\{ u = \sum\, F_{j,k}\,
z^j \overline{z}^k \big\}$
can be brought, by a rigid biholomorphic transformation
fixing the origin,
to:
\[
u
\,=\,
z\overline{z}
+
\Big[
\frac{F_{2,2}\,F_{1,1}-F_{2,1}\,F_{1,2}}{F_{1,1}^3}
\Big]\,
z^2\overline{z}^2
+
z^2\overline{z}^3
\big(\cdots\big)
+
z^3\overline{z}^2
\big(\cdots\big).
\eqno\qed
\]
\end{Proposition}

In other words:
\reqnomode\usetagform{EngelLie}
\begin{align}
0
&
\,=\,
F_{j,0}
\,=\,
F_{0,k}
\tag{(j\,\geqslant\,1,\,\,k\,\geqslant\,1),}
\\
1
&
\,=\,
F_{1,1},
\notag
\\
0
&
\,=\,
F_{j,1}
\,=\,
F_{1,k}
\tag{(j\,\geqslant\,2,\,\,k\,\geqslant\,2).}
\end{align}

Can one normalize the graphing function $F$ further?
For instance, can one annihilate some other $F_{j,k}$?
Not much freedom is left, as states the next

\begin{Lemma}
If two rigid hypersurfaces in $\C^2$ having the form:
\[
u
\,=\,
z\overline{z}
+
\smallsum{j,k\geqslant 2}\,
F_{j,k}\,
z^j\overline{z}^k
\ \ \ \ \ \ \ \ \ \ \ \ \ \ \ \ \ \ \ \
\text{and}
\ \ \ \ \ \ \ \ \ \ \ \ \ \ \ \ \ \ \ \
u'
\,=\,
z'\overline{z}'
+
\smallsum{j,k\geqslant 2}\,
F_{j,k}'\,
{z'}^j{\overline{z}'}^k,
\]
are equivalent through a rigid
biholomorphism fixing the origin, then
there exist $\rho \in \R_+^\ast$ and $\varphi \in \R$ such that:
\[
z'
\,=\,
\rho^{1/2}\,e^{i\varphi}\,z,
\ \ \ \ \ \ \ \ \ \ \ \ \ \ \ \ \ \ \ \ \ \ \ \ \ \
w'
\,=\,
\rho\,w.
\]
\end{Lemma}

In particular, this shows that the group of rigid
transformations fixing the origin
of the Heisenberg sphere 
$\{ u = z \overline{z} \}$ is $2$-dimensional,
generated by these obvious rotation\big/dilation 
commuting transformations (solution of the exercise).

\proof
Write as above $(z', w') = \big( f(z),\, \rho\,w + g(z) \big)$,
with $f(0) = 0 = g(0)$.
The fundamental 
equation~({\ref{fundamental-equation-before-weighting}}) reads:
\[
\rho\,F(z,\overline{z})
+
{\textstyle{\frac{1}{2}}}\,
g(z)
+
{\textstyle{\frac{1}{2}}}\,
\overline{g}(\overline{z})
\,\,\equiv\,\,
F'
\big(
f(z),\overline{f}(\overline{z})
\big).
\]
Put $\overline{z} := 0$, get $\overline{g} (\overline{z}) \equiv 0$.
Thus:
\[
\rho\,
\big(
z\overline{z}
+
z^2\overline{z}^2
(\cdots)
\big)
\,\equiv\,
f(z)\overline{f}(\overline{z})
+
f(z)^2\overline{f}(\overline{z})^2\,
\big(
\cdots
\big),
\]
and using $f(z) = {\rm O}(z)$:
\[
\rho\,z\overline{z}
\,\equiv\,
f(z)\overline{f}(\overline{z})
+
z^2\overline{z}^2
\big(\cdots\big).
\]
Invertibility of the Jacobian yields $f_z(0) \neq 0$.
Apply $\partial_{\overline{z}} \big\vert_0$ and get:
\[
\rho\,z
\,\equiv\,
f(z)\,
\overline{f}'(0),
\]
so $f(z) = \lambda \, z$ for some $\lambda \in \C^\ast$.
Lastly, $\rho = \lambda \overline{\lambda}$, 
which concludes.
\endproof

\begin{Corollary}
\label{Cor-toy-Moser-C2}
Two rigid hypersurfaces in $\C^2$:
\[
u
\,=\,
z\overline{z}
+
\smallsum{j,k\geqslant 2}\,
F_{j,k}\,
z^j\overline{z}^k
\ \ \ \ \ \ \ \ \ \ \ \ \ \ \ \ \ \ \ \
\text{and}
\ \ \ \ \ \ \ \ \ \ \ \ \ \ \ \ \ \ \ \
u'
\,=\,
z'\overline{z}'
+
\smallsum{j,k\geqslant 2}\,
F_{j,k}'\,
{z'}^j{\overline{z}'}^k,
\]
are rigidly biholomorphically equivalent
if and only if 
there exist $\rho \in \R_+^\ast$ and $\varphi \in \R$ such that:
\[
F_{j,k}
\,=\,
\rho^{\frac{j+k-2}{2}}\,\,
e^{i\,\varphi\,(j-k)}\,\,
F_{j,k}'
\eqno
{\scriptstyle{(j\,\geqslant\,2,\,\,k\,\geqslant\,2)}}.
\eqno\qed
\]
\end{Corollary}

At any point $(z_0, w_0) \in M$ close to the origin,
all these results are also valid, and
using the recentered 
holomorphic coordinates $z-z_0$ and $w-w_0$, one obtains:
\[
u-u_0
\,=\,
(z-z_0)\,
\big(\overline{z}-\overline{z}_0\big)
+
\frac{4\,F_{zz\overline{z}\overline{z}}(z_0)\,
F_{z\overline{z}}(z_0)
-
2\,F_{zz\overline{z}}(z_0)\,
2\,F_{z\overline{z}\overline{z}}(z_0)}{
F_{z\overline{z}}(z_0)^3}\,
(z-z_0)^2\,
\big(
\overline{z}-\overline{z}_0
\big)^2
+
\cdots.
\]
The $(2,2)$-coefficient at various points $z_0$ is,
up to a power of $F_{z\overline{z}}$ in the denominator,
exactly equal to the relative invariant function
$\Raux$ found in~({\ref{R-invariant}}) by 
applying Cartan's method.

According to Lie's principle of thought
({\cite[Chap.~1]{Lie-Merker-2015}),
a relative invariant is assumed to be either
identically zero, or nowhere zero, after
restriction to an appropriate open subset. 
Since Proposition~{\ref{Assertion-M-flat-R-zero}} already
understood the branch $\Raux \equiv 0$,
it remains only to treat the branch $\Raux \neq 0$.
This is left as an exercise. 

\Section{\bf Two Invariant Determinants 
for Hypersurfaces $M^5 \subset \C^3$}
\label{two-invariant-determinants}
\HEAD{{\ref{two-invariant-determinants}}.~{\sf Two Invariant 
Determinants for Hypersurfaces $M^5 \subset \C^3$}
}{
Zhangchi {\sc Chen}, Wei Guo {\sc Foo}, Joël {\sc Merker}, 
The Anh {\sc Ta}}

Consider a rigid biholomorphism:
\[
H
\colon\ \ \
(z,\zeta,w)
\,\,\,\longmapsto\,\,\,
\Big(
f(z,\zeta),\,
g(z,\zeta),\,
\rho\,w+h(z,\zeta)
\Big)
\,=:\,
\big(z',\zeta',w'\big)
\eqno
{\scriptstyle{(\rho\,\in\,\R^\ast)}},
\]
hence with Jacobian $f_zg_\zeta -f_\zeta g_z \neq 0$, 
between two rigid $\mathcal{C}^\omega$ hypersurfaces:
\[
w
\,=\,
-\,\overline{w}
+
2\,F\big(z,\zeta,\overline{z},\overline{\zeta}\big)
\,=:\,
Q
\ \ \ \ \ \ \ \ \ \ \ \ \ \ \ \ \ \ \ \
\text{and}
\ \ \ \ \ \ \ \ \ \ \ \ \ \ \ \ \ \ \ \
w'
\,=\,
-\,\overline{w}'
+
2\,F'\big(z',\zeta',\overline{z}',\overline{\zeta}'\big)
\,=:\,
Q'.
\]
Plugging the three components of $H$ in the target equation:
\[
\rho\,w
+
h(z,\zeta)
+
\rho\,\overline{w}
+
\overline{h}(\overline{z},\overline{\zeta})
\,\,=\,\,
2\,F'
\Big(
f(z,\zeta),\,
g(z,\zeta),\,
\overline{f}(\overline{z},\overline{\zeta}),\,
\overline{g}(\overline{z},\overline{\zeta})
\Big),
\]
and replacing $w+\overline{w} = 2\,F$, one receives 
the {\sl fundamental equation} expressing $H(M) \subset M'$:
\[
2\,\rho\,F\big(z,\zeta,\overline{z},\overline{\zeta}\big)
+
h(z,\zeta)
+
\overline{h}(\overline{z},\overline{\zeta})
\,\,\equiv\,\,
2\,F'
\Big(
f(z,\zeta),\,
g(z,\zeta),\,
\overline{f}(\overline{z},\overline{\zeta}),\,
\overline{g}(\overline{z},\overline{\zeta})
\Big).
\]

By differentiating it (exercise! use a computer!), 
one expresses as follows the invariancy of the Levi determinant
defined for general biholomorphisms~{\cite{Merker-Nurowski-2019}} as:
\[
\left\vert\!
\begin{array}{ccc}
Q_{\overline{z}} & Q_{\overline{\zeta}} 
& 
Q_{\overline{w}}
\\
Q_{z\overline{z}} & Q_{z\overline{\zeta}} 
& 
Q_{z\overline{w}}
\\
Q_{\zeta\overline{z}} & Q_{\zeta\overline{\zeta}} 
& 
Q_{\zeta\overline{w}}
\end{array}
\!\right\vert
\,=\,
2^2\,
\left\vert\!
\begin{array}{ccc}
F_{\overline{z}} & F_{\overline{\zeta}} 
& 
-1
\\
F_{z\overline{z}} & F_{z\overline{\zeta}} 
& 
0
\\
F_{\zeta\overline{z}} & F_{\zeta\overline{\zeta}} 
& 
0
\end{array}
\!\right\vert.
\]

\begin{Proposition}
\label{Prp-Levi-determinant}
Through any rigid biholomorphism:
\[
\left\vert\!
\begin{array}{cc}
F_{z'\overline{z}'}'
&
F_{z'\overline{\zeta}'}'
\\
F_{\zeta'\overline{z}'}'
&
F_{\zeta'\overline{\zeta}'}'
\end{array}
\!\right\vert
\,\,=\,\,
\frac{\rho^2}{
\left\vert\!
\begin{array}{cc}
f_z & f_\zeta
\\
g_z & g_\zeta
\end{array}
\!\right\vert\,\,
\left\vert\!
\begin{array}{cc}
\overline{f}_{\overline{z}} & \overline{f}_{\overline{\zeta}}
\\
\overline{g}_{\overline{z}} & \overline{g}_{\overline{\zeta}}
\end{array}
\!\right\vert}\,\,
\left\vert\!
\begin{array}{cc}
F_{z\overline{z}}
&
F_{z\overline{\zeta}}
\\
F_{\zeta\overline{z}}
&
F_{\zeta\overline{\zeta}}
\end{array}
\!\right\vert.
\eqno\qed
\]
\end{Proposition}

Consequently, the property that the Levi form is of constant
rank $1$ is biholomorphically invariant.
The $2$-nondegeneracy property~{\cite{Merker-Nurowski-2019}} 
then expresses as the nonvanishing of: 
\[
\left\vert\!
\begin{array}{ccc}
Q_{\overline{z}} & Q_{\overline{\zeta}} 
& 
Q_{\overline{w}}
\\
Q_{z\overline{z}} & Q_{z\overline{\zeta}} 
& 
Q_{z\overline{w}}
\\
Q_{zz\overline{z}} & Q_{zz\overline{\zeta}} 
& 
Q_{zz\overline{w}}
\end{array}
\!\right\vert
\,=\,
2^2\,
\left\vert\!
\begin{array}{ccc}
F_{\overline{z}} & F_{\overline{\zeta}} 
& 
-1
\\
F_{z\overline{z}} & F_{z\overline{\zeta}} 
& 
0
\\
F_{zz\overline{z}} & F_{zz\overline{\zeta}} 
& 
0
\end{array}
\!\right\vert.
\]

\begin{Proposition}
\label{Prp-2ndg-determinant}
When the Levi form is of constant rank $1$, 
through any rigid biholomorphism:

\[
\left\vert\!
\begin{array}{cc}
F_{z'\overline{z}'}'
&
F_{z'\overline{\zeta}'}'
\\
F_{z'z'\overline{z}'}'
&
F_{z'z'\overline{\zeta}'}'
\end{array}
\!\right\vert
\,\,=\,\,
\frac{\rho^2\,\,
\big(
g_\zeta\,F_{z\overline{z}}
-
g_z\,F_{\zeta\overline{z}}
\big)^3}{
\left\vert\!
\begin{array}{cc}
f_z & f_\zeta
\\
g_z & g_\zeta
\end{array}
\!\right\vert^{3}\,\,
\left\vert\!
\begin{array}{cc}
\overline{f}_{\overline{z}} & \overline{f}_{\overline{\zeta}}
\\
\overline{g}_{\overline{z}} & \overline{g}_{\overline{\zeta}}
\end{array}
\!\right\vert}\,\,
\left\vert\!
\begin{array}{cc}
F_{z\overline{z}}
&
F_{z\overline{\zeta}}
\\
F_{zz\overline{z}}
&
F_{zz\overline{\zeta}}
\end{array}
\!\right\vert.
\eqno\qed
\]

\end{Proposition}

Recall that we denote the class of 
(local) hypersurfaces $M^5 \subset \C^3$
passing by the origin $0 \in M$ 
that are $2$-nondegenerate and whose Levi form
has constant rank $1$ as:
\[
\mathfrak{C}_{2,1}.
\]

\Section{\bf Rigid Infinitesimal CR Automorphisms
\\ 
of the Gaussier-Merker Model}
\label{rigid-infinitesimal-CR-automorphisms-GM-model}
\HEAD{{\ref{rigid-infinitesimal-CR-automorphisms-GM-model}}.~{\sf 
Rigid Infinitesimal CR Automorphisms of the Gaussier-Merker Model}
}{
Zhangchi {\sc Chen}, Wei Guo {\sc Foo}, Joël {\sc Merker}, 
The Anh {\sc Ta}}

The appropriate model $M_{\sf LC}$
is rigid and was set up by Gaussier-Merker
in~{\cite{Gaussier-Merker-2003}} and
Fels-Kaup in~{\cite{Fels-Kaup-2007}}:
\[
M_{\sf LC}
\colon
\ \ \ \ \
u
\,=\,
\frac{z\overline{z}+\frac{1}{2}\,z^2\overline{\zeta}
+\frac{1}{2}\overline{z}^2\zeta}{1-\zeta\overline{\zeta}}
\,\,=:\,\,
\maux\big(z,\zeta,\overline{z},\overline{\zeta}\big).
\]
It is a locally graphed representation of the
tube in $\C^3$ over the future light cone in $\R^3$.
The $10$-dimensional simple Lie algebra of its infinitesimal 
CR automorphisms:
\[
\mathfrak{g}
\,:=\,
\mathfrak{aut}_{CR}\big(M_{\sf LC}\big)
\,\cong\,
\mathfrak{so}_{2,3}(\R),
\]
has $10$ natural generators
$X_1, \dots, X_{10}$, which are $(1,0)$ vector fields 
having holomorphic coefficients with
$X_\sigma + \overline{X}_\sigma$ tangent to $M_{\sf LC}$.
Assigning weights to variables, to vector fields, and 
the same weights to their conjugates:
\leqnomode\usetagform{default}
\begin{align}
\label{weighting-z-zeta-w}
[z]
\,:=\,
1
\ \ \ \ \ \ \ \ 
[\zeta]
\,:=\,
0,
\ \ \ \ \ \ \ \ 
[w]
\,:=\,
2
\ \ \ \ \ \ \ \ 
\big[\partial_z\big]
\,:=\,
-\,1
\ \ \ \ \ \ \ \ 
\big[\partial_\zeta\big]
\,:=\,
0
\ \ \ \ \ \ \ \ 
\big[\partial_w\big]
\,:=\,
-\,2,
\end{align}
this Lie algebra of vector fields
isomorphic to $\mathfrak{so}_{2,3}(\R)$ can be graded as:
\[
\mathfrak{g}
\,=\,
\mathfrak{g}_{-2}
\oplus
\mathfrak{g}_{-1}
\oplus
\mathfrak{g}_0
\oplus
\mathfrak{g}_1
\oplus
\mathfrak{g}_2,
\]
where, as shown in~{\cite{Gaussier-Merker-2003,
Foo-Merker-Ta-2019}}:
\[
\aligned
\mathfrak{g}_{-2}
&
\,:=\,
\Span\,
\big\{
i\,\partial_w
\big\},
\\
\mathfrak{g}_{-1}
&
\,:=\,
\Span\,
\big\{
(\zeta-1)\,\partial_z
-
2z\,\partial_w,
\ \ \
(i+i\zeta)\,\partial_z
-
2iz\,\partial_w
\big\},
\endaligned
\]
where $\mathfrak{g}_0 = \mathfrak{g}_0^{\sf trans} 
\oplus \mathfrak{g}_0^{\sf iso}$:
\[
\aligned
\mathfrak{g}_0^{\sf trans}
&
\,:=\,
\Span\,
\Big\{
z\zeta\,\partial_z
+
(\zeta^2-1)\,\partial_\zeta
-
z^2\,\partial_w,
\ \ \
iz\zeta\,\partial_z
+
(i+i\zeta^2)\,\partial_\zeta
-
iz^2\,\partial_w
\Big\},
\\
\mathfrak{g}_0^{\sf iso}
&
\,:=\,
\Span\,
\big\{
z\,\partial_z
+
2w\,\partial_w,
\ \ \ 
iz\,\partial_z
+
2i\zeta\,\partial_\zeta
\big\},
\endaligned
\]
while:
\[
\aligned
\mathfrak{g}_1
&
\,:=\,
\Span\,
\big\{
\big(z^2-\zeta w-w)\,\partial_z
+
\big(2z\zeta+2z\big)\,\partial_\zeta
+
2zw\,\partial_w,
\\
&
\ \ \ \ \ \ \ \ \ \ \ \ \ \ \ \ \ \ \ \
\ \ \ \ \ 
\big(-iz^2+i\zeta w-iw\big)\,\partial_z
+
\big(-2iz\zeta+2iz\big)\,\partial_\zeta
-
2izw\,\partial_w
\big\},
\\
\mathfrak{g}_2
&
\,:=\,
\Span\,
\big\{
izw\,\partial_z
-
iz^2\,\partial_\zeta
+
iw^2\,\partial_w
\big\}.
\endaligned
\]

Calling these $X_1, \dots, X_{10}$ in order of appearance,
the five $X_\sigma + \overline{X}_\sigma$ for $\sigma = 1, 2, 3, 4, 5$
span $TM^5$ while those for $\sigma = 6, 7, 8, 9, 10$
generate the isotropy subgroup of the origin.

\Section{\bf Prenormalization}
\label{prenormalization}
\HEAD{{\ref{prenormalization}}.~{\sf Prenormalization}
}{
Zhangchi {\sc Chen}, Wei Guo {\sc Foo}, Joël {\sc Merker}, 
The Anh {\sc Ta}}

In coordinates $(z, \zeta, w) \in \C^3$ with $w = u + i\, v$, 
consider a local $\mathcal{C}^\omega$ rigid hypersurface $M^5
\subset \C^3$ graphed as $u = F(z, \zeta, \overline{z}, 
\overline{\zeta})$ passing through the origin. 
Expand $\sum_{a+b+c+d\geqslant 1}\,
F_{a,b,c,d}\, z^a \zeta^b \overline{z}^c \overline{\zeta}^d$,
and define by conjugating only coefficients:
\[
\overline{F}\big(z,\zeta,\overline{z},\overline{\zeta}\big)
\,:=\,
\sum_{a+b+c+d\geqslant 1}\,
\overline{F}_{a,b,c,d}\,
z^a\zeta^b\overline{z}^c\overline{\zeta}^d.
\]
The reality $\overline{u} = u$ forces 
$\overline{F (z,\zeta, \overline{z}, \overline{\zeta})} = 
F(z,\zeta, \overline{z}, \overline{\zeta})$ which becomes:
\[
\overline{F}
\big(
\overline{z},\overline{\zeta},z,\zeta
\big)
\,\equiv\,
F\big(
z,\zeta,\overline{z},\overline{\zeta}
\big).
\]
The $4$ independent derivations $\partial_z$, $\partial_\zeta$,
$\partial_{\overline{z}}$, $\partial_{\overline{\zeta}}$
commute. Applying $\frac{1}{a!}\partial_z^a\,
\frac{1}{b!}\partial_\zeta^b \frac{1}{c!}\partial_{\overline{z}}^c
\frac{1}{d!}\partial_{\overline{\zeta}}^d$ at the origin
$(0,0,0,0)$, it comes:
\[
\overline{F}_{c,d,a,b}
\,=\,
F_{a,b,c,d}.
\]
With $\chi(z, \zeta) := F(z,\zeta,0,0)$ which is holomorphic,
setting $w' := w - 2\, \chi(z,\zeta)$, we get:
\[
{\textstyle{\frac{w'+\overline{w}'}{2}}}
\,=\,
u'
\,=\,
F\big(z,\zeta,\overline{z},\overline{\zeta}\big)
-
\chi(z,\zeta)
-
\overline{\chi}
\big(
\overline{z},\overline{\zeta}
\big)
\,\,=:\,\,
F'\big(z,\zeta,\overline{z},\overline{\zeta}\big),
\]
with now $0 \equiv F'(z,\zeta,0,0) \equiv F'(0,0,\overline{z},
\overline{\zeta})$.

By ${\rm O}_x(3)$, we mean a
(remainder) function equal to $x^3(\cdots)$,
where $(\cdots)$ is any function of one or several variables.
By ${\rm O}_{x,y}(2)$, we mean $x^2 (\cdots) + xy (\cdots) + 
y^2(\cdots)$, and so on.

\begin{Proposition}
\label{Prp-prenormalization}
After a rigid biholomorphism, an $M \in \mathfrak{C}_{2,1}$ satisfies:
\[
F\big(z,\zeta,\overline{z},0\big)
\,=\,
z\overline{z}
+
{\textstyle{\frac{1}{2}}}\,
\zeta\overline{z}^2
+
{\rm O}_{\overline{z}}(3).
\] 
\end{Proposition}

Employing the letter $\mathcal{R}$ for unspecified functions, 
this amounts to:
\leqnomode\usetagform{default}
\begin{align}
\label{F-remainders-R-R}
F
\big(z,\zeta,\overline{z},\overline{\zeta}\big)
\,\,=\,\,
z\overline{z}
+
{\textstyle{\frac{1}{2}}}\,
\zeta\overline{z}^2
+
\overline{z}^3\,
\mathcal{R}
\big(z,\zeta,\overline{z}\big)
+
\overline{\zeta}\,
\mathcal{R}
\big(z,\zeta,\overline{z},\overline{\zeta}\big).
\end{align}
We will use without mention: 
\[
\mathcal{R}
\big(z,\zeta,\overline{z},\overline{\zeta}\big)
\,\,=\,\,
\mathcal{R}
\big(z,\zeta,\overline{z}\big)
+
\overline{\zeta}\,
\mathcal{R}
\big(z,\zeta,\overline{z},\overline{\zeta}\big).
\]

\proof
We will perform rigid biholomorphisms of the form
$z' = z'(z,\zeta)$, $\zeta' = \zeta'(z,\zeta)$, $w' = w$
fixing $0$. They transform $u = F(z,\zeta, \overline{z},
\overline{\zeta})$ into $u' = F' (z', \zeta', \overline{z}',
\overline{\zeta}')$ with:
\[
F'
\big(
z',\zeta',\overline{z}',\overline{\zeta}'
\big)
\,:=\,
F
\big(
z(z',\zeta'),\,
\zeta(z',\zeta'),\,\,
\overline{z}(\overline{z}',\overline{\zeta}'),\,
\overline{\zeta}(\overline{z}',\overline{\zeta}')
\big),
\]
hence they conserves $F'(z', \zeta', 0, 0) \equiv 0$.

The Levi form being of rank $1$ at $0$, we may assume:
\[
u
\,=\,
z\overline{z}
+
{\rm O}_3
\big(z,\zeta,\overline{z},\overline{\zeta}\big).
\]

\begin{Assertion}
\label{Ass-z-bar-power-1}
After a rigid biholomorphism fixing $0$:
\[
F
\,=\,
z\overline{z}
+
\overline{z}^2\,
\mathcal{R}
+
\overline{\zeta}\,
\mathcal{R}.
\]
\end{Assertion}

\proof
We can decompose:
\[
F\big(z,\zeta,\overline{z},\overline{\zeta}\big)
\,\,=\,\,
F\big(z,\zeta,\overline{z},0\big)
+
\overline{\zeta}\,\mathcal{R}
\,\,=\,\,
\overline{z}
\big(z+\chi(z,\zeta)\big)
+
\overline{z}^2\,
\mathcal{R}
+
\overline{\zeta}\,\mathcal{R},
\]
with $\chi = {\rm O}(2)$. Then:
\[
F
\,=\,
\big(z+\chi\big)\,
\big(\overline{z}+\overline{\chi}\big)
-
z\,\overline{\chi}
-
\chi\,\overline{\chi}
+
\overline{z}^2\,\mathcal{R}
+
\overline{\zeta}\,\mathcal{R}.
\]
But $\overline{\chi} = \overline{z}^2\, \mathcal{R}(\overline{z})
+ \overline{\zeta}\, \mathcal{R}
(\overline{z}, \overline{\zeta})$ is absorbable, hence:
\[
F
\,=\,
\big(z+\chi\big)\,
\big(\overline{z}+\overline{\chi}\big)
+
\overline{z}^2\,\mathcal{R}
+
\overline{\zeta}\,\mathcal{R}.
\]

Thus, we perform the rigid biholomorphism $z' := z + \chi(z,\zeta)$,
$\zeta' := \zeta$, with inverse:
\[
z
\,=\,
z'
+
{\rm O}_{z',\zeta'}(2)
\,\,=\,\,
z'
+
{z'}^2\,\mathcal{R}'
+
\zeta'\,\mathcal{R}'.
\]
Hence $\overline{z}^2 = {\overline{z}'}^2 \mathcal{R}' 
+ \overline{\zeta}' \mathcal{R}'$, and lastly:
\[
F'\big(z',\zeta',\overline{z}',\overline{\zeta}'\big)
\,=\,
z'\overline{z}'
+
{\overline{z}'}^2\,\mathcal{R}'
+
\overline{\zeta}'\,\mathcal{R}'.
\qedhere
\]
\endproof

Next, dropping primes, specifying 3\textsuperscript{rd}
order (real) terms $P = P_3$ in $F = z \overline{z} + 
P_3 + {\rm O}_{z, \zeta, \overline{z}, \overline{\zeta}}(4)$,
let us inspect the Levi determinant:
\[
0
\,\equiv\,
\left\vert\!
\begin{array}{cc}
1+P_{z\overline{z}}+{\rm O}_2
&
P_{\zeta\overline{z}}+{\rm O}_2
\\
P_{z\overline{\zeta}}+{\rm O}_2
&
P_{\zeta\overline{\zeta}}
+
{\rm O}_2
\end{array}
\!\right\vert,
\ \ \ \ \ \ \ \ \ \ \ \ 
\text{whence}
\ \ \ \ \ \ 
0
\,\equiv\,
P_{\zeta\overline{\zeta}},
\]
{\em i.e.} $P$ is harmonic with respect to $\zeta$ when $z$, 
$\overline{z}$ are seen as constants. Thus taking account of
$0 \equiv P(z,\zeta, 0, 0)$:
\[
P
\,=\,
a\,z^2\overline{z}
+
\overline{a}\,z\overline{z}^2
+
\zeta\,
\big(
b\,z\overline{z}
+
c\,\overline{z}^2
\big)
+
\overline{\zeta}\,
\big(
\overline{b}\,z\overline{z}
+
\overline{c}\,z^2
\big)
+
\zeta^2\,
\big(
d\,\overline{z}
\big)
+
\overline{\zeta}^2\,
\big(
\overline{d}\,z
\big).
\]
But Assertion~{\ref{Ass-z-bar-power-1}} forces 
$a = 0$, $b = 0$, $d = 0$, whence:
\[
u
\,=\,
z\overline{z}
+
c\,\zeta\,\overline{z}^2
+
\overline{c}\,\overline{\zeta}z^2
+
{\rm O}_{z,\zeta,\overline{z},\overline{\zeta}}(4).
\]

From Proposition~{\ref{Prp-2ndg-determinant}},
we know that $c \neq 0$, hence $c\, \zeta =:
\frac{1}{2}\, \zeta'$ conducts to:
\leqnomode\usetagform{default}
\begin{align}
\label{3rd-order-normalization}
u
\,=\,
z\overline{z}
+
{\textstyle{\frac{1}{2}}}\,
z^2\overline{\zeta}
+
{\textstyle{\frac{1}{2}}}\,
\overline{z}^2\zeta
+
{\rm O}_{z,\zeta,\overline{z},\overline{\zeta}}(4)
\,\,=\,\,
z\overline{z}
+
\overline{z}^2\,\mathcal{R}
+
\overline{\zeta}\,\mathcal{R}.
\end{align}

Next, let us look at 4\textsuperscript{th} 
order terms which depend only on $(z,
\overline{z})$, especially at the monomial $e\, z^2 \overline{z}^2$
with $e := F_{2,0,2,0} \in \R$. We can make $e = 0$ thanks to $\zeta'
:= \zeta + e\, z^2$:
\[
u
\,=\,
z\overline{z}
+
{\textstyle{\frac{1}{2}}}\,
\big(
\zeta
+
e\,z^2
\big)\,
\overline{z}^2
+
{\textstyle{\frac{1}{2}}}\,
\big(
\overline{\zeta}
+
e\,\overline{z}^2
\big)\,
z^2
+
\overline{z}^2\mathcal{R}
+
\overline{\zeta}\mathcal{R}.
\]
So we can assume $F_{2,0,2,0} = 0$. We then write:
\[
u
\,=\,
z\overline{z}
+
{\textstyle{\frac{1}{2}}}\,
\overline{z}^2\,
S\big(z,\zeta,\overline{z}\big)
+
\overline{\zeta}\,
\mathcal{R}\big(z,\zeta,\overline{z},\overline{\zeta}\big),
\]
with $S = \zeta + 
{\rm O}_{z,\zeta,\overline{z}}(2)$ and with {\em no} $z^2$ monomial
in the remainder. Hence with some
function $\tau(z)$ which {\em is} an ${\rm O}_z(3)$, 
and with some function $\omega(z, \zeta) = {\rm O}_{z, \zeta}(1)$,
we devise which biholomorphism to perform:
\[
\aligned
u
&
\,=\,
z\overline{z}
+
{\textstyle{\frac{1}{2}}}\,
\overline{z}^2\,
\big(
\zeta
+
\tau(z)
+
\zeta\,\omega(z,\zeta)
+
\overline{z}\,
\theta(z,\zeta,\overline{z})
\big)
+
\overline{\zeta}\,\mathcal{R}
\\
&
\,=\,
z\overline{z}
+
{\textstyle{\frac{1}{2}}}\,
\overline{z}^2\,
\big(
\underbrace{
\zeta
+
\tau(z)
+
\zeta\,\omega(z,\zeta)}_{=:\,\,\zeta',\,\,\,\text{\rm while}\,\,
z\,=:\,z'}
\big)
+
\overline{z}^3\,
\mathcal{R}
+
\overline{\zeta}\,
\mathcal{R}.
\endaligned
\]

\begin{Assertion}
The inverse $\zeta = \zeta' + {\rm O}(2) = \tau'(z') + 
\zeta' \big[ 1 + \omega'(z',\zeta') \big]$ also satisfies $\tau'(z') = 
{\rm O}_{z'}(3)$.
\end{Assertion}

\proof
Indeed, by definition:
\[
\zeta
\,\equiv\,
\tau'(z)
+
\big[
\tau(z)
+
\zeta\,\big(1+\omega(z,\zeta)\big)
\big]\,
\big[
1
+
\omega'
\big(
z,\,
\tau(z)+\zeta\,\big(1+\omega(z,\zeta)\big)
\big)
\big],
\]
and it suffices to put 
$\zeta := 0$ to get a concluding relation which even
shows that $\ord_0 \tau = \ord_0 \tau'$:
\[
0
\,\equiv\,
\tau'(z)
+
\tau(z)\,
\big[
1
+
\omega'\big(z,\tau(z)\big)
\big].
\qedhere
\]
\endproof

All this enables to reach the goal~({\ref{F-remainders-R-R}})
since $\overline{\tau}'(\overline{z}')$ is absorbable
in ${\overline{z}'}^3 \mathcal{R}'$:
\[
u
\,=\,
z'\overline{z}'
+
{\textstyle{\frac{1}{2}}}\,
{\overline{z}'}^2\,
\zeta'
+
{\overline{z}'}^3\,
\mathcal{R}'
+
\Big(
\overline{\zeta}'
+
\overline{\tau}'(\overline{z}')
+
\overline{\zeta}'\,
\overline{\omega}'
(\overline{z}',\overline{\zeta}')
\Big)\,
\mathcal{R}'.
\qedhere
\]
\endproof

Coordinates like in Proposition~{\ref{Prp-prenormalization}}
will be called {\sl prenormalized}. Equivalently
(exercise):
\reqnomode\usetagform{EngelLie}
\begin{align}
0
&
\,=\,
F_{a,b,0,0}
\,=\,
F_{0,0,c,d},
\notag
\\
0
&
\,=\,
F_{a,b,1,0}
\,=\,
F_{1,0,c,d},
\notag
\\
0
&
\,=\,
F_{a,b,2,0}
\,=\,
F_{2,0,c,d},
\notag
\end{align}
with only three exceptions 
$F_{1,0,1,0} = 1$ and $F_{2,0,0,1} = \frac{1}{2} = F_{0,1,2,0}$.
During the proof, in~({\ref{3rd-order-normalization}}), 
we obtained simultaneously:
\leqnomode\usetagform{default}
\begin{align}
\label{simultaneous-prenormalizations-F}
\boxed{\,\,
u
\,=\,
F
\,=\,
z\overline{z}
+
{\textstyle{\frac{1}{2}}}\,
\overline{z}^2\zeta
+
{\rm O}_{\overline{z}}(3)
+
{\rm O}_{\overline{\zeta}}(1)
\,\,=\,\,
z\overline{z}
+
{\textstyle{\frac{1}{2}}}\,
z^2\overline{\zeta}
+
{\textstyle{\frac{1}{2}}}\,
\overline{z}^2\zeta
+
{\rm O}_{z,\zeta,\overline{z},\overline{\zeta}}(4).\,\,}
\end{align}

Now, recall that the Gaussier-Merker model
is homogeneous of degree $2$ in $z$, 
$\overline{z}$, when $\zeta$, $\overline{\zeta}$ are
treated as constants:
\[
u
\,=\,
\frac{z\overline{z}+\frac{1}{2}\,z^2\overline{\zeta}
+\frac{1}{2}\overline{z}^2\zeta}{1-\zeta\overline{\zeta}}
\,\,=:\,\,
\maux\big(z,\zeta,\overline{z},\overline{\zeta}\big).
\]
A general $M \in \mathfrak{C}_{2,1}$ is just a perturbation of it:
\[
u
\,=\,
F
\,=\,
\maux
+
G,
\ \ \ \ \ \ \ \ \ \ \ \ \ \ \ \ \ \ \ \
\text{with}
\ \ \ \ \
G
\,:=\,
F-\maux
\,=\,
{\rm O}_{z,\zeta,\overline{z},\overline{\zeta}}(4).
\]

\begin{Proposition}
\label{Prp-G-O-z-3}
In prenormalized coordinates, one has
$G = {\rm O}_{z, \overline{z}}(3)$.
\end{Proposition}

\proof
Expand:
\[
\aligned
\maux
&
\,=\,
z\overline{z}\,
\sum_{i\geqslant 0}\,
\zeta^i\overline{\zeta}^i
+
{\textstyle{\frac{1}{2}}}\,
z^2\,
\sum_{i\geqslant 0}\,
\zeta^i\overline{\zeta}^{i+1}
+
{\textstyle{\frac{1}{2}}}\,
\overline{z}^2\,
\sum_{i\geqslant 0}\,
\zeta^{i+1}\overline{\zeta}^i
\,\,=\,\,
z\overline{z}
+
{\textstyle{\frac{1}{2}}}\,
z^2\overline{\zeta}
+
{\textstyle{\frac{1}{2}}}\,
\overline{z}^2\zeta
+
{\rm O}_{z,\zeta,\overline{z},\overline{\zeta}}(4),
\\
G
&
\,=\,
\sum_{k\geqslant 4}\,
\sum_{a+b+c+d=k}\,
G_{a,b,c,d}\,
z^a\zeta^b\overline{z}^c\overline{\zeta}^d
\,\,=:\,\,
\sum_{k\geqslant 4}\,
G^k.
\endaligned
\]
Of course, $F^k = \maux^k + G^k$, with $G^2 = G^3 = 0$.

\begin{Assertion}
For every $k \geqslant 2$, one has $G^k = {\rm O}_{z,
\overline{z}}(3)$.
\end{Assertion}

\proof
For some $k \geqslant 4$, assume by induction
that $G^2, G^3, \dots, G^{k-1}$
are ${\rm O}_{z,\overline{z}}(3)$, whence:
\[
G_{z\overline{z}}^\ell
\,=\,
{\rm O}_{z,\overline{z}}(1),
\ \ \ \ \ \ \ \ \ \ \ \ \ \
G_{\zeta\overline{z}}^\ell
\,=\,
{\rm O}_{z,\overline{z}}(2)
\,=\,
G_{z\overline{\zeta}}^\ell,
\ \ \ \ \ \ \ \ \ \ \ \ \ \
G_{\zeta\overline{\zeta}}^\ell
\,=\,
{\rm O}_{z,\overline{z}}(3)
\eqno
{\scriptstyle{(1\,\leqslant\,\ell\,\leqslant\,k-1)}}.
\]

Next, insert $F = \sum_{i\geqslant 2}\, F^i$ in the Levi determinant:
\[
0
\,\equiv\,
\left\vert\!
\begin{array}{cc}
\smallsum{i}\,F_{z\overline{z}}^i
&
\smallsum{j}\,F_{\zeta\overline{z}}^j
\\
\smallsum{i}\,F_{z\overline{\zeta}}^i
&
\smallsum{j}\,F_{\zeta\overline{\zeta}}^j
\end{array}
\!\right\vert
\,\,=\,\,
\sum_{\ell\geqslant 4}\,
\bigg(
\sum_{i+j=\ell
\atop
i,j\geqslant 2}\,
\Big(
F_{z\overline{z}}^i\,
F_{\zeta\overline{\zeta}}^j
-
F_{z\overline{\zeta}}^i\,
F_{\zeta\overline{z}}^j
\Big)
\bigg).
\]
Behind $\sum_\ell$, all terms are of constant homogeneous order
$i - 2 + j-2 = \ell - 4$, hence $0 \equiv \sum_{i+j=\ell}$ for 
each $\ell \geqslant 4$. Take $\ell := k+2$ and expand:
\[
\aligned
0
&
\,\equiv\,
F_{z\overline{z}}^2\,
\boxed{F_{\zeta\overline{\zeta}}^k}
+
\sum_{3\leqslant i\leqslant k-1}\,
F_{z\overline{z}}^i\,
F_{\zeta\overline{\zeta}}^{k+2-i}
+
F_{z\overline{z}}^k\,
\zero{F_{\zeta\overline{\zeta}}^2}
\,-
\\
&
\ \ \ \ \
-\,
\zero{F_{z\overline{\zeta}}^2}\,
F_{\zeta\overline{z}}^k
-
\sum_{3\leqslant i\leqslant k-1}\,
F_{z\overline{\zeta}}^i\,
F_{\zeta\overline{z}}^{k+2-i}
-
F_{z\overline{\zeta}}^k\,
\zero{F_{\zeta\overline{z}}^2}.
\endaligned
\]
Observe from~({\ref{simultaneous-prenormalizations-F}}) that
$1 \equiv F_{z\overline{z}}^2$ while $0 \equiv F_{\zeta
\overline{\zeta}}^2 \equiv F_{z\overline{\zeta}}^2 \equiv
F_{\zeta\overline{z}}^2$.
Of course, Levi determinant vanishing holds for $F := \maux$:
\[
\aligned
0
&
\,\equiv\,
\maux_{z\overline{z}}^2\,
\maux_{\zeta\overline{\zeta}}^k
+
\sum_{3\leqslant i\leqslant k-1}\,
\maux_{z\overline{z}}^i\,
\maux_{\zeta\overline{\zeta}}^{k+2-i}
+
\maux_{z\overline{z}}^k\,
\zero{\maux_{\zeta\overline{\zeta}}^2}
\,-
\\
&
\ \ \ \ \
-\,
\zero{\maux_{z\overline{\zeta}}^2}\,
\maux_{\zeta\overline{z}}^k
-
\sum_{3\leqslant i\leqslant k-1}\,
\maux_{z\overline{\zeta}}^i\,
\maux_{\zeta\overline{z}}^{k+2-i}
-
\maux_{z\overline{\zeta}}^k\,
\zero{\maux_{\zeta\overline{z}}^2}.
\endaligned
\]

Substituting the boxed term $F_{\zeta \overline{\zeta}}^k$ 
with $\maux_{\zeta\overline{\zeta}}^k + G_{\zeta 
\overline{\zeta}}^k$, solving for $G_{\zeta \overline{\zeta}}^k$,
substituting as well the other
$F_{\bdot\bdot}^\ell = \maux_{\bdot\bdot}^\ell 
+ G_{\bdot \bdot}^\ell$, and subtracting, we obtain:
\[
\aligned
-\,G_{\zeta\overline{\zeta}}^k
&
\,\equiv\,
\sum_{3\leqslant i\leqslant k-1}\,
\bigg(
\maux_{z\overline{z}}^i\,
G_{\zeta\overline{\zeta}}^{k+2-i}
+
G_{z\overline{z}}^i\,
\maux_{\zeta\overline{\zeta}}^{k+2-i}
+
G_{z\overline{z}}^i\,
G_{\zeta\overline{\zeta}}^{k+2-i}
\bigg)
\,-
\\
&
\ \ \ \ \ 
-\,
\sum_{3\leqslant i\leqslant k-1}\,
\bigg(
\maux_{z\overline{\zeta}}^i\,
G_{\zeta\overline{z}}^{k+2-i}
+
G_{z\overline{\zeta}}^i\,
\maux_{\zeta\overline{z}}^{k+2-i}
+
G_{z\overline{\zeta}}^i\,
G_{\zeta\overline{z}}^{k+2-i}
\bigg).
\endaligned
\]
Since we also have $3 \leqslant k+2-i \leqslant k-1$,
induction applies to all six products to get
$G_{\zeta\overline{\zeta}}^k = {\rm O}_{z ,\overline{z}}(3)$.

By integration, $G^k = \lambda^k (z, \zeta, \overline{z}) 
+ \overline{\lambda}^k (\overline{z}, \overline{\zeta}, z) 
+ {\rm O}_{z,\overline{z}}(3)$. After absorption
in ${\rm O}_{z, \overline{z}}(3)$, we can 
assume that $\lambda^k$ is of degree $\leqslant 2$ in
$(z, \overline{z})$, hence contains only monomials 
$z^a \zeta^b \overline{z}^c$ with $a + c \leqslant 2$
and $a + b + c = k$. So $b \geqslant
k-2$.

Further, $G^k(z, \zeta, 0, 0) \equiv 0$ 
imposes $\lambda^k (z, \zeta, 0) \equiv 0$.
So $1 \leqslant c \leqslant 2$. 
Consequently, $\lambda^k$ can contain only three monomials:
\[
\lambda^k(z,\zeta,\overline{z})
\,=\,
a\,\overline{z}\zeta^{k-1}
+
b\,z\overline{z}\,\zeta^{k-2}
+
c\,\overline{z}^2\zeta^{k-2}.
\]
Since $k \geqslant 4$, we see that the conjugate
$\overline{\lambda}^k(\overline{z}, \overline{\zeta}, z)$
is multiple of $\overline{\zeta}^{k-2 \geqslant 2}$, hence:
\[
G^k\big(z,\zeta,\overline{z},0\big)
\,=\,
\lambda^k(z,\zeta,\overline{z})
+
\zero{
\overline{\lambda}^k(\overline{z},0,z)}
+
{\rm O}_{z,\overline{z}}(3).
\]

Finally, because the prenormalized
coordinates of Proposition~{\ref{Prp-prenormalization}} 
require $G^k(z,\zeta,\overline{z},0) = {\rm O}_{\overline{z}}(3)$,
we reach $\lambda^k(z,\zeta,\overline{z}) = {\rm O}_{z,
\overline{z}}(3)$, which forces $a = b = c = 0 = \lambda^k$, 
so as asserted $G^k = {\rm O}_{z,\overline{z}}(3)$. 
\endproof

In conclusion, $G = \sum\, G^k = {\rm O}_{z,\overline{z}}(3)$.
\endproof

According to~{\cite{Foo-Merker-Ta-2019}}, 
the Lie group $G$ of rigid CR automorphisms 
of the Gaussier-Merker model $\{ u = \maux \}$
has Lie algebra
$\mathfrak{g}_{-2} \oplus \mathfrak{g}_{-1} \oplus \mathfrak{g}_0$
of dimension $7$, generated by $X_1, \dots, X_7$.
The $2$-dimensional
isotropy subgroup $G_0 \subset G$ of the origin $0 \in \C^3$
has Lie algebra $\mathfrak{g}_0^{\sf iso}$ generated by:
\[
X_6
\,:=\,
z\,\partial_z
+
2w\,\partial_w,
\ \ \ \ \ \ \ \ \ \ \ \ \ \ \ \ \ \ \ \
X_7
\,:=\,
iz\,\partial_z
+
2i\zeta\,\partial_\zeta.
\]
By computing the flows $\exp\big( t\, X_\sigma \big) (z, \zeta, w)$
for $t \in \R$ and $\sigma = 6, 7$, one verifies that
$G_0$ consists of scalings coupled with `rotations':
\[
z'
\,=\,
\rho^{1/2}\,
e^{i\varphi}\,
z,
\ \ \ \ \ \ \ \ \ \ \ \ \ \ \ \ \ \ \ \
\zeta'
\,=\,
e^{2i\varphi}\,
\zeta,
\ \ \ \ \ \ \ \ \ \ \ \ \ \ \ \ \ \ \ \
w'
\,=\,
\rho\,w
\eqno
{\scriptstyle{(\rho\,\in\,\R_+^\ast,\,\,\varphi\,\in\,\R)}}.
\]

Next, any holomorphic function $e = e(z,w)$ decomposes 
in weighted homogeneous terms as:
\[
e(z,w)
\,=\,
\sum_{a,b}\,
e_{a,b}\,z^a\zeta^b
\,\,=\,\,
\sum_{k\geqslant 0}\,
\bigg(
\sum_b\,
e_{k,b}\,\zeta^b
\bigg)\,
z^k
\,\,=:\,\,
\sum_{k\geqslant0}\,
e_k.
\]
Mind notation: for weights, indices $e_k$ are lower case, 
while for orders, as {\em e.g.} in $G^k$ before, they were
upper case. Similarly:
\[
E\big(z,\zeta,\overline{z},\overline{\zeta}\big)
\,=\,
\sum_{k\geqslant0}\,
\bigg(
\sum_{a+c=k}\,
\Big(
\sum_{b,d}\,
E_{a,b,c,d}\,
\zeta^b\overline{\zeta}^d
\Big)\,
z^a\overline{z}^c
\bigg)
\,\,=:\,\,
\sum_{k\geqslant0}\,
E_k.
\]

According to what precedes, we can assume that both the
source $M$ and the target $M'$ rigid hypersurfaces are prenormalized.
Assume therefore that a rigid biholomorphism:
\[
H
\colon\ \ \
(z,\zeta,w)
\,\,\,\longmapsto\,\,\,
\Big(
f(z,\zeta),\,
g(z,\zeta),\,
\rho\,w+h(z,\zeta)
\Big)
\,\,=:\,\,
\big(z',\zeta',w'),
\]
fixing the origin is given between:
\[
\aligned
u
&
\,=\,
F
\,=\,
z\overline{z}
+
{\textstyle{\frac{1}{2}}}\,
\overline{z}^2\zeta
+
{\rm O}_{\overline{z}}(3)
\,=\,
\maux+G
\,=\,
{\textstyle{
\frac{z\overline{z}+\frac{1}{2}\,z^2\overline{\zeta}
+\frac{1}{2}\,\overline{z}^2\zeta}{1-\zeta\overline{\zeta}}}}
+
{\rm O}_{z,\overline{z}}(3),
\\
u'
&
\,=\,
F'
\,=\,
z'\overline{z}'
+
{\textstyle{\frac{1}{2}}}\,
{\overline{z}'}^2\zeta'
+
{\rm O}_{\overline{z}'}(3)
\,=\,
\maux'+G'
\,=\,
{\textstyle{
\frac{z'\overline{z}'+\frac{1}{2}\,{z'}^2\overline{\zeta}'
+\frac{1}{2}\,{\overline{z}'}^2\zeta'}{1-\zeta'\overline{\zeta}'}}}
+
{\rm O}_{z',\overline{z}'}(3).
\endaligned
\]

\begin{Observation}
Scalings and rotations $(z', \zeta', w') \longmapsto 
\big( \rho^{1/2}e^{i\varphi}z',\, e^{2i\varphi}\zeta',\, 
\rho w'\big)$ preserve prenormalizations.\qed
\end{Observation}

Since $T_0^cM = \{w=0\}$ and $T_0^cM' = \{w'=0\}$, and since
$H_\ast T_0^cM = T_0^cM'$, we necessarily have 
$h = {\rm O}_{z,\zeta} (2)$. After the scaling $w' \longmapsto
\frac{1}{\rho}\, w'$, we may therefore 
assume that the last component of $H$
is $w + {\rm O}_{z,\zeta}(2)$.

Let us decompose the components of $H$
in weighted homogeneous parts:
\[
f
\,=\,
f_0+f_1+f_2+f_3+\cdots,
\ \ \ \ \ \ \ \ \ \ \ \ \ \ \ \ \ \ \ \
g
\,=\,
g_0+g_1+g_2+\cdots,
\ \ \ \ \ \ \ \ \ \ \ \ \ \ \ \ \ \ \ \
h
\,=\,
h_0+h_1+h_2+h_3+h_4+\cdots.
\]

Plug in the components of $H$ in the target rigid equation
$\frac{w'+\overline{w}'}{2} = F'(z', \zeta', \overline{z}', 
\overline{\zeta}')$:
\[
w
+
h(z,\zeta)
+
\overline{w}
+
\overline{h}\big(\overline{z},\overline{\zeta}\big)
\,\,=\,\,
2\,
F'
\Big(
f(z,\zeta),\,
g(z,\zeta),\,
\overline{f}\big(\overline{z},\overline{\zeta}\big),\,
\overline{g}\big(\overline{z},\overline{\zeta}\big)
\Big),
\]
and then, substitute $w+\overline{w} = 2\, F$ to get a
{\sl fundamental equation}, holding identically:
\leqnomode\usetagform{default}
\begin{align}
\label{fundamental-equation-before-weighting}
2\,F(z,\zeta,\overline{z},\overline{\zeta})
+
h(z,\zeta)
+
\overline{h}\big(\overline{z},\overline{\zeta}\big)
\,\,\equiv\,\,
2\,F'
\Big(
f(z,\zeta),\,
g(z,\zeta),\,
\overline{f}\big(\overline{z},\overline{\zeta}\big),\,
\overline{g}\big(\overline{z},\overline{\zeta}\big)
\Big).
\end{align}

\begin{Proposition}
\label{Prp-initial-f-g-h}
Possibly after a rotation $(z', \zeta', w') \longmapsto 
(e^{i\varphi}z', e^{2i\varphi}\zeta', w')$, one has:
\[
f
\,=\,
z
+
f_2
+
f_3
+\cdots,
\ \ \ \ \ \ \ \ \ \ \ \ \ \ \ \ \ \ \ \
g
\,=\,
\zeta
+
g_1
+
g_2
+\cdots,
\ \ \ \ \ \ \ \ \ \ \ \ \ \ \ \ \ \ \ \
h
\,=\,
w
+
h_3
+
h_4
+\cdots.
\]
or equivalently: $f_0 = 0$, $f_1 = z$; $g_0 = \zeta$; $h_0 = 0$, 
$h_1 = 0$, $h_2 = w$.
\end{Proposition}

\proof
Recall that $F = \maux + G$, that $\maux = \maux_2$ and that $G = G_3
+ G_4 + \cdots$, with the same about $F' = \maux' + G'$. 
So $F$ and $F'$ have no terms of weights $0$ or $1$.
Of course
$f_0 = f_0(\zeta)$, $g_0 = g_0(\zeta)$, $h_0 = h_0(\zeta)$ depend on
$\zeta$ only.

In~({\ref{fundamental-equation-before-weighting}}), pick terms 
of weight zero:
\[
0
+
h_0(\zeta)
+
\overline{h}_0(\overline{\zeta})
\,\equiv\,
2\,
F'
\big(
f_0(\zeta),\,
g_0(\zeta),\,
\overline{f}_0(\overline{\zeta}),
\overline{g}_0(\overline{\zeta})
\big),
\]
put $\overline{\zeta} := 0$, use $F'(z',\zeta',0,0) \equiv 0$, and
get $h_0 = 0$.

Once again, pick 
in~({\ref{fundamental-equation-before-weighting}})
terms of weight zero using $F' = \maux' + 
{\rm O}_{z',\overline{z}'}(3)$:
\[
0
\,\equiv\,
\frac{f_0(\zeta)\overline{f}_0(\overline{\zeta})
+\frac{1}{2}f_0(\zeta)^2\overline{g}_0(\overline{\zeta})
+\frac{1}{2}\overline{f}_0(\overline{\zeta})g_0(\zeta)}{
1-g_0(\zeta)\overline{g}_0(\overline{\zeta})}
+
{\rm O}_{f_0(\zeta),\overline{f}_0(\overline{\zeta})}(3).
\]
We claim that $f_0(\zeta) \equiv 0$. Otherwise, $f_0 = c\, \zeta^\nu +
{\rm O}_\zeta(\nu+1)$ with $c \neq 0$, but on the right, the monomial
$c\overline{c}\, \zeta^\nu \overline{\zeta}^\nu$ cannot be
killed\,\,---\,\,contradiction. This finishes examination of weight
zero, for it remains only $0 \equiv 0$.

Hence, pass to weight $1$. We claim that $h_1 = 0$. Of course, $f_1 =
z f_1(\zeta)$ and $h_1 = z h_1(\zeta)$. Since
$\maux'$ is weighted homogeneous of degree $2$,
we have $F' = {\rm O}_{z',
\overline{z}'}(2)$, and
we get
from~({\ref{fundamental-equation-before-weighting}})
what forces $h_1 = 0$:
\[
{\rm O}_{z,\overline{z}}(2)
+
z\,h_1(\zeta)
+
\overline{z}\,\overline{h}_1(\overline{\zeta})
\,\equiv\,
{\rm O}_{zf_1(\zeta),\overline{z}\overline{f}_1(\zeta)}(2)
\,\equiv\,
{\rm O}_{z,\overline{z}}(2).
\]

Before passing to weight $2$,
since $f = z f_1(\zeta) + {\rm O}_z(2)$ and $g = g_0(\zeta) + 
z g_1(\zeta) + {\rm O}_z(2)$, the nonzero Jacobian
$\big\vert \begin{smallmatrix}
f_z & f_\zeta \\ g_z & g_\zeta
\end{smallmatrix} \big\vert$ has value at the origin
$\big\vert \begin{smallmatrix}
f_1(0) & 0 \\ g_1(0) & g_0'(0)
\end{smallmatrix} \big\vert$, hence $f_1(0) \neq 0 \neq g_0'(0)$.

Lastly, picking weighted degree $2$ terms 
in~({\ref{fundamental-equation-before-weighting}}), we get:
\[
2\,\maux(z,\zeta,\overline{z},\overline{\zeta})
+
z^2h_2(\zeta)
+
\overline{z}^2\overline{h}_2(\overline{\zeta})
\,\,\equiv\,\,
2\,\maux
\Big(
zf_1(\zeta),\,
g_0(\zeta),\,
\overline{z}\overline{f}_1(\overline{\zeta}),\,
\overline{g}_0(\overline{\zeta})
\Big).
\]
This identity means that the map
$(z, \zeta, w) \longmapsto \big( zf_1(\zeta),\, g_0(\zeta),\,
w+z^2h_2(\zeta) \big)$ is an automorphism of the Gaussier-Merker
model fixing the origin, hence is a rotation, so that
$f_1(\zeta) = e^{i\varphi}$, $g_0(\zeta) = e^{2i\varphi} \zeta$,
$h_2(z,\zeta) \equiv 0$. Post-composing with the 
inverse rotation, we attain the conclusion.
\endproof

\begin{Question}
\label{Q-finite-dimensionality}
{\sl Suppose given two rigid hypersurfaces prenormalized as before:}
\[
\aligned
u
&
\,=\,
F
\,=\,
z\overline{z}
+
{\textstyle{\frac{1}{2}}}\,
\overline{z}^2\zeta
+
{\rm O}_{\overline{z}}(3)
+
{\rm O}_{\overline{\zeta}}(1)
\,=\,
\maux
+
G
\,=\,
{\textstyle{
\frac{z\overline{z}+\frac{1}{2}\,z^2\overline{\zeta}
+\frac{1}{2}\,\overline{z}^2\zeta}{1-\zeta\overline{\zeta}}}}
+
{\rm O}_{z,\overline{z}}(3),
\\
u'
&
\,=\,
F'
\,=\,
z'\overline{z}'
+
{\textstyle{\frac{1}{2}}}\,
{\overline{z}'}^2\zeta'
+
{\rm O}_{\overline{z}'}(3)
+
{\rm O}_{\overline{\zeta}'}(1)
\,=\,
\maux'
+
G'
\,=\,
{\textstyle{
\frac{z'\overline{z}'+\frac{1}{2}\,{z'}^2\overline{\zeta}'
+\frac{1}{2}\,{\overline{z}'}^2\zeta'}{1-\zeta'\overline{\zeta}'}}}
+
{\rm O}_{z',\overline{z}'}(3).
\endaligned
\]
{\sl Is it true that the group of rigid biholomorphisms
at the origin between them:}
\[
(z,\zeta,w)
\,\,\,\longmapsto\,\,\,
\Big(
z
+
f(z,\zeta),\,\,
\zeta
+
g(z,\zeta),\,\,
w
+
h(z,\zeta)
\Big)
\,\,=:\,\,
\big(z',\zeta',w'\big),
\]
where $f = f_2 + f_3 + \cdots$, $g = g_1 + g_2 + \cdots$,
$h = h_3 + h_4 + \cdots$,
{\sl is finite-dimensional\,{\bf ?}}
\end{Question}

Here, the two appearing remainders
${\rm O}_{z,\overline{z}}(3)$ and
${\rm O}_{\overline{z}}(3) + {\rm O}_{\overline{\zeta}}(1)$
are different. By expanding $1 \big/
(1-\zeta\overline{\zeta})$ we see that:
\[
\maux
\,=\,
z\overline{z}
+
{\textstyle{\frac{1}{2}}}\,
\overline{z}^2\zeta
+
{\textstyle{\frac{1}{2}}}\,
z^2\overline{\zeta}
+
\zeta\overline{\zeta}
\big(\cdots\big)
\,\,=\,\,
z\overline{z}
+
{\textstyle{\frac{1}{2}}}\,
\overline{z}^2\zeta
+
{\rm O}_{\overline{\zeta}}(1),
\]
hence by subtraction, we get that $G$ is more than just
an ${\rm O}_{z,\overline{z}}(3)$.

\begin{Observation}
\label{Obs-G-before-synthesis}
The remainder function satisfies $G = {\rm O}_{z, \overline{z}}
(3) = {\rm O}_{\overline{z}}(3) +
{\rm O}_{\overline{\zeta}}(1)$.\qed
\end{Observation}

The synthesis between these two conditions will be
made in Section~{\ref{normal-form}}.

\Section{\bf Weighted Homogeneous Normalizing Biholomorphisms}
\label{weighted-homogeneous-normalizing-biholomorphisms}
\HEAD{{\ref{weighted-homogeneous-normalizing-biholomorphisms}}.~{\sf
Weighted Homogeneous Normalizing Biholomorphisms}
}{
Zhangchi {\sc Chen}, Wei Guo {\sc Foo}, Joël {\sc Merker}, 
The Anh {\sc Ta}}

Now, inspired by Jacobowitz's presentation~{\cite{Jacobowitz-1990}} of
Moser's normal form in $\C^2$, Propositions~{\ref{Prp-G-O-z-3}}
and~{\ref{Prp-initial-f-g-h}} justify to introduce the spaces:
\[
\footnotesize
\aligned
\mathcal{G}
&
\,:=\,
\Big\{
G=G(z,\zeta,\overline{z},\overline{\zeta})
\colon\,
G
=
G_3+G_4+\cdots
\Big\},
\\
\mathcal{D}
&
\,:=\,
\Big\{
\big(z+f(z,\zeta),\,\zeta+g(z,\zeta),\,w+h(z,\zeta)\big)
\colon\,
f=f_2+f_3+\cdots,\,\,
g=g_1+g_2+\cdots,\,\,
h=h_3+h_4+\cdots
\Big\},
\endaligned
\]
where lower indices denote homogeneous components with respect
to the weighting~({\ref{weighting-z-zeta-w}}) defined by:
\[
\big[
z^a \zeta^b \overline{z}^c \overline{\zeta}^d
\big]
\,=\,
a+c.
\]
The goal is to use the `{\sl freedom}' space $\mathcal{D}$ of rigid
biholomorphisms in order to `{\sl normalize}' as much as possible the
remainder $G$ in the
graphed equation $\{ u = \maux + G\big\}$ of any given hypersurface.
Here, $\maux = \frac{z\overline{z} + 
\frac{1}{2} z^2\overline{\zeta} + \frac{1}{2} \overline{z}^2
\zeta}{1-\zeta\overline{\zeta}}$ is homogeneous of weight $2$.

Both $\mathcal{G}$ and $\mathcal{D}$ decompose as direct sums
graded by increasing weights:
\[
\aligned
\mathcal{G}
&
\,=\,
\underset{\nu\geqslant 3}{\medcup}\,
\mathcal{G}_\nu,
&
\ \ \ \ \ \ \ \ \ \ \ \ \ \ \ \ \ \ \ \
\mathcal{G}_\nu
&
\,:=\,
\big\{
G_\nu
\big\},
\\
\mathcal{D}
&
\,=\,
\underset{\nu\geqslant 3}{\medcup}\,
\mathcal{D}_\nu,
&
\ \ \ \ \ \ \ \ \ \ \ \ \ \ \ \ \ \ \ \
\mathcal{D}_\nu
&
\,:=\,
\big\{
\big(
f_{\nu-1},\,
g_{\nu-2},\,
h_\nu
\big)
\big\},
\endaligned
\]
and the (upcoming) justification for the shifts in $\mathcal{D}_\nu$
will be due to two multipliers:
\[
\maux_z
\,=\,
{\textstyle{\frac{\overline{z}+z\overline{\zeta}}{
1-\zeta\overline{\zeta}}}}
\ \ \ \ \
\text{\rm of weight}\,\,
1
\ \ \ \ \ \ \ \ \ \ \ \ \ \ \ \ \ \ \ \
\text{\rm and}
\ \ \ \ \ \ \ \ \ \ \ \ \ \ \ \ \ \ \ \
\maux_\zeta
\,=\,
{\textstyle{\frac{(\overline{z}+z\overline{\zeta})^2}{
2\,(1-\zeta\overline{\zeta})^2}}}
\ \ \ \ \
\text{\rm of weight}\,\,
2.
\]
One can figure out that $G_2 := \maux$
and $G_2' := \maux'$ are already finalized\big/normalized. 
With increasing weights $\nu = 3, 4, 5, \dots$, we shall
perform successive holomorphic rigid transformations
of the shape:
\leqnomode\usetagform{default}
\begin{align}
\label{biholomorphism-nu-tangent-identity}
z'
\,:=\,
z+f_{\nu-1},
\ \ \ \ \ \ \ \ \ \ \ \ \ \ \ \ \ \ \ \
\zeta'
\,:=\,
\zeta+g_{\nu-2},
\ \ \ \ \ \ \ \ \ \ \ \ \ \ \ \ \ \ \ \
w'
\,:=\,
w+h_\nu.
\end{align}

When $\nu \gg 1$ is high, it is intuitively clear that 
such transformations close to the identity
will preserve previously achieved
low order normalizations; to make this claim precise,
let us follow and adapt~{\cite[Chap.~3]{Jacobowitz-1990}}.

For $\mu \geqslant 0$, 
denote by ${\rm O}(\mu)$ power series 
whose monomials $z^a \zeta^b \overline{z}^c \overline{\zeta}^d$
are all of weight $a + c \geqslant \mu$, and introduce the
projection operators:
\[
\pi_\mu
\Big(
\smallsum{a,b,c,d\geqslant 0}\,
T_{a,b,c,d}\,
z^a\zeta^b\overline{z}^c\overline{\zeta}^d
\Big)
\,:=\,
\sum_{a+c\leqslant\mu}\,
\sum_{b,d\geqslant 0}\,
T_{a,b,c,d}\,
z^a\zeta^b\overline{z}^c\overline{\zeta}^d.
\]

\begin{Proposition}
\label{Prp-5-2}
Through any 
biholomorphism~({\ref{biholomorphism-nu-tangent-identity}})
which transforms:
\[
u
\,=\,
\maux
+
G_3+\cdots+G_{\nu-1}
+
G_\nu
+
{\rm O}(\nu+1)
\ \ \ \ \ \ \ 
\text{into}
\ \ \ \ \ \ \
u'
\,=\,
\maux
+
G_3'+\cdots+G_{\nu-1}'
+
G_\nu'
+
{\rm O}'(\nu+1),
\]
homogeneous terms are kept untouched up to order $\leqslant \nu-1$:
\[
G_\mu'\big(z,\zeta,\overline{z},\overline{\zeta}\big)
\,=\,
G_\mu\big(z,\zeta,\overline{z},\overline{\zeta}\big)
\eqno
{\scriptstyle{(3\,\leqslant\,\mu\,\leqslant\,\nu-1)}},
\]
while:
\[
G_\nu'\big(z,\zeta,\overline{z},\overline{\zeta}\big)
\,=\,
G_\nu\big(z,\zeta,\overline{z},\overline{\zeta}\big)
-
2\,\Re\,
\Big\{
{\textstyle{
\frac{\overline{z}+z\overline{\zeta}}{1-\zeta\overline{\zeta}}}}\,
f_{\nu-1}(z,\zeta)
+
{\textstyle{
\frac{(\overline{z}+z\overline{\zeta})^2}{
2(1-\zeta\overline{\zeta})^2}}}\,
g_{\nu-2}(z,\zeta)
-
{\textstyle{\frac{1}{2}}}\,
h_\nu(z,\zeta)
\Big\}.
\]
\end{Proposition}

Thus, by appropriately choosing $(f_{\nu-1}, g_{\nu-2}, h_\nu)$,
we will be able to `kill' many monomials in $G_\nu$, hence
make $G_\nu'$ simpler, or {\sl normalized}. Exercise: verify that 
in fact
$h_\nu \equiv 0$ necessarily,
when $F$ and $F'$ are assumed to be prenormalized.

\proof
As already seen, the fundamental equation,
holding identically, is:
\[
\Re\,
\big(
w+h_\nu
\big)
\,=\,
F(z,\zeta,\overline{z},\overline{\zeta})
+
\Re\,h_\nu
\,\equiv\,
F'\big(
z
+
f_{\nu-1}(z,\zeta),\,
\zeta
+
g_{\nu-2}(z,\zeta),\,
w
+
h_\nu(z,\zeta)
\big).
\]
Decomposing $F = \maux + G$, $F' = \maux' + G'$ and reorganizing, 
it becomes:
\[
\footnotesize
\aligned
\frac{(z+f_{\nu-1})(\overline{z}+\overline{f}_{\nu-1})
+\frac{1}{2}(z+f_{\nu-1})^2(\overline{\zeta}+\overline{g}_{\nu-2})
+\frac{1}{2}(\overline{z}+\overline{f}_{\nu-1})^2(\zeta+g_{\nu-2})
}{
1-(\zeta+g_{\nu-2})(\overline{\zeta}+\overline{g}_{\nu-2})}
-
\frac{z\overline{z}+\frac{1}{2}z^2\overline{\zeta}
+\frac{1}{2}\overline{z}^2\zeta}{1-\zeta\overline{\zeta}}
-
\Re\,h_\nu
\,\,=\,\,
G
-
G'.
\endaligned
\]

A reduction of the left hand side to the same denominator
shows after algebraic simplifications: 
\[
\!\!\!\!\!\!\!\!\!\!\!\!\!\!\!
\footnotesize
\aligned
\frac{(1-\zeta\overline{\zeta})
\big[
z\overline{f}_{\nu-1}
+
\overline{z}f_{\nu-1}
+
\frac{1}{2}
\big(
2zf_{\nu-1}\overline{\zeta}
+
z^2\overline{g}_{\nu-2}
\big)
+
\frac{1}{2}
\big(
2\overline{z}\overline{f}_{\nu-1}\zeta
+
\overline{z}^2g_{\nu-2}
\big)
\big]
+
\big(
\zeta\overline{g}_{\nu-2}
+
\overline{\zeta}g_{\nu-2}
\big)
\big(
z\overline{z}
+
\frac{1}{2}z^2\overline{\zeta}
+
\frac{1}{2}\overline{z}^2\zeta
\big)}{
(1-\zeta\overline{\zeta})\,
\big(
1-\zeta\overline{\zeta}
-
\zeta\overline{g}_{\nu-2}
-
\overline{\zeta}g_{\nu-2}
-
g_{\nu-2}\overline{g}_{\nu-2}
\big)}
-
\Re\,h_\nu.
\endaligned
\]
that this left-hand side 
is ${\rm O}(\nu)$, hence has zero $\pi_{\nu-1} (
\centersmallbullet) = 0$. Moreover, its homogeneous degree $\nu$
part is obtained by taking only weighted degree zero terms
in the denominator, namely
$\frac{\sf numerator}{(1-\zeta\overline{\zeta})^2} - \Re\, h_\nu$, 
and one recognizes\big/reconstitutes
$\maux_z$, $\maux_\zeta$ as 
homogeneous multipliers of weights $1$, $2$:
\[
\pi_\nu
\Big(
\maux'
-
\maux
-
\Re\,h_\nu
\Big)
\,=\,
2\,\Re\,
\Big\{
{\textstyle{
\frac{\overline{z}+z\overline{\zeta}}{1-\zeta\overline{\zeta}}}}\,
f_{\nu-1}(z,\zeta)
+
{\textstyle{
\frac{(\overline{z}+z\overline{\zeta})^2}{
2(1-\zeta\overline{\zeta})^2}}}\,
g_{\nu-2}(z,\zeta)
-
{\textstyle{\frac{1}{2}}}\,
h_\nu(z,\zeta)
\Big\}.
\]

It remains to treat $\pi_\nu (\centersmallbullet)$
of the right-hand side:
\[
\smallsum{3\leqslant\mu\leqslant\nu}\,
G_\mu(z,\zeta,\overline{z},\overline{\zeta})
-
\pi_\nu
\Big(
\smallsum{3\leqslant\mu\leqslant\nu}\,
G_\mu'
\big(
z+f_{\nu-1},\,
\zeta+g_{\nu-2},\,
\overline{z}+\overline{f}_{\nu-1},\,
\overline{\zeta}+\overline{g}_{\nu-2}
\big)
\Big).
\]

\begin{Assertion}
For each $3 \leqslant \mu \leqslant \nu$:
\[
\pi_\nu
\Big(
G_\mu'
\big(
z+f_{\nu-1},\,
\zeta+g_{\nu-2},\,
\overline{z}+\overline{f}_{\nu-1},\,
\overline{\zeta}+\overline{g}_{\nu-2}
\big)
\Big)
\,\,=\,\,
G_\mu'
\big(z,\zeta,\overline{z},\overline{\zeta}\big).
\]
\end{Assertion}

\proof
All possible monomials in $G_\mu'$ with $a + c = \mu \geqslant 3$
after binomial expansion:
\[
\!\!\!\!\!\!\!\!\!\!\!\!\!\!\!
\scriptsize
\aligned
\big(z+f_{\nu-1}\big)^a
\big(\zeta+g_{\nu-2}\big)^b
\big(\overline{z}+\overline{f}_{\nu-1}\big)^c
\big(\overline{\zeta}+\overline{g}_{\nu-2}\big)^d
&
\,\,=\,\,
\big(
z^a
+
{\rm O}(a-1+\nu-1)
\big)
\big(
\zeta^b
+
{\rm O}(\nu-2)
\big)
\big(
\overline{z}^c
+
{\rm O}(c-1+\nu-1)
\big)
\big(
\overline{\zeta}^d
+
{\rm O}(\nu-2)
\big)
\notag
\\
&
\,=\,
z^a\zeta^b\overline{z}^c\overline{\zeta}^d
+
{\rm O}(a+c-2+\nu),
\endaligned
\]
have the simple projection $\pi_\nu (\centersmallbullet) = z^a \zeta^b
\overline{z}^c \overline{\zeta}^d$ since $a + c - 2 + \nu \geqslant 1
+ \nu$.
\endproof

We therefore obtain an identity in which all
arguments are $(z, \zeta, \overline{z}, \overline{\zeta})$:
\[
2\,\Re\,
\Big\{
{\textstyle{
\frac{\overline{z}+z\overline{\zeta}}{1-\zeta\overline{\zeta}}}}\,
f_{\nu-1}
+
{\textstyle{
\frac{(\overline{z}+z\overline{\zeta})^2}{
2(1-\zeta\overline{\zeta})^2}}}\,
g_{\nu-2}
-
{\textstyle{\frac{1}{2}}}\,
h_\nu
\Big\}
\,\,\equiv\,\,
\smallsum{3\leqslant\mu\leqslant\nu-1}\,
\Big(
\zero{
G_\mu
-
G_\mu'}
\Big)
+
G_\nu
-
G_\nu'.
\]
Applying $\pi_{\nu-1}$ annihilates both the left-hand side
and $G_\nu - G_\nu'$, whence $G_\mu = G_\mu'$ for
$3 \leqslant \mu \leqslant \nu-1$, which concludes.
\endproof

\Section{\bf Normal Form}
\label{normal-form}
\HEAD{{\ref{normal-form}}.~{\sf Normal Form}
}{
Zhangchi {\sc Chen}, Wei Guo {\sc Foo}, Joël {\sc Merker}, 
The Anh {\sc Ta}}

The assumption that the Levi form is of constant rank $1$:
\[
F_{z\overline{z}}
\,\neq\,
0
\,\equiv\,
F_{z\overline{z}}\,
F_{\zeta\overline{\zeta}}
-
F_{\zeta\overline{z}}\,
F_{z\overline{\zeta}},
\]
enables to solve identically as functions
of $(z,\zeta,\overline{z},\overline{\zeta})$:
\[
F_{\zeta\overline{\zeta}}
\,\equiv\,
\frac{
F_{\zeta\overline{z}}\,
F_{z\overline{\zeta}}
}{
F_{z\overline{z}}}.
\]
By successively differentiating this identity and performing
replacements, we get formulas.

\begin{Lemma}
\label{Lm-determination-other-coefficients}
For every jet multiindex
$(a, b, c, d) \in \N^4$ with $b \geqslant 1$ and
$d \geqslant 1$, abbreviating $n := a + b + c + d$,
there exists a polynomomial
$P_{a,b,c,d}$ in its arguments and an
integer $\NN_{a,b,c,d} \geqslant 1$ such that:
\reqnomode\usetagform{EngelLie}
\begin{align}
\footnotesize
\aligned
F_{z^a\zeta^b\overline{z}^c\overline{\zeta}^d}
\,\equiv\,
\frac{1}{\big(F_{z\overline{z}}\big)^{\NN_{a,b,c,d}}}\,\,
P_{a,b,c,d}
\Big(
\big\{
F_{z^{a'}\overline{z}^{c'}}
\big\}_{a'+c'\leqslant n},\,\,
\big\{
F_{z^{a'}\zeta^{b'}\overline{z}^{c'}}
\big\}_{a'+b'+c'\leqslant n}^{b'\geqslant 1},\,\,
\big\{
F_{z^{a'}\overline{z}^{c'}\overline{\zeta}^{d'}}
\big\}_{a'+c'+d'\leqslant n}^{d'\geqslant 1}
\Big).
\endaligned
\tag{\qed}
\end{align}
\end{Lemma}

In other words, the Levi rank $1$ assumption implies that
all Taylor coefficients at the origin of $\sum_{a,b,c,d}\,
F_{a,b,c,d}\, z^a \zeta^b \overline{z}^c \overline{\zeta}^d$
for which $b \geqslant 1$ and $d \geqslant 1$ are
determined by the free Taylor coefficients:
\[
\big\{
F_{a,0,c,0}
\big\}_{a\geqslant 0,\,c\geqslant 0}
\,\medcup\,
\big\{
F_{a,b,c,0}
\big\}_{a\geqslant 0,\,b\geqslant 1,\,c\geqslant 0}
\,\medcup\,
\big\{
F_{a,0,c,d}
\big\}_{a\geqslant 0,\,c\geqslant 0,\,d\geqslant 1}.
\]

In subsequent computations, we will therefore normalize only these
free (independent) Taylor coefficients at the origin, while those
(dependent) attached to monomials that are multiple of $\zeta
\overline{\zeta}$ will then be automatically determined by the
formulas of Lemma~{\ref{Lm-determination-other-coefficients}}.

As promised, we can now explore
Observation~{\ref{Obs-G-before-synthesis}} further.  What precedes
shows that it is best appropriate to expand $G$ with respect to
$(\zeta, \overline{\zeta})$:
\[
\aligned
G
\,=\,
\sum_{a,c\geqslant 0}\,
G_{a,0,c,0}\,
z^a\overline{z}^c
+
\sum_{b\geqslant 1}\,
\zeta^b\,
\Big(
\sum_{a,c\geqslant 0}\,
G_{a,b,c,0}\,
z^a\overline{z}^c
\Big)
&
+
\sum_{d\geqslant 1}\,
\overline{\zeta}^d\,
\Big(
\sum_{a,c\geqslant 0}\,
G_{a,0,c,d}\,
z^a\overline{z}^c
\Big)
\\
&
+
\sum_{b,d\geqslant 1}\,
\sum_{a,c\geqslant 0}\,
G_{a,b,c,d}\,
z^a\zeta^b\overline{z}^c\overline{\zeta}^d.
\endaligned
\]
The last quadruple sum gathers all dependent jets. We will
abbreviate this remainder as $\zeta \overline{\zeta} ( \cdots )$.
With different notations, we can therefore write:
\[
G
\,=\,
a(z,\overline{z})
+
\sum_{k\geqslant 0}\,
\zeta^{k+1}\,\Pi_k(z,\overline{z})
+
\sum_{k\geqslant 0}\,
\overline{\zeta}^{k+1}\,
\overline{\Pi}_k(\overline{z},z)
+
\zeta\overline{\zeta}
\big(\cdots\big),
\]
with $a (z, \overline{z}) \equiv \overline{a}
(\overline{z}, z)$ real, but no reality constraint
on the $\Pi_k(z, \overline{z})$.

Recall that $G = {\rm O}_{z, \overline{z}}(3)$.
In view of Proposition~{\ref{Prp-5-2}},
we must, for every weight $\nu \geqslant 3$, extract $G_\nu$,
while writing $\zeta^{k+1} = \zeta\, \zeta^k$:
\[
\aligned
G_\nu
&
\,=\,
a_{\nu,0}\,z^\nu
+
a_{\nu-1,1}\,z^{\nu-1}\overline{z}
+\cdots+
a_{1,\nu-1}\,z\overline{z}^{\nu-1}
+
a_{0,\nu}\,\overline{z}^\nu
\,+
\\
&
\ \ \ \ \
+
\sum_{k\geqslant 0}\,
\zeta\,
\zeta^k\,
\Big(
z^\nu\,
\Pi_{k,\nu,0}
+
z^{\nu-1}\overline{z}\,
\Pi_{k,\nu-1,1}
+\cdots+
z\overline{z}^{\nu-1}\,
\Pi_{k,1,\nu-1}
+
\overline{z}^\nu\,
\Pi_{k,0,\nu}
\Big)
\,+
\\
&
\ \ \ \ \
+
\sum_{k\geqslant 0}\,
\overline{\zeta}\,
\overline{\zeta}^k\,
\Big(
\overline{z}^\nu\,
\overline{\Pi}_{k,\nu,0}
+
\overline{z}^{\nu-1}z\,
\overline{\Pi}_{k,\nu-1,1}
+\cdots+
\overline{z}z^{\nu-1}\,
\overline{\Pi}_{k,1,\nu-1}
+
z^\nu\,\overline{\Pi}_{k,0,\nu}
\Big)
\,+
\\
&
\ \ \ \ \
+
\zeta\overline{\zeta}\,
\big(\cdots\big).
\endaligned
\]

To reorganize all this in powers of $(z, \overline{z})$, let us
introduce the two collections for all $0 \leqslant \mu \leqslant \nu$
of (anti)holomorphic functions (mind the inversion $\nu - \mu
\longleftrightarrow \mu$ at the end):
\[
B_{\nu-\mu,\mu}(\zeta)
\,:=\,
\smallsum{k\geqslant 0}\,
\zeta^k\,
\Pi_{k,\nu-\mu,\mu}
\ \ \ \ \ \ \ \ \ \ \ \ \ \ \ \ \ \ \ \
\text{and}
\ \ \ \ \ \ \ \ \ \ \ \ \ \ \ \ \ \ \ \
\overline{C}_{\nu-\mu,\mu}(\overline{\zeta})
\,:=\,
\sum_{k\geqslant 0}\,
\overline{\zeta}^k\,
\overline{\Pi}_{k,\mu,\nu-\mu}.
\]
The definition of
these $B_{\centersmallbullet, \centersmallbullet}$ and
$\overline{C}_{\centersmallbullet, \centersmallbullet}$
enables us to
emphasize that the obtained functions
$\zeta\, B_{\centersmallbullet, \centersmallbullet} (\zeta)$ and
$\overline{\zeta}\,
\overline{C}_{\centersmallbullet, \centersmallbullet}
(\overline{\zeta})$
vanish
when either $\zeta := 0$ or $\overline{\zeta} := 0$,
and we therefore obtain, 
taking also account of the fact that $G_\nu$ is real:
\[
\aligned
G_\nu
&
\,=\,
z^\nu\,
\Big(
a_{\nu,0}
+
\zeta\,B_{\nu,0}(\zeta)
+
\overline{\zeta}\,\overline{C}_{\nu,0}(\overline{\zeta})
\Big)
+
z^{\nu-1}\overline{z}\,
\Big(
a_{\nu-1,1}
+
\zeta\,B_{\nu-1,1}(\zeta)
+
\overline{\zeta}\,
\overline{C}_{\nu-1,1}(\overline{\zeta})
\Big)
+
\\
&
\ \ \ \ \
+
\cdots\cdots\cdots\cdots\cdots\cdots\cdots\cdots\cdots
\cdots\cdots\cdots\cdots\cdots\cdots\cdots\cdots\cdots
\cdots\cdots\cdots\cdots\cdots\cdots\cdots\cdots\cdot
+
\\
&
\ \ \ \ \
+
z\overline{z}^{\nu-1}\,
\Big(
\overline{a}_{\nu-1,1}
+
\overline{\zeta}\,
\overline{B}_{\nu-1,1}(\overline{\zeta})
+
\zeta\,C_{\nu-1,1}(\zeta)
\Big)
+
\overline{z}^\nu\,
\Big(
\overline{a}_{\nu,0}
+
\overline{\zeta}\,
\overline{B}_{\nu,0}(\overline{\zeta})
+
\zeta\,C_{\nu,0}(\zeta)
\Big)
+
\zeta\overline{\zeta}\,
\big(\cdots\big).
\endaligned
\]

Of course, all these weighted homogeneous functions $G_\nu$
automatically satisfy $G_\nu = {\rm O}_{z, \overline{z}}(3)$,
since $\nu \geqslant 3$ thanks to
Proposition~{\ref{Prp-G-O-z-3}}.
Now, Observation~{\ref{Obs-G-before-synthesis}}
also requires that they satisfy, since they are real:
\leqnomode\usetagform{default}
\begin{align}
\label{G-nu-O-3-O-1}
G_\nu
\,=\,
{\rm O}_{\overline{z}}(3)
+
{\rm O}_{\overline{\zeta}}(1)
\,=\,
{\rm O}_{z}(3)
+
{\rm O}_{\zeta}(1).
\end{align}

\begin{Lemma}
\label{Lm-synthesis-prenormalization}
For each weight $\nu \geqslant 5$, the function $G_\nu$
satisfies~({\ref{G-nu-O-3-O-1}}) if and only if it is of the form:
\[
\aligned
G_\nu
&
\,=\,
\ \ \ \ \ \ \ \ \
z^\nu\,
\Big(
\ \ \ \
0
\ \ \ \
+
\ \ \ \ \ \ \ \ 0 \ \ \ \ \ \ \ \
+
\ \ \
\overline{\zeta}\,
\overline{C}_{\nu,0}(\overline{\zeta})
\ \ \
\Big)
\\
&
\ \ \ \ \ \
+
z^{\nu-1}\overline{z}\,
\Big(
\ \ \ \
0
\ \ \ \
+
\ \ \ \ \ \ \ \ 0 \ \ \ \ \ \ \ \
+
\overline{\zeta}\,
\overline{C}_{\nu-1,1}(\overline{\zeta})
\Big)
\\
&
\ \ \ \ \ 
+
z^{\nu-2}\overline{z}^2\,
\Big(
\ \ \ \
0
\ \ \ \
+
\ \ \ \ \ \ \ \ 0 \ \ \ \ \ \ \ \
+
\overline{\zeta}\,
\overline{C}_{\nu-2,2}(\overline{\zeta})
\Big)
\\
&
\ \ \ \ \ 
+
z^{\nu-3}\overline{z}^3\,
\Big(
a_{\nu-3,3}
+
\zeta\,B_{\nu-3,3}(\zeta)
+
\overline{\zeta}\,
\overline{C}_{\nu-3,3}(\overline{\zeta})
\Big)
\\
&
\ \ \ \ \
+
\cdots\cdots\cdots\cdots\cdots\cdots\cdots\cdots\cdots\cdots\cdots
\cdots\cdots
+
\\
&
\ \ \ \ \ 
+
z^3\overline{z}^{\nu-3}\,
\Big(
\overline{a}_{\nu-3,3}
+
\zeta\,C_{\nu-3,3}(\zeta)
+
\overline{\zeta}\,\overline{B}_{\nu-3,3}(\overline{\zeta})
\Big)
\\
&
\ \ \ \ \ 
+
z^2\overline{z}^{\nu-2}\,
\Big(
\ \ \ \ 
0
\ \ \ \ 
+
\zeta\,C_{\nu-2,2}(\zeta)
+
\ \ \ \ \ \ \ \ 
0
\ \ \ \ \ \ \ \ 
\Big)
\\
&
\ \ \ \ \ 
+
z^1\overline{z}^{\nu-1}\,
\Big(
\ \ \ \ 
0
\ \ \ \ 
+
\zeta\,C_{\nu-1,1}(\zeta)
+
\ \ \ \ \ \ \ \ 
0
\ \ \ \ \ \ \ \ 
\Big)
\\
&
\ \ \ \ \ \ \ \ \ \ \ \ 
+
\overline{z}^\nu\,
\Big(
\ \ \ \ 
0
\ \ \ \ 
+
\ \
\zeta\,C_{\nu,0}(\zeta)
\ \
+
\ \ \ \ \ \ \ \ 
0
\ \ \ \ \ \ \ \ 
\Big)
+
\zeta\overline{\zeta}\,
\big(
\cdots
\big).
\endaligned
\]
\end{Lemma}

Just after,
we will treat the two weights $\nu = 3, 4$ separately.

\proof
Putting $\overline{\zeta} := 0$ above, it must hold that:
\[
\aligned
{\rm O}_{\overline{z}}(3)
+
0
\,=\,
G_\nu\big\vert_{\overline{\zeta}=0}
&
\,=\,
\ \ \ \ \ \ \ \ 
z^\nu\,
\Big(
a_{\nu,0}
+
\zeta\,B_{\nu,0}(\zeta)
+
0
\Big)
\\
&
\ \ \ \ \
+
z^{\nu-1}\overline{z}\,
\Big(
a_{\nu-1,1}
+
\zeta\,B_{\nu-1,1}(\zeta)
+
0
\Big)
+
\\
&
\ \ \ \ \
+
z^{\nu-2}\overline{z}^2\,
\Big(
a_{\nu-2,2}
+
\zeta\,B_{\nu-2,2}(\zeta)
+
0
\Big)
\\
&
\ \ \ \ \
+
{\rm O}_{\overline{z}}(3)
+
0.
\endaligned
\]
Thus, all the appearing $a_{\smallbullet, \smallbullet}$
and $B_{\smallbullet, \smallbullet}$ should vanish,
as stated, and the converse is clear.
\endproof

Proceeding similarly, the reader will find for $\nu = 3$
that $G_3$ satisfies~({\ref{G-nu-O-3-O-1}}) if and only if:
\[
\aligned
G_3
&
\,=\,
\ \ \ \ \ 
z^3\,
\Big(
0
+
0
+
\overline{\zeta}\,
\overline{C}_{3,0}(\overline{\zeta})
\Big)
\\
&
\ \ \ \ \
+
z^2\overline{z}\,
\big(
0+0+0
\big)
\\
&
\ \ \ \ \
+
z\overline{z}^2\,
\big(
0+0+0
\big)
\\
&
\ \ \ \ \
+
\overline{z}^3\,
\Big(
0
+
\zeta\,
C_{3,0}(\zeta)
+
0
\Big)
+
\zeta\overline{\zeta}\,
\big(\cdots\big),
\endaligned
\]
as well as:
\[
\aligned
G_4
&
\,=\,
\ \ \ \ \ 
z^4\,
\Big(
0
+
0
+
\overline{\zeta}\,
\overline{C}_{4,0}(\overline{\zeta})
\Big)
\\
&
\ \ \ \ \
+
z^3\overline{z}\,
\Big(
0
+
0
+
\overline{\zeta}\,
\overline{C}_{3,1}(\overline{\zeta})
\Big)
\\
&
\ \ \ \ \
+
z^2\overline{z}^2\,
\big(
0+0+0
\big)
\\
&
\ \ \ \ \
+
z\overline{z}^3\,
\Big(
0
+
\zeta\,
C_{1,3}(\zeta)
+
0
\Big)
\\
&
\ \ \ \ \
+
\overline{z}^4\,
\Big(
0
+
\zeta\,
C_{4,0}(\zeta)
+
0
\Big)
+
\zeta\overline{\zeta}\,
\big(\cdots\big).
\endaligned
\]

Now, consider a rigid biholomorphism $z' = f(z,\zeta)$, $\zeta' =
g(z,\zeta)$, $w' = \rho\, w + h(z,\zeta)$ between two rigid hypersurfaces
$M$ and $M'$.  Of course, as in
Question~{\ref{Q-finite-dimensionality}},
we may assume that both $M$
and $M'$ have already been prenormalized,
and thanks to 
Proposition~{\ref{Prp-initial-f-g-h}}
also that $f = f_2 + f_3 + \cdots$, $g = g_1 + g_2 + \cdots$,
$\rho = 1$, $h = h_3 + h_4 + \cdots$.

The goal is to {\sl
normalize} $M'$ even further, by means of appropriate choices of
$f$, $g$, $h$.

We saw that it is natural to decompose $G = G_3 + G_4 + G_5 +
\cdots$ and $G' = G_3' + G_4' + G_5' + \cdots$ in weighted homogeneous
parts, and we just finished to express what prenormalization means
about these $G_\nu$ and $G_\nu'$.  Proceeding with increasing weights
$\nu = 3, 4, 5, \dots$, we therefore consider biholomorphisms of the
shape $z' = z + f_{\nu-1}$, $\zeta' = \zeta + g_{\nu-2}$, $w' = w +
h_\nu$, and we recall that Proposition~{\ref{Prp-5-2}} showed that:
\[
G_\nu'\big(z,\zeta,\overline{z},\overline{\zeta}\big)
\,=\,
G_\nu\big(z,\zeta,\overline{z},\overline{\zeta}\big)
-
2\,\Re\,
\Big\{
{\textstyle{
\frac{\overline{z}+z\overline{\zeta}}{1-\zeta\overline{\zeta}}}}\,
f_{\nu-1}(z,\zeta)
+
{\textstyle{
\frac{(\overline{z}+z\overline{\zeta})^2}{
2(1-\zeta\overline{\zeta})^2}}}\,
g_{\nu-2}(z,\zeta)
-
{\textstyle{\frac{1}{2}}}\,
h_\nu(z,\zeta)
\Big\}.
\]
The freedom to `{\sl normalize}' $G_\nu'$ even more that
$G_\nu$, namely the term $-2\, \Re\, \{ \cdots \}$,
is parametrized by the complely free choice for the triple
of holomorphic functions $(f_{\nu-1}, g_{\nu-2}, h_\nu)$.
However, prenormalizations should be left untouched.

\begin{Lemma}
\label{LM-nu-geqslant-5}
At every weight level $\nu \geqslant 5$, only the identity
biholomorphic transformation $z' = z$, $\zeta' = \zeta$,
$w' = w$ stabilizes prenormalization in source and target
spaces:
\[
G_\nu(z,\zeta,\overline{z},\overline{\zeta})
\,\,=\,\,
{\rm O}_{\overline{z}}(3)
+
{\rm O}_{\overline{\zeta}}(1)
\,\,=\,\,
G_\nu'\big(z,\zeta,\overline{z},\overline{\zeta}\big),
\]
or equivalently, the `{\sl freedom function}'
respects prenormalization:
\[
{\rm O}_{\overline{z}}(3)
+
{\rm O}_{\overline{\zeta}}(1)
\,=\,
2\,\Re\,
\Big\{
{\textstyle{
\frac{\overline{z}+z\overline{\zeta}}{1-\zeta\overline{\zeta}}}}\,
f_{\nu-1}(z,\zeta)
+
{\textstyle{
\frac{(\overline{z}+z\overline{\zeta})^2}{
2(1-\zeta\overline{\zeta})^2}}}\,
g_{\nu-2}(z,\zeta)
-
{\textstyle{\frac{1}{2}}}\,
h_\nu(z,\zeta)
\Big\}
\,\,=:\,\,
\Phi(z,\zeta,\overline{z},\overline{\zeta}),
\]
if and only if $0 = f_{\nu-1} = g_{\nu-2} = h_\nu$.
\end{Lemma}

\proof
It is easy to verify that the
vanishings $G_\nu(z,\zeta,0,0) \equiv 0 \equiv
G_\nu'(z,\zeta,0,0)$, which hold from the very beginning
(of Proposition~{\ref{Prp-prenormalization}})
already suffice
to force $h_\nu(z,\zeta) \equiv 0$.

Next, write:
\[
\aligned
f_{\nu-1}(z,\zeta)
&
\,=\,
z^{\nu-1}\,f(\zeta)
\,=\,
z^{\nu-1}\,
\big(
f_0+f_1\,\zeta+f_2\,\zeta^2
+\cdots
\big),
\\
g_{\nu-2}(z,\zeta)
&
\,=\,
z^{\nu-2}\,g(\zeta)
\,=\,
z^{\nu-2}\,
\big(
g_0+g_1\,\zeta+g_2\,\zeta^2
+\cdots
\big).
\endaligned
\]
The goal is to show that $f(\zeta) \equiv 0$ and
$g(\zeta) \equiv 0$.

Prenormalization being expressed modulo $\zeta \overline{\zeta}
(\cdots)$, when we expand the two denominators
of $\Phi$, we have by luck $\frac{1}{1-\zeta \overline{\zeta}}
\equiv 1$ and $\frac{1}{2\,(1-\zeta\overline{\zeta}^2)} \equiv
\frac{1}{2}$, and hence it suffices to require that:
\[
{\rm O}_{\overline{z}}(3)
+
{\rm O}_{\overline{\zeta}}(1)
\overset{\text{\bf ?}}{\,\,=\,\,}
2\,\Re\,
\Big\{
\big(
\overline{z}
+
z\,\overline{\zeta}
\big)\,
z^{\nu-1}\,
\smallsum{k\geqslant 0}\,
f_k\,\zeta^k
+
{\textstyle{\frac{1}{2}}}\,
\big(
\overline{z}
+
z\,\overline{\zeta}
\big)^2\,
z^{\nu-2}\,
\smallsum{k\geqslant 0}\,
g_k\,\zeta^k
\Big\}.
\]
Using $\nu \geqslant 5$ to guarantee that there is no interference
when extracting the
first three powers $z^\nu$, $z^{\nu-1} \overline{z}$,
$z^{\nu-2} \overline{z}^2$,
let us compute the three relevant terms
of the freedom function:
\[
\footnotesize
\aligned
\Phi(z,\zeta,\overline{z},\overline{\zeta})
&
\,=\,
\big(
\overline{z}+z\overline{\zeta}
\big)\,
z^{\nu-1}\,
\big(
f_0+f_1\,\zeta+f_2\,\zeta^2
+\cdots
\big)
+\big(
{\textstyle{\frac{1}{2}}}\,
\overline{z}^2
+
z\overline{z}\overline{\zeta}
+
{\textstyle{\frac{1}{2}}}\,
z^2\,\overline{\zeta}^2
\big)
z^{\nu-2}\,
\big(
g_0+g_1\,\zeta+g_2\,\zeta^2
+\cdots
\big)
\,+
\\
&
\ \ \ \ \
+
\big(
z+\overline{z}\zeta
\big)\,
\overline{z}^{\nu-1}\,
\big(
\overline{f}_0
+
\overline{f}_1\,\overline{\zeta}
+
\overline{f}_2\,\overline{\zeta}^2
+\cdots
\big)
+
\big(
{\textstyle{\frac{1}{2}}}\,
z^2
+
\overline{z}z\zeta
+
{\textstyle{\frac{1}{2}}}\,
\overline{z}^2\zeta^2
\big)\,
\overline{z}^{\nu-2}\,
\big(
\overline{g}_0
+
\overline{g}_1\,\overline{\zeta}
+
\overline{g}_2\,\overline{\zeta}^2
+\cdots
\big)
\\
&
\,=\,
z^\nu\,
\Big(
f_0\,\overline{\zeta}
+
\zero{
f_1\,\zeta\overline{\zeta}
+
f_2\,\zeta^2\overline{\zeta}
+\cdots}
+
{\textstyle{\frac{1}{2}}}\,
g_0\,
\overline{\zeta}^2
+
\zero{
{\textstyle{\frac{1}{2}}}\,
g_1\,
\zeta\overline{\zeta}^2
+
{\textstyle{\frac{1}{2}}}\,
g_2\,
\zeta^2\overline{\zeta}^2
+\cdots}
\Big)
\\
&
\ \ \ \ \
+
z^{\nu-1}\overline{z}\,
\Big(
f_0
+
f_1\,\zeta
+
f_2\,\zeta^2
+\cdots
+
g_0\,\overline{\zeta}
+
\zero{
g_1\,
\zeta\overline{\zeta}
+
g_2\,
\zeta^2\overline{\zeta}
+\cdots}
\Big)
\\
&
\ \ \ \ \
+
z^{\nu-2}\overline{z}^2\,
\Big(
{\textstyle{\frac{1}{2}}}\,
g_0
+
{\textstyle{\frac{1}{2}}}\,
g_1\,
\zeta
+
{\textstyle{\frac{1}{2}}}\,
g_2\,
\zeta^2
+\cdots
\Big)
\\
&
\ \ \ \ \
+
\overline{z}^3\,
\big(\cdots\big)
+
\zeta\overline{\zeta}\,
\big(\cdots\big).
\endaligned
\]
Since the underlined terms can be absorbed into the remainder
$\zeta \overline{\zeta} ( \cdots )$,
it remains only:
\[
\aligned
\Phi(z,\zeta,\overline{z},\overline{\zeta})
&
\,=\,
{\textstyle{\frac{1}{2}}}\,
z^\nu\,
\big(
2f_0\,
\overline{\zeta}
+
g_0\,\overline{\zeta}^2
\big)
\\
&
\ \ \ \ \
+
z^{\nu-1}\overline{z}\,
\big(
f_0+f_1\,\zeta+f_2\,\zeta^2+\cdots
+
g_0\,\overline{\zeta}
\big)
\\
&
\ \ \ \ \
+
{\textstyle{\frac{1}{2}}}\,
z^{\nu-2}\overline{z}^2\,
\big(
g_0+g_1\,\zeta+g_2\,\zeta^2+\cdots
\big)
\\
&
\ \ \ \ \
+
\overline{z}^3\,
\big(\cdots\big)
+
\zeta\overline{\zeta}\,
\big(\cdots\big).
\endaligned
\]
Putting $\overline{\zeta} := 0$, the result should be an
${\rm O}_{\overline{z}}(3)$, 
hence the first three lines should vanish,
and lines $2$ and $3$ conclude that $f(\zeta) \equiv 0 \equiv
g(\zeta)$, as aimed at.
\endproof

Next, inspect the two remaining weights $\nu = 3, 4$.
For $\nu = 3$,  again modulo $\zeta \overline{\zeta}
(\cdots)$, the freedom function is:
\[
\Phi_3
\,\equiv\,
2\,\Re\,
\Big\{
\big(
\overline{z}
+
z\overline{\zeta}
\big)\,
z^2\,
\big(
f_0+f_1\,\zeta+f_2\,\zeta^2
+\cdots
\big)
+
\big(
{\textstyle{\frac{1}{2}}}\,
\overline{z}^2
+
z\overline{z}\overline{\zeta}
+
{\textstyle{\frac{1}{2}}}\,
z^2\overline{\zeta}^2
\big)\,
z^1\,
\big(
g_0+g_1\,\zeta+g_2\,\zeta^2+\cdots
\big)
\Big\}.
\]

\begin{Assertion}
Prenormalization $\Phi_3 = {\rm O}_{\overline{z}}(3) + 
{\rm O}_{\overline{\zeta}}(1)$ is preserved if and only if:
\[
0
\,=\,
f_0
+
{\textstyle{\frac{1}{2}}}\,
\overline{g}_0,
\ \ \ \ \ \ \ \ 
0
\,=\,
f_1,
\ \ \ \ \ \ \ \ 
0
\,=\,
f_2,
\ \ \ \ \ \ \ \ 
0
\,=\,
\overline{g}_0
+
{\textstyle{\frac{1}{2}}}\,
g_1,
\ \ \ \ \ \ \ \ 
0
\,=\,
g_2,
\ \ \ \ \ \ \ \ 
\dots.
\eqno\qed
\]
\end{Assertion}

Consequently, only $1$ complex constant is free, $f_0$,
in terms of which:
\[
g_0
\,=\,
-\,2\,\overline{f}_0,
\ \ \ \ \ \ \ \ \ \ \ \ \ \ \ \ \ \ \ \ \ \ \ \ \ \
g_1
\,=\,
-\,4\,f_0.
\]
With this, how can one normalize $G_3' = G_3 - \Phi_3$ further?
Still modulo $\zeta \overline{\zeta} ( \cdots )$:
\[
\Phi_3
\,\equiv\,
z^3\,
\big(
f_0\,
\overline{\zeta}
-
\overline{f}_0\,
\overline{\zeta}^2
\big)
+
z^2\overline{z}\,
(0)
+
z\overline{z}^2\,(0)
+
\overline{z}^3\,
\big(
\overline{f}_0\,\zeta
-
f_0\,\zeta^2
\big),
\]
hence:
\[
\aligned
G_{3,0,0,1}'
&
\,=\,
G_{3,0,0,1}
-
f_0,
\\
G_{3,0,0,2}'
&
\,=\,
G_{3,0,0,2}
+
\overline{f}_0.
\endaligned
\]
It is natural to normalize the lowest jet order 
$4 = 3 + 0 + 0 + 1$ coefficient here.

\begin{Assertion}
One can normalize $G_{3,0,0,1}' := 0$ by choosing
$f_0 := G_{3,0,0,1}$.\qed
\end{Assertion}

Once this is done, 
it is easy to see that preserving\big/maintaining the normalization:
\[
G_{3,0,0,1}'
\,=\,
G_{3,0,0,1}
\,=\,
0,
\]
forces $f_0 = 0$ above.

\begin{Assertion}
In prenormalized coordinates which satisfy in addition
$G_{3,0,0,1} = 0$, the coefficient:
\[
G_{3,0,0,2}'
\,=\,
G_{3,0,0,2}
\]
is an invariant (at the origin).\qed
\end{Assertion}

In the next Section~{\ref{invariants-I-0-V-0-Q-0}}, we will show how
to deduce the expression of corresponding invariants at {\em every}
point (not only the origin) of a rigid hypersurface.

After such a normalization, we get:
\[
u
\,=\,
z\overline{z}
+
{\textstyle{\frac{1}{2}}}\,
\overline{z}^2\zeta
+
{\textstyle{\frac{1}{2}}}\,
z^2\overline{\zeta}
+
z\overline{z}\zeta\overline{\zeta}
+
a\,z^2\overline{z}^2
+
{\rm O}_{z,\zeta,\overline{z},\overline{\zeta}}(5),
\]
with, possible, a nonzero real constant $a$,
and possibly, a remainder that is {\em not} prenormalized.

\smallskip

Fortunately, we can apply the process of 
Proposition~{\ref{Prp-prenormalization}}
to prenormalize again the coordinates,
making in particular $a = 0$, 
without perturbing the normalizations obtained up to order $4$ included.

\smallskip

Lastly, treat weight $\nu = 4$. The freedom function 
modulo $\zeta \overline{\zeta} ( \cdots )$, is:
\[
\Phi_4
\,\equiv\,
2\,\Re\,
\Big\{
\big(
\overline{z}
+
z\overline{\zeta}
\big)\,
z^3\,
\big(
f_0+f_1\,\zeta+f_2\,\zeta^2
+\cdots
\big)
+
\big(
{\textstyle{\frac{1}{2}}}\,
\overline{z}^2
+
z\overline{z}\overline{\zeta}
+
{\textstyle{\frac{1}{2}}}\,
z^2\overline{\zeta}^2
\big)\,
z^2\,
\big(
g_0+g_1\,\zeta+g_2\,\zeta^2+\cdots
\big)
\Big\}.
\]

\begin{Assertion}
Prenormalization $\Phi_4 = {\rm O}_{\overline{z}}(3) + 
{\rm O}_{\overline{\zeta}}(1)$ is preserved if and only if:
\[
0
\,=\,
f_0
\,=\,
f_1
\,=\,
f_2
\,=\,
\cdots,
\ \ \ \ \ \ \ \ \ \ \ \ \ \ \ \ \ \ \ \
0
\,=\,
g_0
+
\overline{g}_0
\,=\,
g_1
\,=\,
g_2
\,=\,
\cdots.
\eqno\qed
\]
\end{Assertion}

Thus now, only $1$ {\em real} degree of freedom is left:
\[
g_0
\,=\,
i\,\tau
\eqno
{\scriptstyle{(\tau\,\in\,\R)}}.
\]
With this, how can one normalize $G_4' = G_4 - \Phi_4$ further?
Still modulo $\zeta \overline{\zeta} (\cdots)$:
\[
\Phi_4
\,\equiv\,
z^4\,
\big(
{\textstyle{\frac{i}{2}}}\,
\tau\,
\overline{\zeta}^2
\big)
+
z^3\overline{z}\,
\big(
i\,\tau\,\overline{\zeta}
\big)
+
z^2\overline{z}^2\,(0)
+
z\overline{z}^3\,
\big(
-\,i\,\tau\,\zeta
\big)
+
z^4\,
\big(
-\,
{\textstyle{\frac{i}{2}}}\,
\tau\,\zeta^2
\big),
\]
hence:
\[
\aligned
G_{4,0,0,2}'
&
\,=\,
G_{4,0,0,2}
-
{\textstyle{\frac{i}{2}}}\,
\tau,
\\
G_{3,0,1,1}'
&
\,=\,
G_{3,0,1,1}
-
i\,
\tau,
\\
G_{2,0,2,0}'
&
\,=\,
G_{2,0,0,2}.
\endaligned
\]
The third line shows an invariant.  Notice also that $G_{4,0,0,1}' =
G_{4,0,0,1}$ is an invariant.  We choose to normalize the lowest jet
order $3 + 0 + 1 + 1 = 5$ coefficient here.

\begin{Assertion}
One can normalize $\Im\, G_{3,0,1,1}' := 0$ by choosing
$\tau := \Im\, G_{3,0,1,1}$.\qed
\end{Assertion}

Once this is done, $G_{3,0,1,1}' = G_{3,0,1,1} \in \R$ is an invariant.

\smallskip

Again, we can re-apply the process of 
Proposition~{\ref{Prp-prenormalization}}
to prenormalize the coordinates without touching
the lower order normalizations.

\smallskip

We already saw in Lemma~{\ref{LM-nu-geqslant-5}} that for any weight
$\nu \geqslant 5$, no degree of freedom exists.  Since only $2 + 1 =
3$ real degrees of freedom have been encountered, namely $f_0 \in \C$
in weight $\nu = 3$ and $\Im\, g_0 \in \R$ 
in weight $\nu = 4$, we conclude
that the answer to Question~{\ref{Q-finite-dimensionality}} is
positive.

\smallskip

All this enables us to conclude the present section
by stating results which come from our analysis.

\begin{Theorem}
Every local rigid $\mathcal{C}^\omega$ graphed hypersurface
$M^5 \subset \C^3 \ni (z,\zeta, w = u + i\,v)$
passing through the origin of equation:
\[
u
\,=\,
\sum_{a+b+c+d\geqslant1}\,
F_{a,b,c,d}\,
z^a\zeta^b\overline{z}^c\overline{\zeta}^d,
\]
whose Levi form is of constant
rank $1$ and which is $2$-nondegenerate:
\[
F_{z\overline{z}}
\,\neq\,
0
\,\equiv\,
\left\vert\!
\begin{array}{cc}
F_{z\overline{z}}
&
F_{z\overline{\zeta}}
\\
F_{\zeta\overline{z}}
&
F_{\zeta\overline{\zeta}}
\end{array}
\!\right\vert
\ \ \ \ \ \ \ \ \ \ \ \ \ \ \ \ \ \ \ \
\text{and}
\ \ \ \ \ \ \ \ \ \ \ \ \ \ \ \ \ \ \ \
0
\,\neq\,
\left\vert\!
\begin{array}{cc}
F_{z\overline{z}}
&
F_{z\overline{\zeta}}
\\
F_{zz\overline{z}}
&
F_{zz\overline{\zeta}}
\end{array}
\!\right\vert,
\]
is equivalent, through a local rigid biholomorphism:
\[
(z,\zeta,w)
\,\,\,\longmapsto\,\,\,
\Big(
f(z,\zeta),\,
g(z,\zeta),\,
\rho\,w
+
h(z,\zeta)
\Big)
\,\,=:\,\,
\big(z',\zeta',w'\big)
\eqno
{\scriptstyle{(\rho\,\in\,\R^\ast)}},
\]
to a rigid $\mathcal{C}^\omega$ hypersurface ${M'}^5 \subset {\C'}^3$
which, dropping primes for target coordinates, is a perturbation of
the Gaussier-Merker model\,\,---\,\,homogeneous of order $2$ in $(z,
\overline{z})$\,\,---:
\[
u
\,=\,
\frac{z\overline{z}+\frac{1}{2}\,z^2\overline{\zeta}
+\frac{1}{2}\,\overline{z}^2\zeta}{
1-\zeta\overline{\zeta}}
+
\sum_{a,b,c,d\in\N
\atop
a+c\geqslant 3}\,
G_{a,b,c,d}\,
z^a\zeta^b\overline{z}^c\overline{\zeta}^d,
\]
with a simplified remainder $G$ which:

\smallskip\noindent{\bf (1)}\,
is normalized to be an ${\rm O}_{z, \overline{z}}(3)$;

\smallskip\noindent{\bf (2)}\,
satisfies the prenormalization conditions $G = 
{\rm O}_{\overline{z}}(3) + {\rm O}_{\overline{\zeta}}(1) = 
{\rm O}_z(3) + {\rm O}_\zeta(1)$,
or equivalently:
\[
\aligned
G_{a,b,0,0}
&
\,=\,
0
\,=\,
G_{0,0,c,d},
\\
G_{a,b,1,0}
&
\,=\,
0
\,=\,
G_{1,0,c,d},
\\
G_{a,b,2,0}
&
\,=\,
0
\,=\,
G_{2,0,c,d};
\endaligned
\]

\smallskip\noindent{\bf (3)}\,
satisfies in addition the sporadic normalization conditions:
\reqnomode\usetagform{EngelLie}
\begin{align}
G_{3,0,0,1}
&
\,=\,
0
\,=\,
G_{0,1,3,0},
\notag
\\
\Im\,G_{3,0,1,1}
&
\,=\,
0
\,=\,
\Im\,G_{1,1,3,0}.
\tag{\qed}
\end{align}
\end{Theorem}

There is of course {\em no} uniqueness of a rigid biholomorphic map 
which transfers $M$ to an $M'$ satisfying all these normalization
conditions
{\small\bf (1)}, {\small\bf (2)}, {\small\bf (3)}, 
just because any post-composition with a
dilation-rotation map:
\[
(z',\zeta',w')
\,\,\,\longmapsto\,\,\,
\big(
\rho^{1/2}\,e^{i\varphi}\,z',\,\,
e^{2i\varphi}\,\zeta',\,\,
\rho\,w'
\big)
\,\,=\,\,
(z'',\zeta'',w'')
\eqno
{\scriptstyle{(\rho\,\in\,\R_+^\ast,\,\,
\varphi\,\in\,\R)}},
\]
will transfer $M'$ into an $M'' = \{ u'' = \maux'' + G'' \}$ which
enjoys again the normalization conditions
{\small\bf (1)}, {\small\bf (2)}, {\small\bf (3)}, 
since one obviously has:
\[
G_{a,b,c,d}''\,
\rho^{\frac{a+c-2}{2}}\,
e^{i\varphi(a+2b-c-2d)}
\,=\,
G_{a,b,c,d}'.
\]

Remind that such dilation-rotation maps 
parametrize the $2$-dimensional
isotropy group of the origin
for the Gaussier-Merker model 
$\big\{ u' = \maux(z', \zeta', \overline{z}',
\overline{\zeta}') \big\}$.
Fortunately, an examination of
our analysis above can show that 
these two parameters $\rho$, $\varphi$
are the only ambiguity,
since once one assumes that $f = z + f_2 + f_3 + \cdots$,
with no $\rho^{1/2}\, e^{i\varphi}$ in front of $z$, that
$g = \zeta + g_1 + g_2 + \cdots$, and that
$h = w + h_3 + h_4 + \cdots$,
with no $\rho^{1/2}\, e^{i\varphi}$,
our reasonings showed {\em uniqueness} (exercise)
of the map to normal form.

\smallskip

To finish, let us abbreviate the space of power series $G = G(z,
\zeta, \overline{z}, \overline{\zeta})$ satisfying the normalization
conditions {\small\bf (1)}, {\small\bf (2)}, {\small\bf (3)} as:
\[
\mathfrak{N}_{2,1}.
\]

\begin{Corollary}
Two rigid
$\mathcal{C}^\omega$ 
hypersurfaces $M^5 \subset \C^3$ and ${M'}^5 \subset {\C'}^3$
belonging to $\mathfrak{C}_{2,1}$, both brought into normal form:
\[
\aligned
u
&
\,=\,
\maux
+
G,
&
\ \ \ \ \ \ \ \ \ \ \ \ \ \ \ \ \ \ \ \
G
\,\in\,
\mathfrak{N}_{2,1},
\\
u'
&
\,=\,
\maux'
+
G',
&
\ \ \ \ \ \ \ \ \ \ \ \ \ \ \ \ \ \ \ \
G'
\,\in\,
\mathfrak{N}_{2,1}',
\endaligned
\]
are rigidly biholomorphically equivalent if and only if
there exist two constants $\rho \in \R_+^\ast$,
$\varphi \in \R$, such that for all $a$, $b$, $c$, $d$:
\[
G_{a,b,c,d}
\,=\,
G_{a,b,c,d}'\,
\rho^{\frac{a+c-2}{2}}\,
e^{i\varphi(a+2b-c-2d)}.
\eqno\qed
\]
\end{Corollary}

Granted that hypersurfaces can be put into such a normal form, 
this criterion is quite effective to determine whether
two $M, M' \in \mathfrak{C}_{2,1}$ are rigidly equivalent.

\Section{\bf Finalized Expression of $\Qaux_0$}
\label{finalized-expression-Q0}
\HEAD{{\ref{finalized-expression-Q0}}.~{\sf Finalized Expression 
of $\Qaux_0$}
}{
Zhangchi {\sc Chen}, Wei Guo {\sc Foo}, Joël {\sc Merker}, 
The Anh {\sc Ta}}

In this section, we revisit the secondary invariant $\Qaux_0$. Our
goal is to transform $\Qaux_0$ into a new expression which makes
transparent two interesting features of $\Qaux_0$: that it is
real-valued and of order 5 (not 6 as it was first obtained by Cartan's
method in~\cite{Foo-Merker-Ta-2019}). The calculations in the
following are laborious, and for readers who are only interested in
the finalized expression of $\Qaux_0$, we suggest to use a
mathematical software for symbolic computations to have a quick and
easy check to confirm that the finalized
expression~{~(\ref{final-compact-Q0})} of $\Qaux_0$ indeed agrees with
the expression of $\Qaux_0$ obtained previously
in~\cite{Foo-Merker-Ta-2019}, which will be recalled later in this
section as the formula~{~(\ref{Q0-FMT})} .

\begin{Proposition}
\label{Prp-finalized-expression-Q0}
The secondary invariant $\Qaux_0$ can be brought into the following form
\begin{align}
\label{final-compact-Q0}
    \Qaux_0
&=
    \Baux \, \Iaux_0 
    + 
    \overline{\Baux} \, \overline{\Iaux_0}
    -
    \Baux \, \overline{\Baux}
    +
    \frac{2}{3} \, 
    \Re\bigg\{ \mathcal{L}_1
            \Big[
            \frac{\LbarLbark}{\Lbark} 
            \Big] 
        \bigg\}
    +
    \frac{1}{3} \, 
    \Re\Big( \Lbar 
            (\Paux)
        \Big).
\end{align}

\end{Proposition}

The rest of this section is devoted to the proof of 
Proposition~{\ref{Prp-finalized-expression-Q0}}. 
Let us first recall the formulas of $\Iaux_0, \Vaux_0, \Qaux_0$ 
from~{\cite{Foo-Merker-Ta-2019}}.
\begin{equation}
    \Iaux_0 
=
    -\frac{1}{3}
        \frac{\KLbarLbark}{(\Lbark)^2}
        +
        \frac{1}{3}
        \frac{\KLbark \; \LbarLbark}{(\Lbark)^3}
    +  \frac{2}{3}
        \frac{\LLkbar}{\Lkbar}
        +
        \frac{2}{3}
        \frac{\LLbark}{\Lbark},
\end{equation}

\begin{equation}
    \Vaux_0 
=
    -\frac{1}{3}
    \frac{\LbarLbarLbark}{\Lbark}
    +
    \frac{5}{9}
    \frac{( \LbarLbark )^2}{(\Lbark)^2}
    -
    \frac{1}{9}
    \frac{\LbarLbark \; \barP}{\Lbark}
    +
    \frac{1}{3} \Lbar(\barP)
    -
    \frac{1}{9}
    \barP \, \barP,
\end{equation}
and
\begin{equation}
\label{Q0-FMT}
    \Qaux_0
=
    \frac{1}{2}
    \bigg\{ 
    \Baux \,\Iaux_0
    +
    \Lbar(\Iaux_0)
    -
    \frac{\overline{\Baux} \,\, \Kbar(\Iaux_0)}{ \Lkbar}
    -
    \frac{\mathcal{K}(\Vaux_0)}{\Lbark}
    \bigg\},
\end{equation}
where
\[
   \Baux 
=
    \frac{1}{3}
    \Big(
    \frac{\LbarLbark}{\Lbark} 
    -
    \barP
    \Big)
\;\;\;\;\text{ and }\;\;\;\;
    \overline{\Baux} 
=
    \frac{1}{3}
    \Big(
    \frac{\LLkbar}{\Lkbar} 
    -
    \Paux
    \Big).
\]

For convenience, we will do calculations with  $3\Iaux_0, 9\Vaux_0, 18 \big\vert \Lbark \big\vert^2 \Qaux_0 $ and $3\Baux, 3\overline{\Baux}$.
\begin{equation}
\label{3I0}
    3\Iaux_0 
=
    \frac{\KLbark  \,\, \LbarLbark}{(\Lbark)^3}
    -
    \frac{\KLbarLbark}{(\Lbark)^2}
    +
    2
    \frac{\LLkbar}{\Lkbar}
    +
    2
    \frac{\LLbark}{\Lbark},
\end{equation}
\begin{equation}
\label{9V0}
    9\Vaux_0 
=
    5
    \frac{( \LbarLbark )^2}{(\Lbark)^2}
    -
    \frac{3\LbarLbarLbark  +  \LbarLbark \, \barP}{\Lbark}
    +
    3
    \Lbar(\barP)
    -
    \barP \, \barP,
\end{equation}
\begin{align}
\label{18Q0}
    18  \big\vert \Lbark \big\vert^2  \Qaux_0
=
&\big[
    3\Baux  \, 3\Iaux_0
    +
    3\Lbar(3\Iaux_0)
    \big] \Lkbar  \, \Lbark 
\\ \nonumber
&\;-
    3\overline{\Baux}  \,\, \Kbar(3\Iaux_0)  \, \Lbark 
    -
    \mathcal{K}(9\Vaux_0) \,\Lkbar,
\end{align}
with 
\[
   3\Baux 
=
    \frac{\LbarLbark}{\Lbark} 
    -
    \barP
\;\;\;\;\text{ and }\;\;\;\;
    3\overline{\Baux} 
=
    \frac{\LLkbar}{\Lkbar} 
    -
    \Paux.
\]

In order to transform the expression~{~(\ref{18Q0}}) of $18  \big\vert \Lbark \big\vert^2  \Qaux_0$, we will make use of the following identities.

\begin{Lemma}
\label{relations-of-vector-fields}
We have the following identities:
\begin{enumerate}
    \item $\mathcal{K}(\barP)
            =
            - \Paux \, \Lbark 
            - \LbarLk,
            $
    \item $\mathcal{K} \Lbar (\barP)
            =
            -\Lbark \cdot 2\Re \,
            \big(\Lbar(\Paux) \big)
            -\Paux \, \LbarLbark 
            - \LbarLbarLk,
            $
    \item $\Kbar(\Iaux_0)
            =
            (-2) \, \overline{\Iaux_0} \cdot \Lkbar.
    $
\end{enumerate}
\end{Lemma}

\begin{proof}
The identities (1) and (3) are obtained in Lemma 2.7  and Lemma 10.6 of~{~\cite{Foo-Merker-Ta-2019}}, respectively. 

For the identity (2), we use the relation
$ 
    [\mathcal{K}, \Lbar ] 
= 
    \mathcal{K} \, \Lbar
    -
    \Lbar \, \mathcal{K}
=
    - \Lbark \, \mathcal{L}_1
$
from (2.9) of~{~\cite{Foo-Merker-Ta-2019}} to deduce that
\begin{align*}
        \mathcal{K} \Lbar (\barP)
    &=
        \Lbar \, \mathcal{K}(\barP)
        - 
        \Lbark \, \mathcal{L}_1(\barP)
        \\ 
    &=
        \Lbar 
        \Big[ -\Paux \, \Lbark - \LbarLk  \Big]
        -
        \Lbark \, \mathcal{L}_1(\barP)
        \;\;\;\;\;\; (\text{using (1)})
        \\
    &=
        -\Lbar(\Paux) \, \Lbark 
        -\Paux \, \LbarLbark 
        -\LbarLbarLk 
        -\Lbark \, \mathcal{L}(\barP)
        \\
    &=
        -\Lbark 
        \Big[ \, \Lbar(\Paux) + \mathcal{L}_1(\barP) \, \Big]
        -\Paux \, \LbarLbark 
        - \LbarLbarLk
        \\
    &=
        -\Lbark \cdot 2\Re \,
        \big(\Lbar(\Paux) \big)
        -\Paux \, \LbarLbark 
        - \LbarLbarLk.
\end{align*}
\end{proof}

\begin{proof}[Proof of Proposition~{~\ref{Prp-finalized-expression-Q0}}]
We first substitute the identity (3) of Lemma~{~(\ref{relations-of-vector-fields})} into the term
$
- 3\overline{\Baux} \; \Kbar(3\Iaux_0) \, \Lbark
$ of
$
18 \big\vert \Lbark \big\vert^2 \Qaux_0
$
to obtain
\[
    - 3\overline{\Baux} \; \Kbar(3\Iaux_0) \, \Lbark
=
    - 3\overline{\Baux} (-6 \overline{\Iaux_0} \, \Lkbar) \, \Lbark
=
    2 \cdot 3\overline{\Baux} \cdot 3\overline{\Iaux_0} \; \Lbark \Lkbar,
\]
with which the sum on the right hand side of~{~(\ref{18Q0})} can be rewritten as
\begin{align}
\label{18Q0-first-step}
  18 \big\vert \Lbark \big\vert^2 \Qaux_0
=
    &\Big[
    3\Baux \cdot 3\Iaux_0
    +
    3\overline{\Baux} \cdot 3\overline{\Iaux_0}
    \Big]
    \Lkbar \, \Lbark 
    +
    3\Lbar(3\Iaux_0)
    \Lkbar \, \Lbark 
\\ \nonumber 
&\;+
    3\overline{\Baux} \cdot 3\overline{\Iaux_0} \; \Lkbar \, \Lbark
    -
    \mathcal{K}(9\Vaux_0)\Lkbar.  
\end{align}

Observe that on the right hand side of~{~(\ref{18Q0-first-step})}, the
first term is already real-valued, which hints that we should keep it
untouched until the very end of the proof. We proceed by transforming
the other terms so that the real-valuedness of the sum
in~{~(\ref{18Q0-first-step})} will be transparent. Our strategy is to
look for terms that involve in $\Paux$ and $\barP$ first. From the
expression~{~(\ref{3I0})} of $3\Iaux_0$, one sees that the second term
of the right hand side of~{~(\ref{18Q0-first-step})} doesnot contain
$\Paux$ and $\barP$. Thus, we only need to extract parts involved in
$\Paux$ and $\barP$ from the last two terms of the right hand side
of~{~(\ref{18Q0-first-step})}.

Note that in the expression
$
3\overline{\Baux} \cdot 3\overline{\Iaux_0}
=
\Big(\frac{\LLkbar}{\Lkbar} - \Paux\Big) \cdot 3\overline{\Iaux_0}
=
\frac{\LLkbar}{\Lkbar} \, 3\overline{\Iaux_0} - \Paux \cdot 3\overline{\Iaux_0},
$
the only part involved in $\Paux$ and $\barP$ is 
$
- \Paux \cdot 3\overline{\Iaux_0}. 
$
We will see that by extracting terms involved in $\Paux$ and $\barP$ in 
$
- \mathcal{K}(9\Vaux_0)\Lkbar
$, which is 
\begin{equation}
 -
\mathcal{K}
\bigg\{ 
    -
    \frac{\LbarLbark \, \barP}{\Lbark}
    +
    3\Lbar(\barP)
    -
    \barP \, \barP
\bigg\}  \, \Lkbar,   
\end{equation}
we will obtain a conjugate of 
$
-\Paux \cdot 3\overline{\Iaux_0} \cdot \Lkbar \, \Lbark
$. Indeed, let us expand 
\begin{flushleft}
$   -
    \mathcal{K}
        \bigg\{ 
        -
        \frac{\LbarLbark \, \barP}{\Lbark}
        +
        3\Lbar(\barP)
        -
        \barP \, \barP
        \bigg\}
     \, \Lkbar
$
\end{flushleft}

\begin{flushleft}
$   =
    \bigg\{ 
        \mathcal{K} \, 
        \Big[
        \frac{\LbarLbark\, \barP}{\Lbark}
        \Big]
        -
        3\, \KLbar(\barP)
        +
        2\barP \, \mathcal{K}(\barP)
    \bigg\} 
    \, \Lkbar
$
\end{flushleft}

\begin{flushleft}
$   =
    \bigg\{ 
        \mathcal{K} \, 
        \Big[
        \frac{\LbarLbark}{\Lbark}
        \Big] 
        \, \barP
        +
        \frac{\LbarLbark}{\Lbark} \, \mathcal{K}(\barP)
        -
        3 \, \KLbar(\barP)
        +
        2\barP \, \mathcal{K}(\barP)
    \bigg\}
 \, \Lkbar
$
\end{flushleft}

\begin{flushleft}
$=
    \bigg\{ 
        \Big[
        \frac{\KLbarLbark}{\Lbark}
        -
        \frac{\LbarLbark  \, \KLbark}{(\Lbark)^2}
        \Big]
        \, \barP
        + \,
        \frac{\LbarLbark}{\Lbark} \, 
        \Big[ 
        - \Paux \, \Lbark -\LbarLk 
        \Big]
$
\end{flushleft}

\begin{flushleft}
$\;\;\;-
        3 \, 
        \Big[ 
        - \, \Lbark \cdot 2\Re\big(\Lbar(\Paux)\big)
        -
        \Paux\,\LbarLbark
        -
        \LbarLbarLk
        \Big]
$
\end{flushleft}

\begin{flushleft}
$\;\;\;+ 
        2
        \barP \, 
        \Big[ 
             -\Paux\,\Lbark - \LbarLk 
        \Big]
\bigg\}  \, \Lkbar
\;\; (\text{using (1) and (2) of Lemma~{~\ref{relations-of-vector-fields}}})
$
\end{flushleft}

\begin{flushleft}
$=
    \bigg\{ 
        \barP \,       
        \Big[
        -
        \frac{\KLbark \, \LbarLbark}{(\Lbark)^2}
        +
        \frac{\KLbarLbark}{\Lbark}
        \Big]
        -
        \Paux \, \LbarLbark
$
\end{flushleft}

\begin{flushleft}
$\;\;\;-
        \frac{\LbarLbark \, \LbarLk}{\Lbark}
        +  
        6 \, \Lbark \cdot \,\Re\Big( \Lbar(\Paux) \Big)
        +
        3 \, \Paux \,\LbarLbark
$
\end{flushleft}

\begin{flushleft}
$\;\;\;+
        3\, \LbarLbarLk
        - 
        2\, \barP \,\Paux \,\Lbark 
        - 
        2\, \barP \, \LbarLk 
\bigg\}  \, \Lkbar.
$
\end{flushleft}

At this point, we extract 
$
- \, \barP  \, \Lbark  \, \Lkbar \cdot 3\Iaux_0
$
to obtain

\begin{flushleft}
$-
    \mathcal{K}
    \bigg\{ 
    - \,
    \frac{\LbarLbark \, \barP}{\Lbark}
    +
    3\Lbar(\barP)
    -
    \barP \, \barP
    \bigg\}  \, \Lkbar
$
\end{flushleft}

\begin{flushleft}
$= 
    - 
    \barP  \, 
    \Lbark \, \Lkbar
\bigg\{ 
    \frac{\KLbark \, \LbarLbark}{(\Lbark)^3}
    -
    \frac{\KLbarLbark}{(\Lbark)^2}
    + 
    2 \frac{\LLbark}{\Lbark}
    +
    2 \frac{\LLkbar}{\Lkbar}
\bigg\}$
\end{flushleft}

\begin{flushleft}
$\;\;\;+ 
    2  
    \barP  \, 
    \LLbark \, \Lkbar
    +
    2  
    \barP  \, 
    \LLkbar \, \Lbark
    -
    \Paux \, 
    \LbarLbark \, \Lkbar
$
\end{flushleft}

\begin{flushleft}
$\;\;\;+
    \frac{\LbarLbark \, \LbarLk \, \Lkbar}{\Lbark}
    +
    6 \,
    \Lbark \, \Lkbar \cdot
    \Re\Big(\Lbar(\Paux)\Big)
    +
    3  
    \Paux \,\LbarLbark \, \Lkbar$
\end{flushleft}

\begin{flushleft}
$\;\;\;+
    3 
    \LbarLbarLk \,\Lkbar
    -
    2 
    \barP \,\Paux \,
    \Lbark \, \Lkbar
    - 
    2 
    \barP \,\LbarLk \,\Lkbar
$
\end{flushleft}

\begin{flushleft}
$=
    - \barP  \, 
    \Lbark \, \Lkbar \cdot 3\Iaux_0
    -
    \frac{\LbarLbark \, \LbarLk \, \Lkbar}{\Lbark}
    +
    3 \LbarLbarLk \,\Lkbar$
\end{flushleft}

\begin{flushleft}
$\;\;\;+ 
    2
    \Re\Big(\frac{\LbarLbark}{\Lbark}\, \Paux\Big) 
    \,\Lbark \, \Lkbar
    +
    6 
    \Lbark \, \Lkbar\,
    \Re\Big(\Lbar(\Paux)\Big)
    -
    2\, \barP\,\Paux\, \Lbark \, \Lkbar
$
\end{flushleft}

\begin{flushleft}
$=
    -
    \barP  \, 
    \Lbark  \, \Lkbar \cdot 3\Iaux_0
    -
    \frac{\LbarLbark \, \LbarLk \, \Lkbar}{\Lbark}
    +
    3  \LbarLbarLk \,\Lkbar
$
\end{flushleft}

\begin{flushleft}
$\;\;\;+ 
    2
    \Re\Big(\frac{\LbarLbark}{\Lbark}\, \Paux\Big) 
    \,
    \big\vert \Lbark \big\vert^2
    + 
    6 \,
    \big\vert \Lbark\big\vert^2
    \,\Re\Big(\Lbar(\Paux)\Big)
    -
    2
    \big\vert \Paux \big\vert^2 \, 
    \big\vert \Lbark \big\vert^2,
$
\end{flushleft}
whose last 3 terms are real-valued.

Now, we substitute the just obtained expansion
\begin{flushleft}
$(-1)
    \mathcal{K}
\Big\{ 
    (-1) \,
    \frac{\LbarLbark. \, \barP}{\Lbark}
    +
    3\Lbar(\barP)
    -
    \barP \, \barP
\Big\} . \, \Lkbar
$
\end{flushleft}

\begin{flushleft}
$=
    (-1) \, 
    \barP  \, 
    \Lbark  \, \Lkbar.\, 3\Iaux_0
    + 
    (-1) \, 
    \frac{\LbarLbark \, \LbarLk \, \Lkbar}{\Lbark}
    +
    3 \, \LbarLbarLk \,\Lkbar
$
\end{flushleft}

\begin{flushleft}
$\;\;\;+ 
    2\,
    \Re(\frac{\LbarLbark}{\Lbark}\Paux) 
    .\,
    \big\vert \Lbark \big\vert^2
    + 
    6 \,
    \big\vert \Lbark\big\vert^2
    .\,\Re(\Lbar(\Paux))
    +
    (-2) \,
    \big\vert \Paux \big\vert^2.\, 
    \big\vert \Lbark \big\vert^2
$
\end{flushleft}
back into the expression~{~(\ref{18Q0-first-step})} of $18 \big\vert
\Lbark \big\vert^2 \Qaux_0$ to obtain
\begin{flushleft}
$18 
    \big\vert \Lbark \big\vert^2 
    \Qaux_0
= \big[
    3\Baux \, 3\Iaux_0
    +
    3\overline{\Baux} \, 3\overline{\Iaux_0}
\big]
    \Lkbar \, \Lbark 
+
    3\Lbar(3\Iaux_0)
    \Lkbar \, \Lbark 
$
\end{flushleft}

\begin{flushleft}
$\;\;\;\;\;\;\;\;\;\;\;\;\;\;\;\;\;\;\;\;\;\;\;\;\;\;-
    \Paux \, 3\overline{\Iaux_0} \, \Lkbar \, \Lbark
+
    \frac{\LLkbar}{\Lkbar}\, 3\overline{\Iaux_0} \, \Lkbar \, \Lbark
$
\end{flushleft}

\begin{flushleft}
$\;\;\;\;\;\;\;\;\;\;\;\;\;\;\;\;\;\;\;\;\;\;\;\;\;\;-
    \mathcal{K}
    \Big\{ 
    5
    \frac{(\LbarLbark)^2}{(\Lbark)^2}
    -
    3
    \frac{\LbarLbarLbark}{\Lbark}
    \Big\}  \, \Lkbar
$
\end{flushleft}

\begin{flushleft}
$\;\;\;\;\;\;\;\;\;\;\;\;\;\;\;\;\;\;\;\;\;\;\;\;\;\;-
    \mathcal{K}
    \Big\{ 
    -
    \frac{\LbarLbark. \, \barP}{\Lbark}
    +
    3\Lbar(\barP)
    -
    \barP \, \barP
    \Big\}  \, \Lkbar
$
\end{flushleft}

\begin{flushleft}
$\;\;\;\;\;\;\;\;\;\;\;\;\;\;\;\;\;\;\;\;\;\;\;
=
    \big[
    3\Baux \, 3\Iaux_0
    +
    3\overline{\Baux} \, 3\overline{\Iaux_0}
    \big]
    \Lkbar \, \Lbark 
+
    3\Lbar(3\Iaux_0)
    \Lkbar \, \Lbark 
$
\end{flushleft}

\begin{flushleft}
$\;\;\;\;\;\;\;\;\;\;\;\;\;\;\;\;\;\;\;\;\;\;\;\;\;\;\;-
    \Paux \, 3\overline{\Iaux_0} \, \Lkbar \, \Lbark
+
    \frac{\LLkbar}{\Lkbar}\, 3\overline{\Iaux_0} \, \Lkbar \, \Lbark
 $
\end{flushleft}   

\begin{flushleft}
$\;\;\;\;\;\;\;\;\;\;\;\;\;\;\;\;\;\;\;\;\;\;\;\;\;\;\;-
    \mathcal{K}
    \Big\{ 
    5
    \frac{(\LbarLbark)^2}{(\Lbark)^2}
    -
    3
    \frac{\LbarLbarLbark}{\Lbark}
    \Big\} \, \Lkbar
$
\end{flushleft}

\begin{flushleft}
$\;\;\;\;\;\;\;\;\;\;\;\;\;\;\;\;\;\;\;\;\;\;\;\;\;\;\;-
     \barP \, 
    \Lbark \, \Lkbar\, 3\Iaux_0
    -
    \frac{\LbarLbark \, \LbarLk \, \Lkbar}{\Lbark}
$

\end{flushleft}

\begin{flushleft}
$\;\;\;\;\;\;\;\;\;\;\;\;\;\;\;\;\;\;\;\;\;\;\;\;\;\;\;+
    3 \, \LbarLbarLk \,\Lkbar
    - 
    2 \,
    \Re\Big(\frac{\LbarLbark}{\Lbark}\Paux \Big) 
    \,
    \big\vert \Lbark \big\vert^2
$
\end{flushleft}

\begin{flushleft}
$\;\;\;\;\;\;\;\;\;\;\;\;\;\;\;\;\;\;\;\;\;\;\;\;\;\;\;+ 
    6 \,\big\vert \Lbark\big\vert^2
    \,
    \Re\Big(\Lbar(\Paux)\Big)
    -
    2 \,
    \big\vert \Paux \big\vert^2 \, \big\vert \Lbark \big\vert^2,
$
\end{flushleft}
which after rearranging gives
\begin{equation}
\label{18Q0-second-step}
\begin{split}
    18 \big\vert \Lbark \big\vert^2 \Qaux_0
&=
    \big[
    3\Baux \, 3\Iaux_0
    +
    3\overline{\Baux} \, 3\overline{\Iaux_0}
    \big]
    \Lkbar \, \Lbark 
\\ 
&\;\;\;\;\;-
    \big[ 3\overline{\Iaux_0} \, \Paux
        + 3\Iaux_0 \, \barP
    \big]  \, \Lkbar \, \Lbark
\\ 
&\;\;\;\;\;+
    3\Lbar(3\Iaux_0)
    \Lkbar \, \Lbark 
+
    \frac{\LLkbar}{\Lkbar}\, 3\overline{\Iaux_0} \, \Lkbar \, \Lbark
\\ 
&\;\;\;\;\;-
    \mathcal{K}
    \Big\{ 
    5
    \frac{(\LbarLbark)^2}{(\Lbark)^2}
    -
    3
    \frac{\LbarLbarLbark}{\Lbark}
    \Big\}  \, \Lkbar
\\ 
&\;\;\;\;\;-
         \frac{\LbarLbark \, \LbarLk \, \Lkbar}{\Lbark}
        +
        3 \, 
        \LbarLbarLk \,\Lkbar
\\ 
&\;\;\;\;\;- 
        2\,\Re\Big(\frac{\LbarLbark}{\Lbark}\Paux\Big) \,\big\vert \Lbark \big\vert^2
        +
        6 \,\big\vert \Lbark\big\vert^2
        \,\Re\Big(\Lbar(\Paux)\Big)
\\ 
&\;\;\;\;\;-
        2 \,
        \big\vert \Paux \big\vert^2\, \big\vert \Lbark \big\vert^2.
\end{split}    
\end{equation}

Next, we want to extract a conjugate of
$
\frac{\LLkbar}{\Lkbar} \, 3\overline{\Iaux_0} \, 
\Lkbar \, \Lbark
$
and a copy of 
$
- 3\Lbar(3\Iaux_0) \, \Lbark \, \Lkbar 
$
from
$-
    \mathcal{K}
    \Big\{ 
    5
    \frac{(\LbarLbark)^2}{(\Lbark)^2}
    -
    3
    \frac{\LbarLbarLbark}{\Lbark}
    \Big\}  \, \Lkbar.
$

We first expand
\begin{equation}
\label{3I0-LLkbar}
\frac{\LLkbar}{\Lkbar}\, 3\overline{\Iaux_0} \, \Lkbar \, \Lbark
 =
    3\overline{\Iaux_0} \,\LLkbar \, \Lbark     
\end{equation}
\begin{equation*}
    =
    \frac{\overline{\mathcal{K}}\Lkbar \; (\LLkbar)^2 \;\Lbark}{(\Lkbar)^3}
    -
    \frac{\overline{\mathcal{K}}\LLkbar \; \LLkbar \; \Lbark}{(\Lkbar)^2}    
\end{equation*} 

\begin{equation*}
    +
    2\; \frac{\Lbar\Lkbar \; \LLkbar \; \Lbark}{\Lkbar}
    +
    2 \; \LbarLbark \; \LLkbar,   
\end{equation*}
and
\begin{equation}
\label{Lbar-3I0}
 -3\; \Lbar(3\Iaux_0) \; \Lbark \; \Lkbar    
\end{equation}
\begin{equation*}
  =
    9 \;
    \frac{ \KLbark \; (\LbarLbark)^2 \; \Lkbar}{(\Lbark)^3}
    -
    9 \;
    \frac{ \KLbarLbark \; \LbarLbark \; \Lkbar}{(\Lbark)^2}  
\end{equation*}
\begin{equation*}
    -
    3 \;
    \frac{ \KLbark \; \LbarLbarLbark \; \Lkbar}{(\Lbark)^2}
+
    3 \;
    \frac{ \mathcal{K}\LbarLbarLbark \; \Lkbar}{ \Lbark} 
\end{equation*}
\begin{equation*}
    -
    3 \;
    \frac{\LLbark \; \LbarLbark \; \Lkbar}{\Lbark}
+
    3 \; \mathcal{L}_1\LbarLbark \; \Lkbar
\end{equation*}
\begin{equation*}
    \;\;\;\;\;\;\;\;\;\;\;\;\;\;\;\;+
        12 \; \big\vert \Lbark \big\vert^2
            \; \Re\Big( \frac{\LLbark \; \LbarLbark}{(\Lbark)^2} \Big)
    -
        12 \; \big\vert \Lbark \big\vert^2
            \; \Re\Big( \frac{\mathcal{L}_1\LbarLbark}{\Lbark} \Big).
\end{equation*}

We now use the expansions~{~(\ref{3I0-LLkbar})} and~{~(\ref{Lbar-3I0})} to expand 
$-
    \mathcal{K}
    \Big\{ 
    5
    \frac{(\LbarLbark)^2}{(\Lbark)^2}
    -
    3
    \frac{\LbarLbarLbark}{\Lbark}
    \Big\} \, \Lkbar
$
as follows.
\begin{equation} 
\label{-1K5}
-
    \mathcal{K}
    \Big\{ 
    5 \;
    \frac{(\LbarLbark)^2}{(\Lbark)^2}
    -
    3 \;
    \frac{\LbarLbarLbark}{\Lbark}
    \Big\}  \; \Lkbar
\end{equation}

\begin{flushleft}
$=  \bigg\{
    \Big[ 
    -10 \;
    \frac{\LbarLbark \; \KLbarLbark}{(\Lbark)^2}
    +
    10 \;
    \frac{(\LbarLbark)^2 \; \KLbark}{(\Lbark)^3}
    \Big] 
$
\end{flushleft}

\begin{flushleft}
$\;\;\;\;\;\;+
    \Big[ 
    3 \;
    \frac{\mathcal{K}\LbarLbarLbark}{\Lbark}
    -
    3 \;
    \frac{\LbarLbarLbark \; \KLbark}{(\Lbark)^2}
    \Big] 
    \bigg\}\, \Lkbar
$
\end{flushleft}

\begin{flushleft}
$= 
    10 \;
    \frac{\KLbark \; (\LbarLbark)^2 \; \Lkbar}{(\Lbark)^3}
    -
    10 \;
    \frac{\LbarLbark \; \KLbarLbark  \; \Lkbar}{(\Lbark)^2}
$
\end{flushleft}

\begin{flushleft}
$\;\;\;\;\;\;-
    3 \;
    \frac{\LbarLbarLbark \; \KLbark \; \Lkbar}{(\Lbark)^2}
+
    3 \;
    \frac{\mathcal{K}\LbarLbarLbark}{\Lbark}
$
\end{flushleft}

\begin{flushleft}
$=
    \overline{3\overline{\Iaux_0} \; \LLkbar \; \Lbark  }
+
    9 \;
    \frac{\KLbark \; (\LbarLbark)^2 \; \Lkbar}{(\Lbark)^3}
-
    9 \;
    \frac{\KLbarLbark \; \LbarLbark \; \Lkbar}{(\Lbark)^2}
$
\end{flushleft}

\begin{flushleft}
$\;\;\;\;\;\;-
    3 \;
    \frac{\KLbark \; \LbarLbarLbark \; \Lkbar}{(\Lbark)^2}
+
    \frac{3 \; \KLbar\LbarLbark \; \Lkbar  -2 \; \LLbark \; \LbarLbark \; \Lkbar}{\Lbark}
$
\end{flushleft}

\begin{flushleft}
$\;\;\;\;\;\;-
    2 \; \LbarLbark \; \LLkbar
$
\end{flushleft}

\begin{flushleft}
$
= 
    \overline{3\overline{\Iaux_0} \; \LLkbar \; \Lbark  }
    -
    3 \Lbar(3\Iaux_0) \; \Lbark \; \Lkbar
    +
    \frac{\LbarLbark \; \LbarLk \; \Lkbar}{\Lbark}
$
\end{flushleft}

\begin{flushleft}
$\;\;\;\;\;\;-
        3 \; \LbarLbarLk \; \Lkbar
-
        12 \; \big\vert \Lbark \big\vert^2
            \; \Re\Big( \frac{\LLbark \; \LbarLbark}{(\Lbark)^2} \Big)
$
\end{flushleft}

\begin{flushleft}
$\;\;\;\;\;\;+
        12 \; \big\vert \Lbark \big\vert^2
            \; \Re\Big( \frac{\mathcal{L}_1\LbarLbark}{\Lbark} \Big)
-
        2 \; \big\vert \LbarLbark \big\vert^2.
$
\end{flushleft}

Substituting the expansion~{~(\ref{-1K5})} into the right hand side of~{~(\ref{18Q0-second-step})} leads to
\begin{equation}
\label{18Q0-third-step}
    \begin{split}
   18 \big\vert \Lbark \big\vert^2 \Qaux_0
&= 
    \Big[ 
    3\Baux \; 3\Iaux_0 
    + 
    3\overline{\Baux} \; 3\overline{\Iaux_0}
    \Big] \; \Lbark \; \Lkbar
\\ 
&\;\;\;\;-
    \Big[ 
    3\overline{\Iaux_0} \; \Paux
    + 
    3\Iaux_0 \; \barP
    \Big] \; \Lbark \; \Lkbar
\\ 
&\;\;\;\;+
        3\overline{\Iaux_0} \; \LLkbar \; \Lbark
        +
        \overline{3\overline{\Iaux_0} \; \LLkbar \; \Lbark  }
\\ 
&\;\;\;\;-
        12 \; \big\vert \Lbark \big\vert^2
            \; \Re\Big( \frac{\LLbark \; \LbarLbark}{(\Lbark)^2} \Big)
\\ 
&\;\;\;\;+
        12 \; \big\vert \Lbark \big\vert^2
            \; \Re\Big( \frac{\mathcal{L}_1\LbarLbark}{\Lbark} \Big)
\\ 
&\;\;\;\;+
        2 \; \big\vert \Lbark \big\vert^2 \; \Re\Big( -
        \frac{ \Paux \; \LbarLbark}{\Lbark} + 3\Lbar(\Paux) \Big)
\\
&\;\;\;\;-
        2 \; \big\vert \Lbark \big\vert^2 \; |\Paux|^2
-
        2\; \big\vert \LbarLbark \big\vert^2.    
    \end{split}
\end{equation}

At this point, we can see from the right hand side of~{~(\ref{18Q0-third-step})} that $\Qaux_0$ is real valued and of order 5, but observe that we can contract more terms into 
$ 
\big\vert \Lbark \big\vert^2 \, 
3\overline{\Baux} \; 3\Baux.
$

Let us expand
\begin{equation}
\label{BbarB}
\big\vert \Lbark \big\vert^2 
3\overline{\Baux} \; 3\Baux
=
 \big\vert \Lbark \big\vert^2 \; |\Paux|^2
 -  
 2 \; \big\vert \Lbark \big\vert^2 \; \Re\Big( \frac{\Paux  \LbarLbark}{\Lbark} \Big)
 +
 \big\vert \LbarLbark \big\vert^2.
\end{equation}

By using the identity~{~(\ref{BbarB})}, we now substitute
$
-2 \, \big\vert \Lbark \big\vert^2 \;
3\overline{\Baux} \; 3\Baux
$
into the expansion~{~(\ref{18Q0-third-step})} of 
$
18 \big\vert \Lbark \big\vert^2 \Qaux_0
$ 
in order to obtain
\begin{equation}
\label{18Q0-fourth-step}
    18 \big\vert \Lbark \big\vert^2 \Qaux_0 = 
\end{equation}

\begin{flushleft}
$= 
    \Big[ 
    3\Baux \; 3\Iaux_0 
    + 
    3\overline{\Baux} \; 3\overline{\Iaux_0}
    \Big] \; \Lbark \; \Lkbar
-
    \Big[ 
    3\overline{\Iaux_0} \; \Paux
    + 
    3\Iaux_0 \; \barP
    \Big] \; \Lbark \; \Lkbar 
$
\end{flushleft}

\begin{flushleft}
$\;\;\;\;+
        3\overline{\Iaux_0} \; \LLkbar \; \Lbark
+
        \overline{3\overline{\Iaux_0} \; \LLkbar \; \Lbark  }
-       
        2 \, \big\vert \Lbark \big\vert^2 \;
        3\overline{\Baux} \; 3\Baux
$
\end{flushleft}

\begin{flushleft}
$\;\;\;\;-
        12 \, \big\vert \Lbark \big\vert^2
            \; \Re\Big( \frac{\LLbark \; \LbarLbark}{(\Lbark)^2} \Big)
+
        12\, \big\vert \Lbark \big\vert^2
            \; \Re\Big( \frac{\mathcal{L}_1\LbarLbark}{\Lbark} \Big)
$
\end{flushleft}

\begin{flushleft}
$\;\;\;\;+
        6 \, \big\vert \Lbark \big\vert^2 \; \Re\Big(- \frac{ \Paux \; \LbarLbark}{\Lbark} + \Lbar(\Paux)\Big)
$
\end{flushleft}

\begin{flushleft}
$=
    \Big[ 
    3\Baux \; 3\Iaux_0 
    + 
    3\overline{\Baux} \; 3\overline{\Iaux_0}
    \Big]\; \Lbark \; \Lkbar
+
    \Big[ 
    3\overline{\Iaux_0} \; \Paux
    + 
    3\Iaux_0 \; \barP
    \Big]\; \Lbark \; \Lkbar
$
\end{flushleft}

\begin{flushleft}
$\;\;\;\;+
        3\overline{\Iaux_0} \; \LLkbar \; \Lbark
        +
        \overline{3\overline{\Iaux_0} \; \LLkbar \; \Lbark  }
        -
        2 \, \big\vert \Lbark \big\vert^2 \;
        3\overline{\Baux} \; 3\Baux
$
\end{flushleft}

\begin{flushleft}
$\;\;\;\;+
        12 \, \big\vert \Lbark \big\vert^2
            \; \Re\bigg\{ \mathcal{L}_1
            \Big[
            \frac{\LbarLbark}{\Lbark} 
            \Big] \bigg\}
+
        6 \, \big\vert \Lbark \big\vert^2 \; \Re\Big( \Lbar (  \Paux ) \Big).
$
\end{flushleft}

At this point, a quick look at the first 4 terms on the right hand side of the expansion~{~(\ref{18Q0-fourth-step})} suggests that we should contract them as follows.
\begin{equation}
\label{contract-first-4-terms}
    \begin{split}
    &\Big[ 
        3\Baux \; 3\Iaux_0 
    + 
        3\overline{\Baux} \; 3\overline{\Iaux_0}
    \Big]\; \Lbark \; \Lkbar
    +
    \Big[ 
    3\overline{\Iaux_0} \; \Paux
    + 
    3\Iaux_0 \; \barP
    \Big] \; \Lbark \; \Lkbar
\\ 
&+
        3\overline{\Iaux_0} \;\LLkbar \; \Lbark
        +
        \overline{3\overline{\Iaux_0} \; \LLkbar \; \Lbark}
        -       
        2\, \big\vert \Lbark \big\vert^2 \;
        3\overline{\Baux} \; 3\Baux    
    \end{split}
\end{equation}

\begin{flushleft}
$\;\;\;\;\;\;\;\;\;\;\;\;=
    \Lbark \; \Lkbar
    \bigg\{
    \Big[ 
    3\Baux \; 3\Iaux_0 
    + 
    3\overline{\Baux} \; 3\overline{\Iaux_0}
    \Big]  
    +
    \Big[ 
    3\overline{\Iaux_0} \; \Paux
    + 
    3\Iaux_0 \; \barP
    \Big]
$
\end{flushleft}   

\begin{flushleft}
$\;\;\;\;\;\;\;\;\;\;\;\;\;\;\;\;\;\;\;\;\;\;\;\;\;\;\;\;\;\;\;\;\;\;\;\;\;\;\;+
    \frac{\LLkbar}{\Lkbar} \; 3\overline{\Iaux_0}
    +
    \frac{\LbarLbark}{\Lbark} \; 3\Iaux_0
    -  
    2 \cdot 
    3\overline{\Baux} \; 3\Baux
    \bigg\}
$
\end{flushleft}

\begin{flushleft}
$\;\;\;\;\;\;\;\;\;\;\;\;=
    \Lbark \;\Lkbar
    \bigg\{
    \Big[ 
    3\Baux \; 3\Iaux_0 
    + 
    3\overline{\Baux} \; 3\overline{\Iaux_0}
    \Big]  
    +
    \Big[ 
    3\overline{\Baux} \; 3\overline{\Iaux_0} 
    + 
    3\Baux \; 3\Iaux_0
    \Big]
-
    2\cdot 
    3\overline{\Baux} \; 3\Baux
    \bigg\}
$
\end{flushleft}

\begin{flushleft}
$\;\;\;\;\;\;\;\;\;\;\;\;=
    2 \, \big\vert \Lbark \big\vert^2 \;
    \bigg\{
     \Big[ 
    3\Baux \; 3\Iaux_0 
    + 
    3\overline{\Baux} \; 3\overline{\Iaux_0}
    \Big]  
    -
    3\overline{\Baux} \; 3\Baux
    \bigg\}.
$
\end{flushleft}

Substituting the contraction~{~(\ref{contract-first-4-terms})} into the right hand side of the expression~{~(\ref{18Q0-fourth-step})} gives 
\begin{equation}
\label{18Q0-final}
    \begin{split}
    18 \big\vert \Lbark \big\vert^2 \Qaux_0
&=
    2 \, \big\vert \Lbark \big\vert^2 \;
    \Big(
    3\Baux  \; 3\Iaux_0 
    + 
    3\overline{\Baux} \; 3\overline{\Iaux_0}
    -
    3\overline{\Baux} \; 3\Baux
    \Big)
\\
&\;\;\;\;+
        12 \, \big\vert \Lbark \big\vert^2
            \; \Re\bigg\{ \mathcal{L}_1
            \Big[
            \frac{\LbarLbark}{\Lbark} 
            \Big] \bigg\}
\\
&\;\;\;\;+
        6 \, \big\vert \Lbark \big\vert^2 \; \Re\Big( \Lbar (\Paux ) \Big).    
    \end{split}
\end{equation}

Finally, simplifying the factor 
$
18 \big\vert \Lbark \big\vert^2
$
on both side of~{~(\ref{18Q0-final})} gives us the desired expression~{~(\ref{final-compact-Q0})} of $\Qaux_0$.
\end{proof}

When we fully expand $\Qaux_0$ from the expression~{~(\ref{final-compact-Q0})}  using the formulas of $\Iaux_0$ and $\Baux$, we arrive at the following long expression of $\Qaux_0$, which only involves in the fundamental functions $\mathbf{k}$ and $\Paux$, and their derivatives:
\begin{align}
   \Qaux_0
&= \;\;\;\;
    \frac{2}{9} \,
    \Re\bigg\{
    \frac{\KLbark \; (\LbarLbark)^2}{(\Lbark)^4}
    \bigg\}
\\ \nonumber
&\;\;\;\;-
    \frac{2}{9} \,
    \Re\bigg\{
    \frac{\KLbarLbark \; \LbarLbark + \KLbark \; \LbarLbark \; \barP}{(\Lbark)^3}
    \bigg\}
\\ \nonumber
&\;\;\;\;+
     \frac{2}{9} \,
     \Re\bigg\{
    \frac{2 \, \LbarLbark \; \LbarLk + \KLbarLbark \; \barP}{(\Lbark)^2}
    \bigg\}  
\\ \nonumber
&\;\;\;\;-       
        \frac{2}{9} \, 
        \Re\bigg\{ \frac{2 \, \LLbark \; \barP + \LbarLbark \; \Paux}{\Lbark} 
        \bigg\}
\\ \nonumber
&\;\;\;\;-
        \frac{1}{9} \, |\Paux|^2
        +
        \frac{1}{3} \,
        \bigg\vert
        \frac{\LbarLbark}{\Lbark}
        \bigg\vert^2 
\\ \nonumber
&\;\;\;\;+
        \frac{2}{3}\, \Re\bigg\{ \mathcal{L}_1
            \Big[
            \frac{\LbarLbark}{\Lbark} 
            \Big] \bigg\}
        +
        \frac{1}{3}\, \Re\Big(\Lbar(\Paux)\Big)
\end{align}

\Section{\bf Caves Beneath a Waterfall}
\label{caves-beneath-waterfall}
\HEAD{{\ref{caves-beneath-waterfall}}.~{\sf Caves Beneath a Waterfall}
}{
Zhangchi {\sc Chen}, Wei Guo {\sc Foo}, Joël {\sc Merker}, 
The Anh {\sc Ta}}

This section displays the technique of calculating differential
invariants under infinite dimensional lie group action. First,
introduce some notations.

\Subsection{Finite dimensional approximations}

\begin{Definition} 
The rigid transformation group of $\bc^{2+1}$ fixing the origin is
denoted by:
\[
RT
\,:=\,
\big\{(z,\zeta,w)\mapsto(z',\zeta',w')=\big(f(z,\zeta),g(z,\zeta),\rho\,w\big)\big\},
\]
where $\rho\in\br^*$ and $f$, $g$ are holomorphic functions near $0\in\bc^2$ with $f(0,0)=g(0,0)=0$ and with invertible Jacobian
\[
\left(\!
\begin{array}{cc}
f_z & f_\zeta
\\
g_z & g_\zeta
\end{array}
\!\right).
\]
\end{Definition}

Multiplications and inversions are induced by compositions and inversions of transformations.\begin{Proposition} $(f,g)$ defines a biholomorphism between neighborhoods of $0\in\bc^2$ if and only if the jacobian matrix is invertible at $0$.
\end{Proposition}
\proof Let us explain only the existence of a {\em formal} inverse. Expand the holomorphic functions $f,g$ as
\begin{align*}
f(z,\zeta)&=\sum\limits_{n=1}^{\infty}\sum\limits_{j=0}^{n} \textstyle\frac{f_{j,n-j}}{j!\,(n-j)!}\,z^j\,\zeta^{n-j},\\
g(z,\zeta)&=\sum\limits_{n=1}^{\infty}\sum\limits_{j=0}^{n}\textstyle\frac{g_{j,n-j}}{j!\,(n-j)!}\,z^j\,\zeta^{n-j}.
\end{align*}
Let us construct progressively the formal inverse, which will be expanded as
\begin{align*}
\tilde{f}(z,\zeta)&=\sum\limits_{n=1}^{\infty}\sum\limits_{j=0}^{n}\textstyle\frac{\tilde{f}_{j,n-j}}{j!\,(n-j)!}\,z^j\,\zeta^{n-j},\\
\tilde{g}(z,\zeta)&=\sum\limits_{n=1}^{\infty}\sum\limits_{j=0}^{n}\textstyle\frac{\tilde{g}_{j,n-j}}{j!\,(n-j)!}\,z^j\,\zeta^{n-j}.
\end{align*}
Then
\begin{align*}
f\big(\tilde{f}(z,\zeta),\tilde{g}(z,\zeta)\big)&\equiv z,\\
g\big(\tilde{f}(z,\zeta),\tilde{g}(z,\zeta)\big)&\equiv \zeta.
\end{align*}
At each degree we get a linear system. For example at degree 1 we have
\[
\left(\!
\begin{array}{cc}
f_{1,0} & f_{0,1}
\\
g_{1,0} & g_{0,1}
\end{array}
\!\right)
\cdot
\left(\!
\begin{array}{cc}
\tilde{f}_{1,0} & \tilde{f}_{0,1}
\\
\tilde{g}_{1,0} & \tilde{g}_{0,1}
\end{array}
\!\right)
=
\left(\!
\begin{array}{cc}
1 & 0
\\
0 & 1
\end{array}
\!\right).
\]
Here $\tilde{f}_{1,0}, \tilde{f}_{0,1}, \tilde{g}_{1,0}, \tilde{g}_{0,1}$ can be uniquely solved thanks to the invertibility of the Jacobian of $(f,g)$.

Suppose by induction, for some $\delta\in\bz_{\geqslant 1}$ that all the coefficients $\tilde{f}_{j,k}$ and $\tilde{g}_{j,k}$ with $j+k\leqslant\delta$ have been already solved as rational functions of $f_{l,n-l}$ and $g_{l,n-l}$ with $n\leqslant \delta$. Then for $j+k=\delta+1$, we expand $f(\tilde{f},\tilde{g})$ and $g(\tilde{f},\tilde{g})$ to degree $\delta+1$ and compare the coefficients of $z^j\zeta^{\delta+1-j}$:
\begin{align*}
0&=\text{Coef}_{z^j\,\zeta^{\delta+1-j}}\big\{\sum\limits_{n=1}^{\delta+1}\sum\limits_{l=0}^{n}\textstyle\frac{f_{l,n-k}}{l!\,(n-l)!}\,\big(\tilde{f}(z,\zeta)\big)^j\,\big(\tilde{g}(z,\zeta)\big)^{n-l}\big\}\\
&=f_{1,0}\,\tilde{f}_{j,\delta+1-j}+f_{0,1}\,\tilde{g}_{j,\delta+1-j}+\text{Coef}_{z^j\,\zeta^{\delta+1-j}}\big\{\sum\limits_{n=2}^{\delta+1}\sum\limits_{l=0}^{n}\textstyle\frac{f_{l,n-k}}{l!\,(n-l)!}\,\big(\tilde{f}(z,\zeta)\big)^j\,\big(\tilde{g}(z,\zeta)\big)^{n-l}\big\},\\
0&=\text{Coef}_{z^j\,\zeta^{\delta+1-j}}\big\{\sum\limits_{n=1}^{\delta+1}\sum\limits_{l=0}^{n}\textstyle\frac{g_{l,n-l}}{l!\,(n-l)!}\,\big(\tilde{f}(z,\zeta)\big)^j\,\big(\tilde{g}(z,\zeta)\big)^{n-l}\big\}\\
&=g_{1,0}\,\tilde{f}_{j,\delta+1-j}+g_{0,1}\,\tilde{g}_{j,\delta+1-j}+\text{Coef}_{z^j\,\zeta^{\delta+1-j}}\big\{\sum\limits_{n=2}^{\delta+1}\sum\limits_{l=0}^{n}\textstyle\frac{g_{l,n-l}}{l!\,(n-l)!}\,\big(\tilde{f}(z,\zeta)\big)^j\,\big(\tilde{g}(z,\zeta)\big)^{n-l}\big\}.
\end{align*}
i.e.
\[
\left(\!
\begin{array}{cc}
f_{1,0} & f_{0,1}
\\
g_{1,0} & g_{0,1}
\end{array}
\!\right)
\cdot
\left(\!
\begin{array}{c}
\tilde{f}_{j,\delta+1-j}
\\
\tilde{g}_{j,\delta+1-j}
\end{array}
\!\right)
+
\left(\!
\begin{array}{c}
\mathscr{R}_1
\\
\mathscr{R}_2
\end{array}
\!\right)
=
\left(\!
\begin{array}{c}
0
\\
0
\end{array}
\!\right),
\]
where $\mathscr{R}_1$ and $\mathscr{R}_2$ are polynomials of $f_{l,n-l}, g_{l,n-l}$ with $n\leqslant \delta+1$ and $\tilde{f}_{p,q}, \tilde{g}_{p,q}$ with $p+q\leqslant \delta$. By inductive assumption $\tilde{f}_{p,q}, \tilde{g}_{p,q}$ are rational functions of $f_{l,n-l}, g_{l,n-l}$ with $n\leqslant \delta$. So $\mathscr{R}_1$ and $\mathscr{R}_2$ are rational functions of $f_{l,n-l}, g_{l,n-l}$ with $n\leqslant \delta+1$. We can solve $\tilde{f}_{j,\delta+1-j}$ and $\tilde{g}_{j,\delta+1-j}$ as rational functions of $f_{l,n-l}, g_{l,n-l}$ with $n\leqslant \delta+1$.\endproof

\begin{Definition} The space of all Levi-rank 1 and 2 non-degenerate CR graphed hypersurfaces passing by the origin in $\bc^3$ is denoted by
\[
\ch:=\big\{u:=Re(w)=F(z,\zeta,\overline{z},\overline{\zeta})\big\}
\]
where
\begin{itemize}
\item (real-valued analytic) $F$ is an analytic and real-valued function in a neigborhood of $(0,0)\in\bc^2$;
\item (passing by the origin) $F(0,0,0,0)=0$;
\item (no harmonic monomials) $\partial_z^a\partial_{\zeta}^bF(0,0,0,0)=0$, for any $a,b\geqslant 0$.
\item (Levi-rank 1) the matrix
\[
\left(\!
\begin{array}{cc}
F_{z\,\overline{z}} & F_{z\,\overline{\zeta}}
\\
F_{\zeta\,\overline{z}} & F_{\zeta\,\overline{\zeta}}
\end{array}
\!\right)
\]
has rank 1 everywhere;
\item (2-non degenerate) the matrix
\[
\left(\!
\begin{array}{cc}
F_{z\,\overline{z}} & F_{z\,\overline{\zeta}}
\\
F_{z\,z\,\overline{z}} & F_{z\,z\,\overline{\zeta}}
\end{array}
\!\right)
\]
is invertible at the origin.
\end{itemize}
\end{Definition}

There is a natural action of the group $RT$ on the space $\ch$: a graphed hypersurface $u=Re(w)=F(z,\zeta,\overline{z},\overline{\zeta})$ is transformed into another hypersurface $u'=Re(w')=F'(z',\zeta',\overline{z'},\overline{\zeta'})$. The expression of $F'$ is obtained by solving the fundamental equation
\[
F'\big(f(z,\zeta),g(z,\zeta),\overline{f(z,\zeta)},\overline{g(z,\zeta)}\big)=\rho\,F(z,\zeta,\overline{z},\overline{\zeta}).
\]
Indeed $F'(z,\zeta,\overline{z},\overline{\zeta})=\rho\,F\big(\tilde{f}(z,\zeta),\tilde{g}(z,\zeta),\overline{\tilde{f}(z,\zeta)},\overline{\tilde{g}(z,\zeta)}\big)$ where $(\tilde{f},\tilde{g})$ is the inverse of $(f,g)$. The inverse transformation brings convenience to obtain the explicit action.

Both the group $RT$ and the space $\ch$ are infinite-dimensional in the sense that they admet infinitely many linearly independent parameters.

For $RT$, any transformation is defined by $\rho\in\br^*$ and two holomorphic functions $f$, $g$ with expansions
\begin{align*}
f(z,\zeta)&=\sum\limits_{n=1}^{\infty}\sum\limits_{j=0}^{n} \textstyle\frac{f_{j,n-j}}{j!\,(n-j)!}\,z^j\,\zeta^{n-j},\\
g(z,\zeta)&=\sum\limits_{n=1}^{\infty}\sum\limits_{j=0}^{n}\textstyle\frac{g_{j,n-j}}{j!\,(n-j)!}\,z^j\,\zeta^{n-j}.
\end{align*}
where $f_{j,k}, g_{j,k}\in \bc$, $f_{1,0}\,g_{0,1}-f_{0,1}\,g_{1,0}\neq0$. The group $RT$ is hence parametrized by $f_{j,k}, g_{j,k}$ and $\rho$.

For $\ch$, any graphed hypersurface admets an expansion
\begin{align*}
u=F(z,\zeta,\overline{z},\overline{\zeta})=\sum\limits_{n=2}^{\infty}\sum\limits_{a+b+c+d=n}\textstyle\frac{F_{a,b,c,d}}{a!\,b!\,c!\,d!}\,z^a\,\zeta^b\,\overline{z}^c\,\overline{\zeta}^d,
\end{align*}
where $F_{a,b,c,d}\in\bc$, $F_{c,d,a,b}=\overline{F_{a,b,c,d}}$, $F_{a,b,0,0}=0$ and conditions of constant Levi-rank 1 and of 2-non degeneracy are satisfied. The space is hence parametrized by $F_{a,b,c,d}$.

But these infinite-dimensional objects have finite dimensional approximations. They can be truncated by degrees in expansions. Then they can be viewed as inverse or projective limits of those finite-dimensional truncations.

\begin{Definition} The $\delta^{th}$ residue group $Res_\delta$ is the subgroup of $RT$ with
\[
f(z,\zeta)=z+O(\delta), \ \ \ \ g(z,\zeta)=\zeta+O(\delta), \ \ \ \ \rho=1.
\]
\end{Definition}

\begin{Proposition} The group $Res_\delta$ is a normal subgroup of $RT$.\qed
\end{Proposition}

\begin{Definition} The $\delta^{th}$ approximation group $RT_{\delta}$ is the quotient group $RT/Res_{\delta+1}$. Each element has a representative
\begin{align*}
f(z,\zeta)&=\sum\limits_{n=1}^{\delta}\sum\limits_{j=0}^{n}\textstyle\frac{f_{j,n-j}}{j!\,(n-j)!}\,z^j\,\zeta^{n-j},\\
g(z,\zeta)&=\sum\limits_{n=1}^{\delta}\sum\limits_{j=0}^{n}\textstyle\frac{g_{j,n-j}}{j!\,(n-j)!}\,z^j\,\zeta^{n-j}.
\end{align*}
The group $RT_{\delta}$ is a finite dimensional Lie group parameterized by $\rho$ and $f_{j,n-j},g_{j,n-j}$ with $n\leqslant \delta$.

\begin{center}
\begin{tabular}{r|r|r|r|r|r|r|r|r}
$\delta$ & 1 & 2 & 3 & 4 & 5 & 6 & 7 & $\delta$ \\\hline
$\dim_{\br}RT_{\delta}$ & 9 & 21 & 37 & 57 & 81 & 109 & 141 & $2\,\delta^2+6\,\delta+1$
\end{tabular}
\end{center}

Its multiplication and inversion are obtained by dropping terms of degree $\geqslant\delta+1$ in the multiplication and inversion of $RT$.
\end{Definition}

\begin{Proposition} 
For any $\delta,\delta'\in\bz_+$ with $\delta>\delta'$ there is a projection $RT_{\delta}\rw RT_{\delta'}$ induced by the injection $Res_\delta\rw Res_{\delta'}$. For any $\delta,\delta',\delta''\in\bz_+$ with $\delta>\delta'>\delta''$ The following diagram commutes.
\[
\xymatrix{
RT_{\delta} 
\ar[r]
\ar[rd]  
& 
RT_{\delta'} 
\ar[d]
\\
& 
RT_{\delta''}.
}
\]
These projections define a projective system
$\{RT_{\delta}\}_{\delta\in\bz_+}$. Projections $\pi_\delta:RT\rw
RT_{\delta}$ are compatible with this system. By the universal
property of the projective limit, there is a morphism
\[
RT\rw
\underset{\underset{\delta}{\longleftarrow}}{\lim}\,
RT_{\delta}
\]
which is indeed an inclusion whose image consists of all convergent power series.
\end{Proposition}

\begin{Definition} For any $\delta\geqslant 2$, the $\delta^{th}$ approximation of $\ch$ is a manifold
\[
\ch_\delta:=\big\{u:=F(z,\zeta,\overline{z},\overline{\zeta})=\sum\limits_{n=2}^{\delta}\sum\limits_{a+b+c+d=n}\textstyle\frac{F_{a,b,c,d}}{a!\,b!\,c!\,d!}\,z^a\,\zeta^b\,\overline{z}^c\,\overline{\zeta}^d\big\},
\]
where
\begin{itemize}
\item (real-valued) $F_{a,b,c,d}=\overline{F_{c,d,a,b}}$ for any $a,b,c,d\geqslant0$ ;
\item (passing by the origin) $F_{0,0,0,0}=0$;
\item (no harmonic monomials) $F_{a,b,0,0}=F_{0,0,c,d}=0$ for any $a,b,c,d\geqslant0$.
\item (2-non-degenerate) the matrix
\[
\left(\!
\begin{array}{cc}
F_{1,0,1,0} & F_{1,0,0,1}
\\
F_{2,0,1,0} & F_{2,0,0,1}
\end{array}
\!\right)
\]
is invertible.
\item (Levi-rank 1 until degree $\delta$) $F_{1,0,1,0}$, $F_{1,0,0,1}=\overline{F_{0,1,1,0}}$ and $F_{0,1,0,1}$ are not simultaneously 0. The complex Hessian of $F(z,\zeta,\overline{z},\overline{\zeta})$ vanishes up to order $\delta-2$, i.e. $F_{z\,\overline{z}}\,F_{\zeta\,\overline{\zeta}}-F_{z\,\overline{\zeta}}\,F_{\zeta\,\overline{z}}=O(\delta-1)$.
\end{itemize}
\end{Definition}

The last condition may look strange, but it is reasonable, as shows the

\begin{Proposition}\label{prop-truncation} A polynomial $F(z,\zeta,\overline{z},\overline{\zeta})=\sum\limits_{n=2}^{\delta}\sum\limits_{a+b+c+d=n}\textstyle\frac{F_{a,b,c,d}}{a!\,b!\,c!\,d!}\,z^a\,\zeta^b\,\overline{z}^c\,\overline{\zeta}^d$ is a degree $\delta$ truncation of a formal power series $\tilde{F}(z,\zeta,\overline{z},\overline{\zeta})$ with $\tilde{F}_{z\,\overline{z}}\,\tilde{F}_{\zeta\,\overline{\zeta}}-\tilde{F}_{z\,\overline{\zeta}}\,\tilde{F}_{\zeta\,\overline{z}}=0$ if and only if $F_{z\,\overline{z}}\,F_{\zeta\,\overline{\zeta}}-F_{z\,\overline{\zeta}}\,F_{\zeta\,\overline{z}}=O(\delta-1)$.
\end{Proposition}

\proof (only if) When calculating the
complex Hessian of a power series
\[
\tilde{F}(z,\zeta,\overline{z},\overline{\zeta})=\sum\limits_{n=2}^{\infty}\sum\limits_{a+b+c+d=n}\textstyle\frac{\tilde{F}_{a,b,c,d}}{a!\,b!\,c!\,d!}\,z^a\,\zeta^b\,\overline{z}^c\,\overline{\zeta}^d,
\]
the $\delta-2$ degree terms of $\tilde{F}_{z\,\overline{z}}\,\tilde{F}_{\zeta\,\overline{\zeta}}-\tilde{F}_{z\,\overline{\zeta}}\,\tilde{F}_{\zeta\,\overline{z}}$ involve only coefficients $\tilde{F}_{a,b,c,d}$ with $a+b+c+d\leqslant\delta$. 

Let $F(z,\zeta,\overline{z},\overline{\zeta})$ be its degree $\delta$ truncation
\[
F(z,\zeta,\overline{z},\overline{\zeta}):=\sum\limits_{n=2}^{\delta}\sum\limits_{a+b+c+d=n}\textstyle\frac{\tilde{F}_{a,b,c,d}}{a!\,b!\,c!\,d!}\,z^a\,\zeta^b\,\overline{z}^c\,\overline{\zeta}^d.
\]
Then $F_{z\,\overline{z}}\,F_{\zeta\,\overline{\zeta}}-F_{z\,\overline{\zeta}}\,F_{\zeta\,\overline{z}}=\tilde{F}_{z\,\overline{z}}\,\tilde{F}_{\zeta\,\overline{\zeta}}-\tilde{F}_{z\,\overline{\zeta}}\,\tilde{F}_{\zeta\,\overline{z}}+O(\delta-1)=O(\delta-1)$.

To prove the (if) part, let us introduce dependent and independent coordinates. The manifolds $\ch$ and $\ch_\delta$ are covered by 3 open subsets: $\{F_{1,0,1,0}\neq0\}$, $\{F_{1,0,0,1}=\overline{F_{0,1,1,0}}\neq0\}$ and $\{F_{0,1,0,1}\neq0\}$. We only treat $F_{1,0,1,0}\neq0$ case because the other two cases can be transformed into this one by changes of coordinates $(z',\zeta')=(z+\zeta,z-\zeta)$ or $(z',\zeta')=(z,\zeta)$ preserving the Levi-rank.

When $F_{1,0,1,0}\neq0$ we have $F_{z,\overline{z}}\neq0$ in a neighborhood of the origin. The Levi-rank 1 condition is now equivalent to
\[
F_{\zeta\,\overline{\zeta}}\equiv \textstyle\frac{F_{z\,\overline{\zeta}}\,F_{\zeta\,\overline{z}}}{F_{z\,\overline{z}}}.
\]
By differentiating both sides, all terms $F_{z^a\,\zeta^b\,\overline{z}^c\,\overline{\zeta}^d}$ with $b\geqslant 1$ and $d\geqslant 1$ can be uniquely expressed as rational functions of $F_{z^{a'}\,\zeta^{b'}\,\overline{z}^{c'}}$ with $a'+b'+c'\leqslant a+b+c+d$ and $F_{z^{a''}\,\overline{z}^{c''}\,\overline{\zeta}^{d''}}$ with $a''+b''+c''\leqslant a+b+c+d$. Moreover, only powers of $F_{z\,\overline{z}}$ appears in the denominators. For example:
\begin{align*}
F_{z\,\zeta,\overline{\zeta}}\equiv \textstyle\frac{F_{z\,\zeta\,\overline{z}}\,F_{z\,\overline{\zeta}}}{F_{z\,\overline{z}}}+\textstyle\frac{F_{z^2\,\overline{\zeta}}\,F_{\zeta\,\overline{z}}}{F_{z\,\overline{z}}}-\textstyle\frac{F_{z^2\,\overline{z}}\,F_{z\,\overline{\zeta}}\,F_{\zeta\,\overline{z}}}{F_{z\,\overline{z}}^2}.
\end{align*}
Taking their values at the origin, the coefficients $F_{a,b,c,d}$ with $b\geqslant 1$ and $d\geqslant 1$ can be uniquely expressed as rational functions of $F_{a',b',c',0}$ with $a'+b'+c'\leqslant a+b+c+d$ and $F_{a'',0,c'',d''}$ with $a''+b''+c''\leqslant a+b+c+d$. Moreover, only powers of $F_{1,0,1,0}$ appear in the denominators. For example:

\begin{align*}
F_{1,1,0,1}=\textstyle\frac{F_{1,1,1,0}\,F_{1,0,0,1}}{F_{1,0,1,0}}+\textstyle\frac{F_{2,0,0,1}\,F_{0,1,1,0}}{F_{1,0,1,0}}-\textstyle\frac{F_{2,0,1,0}\,F_{1,0,0,1}\,F_{0,1,1,0}}{F_{1,0,1,0}^2}.
\end{align*}

\begin{Definition} The coefficient $F_{a,b,c,d}$ will be called {\sl dependent} if $b\geqslant 1$ and $d\geqslant 1$. Otherwise, it will be called {\sl independent}.
\end{Definition}

Elements in the open subset $\{F_{1,0,1,0}\neq0\}$ of $\ch$ and $\ch_\delta$ are uniquely determined by the independent coefficients $F_{a,b,c,d}$ with $b\,d=0$. Since $F$ is real-valued, i.e. $F_{c,d,a,b}=\overline{F_{a,b,c,d}}$, one has
\[
\dim_{\br}\ch_\delta=\#\big\{(a,b,c,d)|a+b\geqslant1,c+d\geqslant1,a+b+c+d\leqslant\delta,b\,d=0\big\}.
\]

\begin{center}
\begin{tabular}{r|r|r|r|r|r|r|r|r}
$\delta$ & 2 & 3 & 4 & 5 & 6 & 7 & 8 & $\delta$ \\\hline
$\dim_{\br}\ch_\delta$ & 3 & 11 & 26 & 50 & 85 & 133 & 196 & $\textstyle\frac{1}{6}(2\,\delta^3+3\,\delta^2-5\,\delta)$
\end{tabular}
\end{center}

To prove the (if) part of Proposition \ref{prop-truncation}, one shall construct a power series $\tilde{F}(z,\zeta,\overline{z},\overline{\zeta})=F(z,\zeta,\overline{z},\overline{\zeta})+\sum\limits_{n=\delta+1}^{\infty}\sum\limits_{a+b+c+d=n}\textstyle\frac{\tilde{F}_{a,b,c,d}}{a!b!c!d!}z^a\zeta^b\overline{z}^c\overline{\zeta}^d$ with $\tilde{F}_{z\,\overline{z}}\,\tilde{F}_{\zeta\,\overline{\zeta}}-\tilde{F}_{z\,\overline{\zeta}}\,\tilde{F}_{\zeta\,\overline{z}}=0$. This can be achieved by taking all the independent coefficients $\tilde{F}_{a,b,c,d}=0$  with $a+b+c+d\geqslant n+1$ and $b\,d=0$ and calculate all the dependent coefficients $\tilde{F}_{a,b,c,d}$ with $b\geqslant 1$ and $d\geqslant 1$ by their rational expressions of the independent ones. \endproof

\begin{Proposition} For any $\delta,\delta'\in\bz_+$ with $\delta>\delta'$ there is a projection $\ch_{\delta}\rw \ch_{\delta'}$ by dropping terms of degree $\geqslant\delta'+1$. For any $\delta,\delta',\delta''\in\bz_+$ with $\delta>\delta'>\delta''$ The following diagram commutes.
\[
\xymatrix{
\ch_{\delta} 
\ar[rd] 
\ar[r] 
&
\ch_{\delta'} 
\ar[d]
\\
& 
\ch_{\delta''}.
}
\]
These projections define a projective system $\{\ch_{\delta}\}_{\delta\in\bz_+}$. Projections $\pi_\delta:\ch\rw \ch_\delta$ are compatible with this system. By the universal property of the projective limit, there is a morphism
\[
\ch
\rw
\underset{\underset{\delta}{\longleftarrow}}{\lim}\,
\ch_\delta.
\]
which is indeed an inclusion.\end{Proposition}

The manifold $\ch_\delta$ is a finite-dimensional manifold parameterized by the independent coefficients $F_{a,b,c,d}$ with $a+b+c+d\leqslant\delta$ and $b\,d=0$. The action of the group $RT$ on $\ch$ induces an action on each manifold $\ch_\delta,\forall \delta\geqslant 0$:
\[
\xymatrix{
\ch 
\ar[r]^{\pi_\delta} 
\ar[d]_{(f,g,\rho)}
&
\ch_{\delta} 
\ar@{.>}[d]
\\
\ch 
\ar[r]^{\pi_\delta}
& 
\ch_{\delta}.
}
\]
More precisely, a polynomial $F(z,\zeta,\overline{z},\overline{\zeta})\in\ch_\delta$ is a degree $\delta$ truncation of a (not unique) convergent power series $\tilde{F}(z,\zeta,\overline{z},\overline{\zeta})\in\ch$, which is transformed to another convergent power series $\tilde{F}'(z,\zeta,\overline{z},\overline{\zeta})$ by the fundamental equation
\begin{align*}
\tilde{F}'(z,\zeta,\overline{z},\overline{\zeta})&=\rho\,\tilde{F}\big(\tilde{f}(z,\zeta),\tilde{g}(z,\zeta),\overline{\tilde{f}(z,\zeta)},\overline{\tilde{g}(z,\zeta)}\big)\\
&=\rho\,\sum\limits_{n=2}^{\delta}\sum\limits_{a+b+c+d=n}\textstyle\frac{F_{a,b,c,d}}{a!b!c!d!}\big(\tilde{f}(z,\zeta)\big)^a\,\big(\tilde{g}(z,\zeta)\big)^b\,\big(\overline{\tilde{f}(z,\zeta)}\big)^c\,\big(\overline{\tilde{g}(z,\zeta)}\big)^d+O(\delta+1).
\end{align*}
The degree $\delta$ truncation of $\tilde{F}'(z,\zeta,\overline{z},\overline{\zeta})$, denoted by $F'(z,\zeta,\overline{z},\overline{\zeta})$, is the image of $F(z,\zeta,\overline{z},\overline{\zeta})$ after the group action. It depends on the coefficients $F_{a,b,c,d}$ with $a+b+c+d\leqslant\delta$ only, hence is independent of the choice of $\tilde{F}(z,\zeta,\overline{z},\overline{\zeta})$. The group action is well-defined.

More precisely
\begin{Proposition} There is a group action of $RT_{\delta-1}$ on $\ch_\delta$. The group action of $RT$ on $\ch_\delta$ factors through the projection $\pi_{\delta-1}:RT\rw RT_{\delta-1}$, i.e. the following diagram commutes:
\[
\xymatrix{
RT\times\ch_\delta 
\ar[r] 
\ar[d] 
&
\ch_{\delta}
\\
RT_{\delta-1}\times\ch_\delta.
\ar[ru]
&
}
\]
\end{Proposition}

\proof When calculating the Taylor coefficients $F'_{a,b,c,d}$ in
\[
\tilde{F}'(z,\zeta,\overline{z},\overline{\zeta})=\sum\limits_{n=2}^{\delta}\textstyle\frac{F'_{a,b,c,d}}{a!b!c!d!}z^a\,\zeta^b\,\overline{z}^c\,\overline{\zeta}^d+O(\delta+1),
\]
we are calculating coefficients of $z^a\,\zeta^b\,\overline{z}^c\,\overline{\zeta}^d$ with $a+b+c+d\leqslant\delta$ from
\[\rho\,\sum\limits_{n=2}^{\delta}\sum\limits_{a+b+c+d=n}\textstyle\frac{F_{a,b,c,d}}{a!b!c!d!}\big(\tilde{f}(z,\zeta)\big)^a\,\big(\tilde{g}(z,\zeta)\big)^b\,\big(\overline{\tilde{f}(z,\zeta)}\big)^c\,\big(\overline{\tilde{g}(z,\zeta)}\big)^d.
\]
Each monomial is a product of at least 2 terms among $\{\tilde{f}(z,\zeta),\tilde{g}(z,\zeta),\overline{\tilde{f}(z,\zeta)},\overline{\tilde{g}(z,\zeta)}\}$. Each term
\begin{align*}
\tilde{f}(z,\zeta)&=\sum\limits_{n=1}^{\infty}\textstyle\frac{\tilde{f}_{j,n-j}}{j!(n-j)!}z^j\,\zeta^{n-j},\\
\tilde{g}(z,\zeta)&=\sum\limits_{n=1}^{\infty}\textstyle\frac{\tilde{g}_{j,n-j}}{j!(n-j)!}z^j\,\zeta^{n-j},
\end{align*}
as a power series of $z,\zeta$ or $\overline{z},\overline{\zeta}$, starts from degree 1. So only $\tilde{f}_{j,n-j}$, $\tilde{g}_{j,n-j}$ and their conjugations with $n\leqslant \delta-1$ contribute to $F'_{a,b,c,d}$ with $a+b+c+d\leqslant \delta$. The group action of $RT_{\delta-1}$ on $\ch_\delta$ can be well-defined and the commutative diagram is satisfied. \endproof

Compare the two tables of dimensions:
\begin{center}
\begin{tabular}{r|r|r|r|r|r|r|r}
$\delta$ & 2 & 3 & 4 & 5 & 6 & 7 & 8  \\\hline
$\dim_{\br}RT_{\delta-1}$ & 9 & 21 & 37 & 57 & 81 & 109 & 141  \\\hline
$\dim_{\br}\ch_\delta$ & 3 & 11 & 26 & 50 & 85 & 133 & 196 
\end{tabular}
\end{center}
The theory of differential invariants of finite-dimensional Lie group actions applies: the orbit dimension of $RT_{\delta-1}$ on $\ch_d$ is at most equal to $\dim_{\br}RT_{\delta-1}$ and the equality is achieved only when the action is locally free. We see immediately that the dimension of transversal, which equals to the number of linearly independent differential invariants up to order $\delta$, is positive when $\delta\geqslant 6$.

The infinite-dimensional Lie group $RT$ can be interpreted as an infinitely long flow of water. The space $\ch$ can be interpreted as an infinitely high valley. At the beginning, water fills the space up. But later on as the waterfall grows wider, water cannot fill the space. Some {\sl caves}, corresponding to the transversal dimension, or {\sl differential invariants}, show up.
\[
\xymatrix{
\ast \ar[d] \ar[rd] & & & &
\\
\ast \ar[d] & \ast \ar[rd] & & &
\\
\ast \ar[d] & \text{inv} & \ast \ar[d] \ar[rd] & &
\\
\ast \ar[d] & \text{inv} & \ast \ar[d] & \ast \ar[rd] &
\\
\ast & \text{inv} & \ast & \text{inv} &\ast
}
\]

\Section{\bf Invariants $I_0$, $V_0$, $Q_0$ at Every Point}
\label{invariants-I-0-V-0-Q-0}
\HEAD{{\ref{invariants-I-0-V-0-Q-0}}.~{\sf Invariants $I_0$, $V_0$, $Q_0$ 
at Every Point}
}{
Zhangchi {\sc Chen}, Wei Guo {\sc Foo}, Joël {\sc Merker}, 
The Anh {\sc Ta}}

Since the $RT$ action on $\ch_\delta$ factors through $\pi_{\delta-1}:RT\rw RT_{\delta-1}$, we have the
\begin{Proposition} A rational function on $\ch_\delta$ is invariant under the $RT$ action if and only if it is invariant under the $RT_{\delta}$ action.\qed
\end{Proposition}

Thus, to calculate differential invariants of order $\delta$ under $RT$ is equivalent to calculate those under the finite-dimensional Lie group $RT_{\delta-1}$. The algorithm goes as follows:
\begin{enumerate}
\item[(1)] Write down how $(f,g,\rho)\in RT_{\delta-1}$ acts on some independent parameters $F_{a,b,c,d}$.
\item[(2)] Choose certain $(f,g,\rho)\in RT_{\delta-1}$ to normalize as many independent parameters $F_{a,b,c,d}$ to 0 or 1 as possible, i.e. $(f,g,\rho)$ send $F_{a,b,c,d}$ to $F^{(1)}_{a,b,c,d}$ and some $F^{(1)}_{a,b,c,d}=0$ or $1$.
\item[(3)] Calculate how the other independent parameters $F^{(1)}_{a,b,c,d}$ are changed under this special $(f,g,\rho)$ action, i.e. express them as rational functions of $F_{a,b,c,d}$, $f_{j,n-j}$, $g_{j,n-j}$ and $\rho$.
\item[(4)] Calculate the "stabilizer", the subgroup $RT_{\delta-1}^{(1)}$ of $RT_{\delta-1}$ which preserves current normalizations.
\item[(5)] Repeat (2) (3) (4) by studying $RT_{\delta-1}^{(1)}$ actions on $F^{(1)}_{a,b,c,d}$, $RT_{\delta-1}^{(2)}$ actions on $F^{(2)}_{a,b,c,d}$ and so on, until no more terms can be normalized, i.e. $RT_{\delta-1}^{(k)}$ fixes all $F^{(k)}_{a,b,c,d}$.
\item[(6)] Express those non-constant $F^{(k)}_{a,b,c,d}$ in terms of $F_{a,b,c,d}$. They are rational functions fixed by $RT_{\delta-1}$, i.e. they are differential invariants of order $\leqslant\delta$.
\end{enumerate}

We fix $\delta=5$ in this section. The goal is to show the existence of order 5 invariants and to compute their explicit expressions.

\Subsection{First normalization: degree 2 terms $=z\,\overline{z}$}
We may assume that $F_{1,0,1,0}\neq 0$. In this case
\begin{align*}
F(z,\zeta,\overline{z},\overline{\zeta})&=F_{1,0,1,0}\,z\,\overline{z}+F_{1,0,0,1}\,z\overline{\zeta}+F_{0,1,1,0}\,\zeta\,\overline{z}+\textstyle\frac{F_{1,0,0,1}\,F_{0,1,1,0}}{F_{1,0,1,0}}\,\zeta\,\overline{\zeta}+O(3)\\
&=F_{1,0,1,0}\,\big(z+\textstyle\frac{F_{0,1,1,0}}{F_{1,0,1,0}}\,\zeta\big)\,\big(\overline{z}+\textstyle\frac{F_{1,0,0,1}}{F_{1,0,1,0}}\,\overline{\zeta}\big)+O(3)\\
&=\underbrace{
\big(F_{1,0,1,0}^{1/2}\,z+\textstyle\frac{F_{0,1,1,0}}{F_{1,0,1,0}^{1/2}}\,\zeta\big)}_{=:z'}
\,
\underbrace{
\big(F_{1,0,1,0}^{1/2}\,\overline{z}+\textstyle\frac{F_{1,0,0,1}}{F_{1,0,1,0}^{1/2}}\,\overline{\zeta}\big)}_{=:\overline{z'}}
\,+O(3).
\end{align*}
After the rigid transformation:
\[
z'=F_{1,0,1,0}^{1/2}\,z+\textstyle\frac{F_{0,1,1,0}}{F_{1,0,1,0}^{1/2}}\,\zeta,  \ \ \ \ \zeta'=\zeta, \ \ \ \ w'=w,
\]
the polynomial $F(z,\zeta,\overline{z},\overline{\zeta})$ becomes $F^{(1)}(z',\zeta',\overline{z'},\overline{\zeta'})=z'\,\overline{z'}+O(3)$. The other independent parameters $F^{(1)}_{a,b,c,d}$ with $a+b\geqslant 1,c+d\geqslant 1,b\,d=0$ can also be uniquely expressed as rational functions of $F_{a,b,c,d}$, by the fundamental equation.

Since all the independent parameters $F^{(1)}_{a,b,c,d}$ have $b\,d=0$ and $F^{(1)}_{c,d,a,b}=\overline{F^{(1)}_{a,b,c,d}}$, it suffices to calculate $F^{(1)}_{a,b,c,0}$ in terms of $F_{a,b,c,d}$. The inverse transformation is
\[
z=\textstyle\frac{1}{F_{1,0,1,0}^{1/2}}\,z'-\textstyle\frac{F_{0,1,1,0}}{F_{1,0,1,0}}\,\zeta',  \ \ \ \ \zeta=\zeta', \ \ \ \ w=w'.
\]
In the fundamental equality
\begin{align*}
\sum\limits_{a,b,c,d}\textstyle\frac{F^{(1)}_{a,b,c,d}}{a!b!c!d!}z'^a\,\zeta'^b\,\overline{z
}^c\,\overline{\zeta'}^d&=
\sum\limits_{a,b,c,d}\textstyle\frac{F_{a,b,c,d}}{a!b!c!d!}z^a\,\zeta^b\,\overline{z}^c\,\overline{\zeta}^d\\
&=\sum\limits_{a,b,c,d}\textstyle\frac{F_{a,b,c,d}}{a!b!c!d!}\big(\textstyle\frac{1}{F_{1,0,1,0}^{1/2}}\,z'-\textstyle\frac{F_{0,1,1,0}}{F_{1,0,1,0}}\,\zeta'\big)^a\,\zeta'^b\,\big(\textstyle\frac{1}{F_{1,0,1,0}^{1/2}}\,\overline{z'}-\textstyle\frac{F_{1,0,0,1}}{F_{1,0,1,0}}\,\overline{\zeta'}\big)^c\,\overline{\zeta'}^d,
\end{align*}
we calculate the coefficient of $z'^a\,\zeta'^b\,\overline{z'}^c$. On the left hand side, it is $F^{(1)}_{a,b,c,0}$. On the right hand side only $F_{j,a+b-j,c,0}$ with $a\leqslant j\leqslant a+b$ contribute. Since
\begin{align*}
&\textstyle\frac{F_{j,a+b-j,c,0}}{j!(a+b-j)!c!}\big(\textstyle\frac{1}{F_{1,0,1,0}^{1/2}}\,z'-\textstyle\frac{F_{0,1,1,0}}{F_{1,0,1,0}}\,\zeta'\big)^a\,\zeta'^{a+b-j}\,\big(\textstyle\frac{1}{F_{1,0,1,0}^{1/2}}\,\overline{z'}-\textstyle\frac{F_{1,0,0,1}}{F_{1,0,1,0}}\,\overline{\zeta'}\big)^c\\
=&\textstyle\frac{F_{j,a+b-j,c,0}}{j!(a+b-j)!c!}\textstyle\frac{j!}{a!(j-a)!}\big(\textstyle\frac{1}{F_{1,0,1,0}^{1/2}}\,z'\big)^a\,\big(-\textstyle\frac{F_{0,1,1,0}}{F_{1,0,1,0}}\,\zeta'\big)^{j-a}\,\zeta'^{a+b-j}\,\big(\textstyle\frac{1}{F_{1,0,1,0}^{1/2}}\,\overline{z'}\big)^c+\text{irrelevant monomials},
\end{align*}
We get
\begin{align*}
F^{(1)}_{a,b,c,0}&=\sum\limits_{j=a}^{a+b}\textstyle\frac{F_{j,a+b-j,c,0}}{a!(j-a)!(a+b-j)!c!}\big(\textstyle\frac{1}{F_{1,0,1,0}^{1/2}}\big)^a\,\big(-\textstyle\frac{F_{0,1,1,0}}{F_{1,0,1,0}}\big)^{j-a}\,\big(\textstyle\frac{1}{F_{1,0,1,0}^{1/2}}\big)^c\\
&=\sum\limits_{j=0}^{b}\textstyle\frac{F_{a+j,b-j,c,0}}{a!j!(b-j)!c!}\big(\textstyle\frac{1}{F_{1,0,1,0}^{1/2}}\big)^{a+c}\,\big(-\textstyle\frac{F_{0,1,1,0}}{F_{1,0,1,0}}\big)^{j}.
\end{align*}

We define $\ch_5^{(1)}:=\{u:=F^{(1)}(z,\zeta,\overline{z},\overline{\zeta})=z\,\overline{z}+O(3)\}$, a codimension 3 submanifold of $\ch_5$ since we have normalized $F^{(1)}_{1,0,1,0}=1$ and $F^{(1)}_{1,0,0,1}=\overline{F^{(1)}_{0,1,1,0}}=0$. So $\dim_\br \ch_5^{(1)}=50-3=47$.

Its stabilizer group $RT_{4}^{(1)}$ consists of $(f,g,\rho)$ such that
\[
f(z,\zeta)=r\,e^{i\theta}\,z+O(2), \ \ \ \  g(z,\zeta)=O(1), \ \ \ \ \rho=r^2,
\]
where $r\in\br_+$, $\theta\in[0,2\pi)$. It is a codimension 3 subgroup of $RT_{4}$, hence $\dim_\br RT_{4}^{(1)}=57-3=54$.

\Subsection{Second normalization: $F^{(2)}_{a,b,1,0}=0$ for $(a,b)\neq(1,0)$.}
Now, we study the group action of $RT_{4}^{(1)}$ on $\ch_5^{(1)}$. Any element in $\ch_5^{(1)}$ has expansion:
\begin{align*}
F^{(1)}(z,\zeta,\overline{z},\overline{\zeta})&=z\,\overline{z}+\overline{z}\,\big(\sum\limits_{2\leqslant  a+b\leqslant 4}\textstyle\frac{F^{(1)}_{a,b,1,0}}{a!b!}z^a\,\zeta^b\big)+z\,\big(\sum\limits_{2\leqslant a+b\leqslant 4}\textstyle\frac{\overline{F^{(1)}_{a,b,1,0}}}{a!b!}\overline{z}^c\,\overline{\zeta}^d\big)+R(z,\zeta,\overline{z},\overline{\zeta})\\
&=\underbrace{\big(z+\sum\limits_{2\leqslant  a+b\leqslant 4}\textstyle\frac{F^{(1)}_{a,b,1,0}}{a!b!}z^a\,\zeta^b\big)}_{=:z'}
\,
\underbrace{\big(\overline{z}+\sum\limits_{2\leqslant  a+b\leqslant 4}\textstyle\frac{\overline{F^{(1)}_{a,b,1,0}}}{a!b!}\overline{z}^a\,\overline{\zeta}^b\big)}_{=:\overline{z'}}+R(z,\zeta,\overline{z},\overline{\zeta})
\end{align*}
whose the remainder $R(z,\zeta,\overline{z},\overline{\zeta})$ contains only terms $z^a\,\zeta^b\,\overline{z}^c\,\overline{z'}^d$ with either $(a,b)$ or $(c,d)$ $\notin \{(0,0),(1,0)\}$. After the rigid transformation in $RT_{4}^{(1)}$:
\[
z'=z+\sum\limits_{2\leqslant  a+b\leqslant 4}\textstyle\frac{F^{(1)}_{a,b,1,0}}{a!b!}z^a\,\zeta^b,  \ \ \ \ \zeta'=\zeta, \ \ \ \ w'=w, \ \ \ \ (*)
\]
the polynomial $F^{(1)}(z,\zeta,\overline{z},\overline{\zeta})$ becomes $F^{(2)}(z',\zeta',\overline{z'},\overline{\zeta'})=z'\,\overline{z'}+R'(z',\zeta',\overline{z'},\overline{\zeta'})$. It remains to show that the remainder $R'(z',\zeta',\overline{z'},\overline{\zeta'})$ contains only terms $z^a\,\zeta^b\,\overline{z}^c\,\overline{z}^d$ with either $(a,b)$ or $(c,d)$ $\notin \{(0,0),(1,0)\}$.

\begin{Lemma} The inverse of $(*)$ in $RT_{4}^{(1)}$ is of the form
\[
z=z'+\sum\limits_{n=2}^4\sum\limits_{j=0}^n\textstyle\frac{\tilde{f}_{j,n-j}}{j!(n-j)!}z^j\,\zeta^{n-j},  \ \ \ \ \zeta=\zeta', \ \ \ \ w=w'.
\]
\end{Lemma}

\proof It suffices to show that $z:=\tilde{f}(z',\zeta')=z'+O_{z',\zeta'}(2)$. From $(*)$
\[
z=z'-\sum\limits_{2\leqslant  a+b\leqslant 4}\textstyle\frac{F^{(1)}_{a,b,1,0}}{a!b!}z^a\,\zeta^b=z'-\sum\limits_{2\leqslant  a+b\leqslant 4}\textstyle\frac{F^{(1)}_{a,b,1,0}}{a!b!}\tilde{f}(z',\zeta')^a\,\zeta'^b=z'+O_{z',\zeta'}(2).\qedhere
\]\endproof

In the remainder $R(z,\zeta,\overline{z},\overline{\zeta})$, each term $z^a\,\zeta^b\,\overline{z}^c\,\overline{z}^d$ is transformed to $\big(z'+O_{z',\zeta'}(2)\big)^a\,\zeta'^b\,\big(\overline{z'}+O_{\overline{z'},\overline{\zeta'}}(2)\big)\,\overline{\zeta'}^d$, whose expansion still contains only terms $z'^a\,\zeta'^b\,\overline{z'}^c\,\overline{\zeta'}^d$ with either $(a,b)$ or $(c,d)\notin\{(0,0),(1,0)\}$. 
\endproof

The terms $F^{(2)}_{a,b,c,0}$ such that $2\leqslant a+b+c\leqslant 5$, $(a,b),(c,0)\notin\{(0,0),(1,0)\}$ can be solved in terms of $F^{(1)}_{a,b,c,d}$:
\begin{align*}\,
F^{(2)}_{0, 1, 2, 0} &= \,\,F^{(1)}_{0, 1, 2, 0}, \\ F^{(2)}_{0, 1, 3, 0} &= \,\,-3\,F^{(1)}_{0, 1, 2, 0}\,F^{(1)}_{1, 0, 2, 0}+F^{(1)}_{0, 1, 3, 0}, \\ F^{(2)}_{0, 1, 4, 0} &= \,\,15\,F^{(1)}_{0, 1, 2, 0}\,(F^{(1)}_{1, 0, 2, 0})^2-4\,F^{(1)}_{0, 1, 2, 0}\,F^{(1)}_{1, 0, 3, 0}-6\,F^{(1)}_{0, 1, 3, 0}\,F^{(1)}_{1, 0, 2, 0}+F^{(1)}_{0, 1, 4, 0}, \\ F^{(2)}_{0, 2, 2, 0} &= \,\,-F^{(1)}_{0, 2, 1, 0}\,F^{(1)}_{1, 0, 2, 0}+F^{(1)}_{0, 2, 2, 0}, \\ F^{(2)}_{0, 2, 3, 0} &= \,\,3\,F^{(1)}_{0, 2, 1, 0}\,(F^{(1)}_{1, 0, 2, 0})^2-F^{(1)}_{0, 2, 1, 0}\,F^{(1)}_{1, 0, 3, 0}-3\,F^{(1)}_{0, 2, 2, 0}\,F^{(1)}_{1, 0, 2, 0}+F^{(1)}_{0, 2, 3, 0}, \\ F^{(2)}_{0, 3, 2, 0} &= \,\,3\,F^{(1)}_{0, 2, 1, 0}\,F^{(1)}_{1, 0, 2, 0}\,F^{(1)}_{1, 1, 1, 0}-3\,F^{(1)}_{0, 2, 1, 0}\,F^{(1)}_{1, 1, 2, 0}-F^{(1)}_{0, 3, 1, 0}\,F^{(1)}_{1, 0, 2, 0}+F^{(1)}_{0, 3, 2, 0}, \\ F^{(2)}_{1, 1, 2, 0} &= \,\,-F^{(1)}_{1, 0, 2, 0}\,F^{(1)}_{1, 1, 1, 0}+F^{(1)}_{1, 1, 2, 0}, \\ F^{(2)}_{1, 1, 3, 0} &= \,\,3\,(F^{(1)}_{1, 0, 2, 0})^2\,F^{(1)}_{1, 1, 1, 0}-3\,F^{(1)}_{1, 0, 2, 0}\,F^{(1)}_{1, 1, 2, 0}-F^{(1)}_{1, 0, 3, 0}\,F^{(1)}_{1, 1, 1, 0}+F^{(1)}_{1, 1, 3, 0}, \\ F^{(2)}_{1, 2, 2, 0} &= \,\,F^{(1)}_{0, 2, 1, 0}\,F^{(1)}_{1, 0, 2, 0}\,F^{(1)}_{2, 0, 1, 0}+2\,F^{(1)}_{1, 0, 2, 0}\,(F^{(1)}_{1, 1, 1, 0})^2-F^{(1)}_{0, 2, 1, 0}\,F^{(1)}_{2, 0, 2, 0}\\
&\,\ \ \ \ -F^{(1)}_{1, 0, 2, 0}\,F^{(1)}_{1, 2, 1, 0}-2\,F^{(1)}_{1, 1, 1, 0}\,F^{(1)}_{1, 1, 2, 0}+F^{(1)}_{1, 2, 2, 0}, \\ F^{(2)}_{2, 0, 2, 0} &= \,\,-F^{(1)}_{1, 0, 2, 0}\,F^{(1)}_{2, 0, 1, 0}+F^{(1)}_{2, 0, 2, 0}, \\ F^{(2)}_{2, 0, 3, 0} &= \,\,3\,(F^{(1)}_{1, 0, 2, 0})^2\,F^{(1)}_{2, 0, 1, 0}-3\,F^{(1)}_{1, 0, 2, 0}\,F^{(1)}_{2, 0, 2, 0}-F^{(1)}_{1, 0, 3, 0}\,F^{(1)}_{2, 0, 1, 0}+F^{(1)}_{2, 0, 3, 0}, \\ F^{(2)}_{2, 1, 2, 0} &= \,\,3\,F^{(1)}_{1, 0, 2, 0}\,F^{(1)}_{1, 1, 1, 0}\,F^{(1)}_{2, 0, 1, 0}-F^{(1)}_{1, 0, 2, 0}\,F^{(1)}_{2, 1, 1, 0}-2\,F^{(1)}_{1, 1, 1, 0}\,F^{(1)}_{2, 0, 2, 0}-F^{(1)}_{1, 1, 2, 0}\,F^{(1)}_{2, 0, 1, 0}+F^{(1)}_{2, 1, 2, 0}, \\ F^{(2)}_{3, 0, 2, 0} &= \,\,3\,F^{(1)}_{1, 0, 2, 0}\,(F^{(1)}_{2, 0, 1, 0})^2-F^{(1)}_{1, 0, 2, 0}\,F^{(1)}_{3, 0, 1, 0}-3\,F^{(1)}_{2, 0, 1, 0}\,F^{(1)}_{2, 0, 2, 0}+F^{(1)}_{3, 0, 2, 0}.
\end{align*}
We define $\ch_5^{(2)}:=\{u:=F^{(2)}(z,\zeta,\overline{z},\overline{\zeta})=z\,\overline{z}+O(3)| F^{(2)}_{a,b,1,0}=0, \forall (a,b)\neq(1,0)\}$, a codimension 24 submanifold of $\ch_5^{(1)}$. So $\dim_\br \ch_5^{(2)}=47-24=23$.

It will be a bit strange to talk about stabilizer group from this step. We in fact need to introduce a new definition of stabilizer. But after the final step, we will recover the stabilizer in the standard sense.
\begin{Definition} For any fixed element $F^{(2)}(z,\zeta,\overline{z},\overline{\zeta})\in\ch_5^{(2)}$, the subset of $RT^{(1)}_{0,4}$ consisting of elements $f,g,\rho$ which send $F^{(2)}(z,\zeta,\overline{z},\overline{\zeta})$ to another element in $\ch_5^{(2)}$, is defined as $RT^{(2)}_{0,4}(F^{(2)})$. It depends on the choice of the original element $F^{(2)}$.
\end{Definition}
The stabilizer $RT_{4}^{(2)}(F^{(2)})$ is a codimension 24 subgroup of $RT_{4}^{(1)}$ hence $\dim_\br RT_{4}^{(2)}(F^{(2)})=54-24=30$. It contains elements $(f,g,\rho)=\big(r\,e^{i\,\theta}\,z+O(2),g,r^2\big)\in RT_{4}^{(1)}$ such that
\begin{align*}
f_{2, 0} &= -r\,e^{i\,\theta}\,F^{(2)}_{2, 0, 0, 1}\,\overline{g_{1, 0}}\,\overline{g_{0,1}}^{-1},\\
f_{3, 0} &= -r\,e^{i\,\theta}\,F^{(2)}_{3, 0, 0, 1}\,\overline{g_{1, 0}}\,\overline{g_{0,1}}^{-1},\\
f_{4, 0} &= -r\,e^{i\,\theta}\,F^{(2)}_{4, 0, 0, 1}\,\overline{g_{1, 0}}\,\overline{g_{0,1}}^{-1},\\
f_{0, 2} &= 0, f_{1, 1} = 0,  f_{0, 3} = 0, f_{1, 2} = 0, f_{2, 1} = 0,  f_{0, 4} = 0,  f_{1, 3} = 0, f_{2, 2} = 0, f_{3, 1} = 0,
\end{align*}
which are in total 12 conditions on complex coefficients.

\Subsection{Third normalization: $F^{(3)}_{2,0,0,1}=\overline{F^{(3)}_{0,1,2,0}}=1$}
Any element in $\ch_5^{(2)}$ has expansion:
\[
F^{(2)}(z,\zeta,\overline{z},\overline{\zeta})=z\,\overline{z}+\textstyle\frac{F^{(2)}_{2,0,0,1}}{2}z^2\,\overline{\zeta}+\textstyle\frac{\overline{F^{(2)}_{2,0,0,1}}}{2}\overline{z}^2\,\zeta+O(4).
\]
By 2-non-degeneracy $F_{2,0,0,1}^{(2)}\neq0$. So after the rigid transformation:
\[
z'=z,  \ \ \ \  \zeta'=\overline{F^{(2)}_{2,0,0,1}}\,\zeta=F^{(2)}_{0,1,2,0}\,\zeta, \ \ \ \ w'=w,
\]
it becomes a graph $u=F^{(3)}(z,\zeta,\overline{z},\overline{\zeta})=z\,\overline{z}+\frac{1}{2}z^2\,\overline{\zeta}+\frac{1}{2}\overline{z}^2\,\zeta+O(4)$.

The relations are $F^{(3)}_{a,b,c,0}=F^{(2)}_{a,b,c,0}\,(F^{(2)}_{0,1,2,0})^{-b}$.

We define $\ch_5^{(3)}:=\{u:=F^{(3)}(z,\zeta,\overline{z},\overline{\zeta})=z\,\overline{z}+\frac{1}{2}z^2\,\overline{\zeta}+\frac{1}{2}\overline{z}^2\,\zeta+O(4)| F^{(3)}_{a,b,1,0}=0, \forall (a,b)\neq(1,0)\}$, a codimension 2 submanifold of $\ch_5^{(2)}$. So $\dim_\br \ch_5^{(3)}=23-2=21$.

For any fixed element $F^{(3)}\in\ch_5^{(3)}$, there exists some $F^{(2)}\in\ch_5^{(2)}$ whose third normalization is equal to $F^{(3)}$. For example, we can take $F^{(2)}=F^{(3)}$. The stabilizer $RT_{4}^{(3)}(F^{(3)})$ is a codimension 2 subgroup of $RT_{4}^{(2)}(F^{(3)})$. Hence $\dim_\br RT_{4}^{(3)}(F^{(3)})=30-2=28$. It contains elements $(f,g,\rho)\in RT_{4}^{(2)}(F^{(3)})$ satisfying $g_{0,1}=e^{2\,i\,\theta}$, i.e.
\begin{align*}
f(z,\zeta) &= r\,e^{i\theta}\,z-\textstyle\frac{1}{2}\,r\,e^{3\,i\,\theta}\,\overline{g_{1,0}}\,z^2-
\textstyle\frac{1}{6}\,r\,e^{3\,i\,\theta}\,F^{(3)}_{3,0,0,1}\,\overline{g_{1,0}}\,z^3
-\textstyle\frac{1}{24}\,r\,e^{3\,i\,\theta}\,F^{(3)}_{4,0,0,1}\,\overline{g_{1,0}}\,z^4\\
g(z,\zeta) &=g_{1,0}\,z+e^{2\,i\,\theta}\,\zeta +O(2), \ \ \ \ \rho=r^2.
\end{align*}

\Subsection{Fourth normalization: $F^{(4)}_{2,0,2,0}=0$}
Any element in $\ch_5^{(3)}$ has expansion:
\begin{align*}
F^{(3)}(z,\zeta,\overline{z},\overline{\zeta}) &=z\,\overline{z}+\textstyle\frac{1}{2}z^2\,\overline{\zeta}+\textstyle\frac{1}{2}\overline{z}^2\,\zeta+\textstyle\frac{1}{4}F^{(3)}_{2,0,2,0}z^2\,\overline{z}^2+R(z,\zeta,\overline{z},\overline{\zeta})\\
&=z\,\overline{z}+\textstyle\frac{1}{2}z^2\,
\underbrace{\big(\overline{\zeta}+\textstyle\frac{1}{4}F^{(3)}_{2,0,2,0}\,\overline{z}^2\big)}_{\overline{=:\zeta}'}
+\textstyle\frac{1}{2}\overline{z}^2\,
\underbrace{\big(\zeta+\textstyle\frac{1}{4}F^{(3)}_{2,0,2,0}\,z^2\big)}_{=:\zeta'}
+R(z,\zeta,\overline{z},\overline{\zeta}),
\end{align*}
whose remainder $R(z,\zeta,\overline{z},\overline{\zeta})=O(4)$ contains no $z^2\,\overline{z}^2$ term. After the rigid transformation in $RT_{4}^{(3)}$:
\[
z'=z,  \ \ \ \ \zeta'=\zeta+\textstyle\frac{1}{4}F^{(3)}_{2,0,2,0}\,z^2, \ \ \ \ w'=w, \ \ \ \ (**)
\]
the polynomial $F^{(3)}(z,\zeta,\overline{z},\overline{\zeta})$ becomes $F^{(4)}(z',\zeta',\overline{z'},\overline{\zeta'})=z'\,\overline{z'}+\textstyle\frac{1}{2}z'^2\,\overline{\zeta'}+\textstyle\frac{1}{2}\overline{z'}^2\,\zeta'+R'(z',\zeta',\overline{z'},\overline{\zeta'})$. The inverse of $(**)$ is
\[
z=z',  \ \ \ \ \zeta=\zeta'-\textstyle\frac{1}{4}F^{(3)}_{2,0,2,0}\,z'^2, \ \ \ \ w=w'.
\]
So $R'(z',\zeta',\overline{z'},\overline{\zeta'})=R\big(z',\zeta'-\textstyle\frac{1}{4}F^{(3)}_{2,0,2,0}\,z'^2,\overline{z'},\overline{\zeta'}-\textstyle\frac{1}{4}F^{(3)}_{2,0,2,0}\,\overline{z'}^2\big)=O(4)$ without $z'^2\,\overline{z'}^2$ term.

The relations are
\begin{align*}
F^{(4)}_{0,1,3,0} &= F^{(3)}_{0,1,3,0}, \ \ F^{(4)}_{0,2,2,0}= F^{(3)}_{0,2,2,0}, \ \ F^{(4)}_{1,1,2,0} = F^{(3)}_{1,1,2,0}, \ \ F^{(4)}_{0,1,4,0} = F^{(3)}_{0,1,4,0} ,\\
F^{(4)}_{0,2,3,0} &= F^{(3)}_{0,2,3,0}, \ \ F^{(4)}_{0,3,2,0} = F^{(3)}_{0,3,2,0}, \ \ F^{(4)}_{1,2,2,0} = F^{(3)}_{1,2,2,0}, \\
F^{(4)}_{2,1,2,0} &= -\textstyle\frac{1}{2}F^{(3)}_{0,2,2,0}\,F^{(3)}_{2,0,2,0}+F^{(3)}_{2,1,2,0}, \\
F^{(4)}_{3,0,2,0} &= -\textstyle\frac{3}{2}F^{(3)}_{1,1,2,0}\,F^{(3)}_{2,0,2,0}-\textstyle\frac{1}{2}F^{(3)}_{3,0,0,1}\,F^{(3)}_{2,0,2,0}+F^{(3)}_{3,0,2,0}, \\
F^{(4)}_{1,1,3,0} &= -\textstyle\frac{3}{2}F^{(3)}_{2,0,2,0}+F^{(3)}_{1,1,3,0}, \\
F^{(4)}_{2,0,3,0} &= -\textstyle\frac{1}{2}F^{(3)}_{0,1,3,0}\,F^{(3)}_{2,0,2,0}-\textstyle\frac{3}{2}F^{(3)}_{2,0,1,1}\,F^{(3)}_{2,0,2,0}+F^{(3)}_{2,0,3,0}.
\end{align*}

We define $\ch_5^{(4)}$, a codimension 1 submanifold of $\ch_5^{(3)}$ by requiring $F^{(4)}_{2,0,2,0}=0$. So $\dim_\br \ch_5^{(2)}=21-1=20$.

For any fixed element $F^{(4)}\in\ch_5^{(4)}$, the stabilizer $RT_{4}^{(4)}(F^{(4)})$ is a codimension 1 subgroup of some $RT_{4}^{(3)}(F^{(4)})$. Hence $\dim_\br RT_{4}^{(4)}(F^{(3)})=28-1=27$. It contains elements $(f,g,\rho)\in RT_{4}^{(3)}(F^{(4)})$ satisfying 
\begin{align*}
g_{2,0}&=
e^{-2\,i\,\theta}\,F^{(4)}_{0, 2, 2, 0}\,g_{1,0}^2
+e^{6\,i\,\theta}\,F^{(4)}_{2, 0, 0, 2}\,\overline{g_{1,0}}^2
-e^{-4\,i\,\theta}\,g_{0,2}\,g_{1,0}^2
-e^{8\,i\,\theta}\,\overline{g_{0,2}}\,\overline{g_{1,0}}^2\\
&\, \ \ \ -2\,F^{(4)}_{1, 1, 2, 0}\,g_{1,0}
-2\,e^{4\,i\,\theta}\,F^{(4)}_{2, 0, 1, 1}\,\overline{g_{1,0}}
+2\,e^{-2\,i\,\theta}\,g_{1,0}\,g_{1,1}
+2\,e^{6\,i\,\theta}\,\overline{g_{1,0}}\,\overline{g_{1,1}}\\
&\, \ \ \ +3\,e^{2\,i\,\theta}\,g_{1,0}\,\overline{g_{1,0}}
-e^{4\,i\,\theta}\,\overline{g_{2,0}}.
\end{align*}
In other words
\begin{align*}
{\sf Re}\big(e^{-2\,i\,\theta}\,g_{2, 0}\big)={\sf Re}\big\{
&-e^{-4\,i\,\theta}\,F^{(4)}_{0, 2, 2, 0}\,g_{1, 0}^2
-e^{-6\,i\,\theta}\,g_{0,2}\,g_{1,0}^2\\
&-2\,e^{-2\,i\,\theta}\,F^{(4)}_{1, 1, 2, 0}\,g_{1,0}
+2\,e^{-4\,i\,\theta}\,g_{1,0}\,g_{1,1}
+\textstyle\frac{3}{2}\,g_{1,0}\,\overline{g_{1,0}}\big\}.
\end{align*}

\Subsection{Fifth normalization: $F^{(5)}_{a,b,2,0}=0$ for $2\leqslant a+b\leqslant 3$ and $(a,b)\neq(2,0)$}
Any element in $\ch_5^{(4)}$ has expansion:
\begin{align*}
F^{(4)}(z,\zeta,\overline{z},\overline{\zeta}) =z\,\overline{z}+\textstyle\frac{1}{2}z^2\,
\underbrace{\big(\overline{\zeta}+\sum\limits_{2\leqslant a+b\leqslant 3}\textstyle\frac{F^{(4)}_{2,0,a,b}}{a!b!}\,\overline{z}^a\,\overline{\zeta}^b\big)}
_{=:\overline{\zeta'}}
+\textstyle\frac{1}{2}\overline{z}^2\,
\underbrace{\big(\zeta+\sum\limits_{2\leqslant a+b\leqslant 3}\textstyle\frac{F^{(4)}_{a,b,2,0}}{a!b!}\,z^a\,\zeta^b\big)}
_{=:\overline{\zeta}}
+R(z,\zeta,\overline{z},\overline{\zeta}),
\end{align*}
whose remainder $R(z,\zeta,\overline{z},\overline{\zeta})=O(4)$ contains no $z^a\,\zeta^b\,\overline{z}^2$ term for any $2\leqslant a+b\leqslant 3$. After the rigid transformation in $RT_{4}^{(4)}(F^{(4)})$:
\[
z'=z,  \ \ \ \ \zeta'=\zeta+\sum\limits_{2\leqslant a+b\leqslant 4}\textstyle\frac{F^{(4)}_{a,b,2,0}}{a!b!}\,z^a\,\zeta^b, \ \ \ \ w'=w, \ \ \ \ (***)
\]
the polynomial $F^{(4)}(z,\zeta,\overline{z},\overline{\zeta})$ becomes $F^{(5)}(z',\zeta',\overline{z'},\overline{\zeta'})=z'\,\overline{z'}+\textstyle\frac{1}{2}z'^2\,\overline{\zeta'}+\textstyle\frac{1}{2}\overline{z'}^2\,\zeta'+R'(z',\zeta',\overline{z'},\overline{\zeta'})$. The inverse of $(***)$ is
\[
z=z',  \ \ \ \ \zeta=\zeta'+O_{z',\zeta'}(2), \ \ \ \ w=w'.
\]
So $R'(z',\zeta',\overline{z'},\overline{\zeta'})=R\big(z',\zeta'+O_{z',\zeta'}(2),\overline{z'},\overline{\zeta'}+O_{\overline{z'},\overline{\zeta'}}(2)\big)=O(4)$ without $z'^a\,\zeta'^b\,\overline{z'}^2$ terms for any $2\leqslant a+b\leqslant 3$.

The relations are
\begin{align*}
F^{(5)}_{0, 1, 3, 0} &= F^{(4)}_{0, 1, 3, 0}, & F^{(5)}_{0, 1, 4, 0} &= F^{(4)}_{0, 1, 4, 0}, \\ F^{(5)}_{0, 2, 3, 0} &= -2\,F^{(4)}_{0, 1, 3, 0}\,F^{(4)}_{0, 2, 2, 0}+F^{(4)}_{0, 2, 3, 0}, & F^{(5)}_{1, 1, 3, 0} &= -2\,F^{(4)}_{0, 1, 3, 0}\,F^{(4)}_{1, 1, 2, 0}+F^{(4)}_{1, 1, 3, 0}.
\end{align*}

We define $\ch_5^{(5)}$ a codimension 12 submanifold of $\ch_5^{(4)}$ where $F^{(5)}_{a,b,2,0}=0$ for
\[
(a,b)\in\{(1,1),(0,2),(3,0),(2,1),(1,2),(0,3)\}.
\]
So $\dim_\br \ch_5^{(5)}=20-12=8$.

For any fixed element $F^{(5)}\in\ch_5^{(5)}$, the stabilizer $RT_{4}^{(5)}(F^{(5)})$ is a codimension 12 subgroup of some $RT_{4}^{(4)}(F^{(5)})$. Hence $\dim_\br RT_{4}^{(5)}(F^{(5)})=27-12=15$. It contains element $(f,g,\rho)\in RT_{4}^{(4)}(F^{(5)})$ satisfying
\begin{align*}
g_{0,2}&=0, \ \ \ \ g_{1,1}=-2\,e^{4\,i\,\theta}\,\overline{g_{1,0}}, \ \ \ \ g_{0,3}=0, \ \ \ \ g_{1,2}=0, \\
g_{2,1}&=2\,e^{6\,i\,\theta}\,\overline{g_{1,0}}^2-2\,e^{4\,i\,\theta}\,F^{(5)}_{3,0,0,1}\,\overline{g_{1,0}},\\
g_{3,0}&=
-5\,e^{2\,i\,\theta}\,F^{(5)}_{3, 0, 0, 1}\,\overline{g_{1,0}}\,g_{1,0}
+e^{6\,i\,\theta}\,F^{(5)}_{3, 0, 0, 2}\,\overline{g_{1,0}}^2
-2\,e^{4\,i\,\theta}\,F^{(5)}_{3, 0, 1, 1}\,\overline{g_{1,0}}
-e^{4\,i\,\theta}\,F^{(5)}_{3, 0, 0, 1}\,\overline{g_{2,0}}.
\end{align*}
Since $(f,g,\rho)\in RT_{4}^{(4)}(F^{(5)})$ we have
\[
{\sf Re}\big(e^{-2\,i\,\theta}\,g_{2, 0}\big)={\sf Re}\big(-\textstyle\frac{5}{2}\,g_{1,0}\,\overline{g_{1,0}}\big).
\]
Thus $e^{-2\,i\,\theta}\,g_{2, 0}=-\textstyle\frac{5}{2}\,g_{1,0}\,\overline{g_{1,0}}+i\,b_{2,0}$ for some $b_{2,0}\in\br$. So the last equation becomes
\begin{align*}
g_{3,0}=
-\textstyle\frac{5}{2}\,e^{2\,i\,\theta}\,F^{(5)}_{3, 0, 0, 1}\,\overline{g_{1,0}}\,g_{1,0}
+e^{6\,i\,\theta}\,F^{(5)}_{3, 0, 0, 2}\,\overline{g_{1,0}}^2
-2\,e^{4\,i\,\theta}\,F^{(5)}_{3, 0, 1, 1}\,\overline{g_{1,0}}
-i\,e^{2\,i\,\theta}\,F^{(5)}_{3, 0, 0, 1}\,b_{2,0}.
\end{align*}
The stabilizer $RT_{4}^{(5)}(F^{(5)})$ is parametrized by 3 real variables $b_{2,0},r,\theta$ and 6 complex variables $g_{1,0}, g_{j,4-j}$ for $0\leqslant j\leqslant 4$.

\Subsection{Final normalization: $F^{(6)}_{0,1,3,0}=0$ and ${\sf Im}(F^{(6)}_{1,1,3,0})=0$}
Any element in $\ch_5^{(5)}$ has expansion:
\begin{align*}
F^{(5)}(z,\zeta,\overline{z},\overline{\zeta}) &=z\,\overline{z}+\textstyle\frac{1}{2}\,z^2\,\overline{\zeta}+\textstyle\frac{1}{2}\,\zeta\,\overline{z}^2\\
&+\textstyle\frac{1}{6}\,F^{(5)}_{0, 1, 3, 0}\,\zeta\,\overline{z}^3
+\textstyle\frac{1}{6}\,F^{(5)}_{3, 0, 0, 1}\,z^3\,\overline{\zeta}\\
&+\textstyle\frac{1}{6}\,F^{(5)}_{1, 1, 3, 0}\,z\,\zeta\,\overline{z}^3
+\textstyle\frac{1}{6}\,F^{(5)}_{3, 0, 1, 1}\,z^3\,\overline{z}\,\overline{\zeta}
+\textstyle\frac{1}{24}\,F^{(5)}_{0, 1, 4, 0}\,\zeta\,\overline{z}^4
+\textstyle\frac{1}{24}\,F^{(5)}_{4, 0, 0, 1}\,z^4\,\overline{\zeta}\\
&+\textstyle\frac{1}{12}\,F^{(5)}_{0, 2, 3, 0}\,\zeta^2\,\overline{z}^3
+\textstyle\frac{1}{12}\,F^{(5)}_{3, 0, 0, 2}\,z^3\,\overline{\zeta}^2\\
&+\zeta\,\overline{\zeta}\,(\dots).
\end{align*}
We study how $g_{1,0}$ and $b_{2,0}$ act on this object, i.e. we consider an arbitrary $(f,g,\rho)\in RT_{4}^{(5)}(F^{(5)})$ with $r=1$ and $\theta=g_{j,4-j}=0$. They have the form
\begin{align*}
f(z,\zeta) &=z-\textstyle\frac{1}{2}\,\overline{g_{1,0}}\,z^2+O(3),\\
g(z,\zeta) &=g_{1,0}\,z+\zeta +\textstyle\frac{1}{2}\,(-\textstyle\frac{5}{2}\,g_{1,0}\,\overline{g_{1,0}}+i\,b_{2,0})\,z^2+O(3),\\
\rho&=1.
\end{align*}
This transformation sends $F^{(5)}$ to $F'^{(5)}\in\ch_5^{(5)}$ such that
\begin{align*}
F'^{(5)}_{3,0,0,1}&=F^{(5)}_{3,0,0,1}+3\,\overline{g_{1,0}},\\
F'^{(5)}_{3,0,1,1}&=F^{(5)}_{3,0,1,1}-3\,F^{(5)}_{3, 0, 0, 1}\,g_{1, 0}-F^{(5)}_{3, 0, 0, 2}\,\overline{g_{1, 0}}+\textstyle\frac{15}{2}\,g_{1, 0}\,\overline{g_{1, 0}}-3\,i\,b_{2,0}.
\end{align*}
So by a unique choice of $g_{1,0}$ and $b_{2,0}$, namely
\[
g_{1,0}=
-\textstyle\frac{1}{3}F^{(5)}_{0,1,3,0}, \ \ \ \
b_{2,0}=
\textstyle\frac{i}{18}\,\big(F^{(5)}_{0, 2, 3, 0}\,F^{(5)}_{0, 1, 3, 0}
-F^{(5)}_{3, 0, 0, 2}\,F^{(5)}_{3, 0, 0, 1}
+3\,F^{(5)}_{1, 1, 3, 0}
-3\,F^{(5)}_{3, 0, 1, 1}\big),
\]
we can normalize $F'^{(5)}_{3,0,0,1}$ to 0 and $F'^{(5)}_{3,0,1,1}$ to a real number. The polynomial $F^{(5)}(z,\zeta,\overline{z},\overline{\zeta})$ becomes
\begin{align*}
F^{(6)}(z',\zeta',\overline{z'},\overline{\zeta'}) &=z'\,\overline{z'}+\textstyle\frac{1}{2}\,z'^2\,\overline{\zeta'}+\textstyle\frac{1}{2}\,\zeta'\,\overline{z'}^2\\
&\, \ \ \ \ +\textstyle\frac{1}{6}\,F^{(6)}_{1, 1, 3, 0}\,z'\,\zeta'\,\overline{z'}^3
+\textstyle\frac{1}{6}\,F^{(6)}_{3, 0, 1, 1}\,z'^3\,\overline{z'}\,\overline{\zeta'}
+\textstyle\frac{1}{24}\,F^{(6)}_{0, 1, 4, 0}\,\zeta'\,\overline{z'}^4
+\textstyle\frac{1}{24}\,F^{(6)}_{4, 0, 0, 1}\,z'^4\,\overline{\zeta'}\\
&\, \ \ \ \ +\textstyle\frac{1}{12}\,F^{(6)}_{0, 2, 3, 0}\,\zeta'^2\,\overline{z'}^3
+\textstyle\frac{1}{12}\,F^{(6)}_{3, 0, 0, 2}\,z'^3\,\overline{\zeta'}^2\\
&\, \ \ \ \ +\zeta'\,\overline{\zeta'}\,(\dots)\\
&=z'\,\overline{z'}+\textstyle\frac{1}{2}\,z'^2\,\overline{\zeta'}+\textstyle\frac{1}{2}\,\zeta'\,\overline{z'}^2\\
&\, \ \ \ \ +\textstyle\frac{1}{6}\,Q_0\,z'\,\zeta'\,\overline{z'}^3
+\textstyle\frac{1}{6}\,Q_0\,z'^3\,\overline{z'}\,\overline{\zeta'}
+\textstyle\frac{1}{24}\,V_0\,\zeta'\,\overline{z'}^4
+\textstyle\frac{1}{24}\,\overline{V_0}\,z'^4\,\overline{\zeta'}\\
&\, \ \ \ \ +\textstyle\frac{1}{12}\,I_0\,\zeta'^2\,\overline{z'}^3
+\textstyle\frac{1}{12}\,\overline{I_0}\,z'^3\,\overline{\zeta'}^2\\
&\, \ \ \ \ +\zeta'\,\overline{\zeta'}\,(\dots),
\end{align*}
where $I_0:=F^{(6)}_{0, 2, 3, 0}\in\bc$, $V_0:=F^{(6)}_{0, 1,4, 0}\in\bc$, $Q_0:=F^{(6)}_{1, 1, 3, 0}\in\br$.

The relations are
\begin{align*}
I_0 &= F^{(5)}_{0, 2, 3, 0}+2\,F^{(5)}_{3, 0, 0, 1},\\
V_0 &= -\textstyle\frac{5}{3}\,(F^{(5)}_{0, 1, 3, 0})^2+F^{(5)}_{0, 1, 4, 0},\\
Q_0 &= \textstyle\frac{1}{6}\,F^{(5)}_{0, 2, 3, 0}\,F^{(5)}_{0, 1, 3, 0}+\textstyle\frac{1}{2}\,F^{(5)}_{3, 0, 0, 1}\,F^{(5)}_{0, 1, 3, 0}+\textstyle\frac{1}{6}\,F^{(5)}_{3, 0, 0, 2}\,F^{(5)}_{3, 0, 0, 1}+\textstyle\frac{1}{2}\,F^{(5)}_{1, 1, 3, 0}+\textstyle\frac{1}{2}\,F^{(5)}_{3, 0, 1, 1}.
\end{align*}

We define $\cn=\ch_5^{(6)}$ a codimension 3 submanifold of $\ch_5^{(5)}$ by requiring $F^{(6)}_{0,3,1,0}=0$ and $Im(F^{(6)}_{1,1,3,0})=0$.

For any fixed element $F^{(6)}\in\cn$, the stabilizer $RT_{4}^{(6)}(F^{(6)})$ is a codimension 3 subgroup of some $RT_{4}^{(5)}(F^{(6)})$. Hence $\dim_\br RT_{4}^{(6)}(F^{(6)})=15-3=12$. It contains elements $(f,g,\rho)\in RT_{4}^{(5)}(F^{(6)})$ of the form
\begin{align*}
f(z,\zeta)=r\,e^{i\,\theta}\,z, \ \ \ \
g(z,\zeta)=e^{2\,i\,\theta}\,s+O(4), \ \ \ \ \rho=r^2.
\end{align*}
This group sends $I_0,V_0,Q_0$ to $I'_0,V'_0,Q'_0$ with relations
\[
I'_0=r^{-1}\,e^{-i\,\theta}\,I_0, \ \ \ \ V'_0=r^{-2}\,e^{2\,i\,\theta}\,V_0, \ \ \ \ Q'_0=r^{-2}\,Q_0
\]
So if we ignore dilations and rotations $(z',\zeta',w')=(r\,e^{i\,\theta}\,z,e^{2\,i\,\theta}\,\zeta,r^2\,w)$, then $I_0$, $V_0$, $Q_0$ are invariants.

Each $F^{(t)}_{a,b,c,d}$ is a rational function of $F^{(t-1)}_{a',b',c',d'}$ for $t=5,4,3,2$ and each $F^{(1)}_{a,b,c,d}$ is a rational function of $F_{a',b',c',d'}$. By composing these rational functions, one can express $I_0$, $V_0$, $Q_0$ in terms of original coordinates $F_{a,b,c,d}$:
\begin{align*}
I_0&=\frac{\text{52 terms in degree 9}}{F_{1, 0, 1, 0}^{3/2}\,(F_{0, 1, 1, 0}\,F_{1, 0, 2, 0}-F_{0, 1, 2, 0}\,F_{1, 0, 1, 0})^3\,(F_{1, 0, 0, 1}\,F_{2, 0, 1, 0}-F_{1, 0, 1, 0}\,F_{2, 0, 0, 1})},\\
V_0&=\frac{\text{11 terms in degree 4}}{3\,F_{1, 0, 1, 0}\,(F_{0, 1, 1, 0}\,F_{1, 0, 2, 0}-F_{0, 1, 2, 0}\,F_{1, 0, 1, 0})^2},\\
Q_0&=\frac{\text{824 terms in degree 18}}{6\,F_{1, 0, 1, 0}^3\,(F_{0, 1, 1, 0}\,F_{1, 0, 2, 0}-F_{0, 1, 2, 0}\,F_{1, 0, 1, 0})^4\,(F_{1, 0, 0, 1}\,F_{2, 0, 1, 0}-F_{1, 0, 1, 0}\,F_{2, 0, 0, 1})^4}.
\end{align*}

The numerator of $I_0$ is
\[
\scriptstyle
\aligned
&
\, \ \ \ \ F_{0, 1, 1, 0}^3\,F_{1, 0, 0, 1}\,F_{1, 0, 1, 0}^2\,F_{1, 0, 2, 0}\,F_{2, 0, 1, 0}\,F_{2, 0, 3, 0}
-F_{0, 1, 1, 0}^3\,F_{1, 0, 0, 1}\,F_{1, 0, 1, 0}^2\,F_{1, 0, 3, 0}\,F_{2, 0, 1, 0}\,F_{2, 0, 2, 0}
\\
&
+2\,F_{0, 1, 1, 0}^3\,F_{1, 0, 0, 1}\,F_{1, 0, 1, 0}\,F_{1, 0, 2, 0}^3\,F_{3, 0, 1, 0}
-6\,F_{0, 1, 1, 0}^3\,F_{1, 0, 0, 1}\,F_{1, 0, 2, 0}^3\,F_{2, 0, 1, 0}^2
\\
&
-F_{0, 1, 1, 0}^3\,F_{1, 0, 1, 0}^3\,F_{1, 0, 2, 0}\,F_{2, 0, 0, 1}\,F_{2, 0, 3, 0}
+F_{0, 1, 1, 0}^3\,F_{1, 0, 1, 0}^3\,F_{1, 0, 3, 0}\,F_{2, 0, 0, 1}\,F_{2, 0, 2, 0}
\\
&
-2\,F_{0, 1, 1, 0}^3\,F_{1, 0, 1, 0}^2\,F_{1, 0, 2, 0}^3\,F_{3, 0, 0, 1}
+6\,F_{0, 1, 1, 0}^3\,F_{1, 0, 1, 0}\,F_{1, 0, 2, 0}^3\,F_{2, 0, 0, 1}\,F_{2, 0, 1, 0}
\\
&
-F_{0, 1, 1, 0}^2\,F_{0, 1, 2, 0}\,F_{1, 0, 0, 1}\,F_{1, 0, 1, 0}^3\,F_{2, 0, 1, 0}\,F_{2, 0, 3, 0}
-6\,F_{0, 1, 1, 0}^2\,F_{0, 1, 2, 0}\,F_{1, 0, 0, 1}\,F_{1, 0, 1, 0}^2\,F_{1, 0, 2, 0}^2\,F_{3, 0, 1, 0}
\endaligned
\]
\[
\scriptstyle
\aligned
&
+F_{0, 1, 1, 0}^2\,F_{0, 1, 2, 0}\,F_{1, 0, 0, 1}\,F_{1, 0, 1, 0}^2\,F_{1, 0, 3, 0}\,F_{2, 0, 1, 0}^2
+18\,F_{0, 1, 1, 0}^2\,F_{0, 1, 2, 0}\,F_{1, 0, 0, 1}\,F_{1, 0, 1, 0}\,F_{1, 0, 2, 0}^2\,F_{2, 0, 1, 0}^2
\\
&
+F_{0, 1, 1, 0}^2\,F_{0, 1, 2, 0}\,F_{1, 0, 1, 0}^4\,F_{2, 0, 0, 1}\,F_{2, 0, 3, 0}
+6\,F_{0, 1, 1, 0}^2\,F_{0, 1, 2, 0}\,F_{1, 0, 1, 0}^3\,F_{1, 0, 2, 0}^2\,F_{3, 0, 0, 1}
\\
&
-F_{0, 1, 1, 0}^2\,F_{0, 1, 2, 0}\,F_{1, 0, 1, 0}^3\,F_{1, 0, 3, 0}\,F_{2, 0, 0, 1}\,F_{2, 0, 1, 0}
-18\,F_{0, 1, 1, 0}^2\,F_{0, 1, 2, 0}\,F_{1, 0, 1, 0}^2\,F_{1, 0, 2, 0}^2\,F_{2, 0, 0, 1}\,F_{2, 0, 1, 0}
\\
&
+F_{0, 1, 1, 0}^2\,F_{0, 1, 3, 0}\,F_{1, 0, 0, 1}\,F_{1, 0, 1, 0}^3\,F_{2, 0, 1, 0}\,F_{2, 0, 2, 0}
-F_{0, 1, 1, 0}^2\,F_{0, 1, 3, 0}\,F_{1, 0, 0, 1}\,F_{1, 0, 1, 0}^2\,F_{1, 0, 2, 0}\,F_{2, 0, 1, 0}^2
\\
&
-F_{0, 1, 1, 0}^2\,F_{0, 1, 3, 0}\,F_{1, 0, 1, 0}^4\,F_{2, 0, 0, 1}\,F_{2, 0, 2, 0}
+F_{0, 1, 1, 0}^2\,F_{0, 1, 3, 0}\,F_{1, 0, 1, 0}^3\,F_{1, 0, 2, 0}\,F_{2, 0, 0, 1}\,F_{2, 0, 1, 0}
\endaligned
\]
\[
\scriptstyle
\aligned
&
-2\,F_{0, 1, 1, 0}^2\,F_{1, 0, 0, 1}\,F_{1, 0, 1, 0}^3\,F_{1, 0, 2, 0}\,F_{1, 1, 3, 0}\,F_{2, 0, 1, 0}
+2\,F_{0, 1, 1, 0}^2\,F_{1, 0, 0, 1}\,F_{1, 0, 1, 0}^3\,F_{1, 0, 3, 0}\,F_{1, 1, 2, 0}\,F_{2, 0, 1, 0}
\\
&
+2\,F_{0, 1, 1, 0}^2\,F_{1, 0, 1, 0}^4\,F_{1, 0, 2, 0}\,F_{1, 1, 3, 0}\,F_{2, 0, 0, 1}
-2\,F_{0, 1, 1, 0}^2\,F_{1, 0, 1, 0}^4\,F_{1, 0, 3, 0}\,F_{1, 1, 2, 0}\,F_{2, 0, 0, 1}
\\
&
+6\,F_{0, 1, 1, 0}\,F_{0, 1, 2, 0}^2\,F_{1, 0, 0, 1}\,F_{1, 0, 1, 0}^3\,F_{1, 0, 2, 0}\,F_{3, 0, 1, 0}
-18\,F_{0, 1, 1, 0}\,F_{0, 1, 2, 0}^2\,F_{1, 0, 0, 1}\,F_{1, 0, 1, 0}^2\,F_{1, 0, 2, 0}\,F_{2, 0, 1, 0}^2
\\
&
-6\,F_{0, 1, 1, 0}\,F_{0, 1, 2, 0}^2\,F_{1, 0, 1, 0}^4\,F_{1, 0, 2, 0}\,F_{3, 0, 0, 1}
+18\,F_{0, 1, 1, 0}\,F_{0, 1, 2, 0}^2\,F_{1, 0, 1, 0}^3\,F_{1, 0, 2, 0}\,F_{2, 0, 0, 1}\,F_{2, 0, 1, 0}
\\
&
+2\,F_{0, 1, 1, 0}\,F_{0, 1, 2, 0}\,F_{1, 0, 0, 1}\,F_{1, 0, 1, 0}^4\,F_{1, 1, 3, 0}\,F_{2, 0, 1, 0}
-2\,F_{0, 1, 1, 0}\,F_{0, 1, 2, 0}\,F_{1, 0, 0, 1}\,F_{1, 0, 1, 0}^3\,F_{1, 0, 3, 0}\,F_{1, 1, 1, 0}\,F_{2, 0, 1, 0}
\endaligned
\]
\[
\scriptstyle
\aligned
&
-2\,F_{0, 1, 1, 0}\,F_{0, 1, 2, 0}\,F_{1, 0, 1, 0}^5\,F_{1, 1, 3, 0}\,F_{2, 0, 0, 1}
+2\,F_{0, 1, 1, 0}\,F_{0, 1, 2, 0}\,F_{1, 0, 1, 0}^4\,F_{1, 0, 3, 0}\,F_{1, 1, 1, 0}\,F_{2, 0, 0, 1}
\\
&
-2\,F_{0, 1, 1, 0}\,F_{0, 1, 3, 0}\,F_{1, 0, 0, 1}\,F_{1, 0, 1, 0}^4\,F_{1, 1, 2, 0}\,F_{2, 0, 1, 0}
+2\,F_{0, 1, 1, 0}\,F_{0, 1, 3, 0}\,F_{1, 0, 0, 1}\,F_{1, 0, 1, 0}^3\,F_{1, 0, 2, 0}\,F_{1, 1, 1, 0}\,F_{2, 0, 1, 0}
\\
&
+2\,F_{0, 1, 1, 0}\,F_{0, 1, 3, 0}\,F_{1, 0, 1, 0}^5\,F_{1, 1, 2, 0}\,F_{2, 0, 0, 1}
-2\,F_{0, 1, 1, 0}\,F_{0, 1, 3, 0}\,F_{1, 0, 1, 0}^4\,F_{1, 0, 2, 0}\,F_{1, 1, 1, 0}\,F_{2, 0, 0, 1}
\\
&
-F_{0, 1, 1, 0}\,F_{0, 2, 2, 0}\,F_{1, 0, 0, 1}\,F_{1, 0, 1, 0}^4\,F_{1, 0, 3, 0}\,F_{2, 0, 1, 0}
+F_{0, 1, 1, 0}\,F_{0, 2, 2, 0}\,F_{1, 0, 1, 0}^5\,F_{1, 0, 3, 0}\,F_{2, 0, 0, 1}
\\
&
+F_{0, 1, 1, 0}\,F_{0, 2, 3, 0}\,F_{1, 0, 0, 1}\,F_{1, 0, 1, 0}^4\,F_{1, 0, 2, 0}\,F_{2, 0, 1, 0}
-F_{0, 1, 1, 0}\,F_{0, 2, 3, 0}\,F_{1, 0, 1, 0}^5\,F_{1, 0, 2, 0}\,F_{2, 0, 0, 1}
\endaligned
\]
\[
\scriptstyle
\aligned
&
-2\,F_{0, 1, 2, 0}^3\,F_{1, 0, 0, 1}\,F_{1, 0, 1, 0}^4\,F_{3, 0, 1, 0}
+6\,F_{0, 1, 2, 0}^3\,F_{1, 0, 0, 1}\,F_{1, 0, 1, 0}^3\,F_{2, 0, 1, 0}^2
\\
&
+2\,F_{0, 1, 2, 0}^3\,F_{1, 0, 1, 0}^5\,F_{3, 0, 0, 1}
-6\,F_{0, 1, 2, 0}^3\,F_{1, 0, 1, 0}^4\,F_{2, 0, 0, 1}\,F_{2, 0, 1, 0}
\\
&
+F_{0, 1, 2, 0}\,F_{0, 2, 1, 0}\,F_{1, 0, 0, 1}\,F_{1, 0, 1, 0}^4\,F_{1, 0, 3, 0}\,F_{2, 0, 1, 0}
-F_{0, 1, 2, 0}\,F_{0, 2, 1, 0}\,F_{1, 0, 1, 0}^5\,F_{1, 0, 3, 0}\,F_{2, 0, 0, 1}
\\
&
-F_{0, 1, 2, 0}\,F_{0, 2, 3, 0}\,F_{1, 0, 0, 1}\,F_{1, 0, 1, 0}^5\,F_{2, 0, 1, 0}
+F_{0, 1, 2, 0}\,F_{0, 2, 3, 0}\,F_{1, 0, 1, 0}^6\,F_{2, 0, 0, 1}
\\
&
-F_{0, 1, 3, 0}\,F_{0, 2, 1, 0}\,F_{1, 0, 0, 1}\,F_{1, 0, 1, 0}^4\,F_{1, 0, 2, 0}\,F_{2, 0, 1, 0}
+F_{0, 1, 3, 0}\,F_{0, 2, 1, 0}\,F_{1, 0, 1, 0}^5\,F_{1, 0, 2, 0}\,F_{2, 0, 0, 1}
\\
&
+F_{0, 1, 3, 0}\,F_{0, 2, 2, 0}\,F_{1, 0, 0, 1}\,F_{1, 0, 1, 0}^5\,F_{2, 0, 1, 0}
-F_{0, 1, 3, 0}\,F_{0, 2, 2, 0}\,F_{1, 0, 1, 0}^6\,F_{2, 0, 0, 1}.
\endaligned
\]

The numerator of $V_0$ is
\[
\aligned
&
\,\ \ \ \ 3\,F_{0, 1, 1, 0}^2\,F_{1, 0, 2, 0}\,F_{1, 0, 4, 0}
-5\,F_{0, 1, 1, 0}^2\,F_{1, 0, 3, 0}^2
-3\,F_{0, 1, 1, 0}\,F_{0, 1, 2, 0}\,F_{1, 0, 1, 0}\,F_{1, 0, 4, 0}
\\
&
+12\,F_{0, 1, 1, 0}\,F_{0, 1, 2, 0}\,F_{1, 0, 2, 0}\,F_{1, 0, 3, 0}
+10\,F_{0, 1, 1, 0}\,F_{0, 1, 3, 0}\,F_{1, 0, 1, 0}\,F_{1, 0, 3, 0}
-12\,F_{0, 1, 1, 0}\,F_{0, 1, 3, 0}\,F_{1, 0, 2, 0}^2
\\
&
-3\,F_{0, 1, 1, 0}\,F_{0, 1, 4, 0}\,F_{1, 0, 1, 0}\,F_{1, 0, 2, 0}
-12\,F_{0, 1, 2, 0}^2\,F_{1, 0, 1, 0}\,F_{1, 0, 3, 0}
+12\,F_{0, 1, 2, 0}\,F_{0, 1, 3, 0}\,F_{1, 0, 1, 0}\,F_{1, 0, 2, 0}
\\&
+3\,F_{0, 1, 2, 0}\,F_{0, 1, 4, 0}\,F_{1, 0, 1, 0}^2
-5\,F_{0, 1, 3, 0}^2\,F_{1, 0, 1, 0}^2.
\endaligned
\]

We define $\ch_5^{(6)}$ a codimension 3 submanifold of $\ch_5^{(5)}$ by requiring $F^{(6)}_{0,3,1,0}=0$ and $Im(F^{(6)}_{1,1,3,0})=0$.

For any fixed element $F^{(6)}\in\ch_5^{(6)}$, the stabilizer $RT_{4}^{(6)}(F^{(6)})$ is a codimension 3 subgroup of some $RT_{4}^{(5)}(F^{(6)})$. Hence $\dim_\br RT_{4}^{(6)}(F^{(6)})=15-3=12$. It contains elements $(f,g,\rho)\in RT_{4}^{(5)}(F^{(6)})$ of the form
\begin{align*}
f(z,\zeta)=r\,e^{i\,\theta}\,z, \ \ \ \
g(z,\zeta)=e^{2\,i\,\theta}\,\zeta+O(4), \ \ \ \ \rho=r^2.
\end{align*}

Note that this stabilizer group no longer depends on the choice of $F^{(6)}\in\ch_5^{(6)}$. We simply write it as $RT_{4}^{(6)}$.

\subsection{Passing to the infinite dimension}
After these six normalizations, we killed $f_{0,1}$ and $g_{1,0}$. It is a miracle that now we can work directly on the infinite dimensional objects. We define $\ch^{(7)}$ be the subspace of $\ch$ consisting of all power series $u=F^{(7)}(z,\zeta,\overline{z},\overline{\zeta})=\frac{F^{(7)}_{a,b,c,d}}{a!b!c!d!}z^a\,\zeta^b\,\overline{z}^c\,\overline{\zeta}^d$ such that
\begin{itemize}
\item $F^{(7)}_{a,b,1,0}=0$, $\forall (a,b)\neq (1,0)$; $F^{(7)}_{1,0,1,0}=1$;
\item $F^{(7)}_{a,b,2,0}=0$, $\forall (a,b)\neq (0,1)$; $F^{(7)}_{0,1,2,0}=1$;
\item $F^{(7)}_{3,0,0,1}=0$, $F^{(7)}_{3,0,1,1}=F^{(7)}_{1,1,3,0}$.
\end{itemize}
It is both infinitely-dimensional and infinitely-codimensional in $\ch$. But it has a finitely-dimensional stabilizer.

By definition, any element in $\ch^{(7)}$ has its degree 5 truncation in $\ch_5^{(6)}$.
 
\begin{Theorem} Any element $u=F(z,\zeta,\overline{z},\overline{\zeta})$ in $\ch$ can be sent to some element $u=F^{(7)}(z,\zeta,\overline{z},\overline{\zeta})$ in $\ch^{(7)}$ by some (but not unique) element in $RT$. The ambiguity can be controlled in the following sense: any element $(f,g,\rho)\in RT$ sending one element $F^{(7)}\in\ch^{(7)}$ to another $F'^{(7)}\in\ch^{(7)}$ has the form $f(z,\zeta)=r\,e^{i\,\theta}\,z$, $g(z,\zeta)=e^{2\,i\,\theta}\,\zeta$ and $\rho=r^2$.
\end{Theorem}

\proof One shall simply use the six normalizations above with a bit modification: in the second (killing $F_{a,b,1,0}$) and the fifth (killing $F_{a,b,2,0}$) normalization, we normalize for infinitely many $(a,b)$. More precisely, we start from $u=F(z,\zeta,\overline{z},\overline{\zeta})$ in $\ch$. After the six normalizations above we get $u=F^{(6)}(z,\zeta,\overline{z},\overline{\zeta})$ whose degree 5 truncation $\pi_5\big(F^{(6)}(z,\zeta,\overline{z},\overline{\zeta})\big)$ is in $\ch_5^{(6)}$, i.e.
\begin{itemize}
\item $F^{(6)}_{a,b,1,0}=0$, $\forall 2\leqslant a+b\leqslant 4$; $F^{(6)}_{1,0,1,0}=1$;
\item $F^{(6)}_{a,b,2,0}=0$, $\forall 2\leqslant a+b\leqslant 4$; $F^{(6)}_{0,1,2,0}=1$;
\item $F^{(6)}_{3,0,0,1}=0$, $F^{(6)}_{3,0,1,1}=F^{(6)}_{1,1,3,0}$.
\end{itemize}
Then we do 2 more normalizations. First
\[
z'=z+\sum\limits_{a+b\geqslant 5}\textstyle\frac{F^{(6)}_{a,b,1,0}}{a!b!}z^a\,\zeta^b,  \ \ \ \ \zeta'=\zeta, \ \ \ \ w'=w,
\]
gives us $u'=F'(z',\zeta',\overline{z'},\overline{\zeta'})$ with \begin{itemize}
\item $F'_{a,b,1,0}=0$, $\forall a+b\geqslant 2$; $F'_{1,0,1,0}=1$;
\item $F'_{a,b,2,0}=0$, $\forall 2\leqslant a+b\leqslant 4$; $F'_{0,1,2,0}=1$;
\item $F'_{3,0,0,1}=0$, $F'_{3,0,1,1}=F'_{1,1,3,0}$.
\end{itemize}

Then
\[
z''=z',  \ \ \ \ \zeta''=\zeta'+\sum\limits_{a+b\geqslant 5}\textstyle\frac{F'_{a,b,2,0}}{a!b!}\,z^a\,\zeta^b, \ \ \ \ w'=w, 
\]
gives us $u''=F''(z'',\zeta'',\overline{z''},\overline{\zeta''})$ with
\begin{itemize}
\item $F''_{a,b,1,0}=0$, $\forall a+b\geqslant 2$; $F''_{1,0,1,0}=1$;
\item $F''_{a,b,2,0}=0$, $\forall a+b\geqslant 2$; $F''_{0,1,2,0}=1$;
\item $F''_{3,0,0,1}=0$, $F''_{3,0,1,1}=F''_{1,1,3,0}$.
\end{itemize}
So $u''=F''(z'',\zeta'',\overline{z''},\overline{\zeta''})$ is in $\ch^{(7)}$. It is the form we want.

Now suppose that $(f,g,\rho)\in RT$ sends one element $F^{(7)}\in\ch^{(7)}$ to another $F'^{(7)}\in\ch^{(7)}$. In the truncated setting, $\pi_4(f,g,\rho)\in RT_4$ sends $\pi_5(F^{(7)})\in\ch_5^{(6)}$ to $\pi_5(F'^{(7)})\in\ch_5^{(6)}$. So the truncated action $\pi_4(f,g,\rho)$ should be in the stabilizer $RT_4^{(6)}$. That is to say
\begin{align*}
f(z,\zeta)=r\,e^{i\,\theta}\,z+O(5), \ \ \ \
g(z,\zeta)=e^{2\,i\,\theta}\,\zeta+O(4), \ \ \ \ \rho=r^2.
\end{align*}

Recall the fundamental equation
\[
\rho\,F^{(7)}(z,\zeta,\overline{z},\overline{\zeta})=F'^{(7)}\big(f(z,\zeta),g(z,\zeta),\overline{f(z,\zeta)},\overline{g(z,\zeta)}\big).
\]
When we compare the coefficients of $z^j\,\zeta^{n-j}\,\overline{z}$ for any $n\geqslant 2$ and $0\leqslant j\leqslant n$:
\begin{align*}
0&=\text{Coef}_{z^j\,\zeta^{n-j}\,\overline{z}}\big\{F'^{(7)}\big(f(z,\zeta),g(z,\zeta),\overline{f(z,\zeta)},\overline{g(z,\zeta)}\big)\big\}\\
&=\text{Coef}_{z^j\,\zeta^{n-j}\,\overline{z}}\big\{f(z,\zeta)\,\overline{f(z,\zeta)}\big\}+\text{Coef}_{z^j\,\zeta^{n-j}\,\overline{z}}\big\{\sum\limits_{c=0,d=1}(\dots)\,\overline{g(z,\zeta)}^d\big\}\\
& \ \ \ +\text{Coef}_{z^j\,\zeta^{n-j}\,\overline{z}}\big\{\sum\limits_{c+d\geqslant 2}(\dots)\,\overline{f(z,\zeta)}^c\,\overline{g(z,\zeta)}^d\big\}
\end{align*}
The last two terms are 0 because they only contain monomials with $\deg_{\overline{z}}=0$ or $\deg_{\overline{z}}+\deg_{\overline{\zeta}}\geqslant 2$. The first term gives us $0=r\,e^{-i\,\theta}\textstyle{\frac{f_{j,n-j}}{j!(n-j)!}}$. Hence $f(z,\zeta)=r\,e^{i\,\theta}\,z$.

When we compare the coefficients of $z^j\,\zeta^{n-j}\,\overline{z}^2$ for any $n\geqslant 2$ and $0\leqslant j\leqslant n$:

\begin{align*}
0&=\text{Coef}_{z^j\,\zeta^{n-j}\,\overline{z}^2}\big\{F'^{(7)}\big(f(z,\zeta),g(z,\zeta),\overline{f(z,\zeta)},\overline{g(z,\zeta)}\big)\big\}\\
&=\text{Coef}_{z^j\,\zeta^{n-j}\,\overline{z}^2}\big\{f(z,\zeta)\,\overline{f(z,\zeta)}\big\}
+\text{Coef}_{z^j\,\zeta^{n-j}\,\overline{z}^2}\big\{\sum\limits_{c=0,d=1}(\dots)\,\overline{g(z,\zeta)}\big\}\\
& \ \ \ +\text{Coef}_{z^j\,\zeta^{n-j}\,\overline{z}^2}\big\{\textstyle{\frac{1}{2}}g(z,\zeta)\,\overline{f(z,\zeta)}^2\big\}
+\text{Coef}_{z^j\,\zeta^{n-j}\,\overline{z}^2}\big\{\sum\limits_{c=1,d=1}(\dots)\,\overline{f(z,\zeta)}\,\overline{g(z,\zeta)}\big\}\\
& \ \ \ +\text{Coef}_{z^j\,\zeta^{n-j}\,\overline{z}^2}\big\{\sum\limits_{c=0,d=2}(\dots)\,\overline{g(z,\zeta)}^2\big\}
+\text{Coef}_{z^j\,\zeta^{n-j}\,\overline{z}}\big\{\sum\limits_{c+d\geqslant 3}(\dots)\,\overline{f(z,\zeta)}^c\,\overline{g(z,\zeta)}^d\big\}
\end{align*}
Each term, except the third, is 0. The third term gives us $0=\textstyle{\frac{1}{2}}r^2\,\textstyle{\frac{g_{j,n-j}}{j!(n-j)!}}$. Hence $g(z,\zeta)=e^{2\,i\,\theta}\,\zeta$.
\endproof

\subsection{Branches: $I_0\neq0$, $V_0\neq0$ and $I_0\equiv0\equiv V_0$}
To get a normal form under the full rigid transformation group, including rotations and dilations
\[
z'=r\,e^{i\theta}\,z, \ \ \ \ \zeta'=e^{2\,i\,\theta}\,\zeta, \ \ \ \ \rho=r^2,
\]
we should normalize $I_0$ or $V_0$. Such a rotation and a dilation would send $(I_0,V_0,Q_0)$ to $(I'_0,V'_0,Q'_0)$ with
\[
I'_0=r^{-1}\,e^{-i\,\theta}\,I_0, \ \ \ \ V'_0=r^{-2}\,e^{2\,i\,\theta}\,V_0, \ \ \ \ Q'_0=r^{-2}\,Q_0.
\]
We avoid the mixed type and focus on the 3 possible branches:
\begin{itemize}
\item $I_0\neq 0$;
\item $I_0\equiv 0 $ but $V_0\neq 0$;
\item $I_0\equiv 0\equiv V_0$.
\end{itemize}

\subsubsection{Branch $I_0\neq0$}
In this branch we can normalize $I_0$ to 1 by choose $r\,e^{i\,\theta}=I_0$. More precisely, for any surface in $\ch^{(7)}$ graphed by:
\begin{align*}
F^{(7)}(z,\zeta,\overline{z},\overline{\zeta})
&=z\,\overline{z}+\textstyle\frac{1}{2}\,z^2\,\overline{\zeta}+\textstyle\frac{1}{2}\,\zeta\,\overline{z}^2\\
&\, \ \ \ \ +\textstyle\frac{1}{6}\,Q_0\,z\,\zeta\,\overline{z}^3
+\textstyle\frac{1}{6}\,Q_0\,z^3\,\overline{z}\,\overline{\zeta}
+\textstyle\frac{1}{24}\,V_0\,\zeta\,\overline{z}^4
+\textstyle\frac{1}{24}\,\overline{V_0}\,z^4\,\overline{\zeta}\\
&\, \ \ \ \ +\textstyle\frac{1}{12}\,I_0\,\zeta^2\,\overline{z}^3
+\textstyle\frac{1}{12}\,\overline{I_0}\,z^3\,\overline{\zeta}^2\\
&\, \ \ \ \ +\zeta\,\overline{\zeta}\,(\dots)+O(6),
\end{align*}
where $I_0\neq 0$, after the transformation
\[
z'=I_0\,z, \ \ \ \ \zeta'=\frac{I_0^2}{|I_0|^2}\,\zeta, \ \ \ \ ,\rho=|I_0|^2,
\]
the polynomial $F^{(7)}(z,\zeta,\overline{z},\overline{\zeta})$ becomes
\begin{align*}
F^{(8,1)}(z',\zeta',\overline{z'},\overline{\zeta'})
&=z'\,\overline{z'}+\textstyle\frac{1}{2}\,z'^2\,\overline{\zeta'}+\textstyle\frac{1}{2}\,\zeta'\,\overline{z'}^2\\
&\, \ \ \ \ +\textstyle\frac{1}{6}\,invQ_0\,z'\,\zeta'\,\overline{z'}^3
+\textstyle\frac{1}{6}\,invQ_0\,z'^3\,\overline{z'}\,\overline{\zeta'}
+\textstyle\frac{1}{24}\,invV_0\,\zeta'\,\overline{z'}^4
+\textstyle\frac{1}{24}\,\overline{invV_0}\,z'^4\,\overline{\zeta'}\\
&\, \ \ \ \ +\textstyle\frac{1}{12}\,\zeta'^2\,\overline{z'}^3
+\textstyle\frac{1}{12}\,z'^3\,\overline{\zeta'}^2\\
&\, \ \ \ \ +\zeta'\,\overline{\zeta'}\,(\dots)+O(6),
\end{align*}
where
\[
invV_0=\frac{V_0}{\overline{I_0}^2}, \ \ \ \ invQ_0=\frac{Q_0}{|I_0|^2},       
\]

We define $\ch^{(8,1)}$ a codimension 2 submanifold of $\ch^{(7)}$ by requiring $I_0=1$.

For any fixed element $F^{(8,1)}\in\ch^{(8,1)}$, the stabilizer $RT^{(8,1)}$ is the identity.

\subsubsection{Branch $I_0\equiv0$ but $V_0\neq0$}
In this branch we can normalize $V_0$ to 1 by choose $r^2\,e^{-2\,i\,\theta}=V_0$. This equation has two solutions: $r\,e^{i\,\theta}=\pm x$, where $x^2=\overline{V_0}$ and $\arg(x)\in[0,\pi)$. More precisely, for any surface in $\ch^{(7)}$ graphed by:
\begin{align*}
F^{(7)}(z,\zeta,\overline{z},\overline{\zeta})
&=z\,\overline{z}+\textstyle\frac{1}{2}\,z^2\,\overline{\zeta}+\textstyle\frac{1}{2}\,\zeta\,\overline{z}^2\\
&\, \ \ \ \ +\textstyle\frac{1}{6}\,Q_0\,z\,\zeta\,\overline{z}^3
+\textstyle\frac{1}{6}\,Q_0\,z^3\,\overline{z}\,\overline{\zeta}
+\textstyle\frac{1}{24}\,V_0\,\zeta\,\overline{z}^4
+\textstyle\frac{1}{24}\,\overline{V_0}\,z^4\,\overline{\zeta}\\
&\, \ \ \ \ +\underbrace{\textstyle\frac{1}{12}\,I_0\,\zeta^2\,\overline{z}^3
+\textstyle\frac{1}{12}\,\overline{I_0}\,z^3\,\overline{\zeta}^2}_{=0,\text{ when $I_0\equiv 0$}}\\
&\, \ \ \ \ +\zeta\,\overline{\zeta}\,(\dots)+O(6),
\end{align*}
where $V_0\neq 0$, after the transformation
\[
z'=x\,z, \ \ \ \ \zeta'=\frac{\overline{V_0}}{|V_0|}\,\zeta, \ \ \ \ ,\rho=|V_0|,
\]
the polynomial $F^{(7)}(z,\zeta,\overline{z},\overline{\zeta})$ becomes
\begin{align*}
F^{(8,2)}(z',\zeta',\overline{z'},\overline{\zeta'})
&=z'\,\overline{z'}+\textstyle\frac{1}{2}\,z'^2\,\overline{\zeta'}+\textstyle\frac{1}{2}\,\zeta'\,\overline{z'}^2+\textstyle\frac{1}{6}\,invQ_0\,z'\,\zeta'\,\overline{z'}^3
+\textstyle\frac{1}{6}\,invQ_0\,z'^3\,\overline{z'}\,\overline{\zeta'}\\
& \ \ \ 
+\textstyle\frac{1}{24}\,\zeta'\,\overline{z'}^4
+\textstyle\frac{1}{24}\,z'^4\,\overline{\zeta'}+\zeta'\,\overline{\zeta'}\,(\dots)+O(6),
\end{align*}
where $invQ_0=\frac{Q_0}{|V_0|}$. We define $\ch^{(8,2)}$ a codimension 2 submanifold of $\ch^{(7)}$ by requiring $V_0=1$. For any fixed element $F^{(8,2)}\in\ch^{(8,2)}$, the stabilizer $RT^{(8,2)}$ is a group of two elements: the identity and $(-z,\zeta,1)$.

\subsubsection{Branch $I_0\equiv 0\equiv V_0$}
Since $Q_0$ can be generated by $I_0$, $V_0$ and their differentials, we have $Q_0\equiv 0$. The structure equations degenerate to the model case. The surface is equivalent as the Gaussier-Merker model $u=\frac{z\,\overline{z}+\frac{1}{2}\zeta^2\,\overline{z}+\frac{1}{2}z^2\,\overline{\zeta}}{1-\zeta\,\overline{\zeta}}$.

To conclude, we draw the branches from our root assumption.
\[
\xymatrix{
&&
I_0\,\neq\,0
&&
V_0\,\neq\,0
&&\\
\ar[urr]
\boxed{
F_{z\,\overline{z}}\neq 0\equiv F_{z\,\overline{z}}\,F_{\zeta\,\overline{\zeta}}-F_{\zeta,\overline{z}}\,F_{z\,\overline{\zeta}}}
\ar[rr]
&&
\ar[urr]
I_0\,\equiv\,0
\ar[rr]
&&
V_0\,\equiv\,0
&&
}
\]
where $I_0$ and $V_0$ are relative invariants of order 5.

\begin{Theorem}
Within the branch $I_0 \equiv 0$:
\begin{enumerate}
\item
When $V_0 \equiv 0$, the surface is equivalent to the Gaussier-Merker model $u=\frac{z\,\overline{z}+\frac{1}{2}\zeta^2\,\overline{z}+\frac{1}{2}z^2\,\overline{\zeta}}{1-\zeta\,\overline{\zeta}}$, and conversely;
\item When $V_0 \neq 0$, the surface is, up to $z\mapsto -z$, equivalent to:
\begin{align*}
u&=z\,\overline{z}+\textstyle\frac{1}{2}\,z^2\,\overline{\zeta}+\textstyle\frac{1}{2}\,\zeta\,\overline{z}^2+\textstyle\frac{1}{6}\,\frac{Q_0}{|V_0|}\,z\,\zeta\,\overline{z}^3
+\textstyle\frac{1}{6}\,\frac{Q_0}{|V_0|}\,z^3\,\overline{z}\,\overline{\zeta}
+\textstyle\frac{1}{24}\,\zeta\,\overline{z}^4
+\textstyle\frac{1}{24}\,z^4\,\overline{\zeta}\\
& \ \ \ +\zeta\,\overline{\zeta}\,(\dots)+\textstyle{\sum\limits_{a+b+c+d\geqslant 6,\,b\,d=0}\frac{F_{a,b,c,d}}{a!b!c!d!}}z^a\,\zeta^b\,\overline{z}^c\,\overline{\zeta}^d,
\end{align*}
without any harmonic monomial $z^j\,\zeta^{n-j}$, $\forall n\geqslant 0$, $0\leqslant j\leqslant n$ and any monomial $z^a\,\zeta^b\,\overline{z}^c$, $\forall a+b\geqslant 2$, $c\in\{1,2\}$. Pairs of collection of coefficients:
\[
\frac{Q_0}{|V_0|}, \big\{ F_{a,b,c,d} \big\}_{a+b+c+d\geqslant 6,\,b\,d=0}, \ \ \ \ \ \ \ \
\frac{Q_0}{|V_0|}, \big\{ (-1)^{a+c}\,F_{a,b,c,d} \big\}_{a+b+c+d\geqslant 6,\,b\,d=0}
\]
are in one-to-one correspondence with equivalent classes.
\end{enumerate}
Within the branch $I_0 \neq 0$, the surface is, in a unique way, equivalent to:
\begin{align*}
u&=z\,\overline{z}+\textstyle\frac{1}{2}\,z^2\,\overline{\zeta}+\textstyle\frac{1}{2}\,\zeta\,\overline{z}^2+\textstyle\frac{1}{6}\,\frac{Q_0}{|I_0|^2}\,z\,\zeta\,\overline{z}^3
+\textstyle\frac{1}{6}\,\frac{Q_0}{|I_0|^2}\,z^3\,\overline{z}\,\overline{\zeta}
+\textstyle\frac{1}{24}\,\frac{V_0}{\overline{I_0}^2}\,\zeta\,\overline{z}^4
+\textstyle\frac{1}{24}\,\frac{\overline{V_0}}{I_0^2}\,z^4\,\overline{\zeta}+\textstyle\frac{1}{12}\,\zeta^2\,\overline{z}^3
+\textstyle\frac{1}{12}\,z^3\,\overline{\zeta}^2\\
& \ \ \ +\zeta\,\overline{\zeta}\,(\dots)+\textstyle{\sum\limits_{a+b+c+d\geqslant 6,\,b\,d=0}\frac{F_{a,b,c,d}}{a!b!c!d!}}z^a\,\zeta^b\,\overline{z}^c\,\overline{\zeta}^d,
\end{align*}
without any harmonic monomial $z^j\,\zeta^{n-j}$, $\forall n\geqslant 0$, $0\leqslant j\leqslant n$ and any monomial $z^a\,\zeta^b\,\overline{z}^c$, $\forall a+b\geqslant 2$, $c\in\{1,2\}$. Collections of coefficients: $\frac{V_0}{\overline{I_0}^2}$, $\frac{Q_0}{|I_0|^2}$ and $\big\{ F_{a,b,c,d} \big\}_{a+b+c+d\geqslant 6,\,b\,d=0}$, are in one-to-one correspondence with equivalent classes.
\end{Theorem}



\vfill
\begin{thebibliography}{XL}

{\scriptsize

{\bf\bibitem{Chen-Merker-2019}
{\rm Chen}}, Z.; {\rm Merker}, J.:
{\em On differential invariants of parabolic surfaces},
{\tiny\sf arxiv.org/abs/1908.07867/},
92 pages.

\smallskip

{\bf\bibitem{Doubrov-Komrakov-Rabinovich-1996}
{\rm Doubrov}}, B.; Komrakov, B.; Rabinovich, M.: 
{\em Homogeneous surfaces in the three-dimensional affine geometry}. 
Geometry and topology of submanifolds, 
VIII (Brussels, 1995/Nordfjordeid, 1995), 
168--178, World Sci. Publ., River Edge, NJ, 1996.

\smallskip

{\bf\bibitem{Eastwood-Ezhov-1999}
{\rm Eastwood}}, M.; {\rm Ezhov}, V.:
{\em On affine normal forms and a classification 
of homogeneous surfaces in affine three-space}, 
Geom. Dedicata {\bf 77} (1999), no.~1, 11--69.

\smallskip

{\bf\bibitem{Fels-Kaup-2007}
{\rm Fels}}, M.; {\rm Kaup}, W.:
{\em CR manifolds of dimension $5$: a Lie algebra approach},
J. Reine Angew. Math. {\bf 604} (2007), 47--71.

\smallskip

{\bf\bibitem{Fels-Kaup-2008}
{\rm Fels}}, M.; {\rm Kaup}, W.:
{\em Classification of Levi degenerate homogeneous CR-manifolds in 
dimension $5$},
Acta Math. {\bf 201} (2008), 1--82.

\smallskip

{\bf\bibitem{Foo-Merker-2019}
{\rm Foo}}, W.G.; {\rm Merker}, J.:
{\em Differential $\{e\}$-structures for equivalences
of $2$-nondegenerate Levi rank $1$ hypersurfaces $M^5 \subset \C^3$},
{\tiny\sf arxiv.org/abs/1901.02028/}, 72 pages.

\smallskip

{\bf\bibitem{Foo-Merker-Ta-2019}
{\rm Foo}}, W.G.; {\rm Merker}, J.; {\rm Ta}, T.-A.:
{\em Rigid equivalences of $5$-dimensional $2$-nondegenerate
rigid real hypersurfaces $M^{5}\subset\mathbb{C}^{3}$ 
of constant Levi rank $1$},
{\tiny\sf arxiv.org/abs/1904.02562/},
31 pages.

\smallskip

{\bf\bibitem{Gaussier-Merker-2003}
{\rm Gaussier}}, H.; {\rm Merker}, J.:
{\em A new example of uniformly Levi degenerate hypersurface in 
$\C^3$}, Ark. Mat. {\bf 41} (2003), no.~1, 85--94.
Erratum: {\bf 45} (2007), no.~2, 269--271.

\smallskip

{\bf\bibitem{Isaev-2006}
{\rm Isaev}}, A.:
{\em Analogues of Rossi's map and E. Cartan's classification 
of homogeneous strongly pseudoconvex 3-dimensional hypersurfaces},
J. Lie Theory {\bf 16} (2006), no.~3, 407--426. 

\smallskip

{\bf\bibitem{Isaev-2011}
{\rm Isaev}}, A.:
{\em Spherical tube hypersurfaces},
Lecture Notes in Mathematics, 2020, 
Springer, Heidelberg, 2011, xii+220~pp.

\smallskip

{\bf\bibitem{Isaev-2016}
{\rm Isaev}}, A.:
{\em Affine rigidity of Levi degenerate tube hypersurfaces},
J. Differential Geom. {\bf 104} (2016), no.~1, 111--141.

\smallskip

{\bf\bibitem{Isaev-2016-bis}
{\rm Isaev}}, A.:
{\em On the CR-curvature of Levi degenerate tube hypersurfaces},
{\tiny\sf arxiv.org/abs/1608.02919/}, 11 pages.

\smallskip

{\bf\bibitem{Isaev-2018}
{\rm Isaev}}, A.:
{\em Zero CR-curvature equations for Levi degenerate hypersurfaces
via Pocchiola's invariants}, 
{\tiny\sf arxiv.org/abs/1809.03029/}, 
Ann. Fac. Sci. Toulouse, to appear.

\smallskip

{\bf\bibitem{Isaev-2019}
{\rm Isaev}}, A.:
{\em Rigid Levi degenerate hypersurfaces with vanishing CR-curvature},
{\tiny\sf arxiv.org/abs/1901.03044/}, 10 pages.

\smallskip

{\bf\bibitem{Isaev-Zaitsev-2013}
{\rm Isaev}}, A.; {\rm Zaitsev}, D:
{\em Reduction of five-dimensional uniformly degenerate 
Levi CR structures to absolute parallelisms}, 
J. Anal. {\bf 23} (2013), no.~3, 1571--1605.

\smallskip

{\bf\bibitem{Katok-Hasselblatt-2002}
{\rm Katok}}, B.; {\rm Hasselblatt}, A.:
{\em The development of dynamics in the 20th century 
and the contribution of Jürgen Moser},
Ergodic Theory Dynam. Systems {\bf 22} (2002), no.~5, 1343--1364.

\smallskip

{\bf\bibitem{Kolar-Kossovskiy-2019}
{\rm Kolar}}, M.; {\rm Kossovskiy}, I.:
{\em A complete normal form for everywhere Levi-degenerate 
hypersurfaces in $\C^3$}, 
{\tiny\sf arxiv.org/abs/1905.05629/},
24 pages.

\smallskip

{\bf\bibitem{Jacobowitz-1990}
{\rm Jacobowitz}}, H.:
{\em An introduction to CR structures}, 
Math. Surveys and Monographs, 32, Amer. Math. Soc., Providence, 1990,
x+237~pp.

\smallskip

{\bf\bibitem{Libermann-1955}
{\rm Libermann}}, P.:
{\em Sur les structures presque complexes et autres structures
infinit\'esimales r\'eguli\`eres},
Bull. Soc. Math. France {\bf 83} (1955), 195--224.

\smallskip

{\bf\bibitem{Lie-Merker-2015}
{\rm Lie}}, S. (Author); {\rm Merker}, J. (Editor):
{\em Theory of Transformation Groups I. 
General Properties of Continuous Transformation Groups. 
A Contemporary Approach and Translation}, 
Springer-Verlag, Berlin, Heidelberg, 2015, xv+643~pp.
{\tiny\sf arxiv.org/abs/1003.3202/}

\smallskip

{\bf\bibitem{Loboda-2001}
{\rm Loboda}}, A.V.:
{\em Homogeneous strictly pseudoconvex hypersurfaces in C3 
with two-dimensional isotropy groups},
Sb. Math. {\bf 192} (2001), no.~11-12, 1741--1761.

\smallskip

{\bf\bibitem{Loboda-2003}
{\rm Loboda}}, A.V.:
{\em On the determination of a homogeneous strictly pseudoconvex 
hypersurface from the coefficients of its normal equation},
Math. Notes {\bf 73} (2003), no.~3-4, 419--423.

\smallskip

{\bf\bibitem{Loboda-2013}
{\rm Loboda}}, A.V.:
{\em Affinely homogeneous real hypersurfaces in $\C^2$},
Funct. Anal. Appl. {\bf 47} (2013), no.~2, 113--126.

\smallskip

{\bf\bibitem{Medori-Spiro-2014}
{\rm Medori}}, C.; {\rm Spiro}, A.: 
{\em The equivalence problem for 5-dimensional Levi degenerate 
CR manifolds}, Int. Math. Res. Not. IMRN {\bf 2014},
no.~20, 5602--5647. 

\smallskip

{\bf\bibitem{Medori-Spiro-2015}
{\rm Medori}}, C.; {\rm Spiro}, A.:
{\em Structure equations of Levi degenerate CR hypersurfaces 
of uniform type}. 
Rend. Semin. Mat. Univ. Politec. Torino {\bf 73} (2015), 
no.~1-2, 127--150. 

\smallskip

{\bf\bibitem{Merker-2005}
{\rm Merker}}, J.:
{\em The local geometry of generating submanifolds of
$\C^n$ and the analytic reflection principle}. I. 
(Russian) Sovrem. Mat. Prilozh. {\bf 6}, Kompleks. Anal. (2003), 
3--79; translation in J. Math. Sci. (N.Y.) {\bf 125} (2005), 
no.~6, 751--824.

\smallskip

{\bf\bibitem{Merker-2008}
{\rm Merker}}, J.:
{\em Lie symmetries of partial differential equations
and CR geometry}, Journal of Mathematical
Sciences (N.Y.), {\bf 154} (2008), 817--922.

\smallskip

{\bf\bibitem{Merker-2019}
{\rm Merker}}, J.:
{\em Affine rigidity without integration},
{\tiny\sf arxiv.org/abs/1903.00889/}, 28 pages.

\smallskip

{\bf\bibitem{Merker-Nurowski-2019}
{\rm Merker}}, J.; {\rm Nurowski}, P.:
{\em Equivalences of PDE systems associated to para-CR structures},
June 2019, 47 pages.

\smallskip

{\bf\bibitem{Merker-Pocchiola-2019}
{\rm Merker}}, J.; {\rm Pocchiola}, S.:
{\em Addendum to 
`Explicit absolute parallelism for $2$-nondegenerate real
hypersurfaces $M^5 \subset \C^3$ of constant Levi rank $1$'},
Journal of Geometric Analysis, to appear, 8 pages.

\smallskip

{\bf\bibitem{Merker-Pocchiola-2018}
{\rm Merker}}, J.; {\rm Pocchiola}, S.:
{\em Explicit absolute parallelism for $2$-nondegenerate real
hypersurfaces $M^5 \subset \C^3$ of constant Levi rank $1$},
Journal of Geometric Analysis, 
{\tiny\sf 10.1007/s12220-018-9988-3}, 
42~pp.

\smallskip

{\bf\bibitem{Merker-Pocchiola-Sabzevari-2013-5-CR-II}
{\rm Merker}}, J.; {\rm Pocchiola}, S.; {\rm Sabzevari}, M.:
{\em Equivalences of $5$-dimensional CR manifolds, II: 
General classes 
$\text{\sf I}$,
$\text{\sf II}$, 
$\text{\sf III}_{\text{\sf 1}}$,
$\text{\sf III}_{\text{\sf 2}}$,
$\text{\sf IV}_{\text{\sf 1}}$,
$\text{\sf IV}_{\text{\sf 2}}$},
5 figures, 95 pages, {\tiny\sf arxiv.org/abs/1311.5669/}

\smallskip

{\bf\bibitem{Nurowski-Tafel-1988}
{\rm Nurowski}}, P.; {\rm Tafel}, J.:
{\em Symmetries of Cauchy-Riemann spaces},
Lett. Math. Phys. {\bf 15} (1988), no.~1, 31--38. 

\smallskip

{\bf\bibitem{Olver-1995}
{\rm Olver}}, P.J.:
{\em Equivalence, Invariance and Symmetries}. Cambridge, Cambridge
University Press, 1995, xvi+525~pp.

\smallskip

{\bf\bibitem{Olver-2018}
{\rm Olver}}, P.J.:
{\em Normal forms for submanifolds under group actions},
Symmetries, differential equations and applications, 1--25.
Springer Proc. Math. Stat. {\bf 266}, Springer, Cham, 2018.

\smallskip

{\bf\bibitem{Olver-Pohjanpelto-2008}
{\rm Olver}}, P.J.; Pohjanpelto, J.:
{\em Moving frames for Lie pseudo-groups},
Canad. J. Math. {\bf 60} (2008), no.~6, 1336–1386.

\smallskip

{\bf\bibitem{Olver-Pohjanpelto-2009}
{\rm Olver}}, P.J.; Pohjanpelto, J.:
{\em Differential invariant algebras of Lie pseudo-groups},
Adv. Math. {\bf 222} (2009), no.~5, 1746–1792.

\smallskip

{\bf\bibitem{Pocchiola-2013}
{\rm Pocchiola}}, S.:
{\em Explicit absolute parallelism for 
$2$-nondegenerate real hypersurfaces
$M^5 \subset \C^3$ of constant Levi rank $1$}, 
55 pages, {\tiny\sf arxiv.org/abs/1312.6400/}

\smallskip

{\bf\bibitem{Poincare-1907}
{\rm Poincar\'e}}, H.:
{\em Les fonctions analytiques de deux variables complexes 
et la repr\'esentation conforme},
Rend. Circ. Mat. Palermo {\bf 23} (1907), 185--220.

}

\end{thebibliography}
\end{document}